\documentclass[a4paper,11pt,reqno]{amsart}
\usepackage[utf8]{inputenc}
\usepackage[T1]{fontenc}
\usepackage[mathscr]{eucal}
\usepackage{lmodern}
\usepackage[english]{babel}
\usepackage{microtype}
\usepackage{pdfsync}
\usepackage{import}
\usepackage{hyperref}
\usepackage{upgreek}
\usepackage{enumerate}
\usepackage{xspace}
\usepackage{pifont}
\usepackage{etoolbox}
\usepackage{relsize}
\usepackage{bbm}
\DeclareSymbolFont{myletters}{OML}{ztmcm}{m}{it}
\DeclareMathSymbol{\nicelambda}{\mathord}{myletters}{"15}
\usepackage{amsmath,amssymb,amsfonts,amsthm}
\usepackage{mathtools,accents}
\usepackage{mathrsfs}
\usepackage{aliascnt}
\usepackage{braket}
\usepackage{bm}
\usepackage{esint}
\usepackage{xfrac}    
\usepackage{nicefrac} 
\makeatletter
\g@addto@macro\@floatboxreset\centering
\makeatother
\patchcmd{\abstract}{\scshape\abstractname}{\textbf{\abstractname}}{}{}
\usepackage[x11names, dvipsnames, svgnames]{xcolor}
\definecolor{Burgundy}{RGB}{144,0,32}
\definecolor{Burgundy1}{RGB}{128,0,32}
\definecolor{Burgundy2}{RGB}{158,5,8}
\definecolor{VividBurgundy}{RGB}{159,29,53}
\definecolor{BlueBlue}{RGB}{10,10,120}
\definecolor{link-gray}{gray}{0.5}
\definecolor{ultrablue}{rgb}{0.0,0.0, 5}
\hypersetup{
	linktoc=page,
	colorlinks,
	linkcolor={ultrablue},
	citecolor={ultrablue},
	urlcolor={ultrablue}
}
\usepackage[a4paper,bindingoffset=0in,%
left=1.1in,right=1.1in,top=1.3in,bottom=1.3in,%
footskip=.2in]{geometry}
\usepackage{setspace}
\allowdisplaybreaks
\makeatletter 
\def\l@subsection{\@tocline{2}{0pt}{3pc}{6pc}{}}
\makeatother

\makeatletter
\g@addto@macro\@floatboxreset\centering
\makeatother

\makeatletter
\def\newaliasedtheorem#1[#2]#3{
	\newaliascnt{#1@alt}{#2}
	\newtheorem{#1}[#1@alt]{#3}
	\expandafter\newcommand\csname #1@altname\endcsname{#3}
}
\makeatother
\numberwithin{equation}{section}
\newtheoremstyle{slanted}{\topsep}{\topsep}{\slshape}{}{\bfseries}{.}{.5em}{}
\theoremstyle{plain}
\newtheorem{theorem}{Theorem}[section]
\newaliasedtheorem{proposition}[theorem]{Proposition}
\newaliasedtheorem{lemma}[theorem]{Lemma}
\newaliasedtheorem{corollary}[theorem]{Corollary}
\newaliasedtheorem{counterexample}[theorem]{Counterexample}
\theoremstyle{definition}
\newaliasedtheorem{definition}[theorem]{Definition}
\newaliasedtheorem{question}[theorem]{Question}
\newaliasedtheorem{assumption}[theorem]{Assumption}
\newaliasedtheorem{conjecture}[theorem]{Conjecture}
\theoremstyle{remark}
\newaliasedtheorem{remark}[theorem]{Remark}
\newaliasedtheorem{example}[theorem]{Example}

\newcommand{\T}{\mathcal T}

\newcommand{\Geod}{\textsf{Geod}}
\newcommand{\X}{{\rm X}}
\newcommand{\Y}{{\rm Y}}
\newcommand{\M}{{\rm M}}
\newcommand{\Q}{{\rm Q}}
\newcommand{\B}{{\rm B}}
\newcommand{\LL}{{\rm L}}
\newcommand{\vX}{\mathscr{X}}
\newcommand{\Sing}{\mathsf{S}}

\let\altphi\phi
\let\phi\varphi
\let\varphi\altphi
\let\altphi\undefined

\DeclareMathOperator{\Hess}{Hess}

\newcommand{\di}{\mathop{}\!\mathrm{d}}

\DeclareMathOperator{\supp}{supp}

\DeclareMathOperator{\Lip}{Lip}

\newcommand{\Prob}{\mathscr{P}}

\newcommand{\restr}{\raisebox{-.1908ex}{$\bigr\rvert$}}
\newcommand{\mres}{\mathop{\hbox{\vrule height 6pt width .5pt depth 0pt
			\vrule height .5pt width 6pt depth 0pt}}\nolimits\hspace{0.1pt}}

\newcommand{\ppi}{{\boldsymbol{\uppi}}}
\newcommand{\PP}{{\boldsymbol{\Uppi}}}

\newcommand{\dist}{\mathsf{d}}

\newcommand{\meas}{\mathfrak{m}}

\newcommand{\vol}{{\rm vol}}

\newcommand{\mx}{\mathsf{max}}

\DeclareMathOperator{\CD}{CD}
\DeclareMathOperator{\RCD}{RCD}
\DeclareMathOperator{\BE}{BE}

\newfont{\tmpf}{cmsy10 scaled 2500}

\newcommand{\dif}{\text{d}}
\DeclareMathOperator{\Ric}{Ric}
\DeclareMathOperator{\g}{g}

\def\XXint#1#2#3{{\setbox0=\hbox{$#1{#2#3}{\int}$ }
		\vcenter{\hbox{$#2#3$ }}\kern-.6\wd0}}

\newcommand{\K}{\mathcal{K}}
\newcommand{\N}{\mathcal{N}}
\newcommand{\W}{\mathcal{W}}
\newcommand{\D}{\mathrm{D}}
\newcommand{\Was}{\mathscr{W}}

\DeclareMathOperator{\length}{length}
\DeclareMathOperator{\sgn}{sgn}
\newcommand{\R}{\mathbb R}
\newcommand{\F}{\mathscr F}
\DeclareMathOperator{\diam}{diam}

\DeclareMathOperator*\lowlim{\underline{lim}}
\DeclareMathOperator{\Ent}{Ent}
\DeclareMathOperator*\uplim{\overline{\text{\raisebox{0pt}[1.15\height]{lim}}}}
\newcommand{\plmi}{\scalebox{0.7}{$\pm$}}
\DeclareMathSymbol{\shortminus}{\mathbin}{AMSa}{"39}

\newcommand{\ovs}[1]{\overline{\text{\raisebox{0pt}[1.15\height]{$#1$}}}}
\newcommand{\uphat}[1]{\widehat{\text{\raisebox{0pt}[1.1\height]{$#1$}}}}
\newcommand{\Dom}{\mathsf{Dom}}
\newcommand{\Diam}{\mathcal{D}}
\newcommand{\ai}{\mathsf a}
\newcommand{\bi}{\mathsf b}
\newcommand{\BERic}{\Ric^{\, \N}_{\substack{\\[1pt]\vX}}}
\newcommand{\BEL}{\Delta_{\substack{\\[1pt] \vX}}}
\newcommand{\CHL}{\Delta_{\substack{\\\hspace{1pt} \mathsf{CH}}}}
\newcommand{\CHE}{\mathcal{E}_{\substack{\\ \mathsf{CD}}}}
\newcommand{\BG}{{\Upgamma}}
\newcommand{\Sobol}{\mathcal{W}^{^{1,2}}}
\newcommand{\Lp}{\LL^{\hspace{-2pt}^p}}
\newcommand{\Ltwo}{\LL^{\hspace{-2pt}^2}}
\newcommand{\feig}{\nicelambda_{\hspace{0.5pt}_1}}
\newcommand{\ball}[1]{\mathrm{B}_{_{#1}}}
\newcommand{\e}{\mathsf{e}}
\newcommand{\markov}{\mathsf{P}_{\substack{\\\hspace{-3pt}t}}}
\newcommand{\EE}{\mathrm{E}}
\usepackage{accents}
\newcommand*{\dt}[1]{%
	\accentset{\mbox{\bfseries .}}{#1}}

\newcommand{\bdot}[1]{\dt{\text{\raisebox{0pt}[1.1\height]{$#1$}}}}

\AtBeginDocument{%
	\def\MR#1{}
}

\usepackage[outline]{contour}
\newcommand*{\fancy}[1]{\color{link-gray}\contour{black}{#1}}

\newcommand\fgamma{{\fancy{$\Gamma$}}}
\newcommand\ray{{\fancy{$r$}}}

\date{\today}
\title[Rigidity of spectral gap in non-negatively curved spaces]{The rigidity of sharp spectral gap \\ in non-negatively curved spaces}

\address{Christian Ketterer \newline Institute for Algebra and Geometry \\
	Karlsruhe Institute of Technologie \\
	Englerstr. 2\\
	76131 Karlsruhe, Germany}
\email{\href{mailto:christian.ketterer@kit.edu}{christian.ketterer@kit.edu}}
\address{ Yu Kitabeppu \newline Faculty of Advanced Science and Technology\\
	Kumamoto University\\
	Kumamoto 860-8555, Japan}
\email{\href{mailto:ybeppu@kumamoto-u.ac.jp}{ybeppu@kumamoto-u.ac.jp}}
\address{Sajjad Lakzian\newline Department of Mathematical Sciences\\
	Isfahan University of Technology \\
	 Isfahan Province Postal Code: 8415683111\\
	Iran
}
\email{\href{mailto:slakzian@iut.ac.ir}{slakzian@iut.ac.ir}}

\address{School of Mathematics\\
	Institute for Research in Fundamental Sciences (IPM) \\
	P.O. Box 19395-5746\\
	Iran
}
\email{\href{mailto:slakzian@ipm.ir}{slakzian@ipm.ir}}

\begin{document}
	\maketitle
\begin{center}
	\begin{minipage}{0.3\textwidth}
		\centering
		\bf{\small Christian Ketterer}\\
		\textit{\scriptsize University of Freiburg,}
		\textit{\scriptsize Germany}\\
		\textsf{\scriptsize \email{\href{mailto:christian.ketterer@math.uni-freiburg.de}{christian.ketterer@math.uni-freiburg.de}}}
	\end{minipage}
	\begin{minipage}{0.3\textwidth}
		\centering
		\bf{\small Yu Kitabeppu}\\
		\textit{\scriptsize Kumamoto University,}\\
		\textit{\scriptsize Japan}\\
		\textsf{\scriptsize \email{\href{mailto:ybeppu@kumamoto-u.ac.jp}{ybeppu@kumamoto-u.ac.jp}}}
	\end{minipage}
	\begin{minipage}{0.3\textwidth}
	 \centering
   	\bf	{ \small Sajjad Lakzian}\\
	\textit{\scriptsize  Isfahan University of Technology,}
	\textit{\scriptsize  Iran}\\
	\textsf{\scriptsize  \email{\href{mailto:slakzian@iut.ac.ir}{slakzian@iut.ac.ir}}}\\
\end{minipage}
\end{center}
\vspace{5mm}
%
	\begin{abstract}
		\par We extend the celebrated \emph{rigidity of the sharp first spectral gap under $\Ric \ge0$} to compact infinitesimally Hilbertian spaces with non-negative (weak, also called synthetic) Ricci curvature and bounded (synthetic) dimension i.e. to so-called compact $\RCD\hspace{-1pt}\left(0,\N \right)$ spaces; this is a category of metric measure spaces which in particular includes (Ricci) non-negatively curved Riemannian manifolds, Alexandrov spaces, Ricci limit spaces, Bakry-\'Emery manifolds along with products, certain quotients and measured Gromov-Hausdorff limits of such spaces. In precise terms, \emph{we show in such spaces, $\feig = \nicefrac{\uppi^{\hspace{0.5pt}2}}{\diam^2}$ if and only if the space is one dimensional with a constant density function}. We use new techniques mixing Sobolev theory and singular $1\D$-localization which might also be of independent interest. As a consequence of the rigidity in the singular setting, we also derive almost rigidity results.
	\end{abstract}
	\tableofcontents
\section{Introduction}\label{sec:intro}
\par One of the most prevailing endeavors in modern analysis has been the study of the spectrum of the Laplace operator in various spaces (or domains within) and under different geometric constraints and/or boundary conditions. Kac's famous ``hearing the shape of a drum'' conundrum -- or at least a restatement thereof -- asks under what restrictions, one can hear the shape of a drum; said more precisely, when is it the case that the spectrum of the Laplacian completely determines the underlying space? Spectral rigidity results provide partial answers to this question. 
\par Even though one can construct (locally) nonisometric isospectral (Riemannian) spaces (e.g. see~\cite{Milnor, Ejiri, Szab, Gordon}), estimates on the eigenvalues still carry plenty of geometric and analytic data about the underlying space.

\par Estimates on the first spectral gap are particularly important since they provide valuable estimates on, among other things, curvature bounds for the underlying space (geometric information) and the rate of convergence of diffusion processes to their equilibrium states (analytic information). 

\par In these notes we show the rigidity of the sharp spectral gap in compact $\RCD(0,\N)$ (or equivalently $\RCD^* \hspace{-2pt}\left(0,\N\right)$) spaces; these spaces include (Ricci) non-negatively curved Riemannian manifolds, Alexandrov spaces, Ricci limit spaces, Bakry-\'Emery manifolds along with products, certain quotients and measured Gromov-Hausdorff limits of such spaces. 
\subsection*{A glance into technicalities}
\par We will be considering the first (nonzero) eigenvalue $\feig$ of the Cheeger Laplacian (think of a counterpart for the first Neumann eigenvalue for a Riemannian domain) in a compact $\RCD^* \hspace{-2pt}\left(0,\N\right)$, or equivalently $\RCD\hspace{-1pt}\left(0,\N\right)$, metric measure space $\left(\X, \dist , \meas \right)$. Here, $ \meas$ is assumed to have full support. $\feig$ can equivalently be defined by the classic Rayleigh quotient minimization over locally Lipschitz functions or over the $\Sobol$-Sobolev space. 
\par The goal is to show that if $\feig = \nicefrac{\uppi^{\hspace{0.5pt}2}}{\diam^2}$, where $\diam$ denotes the diameter, then the underlying space must be a circle or a line segment. To clarify the context, we emphasize that here, the $\Sobol$-Sobolev space (see~\cite{AGSRiem, AGSbercd, AGMR}) means the space of functions admitting finite Cheeger-Dirichlet energy (i.e. $\Ltwo$-norm of smallest weak upper gradient is finite). There are various \emph{equivalent} notions of $(1,2)$-Sobolev spaces for metric measure spaces and we will later see a brief discussion about their equivalence. 
\par One of the big hurdles along the way, is to show Lipschitz regularity of $\left|\nabla u\right|$ so that we can apply a strong maximum principle to it. For this one needs to resort to the 1D-localization scheme in order to -- in a suitable sense -- reduce the question to the one dimensional case. 
\par Finishing the proof requires showing that the underlying metric measure space minus a closed singular subset of codimension at least $1$ splits off an interval isometrically; Notice that in the smooth setting and if $n\ge3$, the codimension of this closed singular set is at least $2$ but then aposteriori, the rigidity implies the case $n\ge 2$ cannot occur. 
\par To show the splitting phenomenon, we will construct a harmonic potential from the eigenfunction and use its gradient flow. Notice this splitting phenomenon is standard in Riemannian geometry due to the de Rham's decomposition theorem which plays a key role in the proof of celebrated the Cheeger-Gromoll's splitting theorem. Yet, in this context, we will need a weak  local version of such a splitting result and this will be achieved by a careful analysis of the aforementioned gradient flow.
\subsection*{The work leading up to ours}
\par It is virtually impossible to fully review the huge body of literature on this subject. So, we will only try to highlight the most important ones as they are directly pertinent to the problem considered in this paper; namely, the important results which deal with (sharp) spectral gap estimates under lower bounds on the Ricci curvature and upper bounds on the diameter.
\subsubsection*{\bf \textit{Spectral bounds}}
\par In the case of \emph{positive Ricci curvature}, it was shown in~\cite{Lich} that in compact manifolds, $\Ric \ge \K >0$ implies $\feig \ge \nicefrac{n}{\left(n - 1\right)}\, \K$; later, the rigidity result obtained in~\cite{Obata} proved the equality can only happen if the underlying manifold is the $n$-dimensional spherical space form $\mathbb{S}^{^n}_{\nicefrac{\K}{n-1}}$. This rigidity is now known as Obata's rigidity. 
\par In the much subtler setting of \emph{non-negative Ricci curvature}, it was shown in the celebrated papers~\cite{Li} and~\cite{Li-Yau} that in a compact manifold with non-negative Ricci curvature, $\feig \ge \nicefrac{\uppi^{\hspace{0.5pt}2}}{2\diam^2}$ holds. This estimate is obviously not sharp; for instance, it can be easily verified that $\feig = \nicefrac{\uppi^{\hspace{0.5pt}2}}{\diam^2}$ for the circle in $\R^{^2}$. The sharp first spectral gap (general case) for non-negative Ricci was obtained in~\cite{ZY} by improving upon the previous results in~\cite{Li} and~\cite{Li-Yau}. It was shown in~\cite{ZY} that indeed the sharp estimate $\feig \ge \nicefrac{\uppi^{\hspace{0.5pt}2}}{\diam^2}$ holds. It is worth mentioning that this sharp bound had previously been obtained  for convex Euclidean domains in~\cite{PW}.
\par In the 90's, the collected works in~\cite{Kro,BQ,Chen-Wang1,Chen-Wang2} unified the spectral gap results in the form of spectral comparison with $1$-dimensional model spaces which will be recalled in below. It is worth mentioning that the most inclusive result among the said references, are the ones presented in~\cite{BQ}. As a result, in an $\N$-dimensional compact weighted Riemannian manifold $\M^{^\N}$ with Ricci curvature tensor $\Ric \ge \K$ and $\diam \left(\M^{^\N}\right) \le \Diam$ without boundary (or with convex boundary components), the first nonzero (Neumann) eigenvalue $\feig$ satisfies 
\begin{align}\label{eq:eig-compar}
	\feig \ge {\nicelambda}\hspace{-2pt} \left( \K,\N,\Diam  \right),
\end{align}
where ${\nicelambda} \hspace{-2pt} \left( \K,\N,\Diam  \right)$ is the first (nonzero) eigenvalue of the $1$-dimensional model problem with Neumann boundary conditions
\begin{align}\label{eq:model-1}
	\begin{cases}
	\mathsf{L}_{_{\hspace{1pt}\K,\,\N}} \,v\left(x\right) =\shortminus  \nicelambda \, v\left(x\right)	\quad x \in \left( \shortminus  \nicefrac{\Diam}{2}, \nicefrac{\Diam}{2} \right),  \\[4pt]
		\raggedleft v^{'}\left(\shortminus \nicefrac{\Diam}{2} \right) = v^{'}\left( \nicefrac{\Diam}{2}  \right) = 0, & \\ 
	\end{cases}
\end{align}
in which, the drifted diffusion operator is given by
\begin{align*}
\mathsf{L}_{_{\hspace{1pt}\K,\,\N}} := \nicefrac{\dif^{^{2}}}{\dif x^{2}}\, \shortminus \mathsf{T}_{\hspace{-1pt}_{\K,\,\N}},
\end{align*}
where the \emph{drift term}, $\mathsf{T}_{\hspace{-1pt}_{\K,\,\N}}$, is given by
\begin{align}\label{eq:model-2}
\mathsf{T}_{\hspace{-1pt}_{\K,\,\N}}\left(x\right) := \begin{cases}
	\left(\N \shortminus 1\right)^{\hspace{-3pt}^{\nicefrac{1}{2}}}  \K^{^{\nicefrac{1}{2}}} \ \tan \left(  \left( \N \shortminus 1 \right)^{\hspace{-3pt}^{\shortminus \nicefrac{1}{2}}}\K^{^{\nicefrac{1}{2}}} \, x \right)&\text{if } \K > 0  , 1 < \N < \infty, \\[4pt]
		\shortminus \left( \N \shortminus 1 \right)^{\hspace{-3pt}^{\nicefrac{1}{2}}}\left| \K \right|^{\hspace{-1pt}^{\nicefrac{1}{2}}}  \ \tanh \left(\left( \N \shortminus 1 \right)^{\hspace{-2pt}^{\shortminus \nicefrac{1}{2}}}\left| \K \right|^{\hspace{-1pt}^{\nicefrac{1}{2}}}\, x \right)&\text{if } \K < 0  , 1 < \N < \infty,\\[4pt]
		0 &\text{if } \K = 0  , 1 < \N <\infty ,\\[4pt]
		\K x &\text{if } \N = \infty.
\end{cases} 
\end{align}
\par In the Riemannian setting, this comparison result has since been generalized to the $n$-dimensional weighted manifolds $\M^{^n}$ (without boundary or with convex boundary components) with the \emph{$\N$-Bakry-\'Emery Ricci tensor $\Ric^{\,\N} \ge \K$ ($n \le \N$)} in~\cite{BQ,AC,Mil-1}. Indeed, by~\cite{BQ}, this comparison holds for the first \emph{real} eigenvalue of a smooth elliptic operator $L=\Delta +  B$ and under Bakry-\'Emery conditions  $\BE\hspace{-1pt}\left( \K,\N \right)$ (in~\cite{BQ}, these are called ``$CD\hspace{-1pt}\left( \K,\N \right)$'' conditions not to be mistaken with the (not unrelated) $\CD\hspace{-1pt}\left( \K,\N\right)$ conditions we will encounter later which are defined via $\left(\K,\N \right)$-convexity of the entropy functional along $\Ltwo$-optimal transportation of probability measures). 
\par In the setting of Alexandrov spaces, the said comparison estimates for $\feig$ have been obtained in  \emph{$\N$-dimensional  Alexandrov spaces without boundary} in~\cite{QZZ}.
\par The generalization of the above comparison results to  $n$-dimensional Finsler manifolds (without boundary or with a convex one) with \emph{weighted Ricci bound $\Ric_{_\N} \ge \K$ ($n \le \N$)} can be found in~\cite{WX}.
\par In the context of metric measure spaces satisfying reduced weak Riemannian Ricci curvature bounds, also known as $\RCD^* \hspace{-2pt}\left( \K,\N \right)$ (RCD: Riemannian curvature dimension) conditions, the above spectral comparison estimates have recently been obtained in~\cite{KeObata,Jiang-Zhang,CM-1,CM-2}. 
\par So far, the most general comparison result in the setting of metric measure spaces is the generalization of the spectral comparison estimates proven in~\cite{CM-1,CM-2} for the $p\hspace{0.5pt}$-spectral gap $\nicelambda^{\hspace{-2pt}^{1,p}}_{_{\left(\X, \dist , \meas   \right)}}$ for essentially non-branching $\mathrm{CD}^* \hspace{-2pt}\left( \K,\N \right)$ spaces $\left(\X,\dist , \meas\right)$; and proven for $\mathrm{MCP}^*\hspace{-2pt}\left( \K,\N \right)$ spaces in~\cite{BX}. Notice that in an $\mathrm{MCP}^* \hspace{-2pt}\left( \K,\N \right)$ space, the $p\hspace{0.5pt}$-spectral gap can merely be defined as the infimum of the $\Lp$-Rayleigh quotient and no Laplacian can be used as of yet. For $p=2$, this comparison result coincides with the one we stated earlier. Recall that when the underlying space is $\RCD\hspace{-1pt}\left( \K,\N \right)$, the $p\hspace{1pt}$-spectral gap $\nicelambda^{\hspace{-2pt}^{1,p}}_{_{\left(\X, \dist , \meas   \right)}}$ coincides with the first nonzero eigenvalue of the $p\hspace{1pt}$-Laplacian. The proofs in~\cite{CM-1,CM-2} cleverly use the generalization of $1\D\hspace{0.5pt}$-localization of curvature dimension conditions to metric measure spaces (developed earlier in the smooth setting in~\cite{klartagneedle}) hence are different in nature from the techniques used by other authors. In this article, we also use $1\D\hspace{0.5pt}$-localization in an essential way. 
\subsubsection*{\bf \textit{Rigidity results}}
\par  Obata's rigidity theorem in~\cite{Obata} is the first result that characterizes the equality case of the sharp spectral gap estimates (when $\K>0$). In general (i.e. when the underlying space is not necessarily smooth), for $\K>0$, the rigidity result states that $\feig = {\nicelambda} \left( \K,n, \mathcal{D} \right)$ if and only if the underlying space is a \emph{spherical suspension}; see~\cite{KeObata,Jiang-Zhang,CM-2,GKK}.
\par For compact Riemannian manifolds $\M^{^n}$ with $\Ric \ge 0$ (the case $\K=0$), the characterization of the equality case of the sharp spectral gap estimates is obtained in~\cite{Hang-Wang}. In~\cite{Hang-Wang}, it was shown that in a closed manifold $\M^{^n}$, $\feig = \nicefrac{\uppi^{\hspace{0.5pt}2}}{\diam^2}$ holds if and only if  $n=1$ and $\M^{^1}$ is a circle of perimeter $2\uppi$; this is proven in~\cite{Hang-Wang} by showing the equality implies the space is (minus two points) one or two copies of the open interval $\left( \shortminus \nicefrac{\uppi}{2} , \nicefrac{\uppi}{2} \right)$; this is the desired splitting phenomenon that we wish to verify in a much more general context.
\par In Finsler structures with non-negative weighted Ricci curvature, $\Ric_{_\N} \ge 0$, such a rigidity result has recently been proven in~\cite{Xia-rig}. 
\par In the Riemannian setting, the spectral comparison for $p\hspace{1pt}$-Laplacian and the equality case have been considered, among other places, in~\cite{Valtorta,NV}.
\par It is worth mentioning that an interesting rigidity result for $\nicelambda^{\hspace{-2pt}^{1,p}}_{_{\left(\X, \dist , \meas   \right)}}$ (not as an eigenvalue of a $p\hspace{1pt}$-Laplacian but only as the minimum of $\Lp$-Rayleigh quotient), in $\mathrm{MCP}\hspace{-1pt}\left(\K,\N\right)$ spaces and for $\K \le 0$, has been obtained in~\cite{BX}. For some exciting recent work about sharp and rigid inequalities for $\nicelambda^{\hspace{-2pt}^{1,p}}_{_{\left(\X, \dist , \meas   \right)}}$ on bounded domains, we refer the reader to~\cite{MV, MS, APPS}. 
\subsection*{Main result}
\par In these notes, we will mainly attend to the equality case of the sharp spectral gap estimates in the general setting of metric measure spaces satisfying $\RCD^* \hspace{-2pt}\left( 0,\N \right)$ or equivalently $\RCD\hspace{-1pt}\left( 0,\N \right)$ conditions. These spaces could be nonsmooth yet are, in many ways, reminiscent of Riemannian manifolds with non-negative Ricci curvature and dimension less than or equal to $\N$. 
\par Let
\begin{align*}
	\feig&:=\inf\left\{ \left\| u \right\|^{^{-2}}_{\substack{\\\Ltwo}}\, \left\| \nabla u \right\|^{^{2}}_{\substack{\\\Ltwo}} \;\; \text{\textbrokenbar} \;\; u\in \Lip\left( \X \right)\cap \Ltwo\left(\meas\right) \smallsetminus \left\{0\right\}\; \text{and}\; \scalebox{1.2}{$\int$}_{\hspace{-5pt}_{\X}} u\, \dif\meas = 0\right\}\notag \\
	&= \inf \left\{ \left\| \nabla u \right\|^{^{2}}_{\substack{\\\Ltwo}}\;\; \text{\textbrokenbar} \;\; u \in \Lip\left(\X\right)\cap \Ltwo\left(\meas\right) \smallsetminus \left\{1\right\} \; \text{and}\; \ \left\| u \right\|_{\substack{\\\Ltwo}} = 1 \right\}. 
\end{align*}
Here $\left\|\nabla u\right\|_{\substack{\\\Ltwo}} $ is a notation for twice the {\it Cheeger-Dirichlet energy} $\CHE\left( u\right)$ of $u$ on $\left( \X,\dist, \meas \right)$ and also is the $\Ltwo$-norm of $|\nabla u|$ once one makes sense of the norms of gradients. 
Our main theorem is the following generalization of the rigidity statement for the equality case of the sharp spectral gap. Let $\diam$ denote the diameter of the whole space $\X$. 
\begin{theorem}\label{thm:main-YZ}
	Suppose $\left(\X, \dist,\meas \right)$ is a compact $\RCD\hspace{-1pt}\left(0,\N\right)$ space with $\supp \left( \meas \right)  = \X$.  $\feig = \nicefrac{\uppi^{\hspace{0.5pt}2}}{\diam^2}$ if and only if $\X$ is either a weighted circle or a weighted line segment; In either cases, the space is equipped with a constant weight function i.e. $\meas = c \mathcal{H}^1$ (in other words, $X$ is a non-collapsed 1D RCD space). 
\end{theorem}

\par A special case of this is the rigidity for the Bakry-\'Emery Ricci tensor. In~\textsection\thinspace\ref{app:Drif-Laplacian}, we have discussed this along with an alternative proof that perhaps will b more favorable to those who are more differential-geometrically inclined.
\subsection*{Road map to the proof}
\par As is the nature of the field, since we are only allowed to use the tools currently available in the setting of $\RCD\hspace{-1pt}\left(\K,\N\right)$ spaces, we inevitably have to adapt many arguments to the setting of metric measure spaces; that is to say, eventually all the arguments have to only depend on the distance $\dist$ and measure $\meas$. This means the proof becomes substantially more complicated and lengthier, and is not a derivative or direct modification of smooth techniques; this is substantiated by the novel  techniques needed to be developed such as a mixture of  Sobolev theory and $1\D$ localization which is on its own, of interest to experts. 
\par Some of the tools we need to employ, are recent and new technology in metric measure spaces with lower weak Ricci curvature bounds, the most noteworthy among which, are the 1D localization scheme (see~\cite{Fabio-decomp, CM-1, CM-2}), a localized version of nonsmooth splitting theorem (see~\cite{Gsplit}), calculus of the square field operators and of the measure-valued Laplacian (see~\cite{Sav-1}) and a new second variation formula (see~\cite{GT-2}). 
\par To provide a road map, in below, we highlight the main steps of the proof of our main result, Theorem~\ref{thm:main-YZ}.
\begin{itemize}
	\item[\ding{202}] By rescaling the metric $\dist$ we can assume $\diam \left( \X \right) = \uppi$ and $\feig = 1$. Then, we pick an eigenfunction $u$ associated to $\feig = 1$. By using the heat flow and a Sobolev-to-Lipschitz property, we will argue that $u$ has a Lipschitz representative again denoted by $u$.
	\smallskip
	\item[\ding{203}] Applying the $1\D$-localization scheme, we will reduce some of the analysis to the 1D case; in particular, we prove constancy of a crucial quantity $\upalpha = u^2 + \left| \nabla u  \right|^{^2}$ (the density of the Sobolev norm) along transport geodesics in the so-called transport set.
	\smallskip
	\item[\ding{204}]  Taking advantage of the calculus of measured-valued Laplacian along with some $\BG$-calculus, we show $\Sobol$-regularity of $\left| \nabla u \right|^{^2}$. 
	\smallskip
	\item[\ding{205}] By means of the first variation formula (see~\cite{GT}) along with some $\BG$-calculus computations and by using a Sobolev-to-Lipschitz property, we show $\upalpha := u^2 + \left| \nabla u \right|^{^2} \in \Sobol$ is subharmonic and that \emph{$\left| \nabla u \right|^{^2}$ indeed has a Lipschitz representative}. 
	\smallskip
	\item[\ding{206}] Using the strong maximum principle due to~\cite{Gigli-Rigoni} and since $\upalpha$ (which has been proven to be Lipschitz in the previous steps) is subharmonic, $\X$ is compact and $\upalpha$ achieves an interior maximum, we deduce $\upalpha$ is identically equal to $1$. This implies the weak harmonicity of $f := \sin^{\hspace{-1pt}^{-1}} \circ \ u$ and indeed, the stronger  fact that the contracted Hessian of $f$ vanishes; it also implies the identity $\left| \nabla f \right| \equiv 1$. 
	\smallskip
	\item[\ding{207}] With arguments that are similar in essence to, and are inspired by the nonsmooth splitting theorem (see~\cite{Gsplit,Goverview}) and also by using some analytic techniques from~\cite{GKK}, we will prove that the \emph{harmonicity} of $f$ results in the gradient flow of $f$ being measure-preserving and that roughly speaking, $f$ being ``Hessian-free'' implies the gradient flow of $f$ is distance-preserving.
	\smallskip
	\item[\ding{208}]  Defining a suitable quotient space $\left( \widetilde{\Y}, \dist_{\substack{\\ \widetilde Y}} , \meas_{\substack{\\ \widetilde{Y} }} \right)$ and the singular sets $\Sing_{\substack{\\\pm} } := u^{-1}\left( \plmi 1  \right) = f^{^{-1}}\left( \plmi \nicefrac{\uppi}{2} \right)$, we will demonstrate that the metric completion of the regular set $\widetilde{\X} := {\sf MetCompl} \left( \X \smallsetminus \left( \Sing_{\substack{\\-}}\, \dot{\sqcup}\, \Sing_{\substack{\\+}} \right) \right)$ splits isometrically as
	\begin{align*}
	\left( \widetilde{\X}, \dist_{\tilde{\X}}, \meas_{\tilde{\X}}  \right) \stackrel{\mathsf{isom}}{\cong}  \left(  \widetilde{Y} \times \left[ \shortminus \nicefrac{\uppi}{2},  \nicefrac{\uppi}{2}  \right] ,  \dist_{\substack{\\ \widetilde \Y}}  \oplus \dist_{_\textsf{Euc}} , \meas_{\substack{\\ \widetilde \Y}}  \otimes \mathscr{L}^{^1}\mres_{\left[ \shortminus \nicefrac{\uppi}{2} , \nicefrac{\uppi}{2}  \right]}  \right),
	\end{align*}
	and that $\widetilde \Y$ in fact satisfies the $\RCD\hspace{-1pt}\left(0,\N-1\right)$ conditions. Note the distances used are induced ones. 
	\smallskip
	\item[\ding{209}] To finish the proof, we will then show that either $\widetilde \X$ admits nonempty \emph{one dimensional regular set} $\mathcal{R}_1\neq \emptyset$ (not to be confused with the regular set $\X \smallsetminus \Sing$) or it is isometric to $\X$; this is achieved by showing that the the regular set $\X \smallsetminus \Sing$ admits a foliation by non-branching so-called horizontal geodesics. The former will lead to $\X$ being one dimensional (by the characterization of low-dimensional RCD spaces~\cite{KL-2}) and the former will lead to a contradiction by using the metric \emph{Pythagorean theorem} for the aforementioned isometric splitting.  This means -- again by~\cite{KL-2} --  $\X$ has to be a circle (with constant weight) or a line segment with a $\left(\K, \N\right)$-convex weight function. In the latter case, a simple calculation shows that since the sharp lower bound on the eigenvalue is achieved, the weight function must be constant. 
\end{itemize}

\subsubsection*{ \small \bf \textit{Terminology}} Sometimes to wit, we say a space is spectrally-extremal whenever it has non-negative Ricci curvature (in the sense which will be clear from the context) and that the bottom of the real (Neumann whenever applicable) spectrum of the associated diffusion operator (which will be known from the context) satisfies $\feig = \nicefrac{\uppi^{\hspace{0.5pt}2}}{\diam^2}$. So our main results could be summarized as \emph{spectrally-extremal spaces are $1$-dimensional and with constant weight}. 
\subsection*{Main Consequences}
\par Let us briefly touch upon corollaries of Theorem~\ref{thm:main-YZ}. 
\begin{corollary}[Almost rigidity]\label{cor:almost-rigidity} For every $\varepsilon>0$ and $\N\geq 1$ there exists $\updelta>0$ such that the following holds.
\par For all compact Riemannian manifolds $\left(\M^{^n},\g \right)$ with $n \leq \N$, $\Ric_{\g} \geq \shortminus \updelta$, $\feig\leq \uppi^2+\updelta$ and $\diam\leq 1$, it holds
	\begin{align*}\text{$\dist_{\mathsf{GH}}\left(  \left( \M^{^n},\dist_{\g} \right), \nicefrac{1}{\uppi} \ \mathbb S^{^1} \right)\leq \varepsilon$ \quad  \text{or} \quad  $\dist_{\mathsf{GH}}\left( \left( \M^{^n},\dist_{\g} \right), \left[0,1\right] \right)\leq \varepsilon$.}
	\end{align*}
\end{corollary}
\begin{proof}[\footnotesize \textbf{Proof}]
Suppose not. Then, for some $\varepsilon>0$ and $\N\geq 1$ fixed, there exists a sequence $\left(\M_i^{^n},\g_{\substack{\\i}} \right)$ defying the conclusion. Passing to a converging subsequence, one gets a $(0,n)$-Ricci limit space,$\left(\X, \dist_{\substack{\\\X}}\right)$, with $\diam \le 1$ and $\feig \le \uppi^2$ and with
\begin{align}\label{eq:almost-rig}
	\text{$\dist_{\textsf{GH}}\left( \left(\X, \dist_{\substack{\\\X}}\right), \nicefrac{1}{\uppi} \ \mathbb S^{^1} \right)\geq \varepsilon$ \quad  \text{or} \quad  $\dist_{\textsf{GH}}\big( \left(\X, \dist_{\substack{\\\X}}\right) , \left[0,1\right] \big)\geq \varepsilon$.}
\end{align}
This Ricci limit space $\X$ must be one dimensional by Theorem~\ref{thm:main-YZ} since it achieves the sharp gap; note that stability of the spectrum under Gromov-Hausdorff convergence was proven in \cite{GMS, CC3}.  In the light of our rigidity result, this is in contradiction with \eqref{eq:almost-rig}. 
\end{proof}
\begin{corollary}[Stronger rigidity in the closed noncollapsed case]\label{cor:diff1} For every $\varepsilon>0$ and $v>0$, there exists $\updelta>0$ such that the following holds.
\par Let $\left( \M^{^n},\g \right)$ be a compact Riemannian manifold without boundary such that $\Ric_{\g}\geq \shortminus \updelta$, $\feig \leq \uppi^2 + \updelta$,  $\diam\leq 1$ and $\mathrm{vol}_{\g}\left( \M^{^n} \right)\geq v>0$. Then $n=1$ and $\M^{^1}$ is diffeomorphic to $\nicefrac{1}{\uppi} \ \mathbb S^{^1}$.
\end{corollary}
	\begin{proof}[\footnotesize \textbf{Proof}]
		We first show that the second case in Corollary \ref{cor:almost-rigidity} cannot occur; indeed, otherwise one can find a sequence of  compact Riemmannian manifolds with a uniform lower Ricci curvature bound that converges in Gromov-Hausdorff sense to $\left[0,1\right]$. However, since the sequence is noncollapsed, the limit admits no points that have blow-up tangent cones isometric to $\left[0,\infty\right)$; see~\cite{CC1}. This is a contradiction and hence the first case in Corollary~\ref{cor:almost-rigidity} must occur.  Then, by the well-known stability theorem due to Cheeger and Colding (see~\cite{CC1}), there exists an $\varepsilon$ sufficiently small such that $\M^{^n}$ is diffeomorphic to $\mathbb S^{^1}$ when $\dist_{\textsf{GH}}\left( \left( \M^{^n},\dist_{\g} \right), \nicefrac{1}{\uppi} \ \mathbb S^{^1} \right)\leq \varepsilon$.
\end{proof}
\begin{corollary}[Stronger rigidity in noncollapsed case with boundary]\label{cor:diff} For every $\varepsilon>0$ and $v>0$, there exists $\updelta>0$ such that the following holds.
\par  Let $\left( \M^{^n}, \g \right)$ be a compact Riemannian manifold with boundary which verifies $\RCD\left(\shortminus \updelta,n\right)$, $\feig\leq \uppi^2+\updelta$,  $\diam\leq 1$ and $\vol_{\g}\left(\M^{^n}\right)\geq v>0$. Then $n=1$ and $\M^{^1}$ is diffeomorphic to $\left[0,1\right]$.
\end{corollary}
	\begin{proof}[\footnotesize \textbf{Proof}]
		The corollary follows as before by taking into account that  a noncollapsed Gromov-Hausdorff limit of Riemannian manifolds with boundary that satisfy the condition $\RCD\hspace{-1pt}\left(\shortminus \K,n\right)$ for some $\K>0$ has to have intrinsic boundary by the new developments in~\cite{BNS}.
\end{proof}
The following is a small scale first gap estimate which is a near sharp lower bound for small balls in a Riemannian manifold with negative lower Ricci curvature bound. 
\begin{corollary}[First gap estimate for small balls in a general Riemannian manifold]
	For every $\varepsilon>0$, $v>0$ and $n\in \mathbb N$ with $n\geq 2$, there exists $\updelta>0$ such that the following holds.
	\par Let $\left(\M^{^n},\g\right)$ be a Riemannian manifold such that $\Ric_{\g}\geq \shortminus \left(n \shortminus 1\right)$,  $p\in \M$ and $\vol_{\g}\left(\ball{1}\left(p\right)\right)\geq v >0$. Then for any $\upeta \in \left(0,\updelta\right)$ such that $\ovs{\ball{\upeta}\left(p\right)}$ is geodesically convex, it follows that $\feig \left( \ovs{\ball{\upeta}\left(p\right)} \right)\geq \nicefrac{\uppi^{\hspace{0.5pt}2}}{\mathrm{d}_{p,\upeta}^{\hspace{-1pt}^2}} +\varepsilon$ where  $\mathrm{d}_{p,\upeta}=\diam \left( \ovs{\ball{\upeta}\left(p\right)} \right)$. 
\end{corollary}
\begin{remark}
	Notice the above Corollary is a local result in nature and upon rescaling the metric, we can find an estimate for small balls and any lower bound on their Ricci curvature within them. 
\end{remark}
	\begin{proof}[\footnotesize \textbf{Proof}]
		We argue by contradiction. Let $\left( \M^{^n}_i,\g_{\substack{\\[1pt]i}} \right)$ with $\Ric_{\g_{_i}}\geq \shortminus \left( n \shortminus 1 \right)$. Suppose there exist $p_{\substack{\\[3pt]\hspace{-2pt}i}}\in \M^{^n}_i$ and $\upeta_{\substack{\\[2pt]i}}>0$ such that $\ovs{\ball{\upeta_{\substack{\\[1pt]i}}}\left(p\right)}$ is geodesically convex and $\feig\left( \ovs{\ball{\upeta_{\substack{\\[1pt]i}}}\left(p_{\substack{\\[3pt]\hspace{-2pt}i}}\right)} \right)\rightarrow \nicefrac{\uppi^{\hspace{0.5pt}2}}{\mathrm{d}_i^{^2}}$ where $\mathrm{d}_{\substack{\\[1pt]i}}= \diam \ovs{\ball{\upeta_{\substack{\\[1pt]i}}}\left(p_{\substack{\\[3pt]\hspace{-2pt}i}}\right)}\leq 2\upeta_{\substack{\\[2pt]i}}$. 
		\par  We rescale $\g_{\substack{\\[1pt]i}}$ by $\upeta_{\substack{\\[2pt]i}}^{\hspace{-2pt}^{-2}}$, and obtain a sequence of geodesically convex balls $\left( \ovs{\ball{1}\left( p_{\substack{\\[3pt]\hspace{-2pt}i}} \right)}, \upeta_{\substack{\\[2pt]i}}^{\hspace{-2pt}^{-2}} \g_{\substack{\\[1pt]i}} \right)$ with $\Ric_{\upeta_{\substack{\\[1pt]i}}^{-2}\g_{\substack{\\[1pt]i}}}\geq \upeta_{\substack{\\[1pt]i}}^{^{2}}\left(n \shortminus 1\right)\rightarrow 0$. 
		\par Then based on, by now standard, compactness theorem, a subsequence converges to an $\RCD\left(0,n\right)$ space $\left( \X,\dist,\meas \right)$; by the convergence of the spectrum of $\RCD$ spaces under (pointed) measured GH convergence (see~\cite{GMS}), it follows $\feig\left(\X\right)=\nicefrac{\uppi^{^2}}{\mathrm{d}^{^2}}$ where $\mathrm{d}= \diam \left( \X \right) \leq 2$. 
	\par Hence, $\X$ is diffeomorphic to $\mathbb S^{^1}$ by Corollary \ref{cor:diff}. But for $i$ large enough, $\X$ is also homeomorphic  to $\ovs{\ball{1}\left( p_{\substack{\\[3pt]\hspace{-2pt}i}} \right)}$ due to the Reifenberg theorem for noncollapsing Ricci limit spaces; see~\cite{CC1}. This is a contradiction. 
	\end{proof}
\par Analytic interpretation of Theorem~\ref{thm:main-YZ} is the following estimate on the rate of convergence, to their harmonic states, for diffusion processes generated by $\CHL$.  
\par Here by the rate of convergence, we mean the rate of strong convergence  to the invariant measure of the solution to the corresponding Fokker-Planck equation (also called continuity equation); alternatively, one can consider the rate of convergence to the equilibrium of the corresponding Feller process on continuous functions. This rate is given by the dominant (non-zero) eigenvalue of the Fokker-Planck generator which is dual to the operator which generates the process. Also recall, in the smooth setting, the  Neumann boundary condition translates to the processes being elastically reflected at the boundary.
\begin{corollary}\label{cor:diffusion-rate}
	Let $\left(\X , \dist, \meas\right)$ be a compact $\RCD\hspace{-1pt}\left( 0,\N \right)$ metric measure space. If the Hausdorff dimension $\dim_{\substack{\\\mathsf{Hauss}}}\left( \X \right) > 1$ or if $\X$ is one-dimensional which is not noncollapsed (i.e. $\X$ is not an interval with a nonconstant weight function; see~\cite{KL-2}) then the rate of the convergence to equilibrium of the diffusion process generated by $\CHL$ in $\X$ is \textsf{strictly larger} than the rate of the convergence to equilibrium of the diffusion process (elastically reflected at the boundary) generated by the operator $\nicefrac{\dif^{^{\,2}}}{\dif x^2}$ in $\left[ \shortminus \nicefrac{\uppi}{2},  \nicefrac{\uppi}{2} \right]$ (precisely, the first Neumann gap in $1\D$ model space). 
\end{corollary}
Alternatively, based on Theorem~\ref{thm:main-2}, we get the following.
\begin{corollary}
Let $\vX$ be a nonzero complete vector field on a compact manifold $\M^{^n}$ (without or with convex boundary components). If $\BERic \ge 0$ for some $n < \N<\infty$, then the rate of convergence to equilibrium of the diffusion process (elastically reflected at the boundary) generated by $\BEL := \Delta \shortminus \mathcal{L}_{\substack{\\\mathscr{X}}},$ is \textsf{strictly larger} than the rate of the convergence to equilibrium of the Brownian motion (elastically reflected at the boundary) in $\left[ \shortminus \nicefrac{\uppi}{2},  \nicefrac{\uppi}{2} \right]$.
\end{corollary}
\begin{remark}
	Intuitively, the more the Hausdorff dimension is, the more room is available for the Brownian motion to diffuse; consequently, the faster the said diffusion should reach its steady state. Also existence of a drift could speed up the diffusion process. This intuition is now precisely captured by Corollary~\ref{cor:diffusion-rate} that shows \emph{nontrivial drift and any larger dimension strictly increase the rate of convergence to equilibrium. }
\end{remark}
\addtocontents{toc}{\protect\setcounter{tocdepth}{-1}}
\section*{\small \bf  Organization of the Materials}
\addtocontents{toc}{\protect\setcounter{tocdepth}{1}}
\small
\par The presentation in these notes is organized as follows. In~\textsection\thinspace\ref{sec:intro}, we have set up and motivated the problem at hand and stated our main theorem as well as a few important corollaries along with brief proofs of the corollaries. \textsection\thinspace\ref{sec:prelim} reviews some recent materials about the geometry of metric measure spaces; in this section, we have highlighted the main concepts and the technology that is needed for reading our work.  In ~\textsection\thinspace\ref{sec:rcd-specanal}, we embark on the analysis of the eigenfunctions of the Cheeger Laplacian in the $\RCD$ setting in particular in non-negatively curved such spaces and prove some crucial regularity results. \textsection\thinspace\ref{sec:grad-flow} ushers in the required analysis of the gradient flow of the harmonic potential that is needed for the splitting phenomenon. The splitting phenomenon is discussed and proven in~\textsection\thinspace\ref{sec:isom-split}. In conclusion, \textsection\thinspace\ref{sec:proof-main} presents the culmination of the proof of our main theorem;~\textsection\thinspace\ref{app:Drif-Laplacian} discusses an alternative proof for the special case of the rigidity for sooth measure spaces under non-negative Bakry-\'Emery Ricci tensor (both the boundary-less case and the case with convex boundary); this alternative proof is new and uses less machinery.
\normalsize

\addtocontents{toc}{\protect\setcounter{tocdepth}{-1}}
\section*{\small \bf  Acknowledgements}
\addtocontents{toc}{\protect\setcounter{tocdepth}{1}}
\vspace{-80pt}
\begin{minipage}[t][7cm][b]{0.9\textwidth}
\small
\begin{itemize}
	\renewcommand{\labelitemi}{\scalebox{1.2}{$\centerdot$}}
	\item We would like to express our gratitude to the anonymous referees for their valuable comments and for bringing new relevant references to our attention;\smallskip
	\item We would like to extend our gratitude to Nicola Gigli for providing us with insightful remarks and for his interest in our work; \smallskip
	\item CK worked on this project while he was a postdoc in the Department of Mathematics at the University of Toronto. CK was funded by the Deutsche Forschungsgemeinschaft (DFG) - Projektnummer 396662902, ``Synthetische Kr\"ummungsschranken durch Methoden des optimal Transports'';\smallskip
	\item YK acknowledges the supports of the Grand-in-Aid for Young Scientists of 18K13412; 
	\item SL acknowledges partial support from IPM, Grant No. 1400460424.
\end{itemize}
\end{minipage}
	\normalsize
\section{Preliminaries}\label{sec:prelim}
\subsection{Some details on curvature dimension conditions using optimal transport theory}  
\par Let $\left( \X,\dist,\meas \right)$ be a complete separable geodesic metric measure space and $\mathscr{P}\left(\X\right)$, the set of all Borel probability measures. For $p\ge1$, we denote by $\Prob_{\hspace{-5pt}_p}\left(\X\right)$, the set of all Borel probability measures with finite $p\hspace{1pt}$-th moments. 
\par For any $\upmu_{\hspace{1pt}_0}\,,\upmu_{\hspace{1pt}_1}\in \Prob_{\hspace{-5pt}_p} \left(\X\right)$, the $\Lp$-Wasserstein distance is defined as\\
\begin{align} \label{eq:L2W}
\Was_{\hspace{-3pt}_p}\left( \upmu_{\hspace{1pt}_0},\upmu_{\hspace{1pt}_1} \right):=\inf\left\{ \left\| \dist  \right\|^{^p}_{\substack{\\[1pt]\Lp \left( \X^{^2}, \, \mathfrak{q} \right)}}\; \text{\textbrokenbar}\;\; {\frak q}\text{ is a coupling of $\upmu_{\hspace{1pt}_0}$ and $\upmu_{\hspace{1pt}_1}$} \right\}^{\nicefrac{1}{p}}.  
\end{align}
\par For every complete separable geodesic space $\left( \X,\dist \right)$, the $\Lp$-Wasserstein space
\begin{align*}
\left( \Prob_{\hspace{-5pt}_p} (\X), \Was_{\hspace{-3pt}_p} \right) =:\Prob_{\hspace{-5pt}_p} \left( \X,\dist \right),
\end{align*}
is also a complete separable geodesic space.
The space $\Prob_{\hspace{-5pt}_2}\left( \X,\dist,\meas \right)$ is defined as the subspace of measures in $\Prob_{\hspace{-5pt}_2}\left( \X,\dist \right)$ that are absolutely continuous w.r.t. $\meas$.
\par We denote the space of all constant speed geodesics from $\left[0,1\right]$ to $\left( \X,\dist \right)$ equipped with the sup-norm by $\Geod\left(\X\right)$. The evaluation map for each $t\in\left[0,1\right]$ is denoted by $\e_{\substack{\\\hspace{0.5pt}t}}: \Geod\left(\X\right)\rightarrow \X$. 
\par For $p>1$, it is by now standard that any geodesic $\left( \upmu_{\substack{\\[1pt]\hspace{1pt}t}} \right)_{t\in\left[0,1\right]}\subset \Geod\left( \Prob_{\hspace{-5pt}_p}\left(\X\right) \right)$ can be lifted to a measure $\PP\in\mathscr{P}\left( \Geod\left(\X\right) \right)$ so that $\left( \e_{\substack{\\\hspace{0.5pt}t}} \right)_{*}\PP=\upmu_{\substack{\\[1pt]\hspace{1pt}t}}$ for all $t\in\left[0,1\right]$. Given two probability measures $\upmu_{\hspace{1pt}_0},\upmu_{\hspace{1pt}_1}\in \Prob_{\hspace{-5pt}_p}\left(\X\right)$, we denote by $\mathsf{OptGeod}_p\left( \upmu_{\hspace{1pt}_0},\upmu_{\hspace{1pt}_1} \right)$, the space of all probability measures $\PP\in\Prob\left( \Geod\left(\X\right) \right)$ for which, $\left( \e_{_0}, \e_{_1} \right)_{*}\PP$ is an optimal coupling between $\upmu_{\hspace{1pt}_0}$ and $\upmu_{\hspace{1pt}_1}$. $\PP$ is called an optimal {\it dynamical} coupling.
\par For $\left( \X, \dist, \meas \right)$, the dimension-free curvature conditions $\CD\hspace{-1pt}\left(\K,\infty\right)$, $\K\in \mathbb{R}$, are given by the weak $\K$-convexity of the relative entropy. Namely, we say $\left( \X, \dist, \meas \right)$ satisfies the $\CD\hspace{-1pt}\left(\K,\infty\right)$ condition whenever for any two probability measures $\upnu_{\hspace{-1pt}_0}$ and $\upnu_{\hspace{-1pt}_1}$ in $\Prob_{\hspace{-5pt}_2}\left( \X,\dist,\meas \right) $ with finite entropies, there exists a $\Was_{\hspace{-3pt}_2}$-geodesic $\left( \upnu_{\substack{\\[1pt]\hspace{-0.5pt}t}} \right)_{t\in \left[0,1\right]}$ in $\Prob_{\hspace{-5pt}_2}  \left( \X, \dist, \meas \right)$  (consisting of probability measures with finite entropies) such that
\begin{align}\label{convexitcd}
	\Ent\left(   \upnu_{\substack{\\[1pt]\hspace{-0.5pt}t}} \,|\, \meas   \right) \le \left(1 \shortminus t\right)\Ent\left(  \upnu_{\hspace{-1pt}_0}\,|\, \meas   \right)  + t \Ent\left( \upnu_{\hspace{-1pt}_1}\,|\, \meas  \right) \shortminus \nicefrac{1}{2} \, \K \,t\left(1 \shortminus t\right)\Was_{\hspace{-3pt}_2}^{^2}\left(  \upnu_{\hspace{-1pt}_0} ,  \upnu_{\hspace{-1pt}_1} \right),
\end{align}
holds; here, ``$\Ent$'' is the relative entropy. 
\par A metric measure space $\left( \X,\dist,\meas \right)$ satisfies the {\it strong} $\CD\hspace{-1pt}\left(\K,\infty\right)$ condition if  \eqref{convexitcd} holds for every geodesic $\upnu_{\substack{\\[1pt]\!t}}$ that connects $\upnu_{\hspace{-1pt}_0}$ and $\upnu_{\hspace{-1pt}_1}$.
For given numbers $\K\in \mathbb{R}$ and $\N\in\left[1,\infty\right)$, consider the distortion coefficients $\upsigma^{^{\left(t\right)}}_{\hspace{-1pt}_{\K,\,\N}}\left( \uptheta \right)$ that are given by 
\begin{align*}
\upsigma^{^{\left(t\right)}}_{\hspace{-1pt}_{\K,\,\N}}\left( \uptheta \right) :=\begin{cases}
	\infty&\text{if }\K\,\uptheta^{^2}\geq \N\uppi^{^2},\\
\nicefrac{	\sin\left( t\uptheta \N^{^{\, \shortminus \nicefrac{1}{2}}}\K^{^{\nicefrac{1}{2}}} \right)}{ \sin\left( \uptheta \N^{^{\, \shortminus \nicefrac{1}{2}}}\K^{^{\nicefrac{1}{2}}}\right)}&\text{if }0<\K\,\uptheta^{^2}<\N\uppi^{^2},\\
	t&\text{if }\K\,\uptheta^{^2}=0,\\
	\nicefrac{\sinh\left( t\uptheta \N^{^{\, \shortminus \nicefrac{1}{2}}}\left|\K\right|^{^{\nicefrac{1}{2}}} \right)}{\sinh\left( \uptheta \N^{^{\,\shortminus \nicefrac{1}{2}}}\left|\K\right|^{^{\nicefrac{1}{2}}} \right)}&\text{if } \K\,\uptheta^{^2}<0.
\end{cases}
\end{align*}
We say $\left(\X, \dist, \meas \right)$ satisfies the $\CD\left( \K,\N \right)$ curvature dimension conditions if for each pair $\upnu_{\hspace{-1pt}_0} ,  \upnu_{\hspace{-1pt}_1}  \in \Prob_{\hspace{-5pt}_2} \left(\X, \dist, \meas \right)$, there exists an optimal coupling $\frak q$ and a geodesic $\left(\upnu_{\substack{\\[1pt]\!t}} \right)_{t\in \left[0,1\right]}$ in $\Prob_{\hspace{-5pt}_2}  \left(\X, \dist, \meas \right)$ joining $\upnu_{\hspace{-1pt}_0}$ to $\upnu_{\hspace{-1pt}_1}$ with densities $\uprho_{\substack{\\[3pt]\hspace{-2pt}t}}$ that satisfy
\begin{align*}
&	\int_{\hspace{-1pt}_{\X}}   \uprho_{\substack{\\[3pt]\hspace{-2pt}t}}^{^{\shortminus \nicefrac{1}{\N^{\hspace{1pt}'}}}} \, \dif\meas \geq \\ &	\int_{\hspace{-2pt}_{\X^{^2}}}   \Big( \uptau^{^{\left(1- t\right)}}_{\hspace{-1pt}_{\K,\,\N^{\hspace{1pt}'}}} \left(  \dist\left( x_{\hspace{-0.5pt}_0} \, , x_{\hspace{-0.5pt}_1}\right) \right) \uprho_{\hspace{-1pt}_0}^{\shortminus \nicefrac{1}{\N^{\hspace{1pt}'}}}\left(x_{\hspace{-0.5pt}_0}\right)  +     \uptau^{^{\left(t\right)}}_{\hspace{-1pt}_{\K,\,\N^{\hspace{1pt}'}}}\ \left(  \dist\left( x_{\hspace{-0.5pt}_0}\, , x_{\hspace{-0.5pt}_1}\right) \right)\uprho_{\hspace{-1pt}_1}^{\shortminus \nicefrac{1}{\N\hspace{0.5pt}'}}\left( x_{\hspace{-0.5pt}_1} \right)     \Big) \, \dif{\frak q}\left( x_{\hspace{-0.5pt}_0}\, , x_{\hspace{-0.5pt}_1} \right)  \notag,
\end{align*}
for all $t\in \left[0,1\right]$ and $\N^{\hspace{1pt}'} \ge \N$.
\par It is known that $\CD\hspace{-1pt}\left( \K,\N \right)$ spaces $\left( \X,\dist,\meas \right)$ for $\N<\infty$ are proper geodesic metric spaces and the measure $\meas$ is $\upsigma$-finite measure; e.g. see \cite{Stmms2}.
\par The reduced curvature dimension conditions, $\CD^* \hspace{-2pt} \left(\K,\N\right)$, introduced in~\cite{BaSt}, are obtained in the same way by replacing the coefficients $\uptau^{^{\left(t\right)}}_{\hspace{-1pt}_{\K,\,\N}}\left( \uptheta \right)$ with $\upsigma^{^{\left(t\right)}}_{\hspace{-1pt}_{\K,\,\N}}\left( \uptheta \right)$. These coincide with the $\CD\hspace{-1pt}\left(\K,\N\right)$ bounds in all the known cases; furthermore, these conditions enjoy local-to-global and tensorization properties for essentially non-branching metric measure spaces.
\par Recently, it is shown in~\cite{cavmil} that $\CD\hspace{-1pt}\left(\K,\N\right)$ and $\CD^* \hspace{-2pt} \left(\K,\N\right)$ are equivalent for essentially non-branching metric measure spaces with finite measures.
\subsection{Essentially non-branching spaces}\label{sec:ess-non-branching}
\begin{definition}[Essentially non-branching property~\cite{RS}]
	A metric measure space $\left(\X, \dist, \meas \right)$ is said to be \emph{essentially non-branching} when any optimal dynamical plan $\PP \in \mathsf{OptGeod}_{_2}\left( \upmu , \upnu \right)$, corresponding to any two absolutely continuous measures $\upmu , \upnu \in\Prob\left(\X\right)$, is concentrated on a non-branching set of geodesics. 
\end{definition}
\begin{theorem}[Essentially non-branching spaces~\cite{RS}]\label{thm:nonbranching}
	Any strongly $\CD\hspace{-1pt}\left(\K,\infty\right)$ space is essentially non-branching. 
\end{theorem} 
\begin{remark}
	Very recently, it has been shown in~\cite{qin} that $\RCD\hspace{-1pt}\left( \K,\N \right)$ spaces with $\N<\infty$ are actually non-branching. However this result will not simplify the proof of our main theorem.
\end{remark}
\subsection{Infinitesimally Hilbertian Spaces}
\par For further details about this section, we refer the reader to~\cite{AGMR, Gigli-1, AGScal}.  Let $\left( \X,\dist,\meas \right)$ be a metric measure space. We denote the set of all Lipschitz functions in $\X$ by $\Lip\left(\X\right)$. For every $f\in \Lip\left(\X\right)$ and $x \in \X$, the local Lipschitz constant at $x$, denoted by $\Lip f\left(x\right)$, is defined by
\begin{align*}
\Lip f\left(x\right):=\uplim_{y\rightarrow x} \; \nicefrac{\left| f\left(x\right)-f\left(y\right) \right|}{\dist\left(x,y\right)},
\end{align*}
when $x$ is not isolated; otherwise, we set $\Lip f \left(x\right):=\infty$. 
\par  The \emph{Cheeger-Dirichlet energy} of a function $f\in \Ltwo\left(\X,\meas\right)$ is defined as 
\begin{align*}
\CHE\left(f\right):= \nicefrac{1}{2} \, \inf\left\{\lowlim_{n\rightarrow\infty}	\int_{\hspace{-1pt}_{\X}}   \left(\Lip f_{\hspace{-1pt}_n} \right)^{^2}\, \dif \meas \;\; \text{\textbrokenbar} \;\; f_{\hspace{-1pt}_n} \in \Lip\left(\X\right),\, f_{\hspace{-1pt}_n} \rightarrow f\text{ in }\Ltwo\right\}.
\end{align*}
\par Set 
\begin{align*}
\Sobol\left(\X\right):= \Dom\left( \CHE \right):=\left\{f\in \Ltwo\left( \X,\meas \right) \; \text{\textbrokenbar} \;\; \CHE\left( f \right)<\infty\right\}.
\end{align*}
 It is known that for any $f\in \Sobol\left( \X \right)$, there exists $\left|\nabla f\right|_{\sf w}\in \Ltwo\left( \X,\meas \right)$ such that 
 \begin{align}\label{some_energy_id}
 	2 \CHE\left( f \right)=\int_{\hspace{-1pt}_{\X}}  \left|\nabla f\right|^{^2}_{\sf w}\, \dif \meas.
 \end{align}
 The function $\left|\nabla f\right|_{\sf w}$ is called minimal weak upper gradient. For a $\CD\hspace{-1pt}\left( \K,\N \right)$ space $\X$ with $\N<\infty$, we have that $\Lip f=\left|\nabla f\right|_{\sf w}$ $\meas$-a.e. for $f\in \Lip\left( \X \right)$. This follows from \cite{Cheeg1} since $\CD\hspace{-1pt}\left( \K,\N \right)$ spaces with finite $\N$, are locally doubling and satisfy a local Poincar\'e inequality.
\subsubsection*{\small \bf \textit{Local Sobolev spaces}}
\begin{remark}\label{rem:locgra}
Let us note the locality property of $\left|\nabla f\right|_{\sf w}$ for $f\in \Sobol\left( \X \right)$. More precisely, if $\Omega\subset \X$ is open with $\meas\left( \partial \Omega \right)=0$, $f\restr_{\ovs{\Omega}}$ admits a minimal weak upper gradient on $\ovs{\Omega}$ equipped with the restricted distance (likewise, the induced intrinsic distance function can be used since in locally geodesic spaces, both distances locally coincide) and with the restricted measure; furthermore, $\left| \nabla \left( f\restr_{\ovs{\Omega}} \right) \right|_{\sf w}= \left|\nabla f\right|_{\sf w}\restr_{\ovs{\Omega}}$ ; see~\cite[Lemma 4.19]{AGSRiem}.
\end{remark}
	Let $\Omega\subset \X$ be open. A Borel function $f:\Omega \rightarrow \R$ belongs to $\Sobol_{\textsf{loc}}\left( \Omega \right)$ provided $f\cdot \upchi\in \Sobol\left( X \right)$ holds for any cutoff Lipschitz function $\upchi: \Omega\rightarrow \R$ with $\supp \left( \upchi \right) \subset \Omega$. Thanks to the locality of the minimal weak upper gradient, one can define $\left|\nabla f\right|_{\sf w}:= \left|\nabla \left( f\cdot \upchi \right)\right|_{\sf w}$ $\meas$-a.e. in the interior of $\upchi\equiv 1$. The space $\Sobol\left( \Omega \right)$ is the set of $f\in \Sobol_{\textsf{loc}}\left( \Omega \right)$ such that $\left|\nabla f\right|_{\sf w}\in \Ltwo\left( \X,\meas \right)$. 
\subsubsection*{\small \bf \textit{Test plans}}
\begin{definition}[Test plans~\cite{AGSRiem}]
	 Let $\ppi\in \Prob \left( \mathcal{C} \left( \left[0,1\right],\X \right) \right)$. We say that $\ppi$ has bounded compression provided there exists a constant $C>0$ such that $
	 \left(\e_{\substack{\\\hspace{0.5pt}t}} \right)_*\ppi\leq C \meas, \ \forall t\in \left[0,1\right].$
	We say $\ppi$ is a test plan if it has bounded compression, is concentrated on $\mathcal{AC}^{\hspace{0.5pt}^2} \left( \left[0,1\right], \X \right)$ -- the set of absolutely continuous curves with metric derivative in $\Ltwo([0,1])$ -- and if
	\begin{align*}
		\int_{\hspace{-1pt}_{\mathsf{Curves}}}  \int_{\hspace{-1pt}_0}^{1}\left|\bdot{ \upxi} \left( t \right)\right|^{^2} \, \dif t \, \dif\ppi\left( \upxi \right)=:\left\| \ppi \right\|_{_2}^{^2}<\infty, \quad \textsf{Curves}:= \mathcal{C} \left( \left[0,1\right],\X \right).
	\end{align*}
	A subset $\mathscr{A}\subset \mathcal{C}\left( \left[0,1\right],\X \right)$ is negligible if $\ppi\left( \mathscr{A} \right)=0$ for all test plans $\ppi$.
\par In $\CD\hspace{-1pt}\left(\K,\N\right)$ spaces, it follows from~\cite{Raj-1} that optimal dynamical couplings are test plans.
\end{definition}
 A function $f\in \Dom\left( \CHE   \right)$ is contained in the {\it Sobolev class} $\mathcal{S}^{\hspace{1pt}^2}_{\textsf{loc}}\left(\X\right)$; see~\cite{Gigli-1,agslipschitz}. Recall, a Borel function $f:X\rightarrow \R$ is contained in $\mathcal{S}^{\hspace{1pt}^2}_{\textsf{loc}}\left(\X\right)$ if  there exists $G\in \Ltwo(\X, \meas)$ and a negligible set $\mathscr{A}\subset \mathcal{C}\left( \left[0,1\right],\X \right)$ such that for any curve $\upgamma:\left[0,1\right]\rightarrow \X$ that is not contained in $\mathscr{A}$, it holds that 
\begin{align*}
\left\vert f\left( \upgamma(s) \right) \shortminus f\left( \upgamma(t) \right)\right\vert \leq \int_{\hspace{-2pt}s}^t G\circ \upgamma(\tau) \left| \bdot \upgamma \right|(\tau) \, \dif \tau\ \quad  \forall s\leq t\in [0,1].
\end{align*}
In this case, $G$ is called a $2$-weak upper gradient, or just a weak upper gradients.
In fact, this allows an alternative equivalent definition of $\CHE$ via the notion of minimal weak upper gradient.

For later purposes, let us state a useful stability property of minimal weak upper gradients.  
\begin{theorem}[Stability of minimal weak upper gradients~{\cite[Theorem 5.3]{agslipschitz}}]\label{th:stagra}
Let $f_{_n}\in \Sobol(X)$ and $f\in \Ltwo(\meas)$ and assume that $f_{_n} \rightarrow f$ pointwise $\meas$-a.e. and that $|\nabla f_{_n}|_{\substack{\\ \sf w}} \rightarrow G$ weakly in $\Ltwo(\meas)$. Then $f\in \Sobol(X)$ and $G\geq |\nabla f|_{\substack{\\ \sf w}}$ $\meas$-a.e.
\end{theorem}
\subsubsection*{\small \bf \textit{Riemannian curvature dimension conditions}}
\par We say that $\left( \X,\dist,\meas \right)$ is \emph{infinitesimally Hilbertian} if the Cheeger-Dirichlet energy is a quadratic form. Infinitesimal Hilbertianity is equivalent to the Sobolev space $\Sobol\left(\X\right)$ -- that is $\Dom\left(  \CHE \right)$ equipped with the norm $\left\| f\right\|^{\hspace{-0.5pt}^2}_{\hspace{-0.5pt}_{1,2}}:=\left\| f\right\|^{\hspace{-0.5pt}^2}_{\hspace{-0.5pt}_{2}}+2\CHE\left(f\right)$ -- being a Hilbert space.
	\begin{definition}\label{def:rcd} Let $\K\in \mathbb{R}$ and $\N\in \left[1,\infty\right]$. An infinitesimally Hilbertian metric measure space $\left( \X,\dist,\meas \right)$ that also satisfies the $\CD\hspace{-1pt}\left(\K,\N\right)$ conditions, is called an $\RCD\hspace{-1pt}\left(\K,\N\right)$ space. 
\end{definition}
\begin{remark} Any $\RCD\hspace{-1pt}\left(\K,\N\right)$ space for $\N\in \left[1,\infty\right)$ is $\RCD\hspace{-1pt}\left(\K,\infty\right)$ and therefore strongly $\CD\hspace{-1pt}\left(\K,\infty\right)$; see~\cite{AGMR}. It follows from Theorem~\ref{thm:nonbranching} that any $\RCD\hspace{-1pt}\left(\K,\N\right)$ space is \emph{essentially non-branching}.
\end{remark}
\begin{remark} Since $\RCD^*\hspace{-2pt}\left(\K,\N\right)$ spaces are essentially non-branching, it follows from \cite{cavmil} that $\RCD\hspace{-1pt}\left(\K,\N\right)$ and $\RCD^*\hspace{-2pt}\left( \K,\N \right)$ are equivalent provided the reference measure $\meas$ is finite.
	Since we consider compact spaces, the $\CD\hspace{-1pt}\left(\K,\N\right)$ or $\CD^*\hspace{-2pt}\left(\K,\N\right)$ condition with $\N<\infty$, implies finite measure. Therefore, we will usually assume the $\RCD\hspace{-1pt}\left(\K,\N\right)$ condition and  any result for $\RCD^*\hspace{-2pt}\left(\K,\N\right)$ spaces  is  at our disposal. 
\end{remark}
\par In an $\RCD\hspace{-1pt}\left(\K,\N\right)$ space with finite $\N$, the Cheeger-Dirichlet energy $\CHE$ (or twice thereof) is a regular \emph{strongly local} Dirichlet form. The set $\Lip\left(\X\right)$ of Lipschitz functions is dense in the Sobolev space $\Sobol\left(\X\right)$; see~\cite{Gigli-1} for definitions. This implies,  by polarization, that the Cheeger-Dirichlet energy $\CHE$ admits a \emph{square field operator} operator $\BG:\Sobol\left(\X\right)\times \Sobol\left(\X\right)\rightarrow \LL^{\hspace{-3pt}^1}\left(\meas\right)$. We will suppress the $\BG$ notation and simply point out that we have the following symmetric bi-linear form
\begin{align*}
\nabla \left(\cdot\right) \boldsymbol{\cdot} \nabla \left(\cdot\right)  :  \Sobol\left(\X\right) \times \Sobol\left(\X\right) \to 	\LL^{\hspace{-3pt}^1}\left(\meas\right). 
\end{align*}
The ``dot'' notation is a reminder of the fact that, in the smooth setting, the square field operator is indeed the same as inner product of gradients. Sometimes, we also use the notation $\langle \nabla \cdot \,, \nabla \cdot\rangle$. 
\subsection{Laplacian in $\RCD\hspace{-1pt}\left(\K,\N\right)$ spaces}
\par  Let $\left( \X,\dist,\meas \right)$ be an $\RCD\hspace{-1pt}\left(\K,\N\right)$ space. 
A function $f \in \Sobol\left(\X\right)$ is in the domain of the {\it Laplace operator} $\Dom\left( {\CHL} \right)$ whenever there exists $h\in \Ltwo\left( \X \right)$ such that 
\begin{align*}
\shortminus \int_{\hspace{-1pt}_{\X}} \left<\nabla{ f} ,  \nabla g   \right> \, \dif\meas = \int_{\hspace{-1pt}_{\X}} g h\, \dif\meas,
\end{align*}
holds for all Lipschitz functions {$g$}. In this case we write $h=: \CHL \left(f\right)$ and call the operator $\CHL$, the \emph{Cheeger Laplacian}.
\par A function $f \in \Sobol\left( \X \right)$ is in the domain of the {\it distributional Laplace operator} $\Dom\left( \mathbf{\Delta} , \Omega \right)$ whenever there exists a {\it Radon functional} $\mathscr{T}\in \Lip\left( \X \right)^{\boldsymbol{*}}$ in the sense of \cite[Definition 2.11]{cavmonlap} such that 
\begin{align*}
\shortminus \int_{\hspace{-1pt}_{\Omega}}  \left<\nabla f ,  \nabla g   \right> \, \dif\meas = \mathscr{T} \left( g \right)=: \mathbf{\Delta}\left(f\right) \left( g \right).
\end{align*}
holds for all Lipschitz functions $g$ with $\supp \left( g \right) \subset \Omega^{\hspace{0.5pt}'} \subset \Omega$ for some open domain $\Omega^{\hspace{0.5pt}'}$; the definition of a Radon functional can be found in \cite[Definition 2.9]{cavmonlap}. Letting $\Omega = \X$, one writes $\Dom \left( \mathbf{\Delta} \right)= \Dom \left( \mathbf{\Delta} , \X \right)$. It is well-known that one can write $\mathscr{T}\in \Lip\left( \X \right)^{\boldsymbol{*}}$ as difference of two Radon measures $\upmu^{\hspace{-1pt}^+}$ and $\upmu^{\hspace{-1pt}^-}$. Hence, if at most one of the two measures $\upmu^{\hspace{-1pt}^\pm}$ attains the value $\infty$, then ${\bf \Delta} \left(f\right)$ is a well-defined signed Radon measure.
\par For $f \in \Dom\left(\mathbf{\Delta}\right)$ such that ${\bf \Delta}\left(f\right)$ is a signed Radon measure, the Radon-Nikodym decomposition of $\mathbf{\Delta}$ w.r.t. $\meas$ will be written as
\begin{align*}
\mathbf{\Delta}\left( f \right) = {\Delta}_{\hspace{0.5pt}\textsf{abs}}\left(f \right) \, \meas + \mathbf{\Delta}_{\hspace{0.5pt}\textsf{sing}}\left(f\right).
\end{align*}
\subsubsection*{\small \bf \textit{Chain rule for the Laplacian}}
\begin{proposition}[Laplacian chain rule~{\cite[Proposition 4.28]{Gigli-1}}]\label{prop:chainrulelaplacian}
	Let $\left(\X,\dist,\meas\right)$ be an infinitesimally Hilbertian metric measure space; let $f\in \Dom\left({\bf \Delta},\Omega \right)\cap \Lip\left(\X\right)$ for an open subset $\Omega\subset \X$;  suppose $\mathcal{J} \subset \mathbb{R}$  is an open subset  with $\meas\left( f^{^{-1}}\left( \mathbb{R}\smallsetminus \mathcal{J}  \right) \right)=0$. Then, for any  $\upvarphi\in \mathcal{C}^{\hspace{0.5pt}^2}\left( \mathcal{J} \right)$, we have $\upvarphi\circ f\in \Dom \left( {\bf \Delta},\Omega \right)$ and 
	\begin{align*}
		{\bf\Delta}\left(\upvarphi\circ f \right)= \upvarphi^{'}\left(f\right){\bf\Delta} f + \upvarphi^{''}\left( f \right)\left|\nabla f\right|^{\hspace{-0.5pt}^2}_{{\sf w}} \, \meas \quad \text{ on }\quad \Omega.
	\end{align*}
\end{proposition}
\subsection{Bakry-\'Emery Curvature dimensions for regular strongly local Dirichlet forms}
\subsubsection*{\small \bf \textit{Setup}}
\par Throughout this section, $\X$ is a Polish topological space equipped with a $\upsigma$-finite measure $\meas$. We will assume $\mathcal{E}$ is a \emph{quasi-regular strongly local symmetric Dirichlet form} (see \cite{Sav-1} for definitions) with domain
\begin{align*}
\mathbb{V}:= \left\{ f \in \Ltwo\left(\X,\meas\right) \; \text{\textbrokenbar} \;\; \mathcal{E}\left(f\right) < \infty   \right\}   \underset{\textsf{dense}}{\subset} \Ltwo\left( \X,\meas  \right),
\end{align*}  
which admits a square field operator $\BG: \mathbb{V} \times \mathbb{V} \to \LL^{\hspace{-3pt}^1}\left(\X,\meas\right)$ which restricted to $\mathbb{V}_{\hspace{-2pt}_\infty} := \mathbb{V} \cap \LL^{\hspace{-2pt}^\infty}$, is uniquely characterized via
\begin{align*}
2\int_{\hspace{-1pt}_{\X}}  \upvarphi \BG \left(  f,g  \right)\, \dif\meas = \mathcal{E}\left( f, g\upvarphi \right) + \mathcal{E}\left( f\upvarphi, g \right) \shortminus  \mathcal{E}\left( fg, \upvarphi \right),
\end{align*} 
for all $f,g,\upvarphi \in \mathbb{V}_{\hspace{-2pt}_\infty}$. In particular, we are interested in the case where $\X$ is $\RCD\hspace{-1pt}\left(\K,\N\right)$ and $\mathcal E=2\CHE$.
\par Furthermore, assume there is a mass-preserving Markov semigroup $\markov$ for $t \ge 0$ generated by $\mathcal E$ which admits an infinitesimal generator $\Delta_{\substack{\\[2pt]\mathcal{E}}}$ with domain
\begin{align*}
\Dom \left( \Delta_{\substack{\\[2pt]\mathcal{E}}} \right)  \underset{\textsf{dense}}{\subset} \mathbb{V}.
\end{align*}
If $\mathcal E=2\CHE$, we have $\mathbb{V}= \Dom\left(\CHE \right)=\Sobol\left( \X \right)$ and the generator $\Delta_{\substack{\\[2pt]\mathcal{E}}}$ in fact coincides with the Cheeger-Laplacian $\CHL$ that was defined before.
\par The domain $\mathbb{V}$ equipped with the Sobolev norm
\begin{align*}
\left\|  \cdot\right \|_{\substack{\\[1pt]\mathbb{V}}} :=\left( \left\|  \cdot\right \|^{^2}_{\substack{\\\Ltwo}}  + \mathcal{E}\left( \cdot, \cdot  \right)  \right)^{\hspace{-3pt}^{\nicefrac{1}{2}}},
\end{align*}
is a Hilbert space. 
\par Define the good set 
\begin{align*}
\mathbb{G}_{\hspace{-0.5pt}_{\infty}}  := \mathbb{V}_{\hspace{-2pt}_\infty} \cap \BG^{^{-1}}\left( \LL^{\hspace{-2pt}^\infty} \right).
\end{align*}
\par In the case of $\mathcal E={2}\CHE$ on $\RCD$ spaces, it follows from \eqref{some_energy_id}  that $\Gamma(\cdot) = |\nabla \cdot|^{^2}_{\sf w}$ and the good set is precisely
\begin{align*}
\mathbb{G}_{\hspace{-0.5pt}_{\infty}}   = \left\{ f\in \Sobol\left( \X \right)\cap \LL^{\hspace{-2pt}^\infty}\left( \meas \right) \; \text{\textbrokenbar} \;\; \left|\nabla  f\right|_{\sf w}\in \LL^{\hspace{-2pt}^\infty} \left( \meas \right)\right\}.
\end{align*}
\subsubsection*{\small \bf \textit{General Bakry-\'Emery conditions}}
\par The weak multilinear Bakry-\'Emery iterated square field operator $\BG_{\hspace{-1pt}_2}$ is defined in~\cite{Sav-1} by 
\begin{align*}
\BG_{\hspace{-1pt}_2}\left( f, g ; \upvarphi    \right):= \nicefrac{1}{2} \int_{\hspace{-1pt}_{\X}}   \BG\left(f,g\right)\Delta_{\substack{\\[2pt]\mathcal{E}}}\left( \upvarphi  \right) \shortminus \big(  \BG\left( \Delta_{\substack{\\[2pt]\mathcal{E}}} \left(f\right),g \right) +  \BG\left( f,\Delta_{\substack{\\[2pt]\mathcal{E}}} \left(g\right) \right)  \big) \upvarphi    \; \dif\meas,
\end{align*}
for triples 
\begin{align*}
	\left( f, g , \upvarphi   \right) \in \Dom_{\hspace{1pt}\substack{\\\mathbb{V}}}\left( \Delta_{\substack{\\[2pt]\mathcal{E}}}  \right) \times \Dom_{\hspace{1pt}\substack{\\\mathbb{V}}}\left( \Delta_{\substack{\\[2pt]\mathcal{E}}}  \right) \times \Dom_{\hspace{1pt}\substack{\\ \LL^{\hspace{-2pt}^\infty}}}\left( \Delta_{\substack{\\[2pt]\mathcal{E}}}  \right) =: \Dom\left( \BG_{\hspace{-1pt}_2} \right),
\end{align*}
 where 
\begin{align*}
\Dom_{\hspace{1pt}\substack{\\\mathbb{V}}}\left( \Delta_{\substack{\\[2pt]\mathcal{E}}} \right) := \Dom\left(  \Delta_{\substack{\\[2pt]\mathcal{E}}}  \right) \cap \Delta_{\substack{\\[2pt]\mathcal{E}}}^{\hspace{-3pt}^{-1}}\left( \mathbb{V}  \right),
\end{align*}
and
\begin{align*}
\Dom_{\hspace{1pt}\substack{\\ \LL^{\hspace{-2pt}^\infty}}}\left( \Delta_{\substack{\\[2pt]\mathcal{E}}} \right) := \Dom\left( \Delta_{\substack{\\[2pt]\mathcal{E}}}  \right) \cap \LL^{\hspace{-2pt}^\infty}\left( \X,\meas \right) \cap \Delta_{\substack{\\[2pt]\mathcal{E}}}^{\hspace{-3pt}^{-1}}\left( \LL^{\hspace{-2pt}^\infty}  \right).
\end{align*}
\par As is customary in this context and by using usual conventions, when $f=g$, one instead writes the abbreviated form
\begin{align*}
 \BG_{\hspace{-1pt}_2}\left( f; \upvarphi  \right) := \int_{\hspace{-1pt}_{\X}}  \nicefrac{1}{2}\, \BG\left(f\right)\Delta_{\substack{\\[2pt]\mathcal{E}}}   \left(\upvarphi\right) \shortminus \BG\left( f, \Delta_{\substack{\\[2pt]\mathcal{E}}}  \left(f\right)  \right) \upvarphi  \; \dif\meas.
\end{align*}
 Similar as in \cite{GKK}, we also define 
	\begin{align*}
		\uphat{\BG_{\hspace{-1pt}_2}}\left( f;\upvarphi \right):= \nicefrac{1}{2}\, \int_{\hspace{-1pt}_{\X}}  \BG\left(f\right) \Delta_{\substack{\\[2pt]\mathcal{E}}} \left(\upvarphi\right) \, \dif\meas + \int_{\hspace{-1pt}_{\X}}  \left(\Delta_{\substack{\\[2pt]\mathcal{E}}} \left(f\right) \right)^{\hspace{-3pt}^2}\upvarphi \, \dif\meas + \int_{\hspace{-1pt}_{\X}}  \BG\left(\upvarphi, f\right) \Delta_{\substack{\\[2pt]\mathcal{E}}} \left(f\right)\, \dif\meas,
	\end{align*}
	for $f\in \Dom\left( \Delta_{\substack{\\[2pt]\mathcal{E}}} \right)$ and $\upvarphi\in \Dom_{\hspace{1pt}\substack{\\ \LL^{\hspace{-2pt}^\infty}}}\left( \Delta_{\substack{\\[2pt]\mathcal{E}}} \right) \cap \mathbb V\cap \BG^{^{-1}}\left( \LL^{\hspace{-2pt}^\infty} \right)$. Integration by parts shows that $\BG_{\hspace{-1pt}_2}$ and $\uphat{\BG_{\hspace{-1pt}_2}}$ coincide on $\Dom \left( \BG_{\hspace{-1pt}_2} \right) \cap \Dom \left(  \uphat{\BG_{\hspace{-1pt}_2}} \right)$.	
\begin{definition}[BE conditions]\label{def:be}
	For given $\K \in \mathbb{R}$ and $\N \in \left[1,\infty\right]$, we say that the infinitesimal generator $ \Delta_{\substack{\\[2pt]\mathcal{E}}}$ satisfies the $\left( \K,\N \right)$-Bakry-\'{E}mery curvature conditions (in short, $\BE\hspace{-0.5pt}\left(\K,\N\right)$ curvature dimension conditions) whenever for every $f \in \Dom_{\hspace{1pt}\substack{\\\mathbb{V}}}\left(\Delta_{\substack{\\[2pt]\mathcal{E}}}\right)$ and $\upvarphi \left(\ge 0\right) \in \Dom_{\hspace{1pt}\substack{\\ \LL^{\hspace{-2pt}^\infty}}}\left( \Delta_{\substack{\\[2pt]\mathcal{E}}} \right)$, there holds
	\begin{align*}
\BG_{\hspace{-1pt}_2}\left(f;\upvarphi\right) \ge \int_{\hspace{-1pt}_{\X}}  \left( \nicefrac{1}{\N} \left( \Delta_{\substack{\\[2pt]\mathcal{E}}} \left(f\right)  \right)^{\hspace{-2pt}^2} +  \K \BG\left(f\right)  \right) \upvarphi \, \dif\meas.
	\end{align*}
\end{definition}
\subsection{Bakry-\'Emery curvature dimension conditions for the Cheeger-Dirichlet energy}\label{sec:BE-cond}
\emph{For the rest of these notes, we set $\Delta:= \CHL$ and $\left| \nabla\, \cdot\, \right|:= \left| \nabla\, \cdot\, \right|_{\sf w}$.}
\par The Bakry-\'Emery curvature dimension conditions $\BE\left( \K,\N \right)$ for the Cheeger Laplacian now amount to 
\begin{align*}
 \BG_{\hspace{-1pt}_2}\left( f;\upvarphi \right) \ge \int_{\hspace{-1pt}_{\X}} \left( \nicefrac{1}{\N} \left( \Delta f  \right)^{\hspace{-1pt}^2} + \K \left| \nabla f \right|^{^2} \right)\upvarphi \, \dif\meas,
\end{align*}
holding for $\left( f,f,\upvarphi \right) \in \Dom\left( \BG_{\hspace{-1pt}_2} \right)$ with $\upvarphi \ge 0 $.
\begin{definition}[Sobolev-to-Lipschitz property~\cite{AGSRiem,Gigli-1}]\label{def_sobtolip}
	A metric measure space satisfies the \textit{Sobolev-to-Lipschitz} property if corresponding to any Sobolev function with bounded minimal weak upper gradient, there exists a Lipschitz function $\bar{f}$ that coincides $\meas$-almost everywhere with $f$ such that $\Lip \bar{f} \geq \left|\nabla f\right|$ in $\meas$-a.e. sense.
\end{definition}
\begin{theorem}[Equivalence of BE and RCD~\cite{AGSbercd, EKS, ams_nonlinear}]\label{th:be}
	Let $\left( \X,\dist,\meas \right)$ be a metric measure space. Assume the growth condition $\meas\left( \ball{r}\left( x_{\hspace{-0.5pt}_0} \right) \right)\leq A\,\mathrm{exp}\left( Br \right)$ for constants $A,B>0$ and $x_{\hspace{-0.5pt}_0} \in \X$.  Then, curvature dimension conditions $\RCD\left( \K,\N \right)$ for $\K\in \mathbb{R}$ and $\N>1$ hold if and only if $\left( \X,\dist,\meas \right)$ is infinitesimally Hilbertian, it satisfies the 
	Sobolev-to-Lipschitz property and it satisfies the Bakry-\'{E}mery curvature dimension conditions $\BE\left( \K,\N \right)$. 
\end{theorem}
\begin{remark}
	The case $\N=\infty$ was proved in~\cite{AGSbercd} and the case $\N<\infty$ in~\cite{EKS}.
	Shortly after the work~\cite{EKS}, an alternative proof for the finite dimensional case - following a completely different strategy - was established in~\cite{ams_nonlinear}.
\par We also note that for a compact metric measure space satisfying any curvature dimension condition the growth condition is always satisfied \cite{Stmms1}.
\end{remark}
\par Let us assume that $\X$ is $\RCD\left(\K,\N\right)$ for $\N<\infty$ and {$\mathcal E=2\CHE$}.
	The important class of functions
	\begin{align*}
		\mathbb{D}_{\hspace{-1pt}_{\infty}} &:= \mathbb{G}_{\hspace{-0.5pt}_{\infty}}  \cap \Dom\left(\Delta\right) \cap \Delta^{\hspace{-2pt}^{-1}}\left(\mathbb{V}\right) \\ &= \left\{f\in \Dom\left( \Delta \right)\cap \LL^{\hspace{-2pt}^\infty}\left(\meas\right) \; \text{\textbrokenbar} \;\; \Delta f\in \Sobol\left(\X\right), \ \left|\nabla f\right| \in \LL^{\hspace{-2pt}^\infty}\left(\meas\right)\right\},
	\end{align*}
is introduced in~\cite{Sav-1}.
\par If $\X$ is a compact $\RCD$ space, then
$\markov  \left( \LL^{\hspace{-2pt}^\infty} \hspace{-2pt}\left(\meas\right) \right)$ is a subset of $\mathbb{D}_{\hspace{-1pt}_{\infty}}$; $\markov  \left( \LL^{\hspace{-2pt}^\infty} \hspace{-2pt}\left(\meas\right) \right)$ is dense in $\Sobol\left(\X\right)$ and in $\Dom_{\substack{\Ltwo\left( \meas\right)}} \left(\Delta\right)$ w.r.t. the Sobolev norm and the graph norm of the operator $\Delta$ respectively; by the 
Sobolev-to-Lipschitz property, we have $\mathbb{D}_{\hspace{-1pt}_{\infty}}\subset \Lip_{\substack{\\{\sf b}}}\left(\X\right)$,  where $\Lip_{\substack{\\{\sf b}}}\left(\X\right)$ is the set of bounded Lipschitz functions.
It is proven in~\cite{Sav-1}  that $f\in \mathbb{D}_{\hspace{-1pt}_{\infty}}$ implies $\left|\nabla f\right|^{^2} \in \Sobol\left(\X\right)$ and the distributional Laplacian ${\bf \Delta} \left(  \left|\nabla f\right|^{^2} \right)$ of $\left|\nabla f\right|^{^2}$ is a signed Radon measure.
Moreover, for $f\in \mathbb{D}_{\hspace{-1pt}_{\infty}}$, one can introduce the \emph{measure-valued iterated square field operator}
\begin{align*}
	{\overset{\star}{\BG}}_{\hspace{-2.5pt}_2} \left(f\right) := \nicefrac{1}{2} \, {\bf\Delta} \left(  \left|\nabla f\right|^{^2} \right) \shortminus \langle \nabla f, \nabla \Delta f  \rangle \, \meas.
\end{align*}
We will write the Radon-Nikodym decomposition of ${\overset{\star}{\BG}}_{\hspace{-2.5pt}_2} \left(f\right) $ w.r.t. $\meas$ as
\begin{align*}
	{\overset{\star}{\BG}}_{\hspace{-2.5pt}_2} \left(f\right)  =  \upgamma^{\hspace{1pt}\textsf{abs}}_{\substack{\\ \hspace{-2pt}_2}} \left(f\right) \, \meas + \boldsymbol{\upgamma}_{\substack{\\ \hspace{-2pt}_2}}^{\hspace{1pt}\textsf{sing}} \left(f\right).
\end{align*}
Based on the definition of ${\overset{\star}{\BG}}_{\hspace{-2.5pt}_2} \left(f\right) $, it is easy to deduce
\begin{align*}
\boldsymbol{\upgamma}_{\substack{\\ \hspace{-2pt}_2}}^{\hspace{1pt}\textsf{sing}}\left(f\right)=\nicefrac{1}{2} \, {\bf \Delta}_{\hspace{1pt}\textsf{sing}}\left(\left|\nabla f\right|^{^2}\right) \ge 0,
\end{align*}\\
 i.e. the negative part of ${\overset{\star}{\BG}}_{\hspace{-2.5pt}_2} \left(f\right)$ is absolutely continuous w.r.t. $\meas$; see~\cite[Lemma 2.6]{Sav-1}.
\par The $\BE\left(\K,\N\right)$ condition gives the following Bochner type inequality 
\begin{align*}
\nicefrac{1}{2} \, {\bf \Delta} \left( \left|\nabla f \right|^{^2} \right) \ge  \langle \nabla f, \nabla \Delta f \rangle \,  \meas +  \K \left|\nabla f\right|^{^2} \,  \meas \quad \forall f\in \mathbb{D}_{\hspace{-1pt}_{\infty}},
\end{align*}
in the sense of measures; since $\boldsymbol{\upgamma}^{\hspace{1pt}\textsf{sing}}_{\hspace{-2pt}_2}\left(f\right)\geq 0$, it holds
\begin{align*}
\nicefrac{1}{2} \, {\Delta}_{\textsf{abs}} \left(  \left|\nabla f\right|^{^2} \right) \ge \langle \nabla  f, \nabla \Delta f  \rangle +  \K\left|\nabla f\right|^{^2} \quad \meas\mbox{-a.e. in } \X. 
\end{align*}
	Another feature of elements ${f}\in \mathbb{D}_{\hspace{-1pt}_{\infty}}$ is that one can define a {\it Hessian} of $f$, $\Hess\left({f}\right)$, via the (Riemannian) formula
\begin{align}\label{eq:hessian}
\Hess\left({f}\right)\left(v,w\right) := \nicefrac{1}{2} \, \Big( \big\langle \nabla v, \nabla \langle \nabla f, \nabla w\rangle \big\rangle + \big\langle \nabla w,\nabla \langle \nabla f,\nabla v\rangle \big\rangle  \shortminus \big\langle \nabla  f, \nabla \langle \nabla v,\nabla w\rangle \big\rangle \Big),
\end{align}
for any  $v,w \in \mathbb{D}_{\hspace{-1pt}_{\infty}}$. 
Note that expression on the right hand side is well-defined by virtue of the Sobolev-to-Lipschitz property. One can even extend the operator $\Hess$ to a much bigger class that one calls $\mathcal{W}^{\hspace{1pt}^{2,2}}\left(\X\right)$. The class $\mathcal{W}^{\hspace{1pt}^{2,2}}\left(\X\right)$ contains $\Dom\left(\Delta\right)\cap \Delta^{\hspace{-2pt}^{-1}}\left(\Ltwo\right)$. Applied to elements of $\mathcal{W}^{\hspace{1pt}^{2,2}}\left(\X\right)$, the Hessian is a tensor object; see \cite{Gigli-2} for the full theory .
A very crucial second order gradient estimates for functions in $\mathbb{D}_{\hspace{-1pt}_{\infty}}$ is obtained, among other places, in~\cite{Sav-1}.
	\begin{theorem}[{Improved Bochner inequalities~\cite{Sav-1,Gigli-2,sturm-tensor}}]\label{thm:improved}
	Let $\left( \X,\dist,\meas \right)$ be an $\RCD\left(\K,\N\right)$ space. 
Then for every $f \in \mathbb{D}_{\hspace{-1pt}_{\infty}}$, it holds
\begin{enumerate}
	\item  \makebox[\linewidth]{\(\begin{aligned}[t]
		{\overset{\star}{\BG}}_{\hspace{-2.5pt}_2} \left(f\right) \shortminus \K\left|\nabla f\right|^{^2}\, \meas&\geq \left\| \Hess\left(f\right)\right\|^{^2}_{_\mathsf{HS}} \, \meas,\label{improved1}
	\end{aligned}		\)} in the sense of measures; and 
	\medskip
	\item \makebox[\linewidth]{\(\begin{aligned}[t]
4 \left( \upgamma^{\hspace{1pt}\mathsf{abs}}_{\hspace{-2pt}_2} \left(f\right) \shortminus \K\left|\nabla f\right|^{^2}\right) \left|\nabla f\right|^{^2}&\geq \left|\nabla \left|\nabla f\right|^{^2} \right|^{^{2}}  \quad  \meas\mbox{-a.e.}\,
\end{aligned}	\)}
\end{enumerate}
where in above, $\left\| \Hess\left(f\right)\right\|_{_\mathsf{HS}}$ is the Hilbert-Schmidt norm of $\Hess\left(f\right)$; e.g. see \cite{Gigli-2}.
\end{theorem}
 Upon multiplying \eqref{improved1} by non-negative $\upvarphi\in \Dom_{\hspace{1pt}\substack{\\ \LL^{\hspace{-2pt}^\infty}}}\left(\Delta\right)$ and integrating, we can rewrite the weak version of \eqref{improved1} as
\begin{align}\label{eq:bochner-modified}
\BG_{\hspace{-1pt}_2}\left(f;\upvarphi\right)\geq \K \int_{\hspace{-1pt}_{\X}}  \left|\nabla f\right|^{^2}\upvarphi \; \dif\meas + \int_{\hspace{-1pt}_{\X}}  \left\| \Hess\left(f\right)\right\|^{^2}_{_\textsf{HS}} \upvarphi \; \dif\meas.
\end{align}
 In particular, for $f\in \Dom\left(\Delta\right)\cap \Delta^{\hspace{-2pt}^{-1}}\left(\Ltwo\right)$ one can choose a sequence $f_{\hspace{-1pt}_n}$ that converges to $f$ in the sense that $f_{\hspace{-1pt}_n}$ and $\Delta \left( f_{\hspace{-1pt}_n} \right) \in \Ltwo$ converge to $f$ and $\Delta \left(f\right)$ resp. in $\Ltwo$; recall $\Delta$ is closed under the graph-norm. 
Then, by choosing $\upvarphi\equiv 1$ we get (see \cite{Gigli-2})
\begin{align*}
\int_{\hspace{-1pt}_{\X}} \left(\Delta \left(f\right) \right)^{\hspace{-1pt}^2}  \shortminus \K\left|\nabla f\right|^{^2}  \, \dif\meas \geq \int_{\hspace{-1pt}_{\X}} \left\|\Hess\left(f\right)\right\|^{^2}_{_\textsf{HS}} \,  \dif\meas,
\end{align*}
where the operator $\Hess\left(f\right)$ is the well-defined \emph{Hessian} of $f\in \Dom\left(\Delta\right) \cap  \Delta^{\hspace{-2pt}^{-1}}\left(\Ltwo\right)$ as mentioned before. More precisely, $f\in \mathcal{H}^{^{2,2}}\left(\X\right)$. The class $\mathcal{H}^{^{2,2}}\left(\X\right)$ is the closure of $ \Dom\left(\Delta\right) \cap  \Delta^{\hspace{-2pt}^{-1}}\left(\Ltwo\right)$ in $\mathcal{W}^{\hspace{1pt}^{2,2}}\left(\X\right)$; c.f. \cite{Gigli-2}. The weak Sobolev space $\mathcal{W}^{\hspace{1pt}^{2,2}}\left(\X\right)$ is defined as functions in $\Sobol$ which admit a $2$-form $\mathsf{Hess}(f)$ as their weak Hessian where the weak Hessian is characterized by satisfaction of the weak formulation
\begin{align*}
	&2 \int_{\hspace{-1pt}_{\X}}   h \, \mathsf{Hess}(f) \left(  \nabla g_{_1}, \nabla g_{_2} \right) \, \dif \meas = \\ & \phantom{cys}\shortminus \int_{\hspace{-1pt}_{\X}}  \Big( \nabla f \boldsymbol{\cdot} \nabla g_{_1} \, \mathsf{div} \left( h \nabla g_{_2} \right) +  \nabla f \boldsymbol{\cdot} \nabla g_{_2}\, \mathsf{div} \left( h \nabla g_{_1}	 \right)  + h \nabla f  \boldsymbol{\cdot}\nabla \left(  \nabla g_{_1}    \boldsymbol{\cdot}  \nabla g_{_2} \right) \Big) \, \dif \meas \\ & \phantom{cys}\phantom{cys}\phantom{cys} \forall h, g_{_1}, g_{_2} \in \Lip\left(\X\right) \cap \mathbb{D}_{\hspace{-1pt}_{\infty}}.
\end{align*}
\par We will not go into more details about $\mathcal{W}^{^{2,2}}\left(\X\right)$ or $\mathcal{H}^{^{2,2}}\left(\X\right)$ and the divergence operator $\mathsf{div}$ or about $2$-forms which are notions involved in the above definition; we only highlight the fact that in order to define these notions on metric measure spaces, one needs to study the normed modules which are nonsmooth counterparts of sections of bundles. Further details can be found in the original source~\cite{Gigli-2}.
\par For the application purposes that we have in mind, it is only important to know that there is a \emph{second variation formula} for $f\in \mathcal{H}^{^{2,2}}\left(\X\right)$; see Theorem \ref{thm:secvar}. 
\par Finally, let us  recall the following chain rule for $\Hess\left(f\right)$ \cite[Proposition 3.3.21]{Gigli-2} that we will need later to show the potential function is Hessian-free.
\begin{proposition}\label{prop:chainrulehessian}
	Let $f\in \mathcal{W}^{^{2,2}}\left(\X\right)$ with $\left|\nabla f\right|\in \LL^{\hspace{-2pt}^\infty}\left(\X, \meas\right)$ and let $\upvarphi\in \mathcal{C}^{^2}\left(\R\right)$ with uniformly bounded first and second derivatives. 
	Then $\upvarphi\circ f\in \mathcal{W}^{^{2,2}}\left(\X\right)$ and 
	\begin{align}\label{eq:hessian-chain}
		\Hess\left( \upvarphi \circ f\right) = \upvarphi^{{''}}\circ f \; \dif f\otimes \dif f + \upvarphi^{{'}}\circ f \; \Hess\left(f\right),
	\end{align}
	where $\dif f$ is the {\it differential} of $f$ that is defined via $\dif f\left(\nabla g\right)= \langle \nabla f, \nabla g\rangle$ for all $g\in \Sobol\left(\X\right)$; see \cite{Gigli-2}. 
\end{proposition}
\begin{remark}
	If $f\in \Dom\left(\Delta\right) \cap  \Delta^{\hspace{-2pt}^{-1}}\left(\Ltwo\right)$, then $\upvarphi \circ f\in \mathcal{H}^{^{2,2}}\left(\X\right)$. Indeed, by Proposition \ref{prop:chainrulelaplacian}, $\upvarphi \circ f\in \Dom\left(\Delta\right) \cap  \Delta^{\hspace{-2pt}^{-1}}\left(\Ltwo\right)$ and therefore $\upvarphi \circ f\in \mathcal{H}^{^{2,2}}\left(\X\right) \subset \mathcal{W}^{^{2,2}}\left(\X\right)$. Moreover, there is at most one Hessian of $u$ (compare with the remarks after in~\cite[Definition 3.3.1]{Gigli-2}) and hence one can compute $\Hess\left(\upvarphi \circ f\right)$ according to the formula~\eqref{eq:hessian-chain}.
\end{remark}
\subsection{Calculus and geometry in $\RCD\hspace{-1pt}\left(\K,\N\right)$ spaces}
\par Below, we will delve deeper into the more intricate analytic and geometric properties of $\RCD^*\hspace{-2pt}\left(\K,\N\right)$ spaces that we will take advantage of in this paper. 
\subsubsection*{\small \bf \textit{Cutoff functions and higher order locality}}
To localize some analytic computations, we will need the existence of suitable cutoff functions on $\RCD$ spaces. The proof of the following lemma can be obtained by arguing as in the proof of~\cite[Lemma 3.1]{MN} or~\cite[Lemma 6.7]{AMS}.
\begin{lemma}[Cutoff functions]\label{lem:cutoff}
		Let $\left(\X,\dist,\meas\right)$ be a compact $\RCD\hspace{-1pt}\left(\K,\N\right)$ space for $\K\in \mathbb R$ and $\N\in \left[1,\infty\right)$. Let $U, V\subset \X$ be open sets such that $\ovs{V}\subset U$. Then, there exists $\upvarphi\in \mathbb{D}_{\hspace{-1pt}_{\infty}}$ with the properties
		\begin{enumerate}
			\item $0\leq \upvarphi  \leq 1$, $\upvarphi \equiv 1$ on $\ovs{V}$ and $\supp\left(\upvarphi\right)\subset U$, 
			\smallskip
			\item[] and \smallskip
			\item $\Delta \upvarphi ,\ \left|\nabla \upvarphi  \right| \in \LL^{\hspace{-2pt}^\infty}\left(\X, \meas\right)$. 
		\end{enumerate}
	\end{lemma}
Now let us briefly touch upon the locality property for $\Hess$. Let $ \Omega$ be an open subset in $\X$ with $\meas \left(  \partial  \Omega \right) = 0$. Denote the corresponding subspace (closure) by $\left( \ovs{ \Omega}, \dist_{\ovs{\Omega}}, \meas\mres_{\ovs{\Omega} }  \right)$ where $\dist_{\ovs{\Omega}}$ is the induced intrinsic metric on $\ovs{\Omega}$. Denote the minimal weak upper gradient in $\ovs{\Omega}$ by $\left|  \ovs{\nabla} \, \cdot \,\right|$. 
\begin{lemma}[Higher order locality]\label{lem:ho-locality}
	For $f \in  \mathbb{D}_{\hspace{-1pt}_{\infty}}$, it holds
\begin{align*}
	\left|\ovs{\nabla} \left(\left|\ovs{\nabla} \left(f\restr_{\ovs{\Omega}}\right)\right|^{^2} \right) \right|= \left|\nabla \left|\nabla f\right|^{^2}\right|\restr_{\ovs{\Omega}}.
\end{align*}	
\end{lemma}
\begin{proof}[\footnotesize \textbf{Proof}]
Recall $\left| \ovs{\nabla} \left( f\restr_{\ovs{\Omega}} \right) \right|= \left|\nabla f\right|\restr_{\ovs{\Omega}}$ holds for $f\in \Sobol\left(\X\right)$ by locality of the minimal weak upper gradient (see Remark \ref{rem:locgra}) . Therefore, if $\left|\nabla f\right|^{^2}\in \Sobol$ (as holds for $f \in  \mathbb{D}_{\hspace{-1pt}_{\infty}}$), one gets the conclusion by iterating the locality. 
\end{proof}
By virtue of locality of minimal weak upper gradient and higher order locality, we will not differentiate between $\ovs{\nabla}$ and $\nabla$ and will use $\nabla$ for both. 
\par The locality of Hessian is then easily inherited from the higher locality as we saw in Lemma~\ref{lem:ho-locality}. Denote the Hessian in $\ovs{\Omega}$ by $\ovs{\Hess}$.
\begin{lemma}\label{lem:hes-local}
		For $f \in  \mathbb{D}_{\hspace{-1pt}_{\infty}}$,
\begin{align*}
	\ovs{\Hess}\left(  f\restr_{\ovs{\Omega}}\right)\left( \left(\nabla v\right)\restr_{\ovs{\Omega}}, \nabla \left(w\right)\restr_{\ovs{\Omega}} \right) = \Hess\left( f \right)\left( \nabla v, \nabla w \right)\restr_{\ovs{\Omega}}.
\end{align*}
\end{lemma}
\begin{proof}[\footnotesize \textbf{Proof}]
	The proof is an easy iteration of the higher-order locality; see the explanation at the very end of \textsection\thinspace\ref{sec:rcd-specanal}. 
\end{proof}
Based on Lemma~\ref{lem:hes-local}, going forward we will not differentiate between $\ovs{\Hess}$ and $\Hess$ and will just write $\Hess$ for both. 
\subsubsection*{\small \bf \textit{First and second variation formulae in the $\RCD$ context}}
\begin{definition}
	Let $f\in \Sobol\left(\X\right)$. We say $\ppi\in \Prob\left(\mathcal{C}\left(\left[0,1\right], \X\right)\right)$ represents  $\nabla f$ if $\ppi$ is a test plan and it holds 
	
	\begin{align}\label{ineq:grad}
	&	\lowlim_{t\downarrow 0} \int_{\hspace{-1pt}_{\mathsf{Curves}}} \nicefrac{1}{t} \, \big( f\circ \upxi\left(t\right) \shortminus f\circ \upxi\left(0\right)  \big) \, \dif\ppi\left(\upxi\right)\geq \\ & \nicefrac{1}{2} \,\int_{\hspace{-1pt}_{\X}}  \left|\nabla f\right|^{^2} \, \dif \left(\e_{_0}\right)_{*}\ppi+ \nicefrac{1}{2}\, \uplim_{t\downarrow 0} \, \nicefrac{1}{t}\, \int_{\hspace{-1pt}_{\mathsf{Curves}}}  \int_{\hspace{-1pt}_0}^{t} \left|\bdot{\upxi}\left(s\right)\right|^{^2} \, \dif s \, \dif{\ppi}\left(\upxi\right)\notag
	\end{align}
in which $\mathsf{Curves} := \mathcal{C} \left( \left[0,1\right],\X \right)$,
\end{definition}
The next theorem can be found in~\cite[Theorem 3.10]{Gigli-1}; also see~\cite{cavmonlap}.
\begin{theorem}[First variation formula~\cite{GT-2}]\label{th:firstvariation}
	Let $f,u \in \Sobol\left(\X\right)$ and consider $\ppi\in \Prob\left(\mathcal{C}\left(\left[0,1\right], \X\right) \right)$ that represents $\nabla u$. Then it holds
	\begin{align*}
		\int_{\hspace{-1pt}_{\X}} \langle \nabla f, \nabla u\rangle \,  \dif \left(\e_{_0}\right)_{*} \ppi= \lim_{t\downarrow 0} \int_{\hspace{-1pt}_\mathsf{Curves}} \nicefrac{1}{t}\, \big( f\circ \upxi\left(t\right) \shortminus f\circ \upxi\left(0\right)   \big) \, \dif\hspace{0.7pt}\ppi\left(  \upxi \right).
	\end{align*}
\end{theorem}
\begin{remark}[First variation for Wasserstein geodesics and the metric Brenier theorem]\label{rem:metbre}
		A special situation where the previous theorem applies, is the case of a Kantorovich potential $\Upphi$ for the pair $ \upmu_{\hspace{1pt}_0} \,, \upmu_{\hspace{1pt}_1} \in \Prob_{\hspace{-5pt}_2}\left(\X,\meas\right)$  and a plan $\PP$ such that $\left( \e_{\substack{\\\hspace{0.5pt}t}} \right)_*{\PP}$ is a Wasserstein geodesic between $\upmu_{\hspace{1pt}_0}$ and $\upmu_{\hspace{1pt}_1}$; in this case, for $f\in \Sobol\left(\X\right)$, the map $t\mapsto \int f \, \dif \upmu_{\substack{\\\hspace{0.5pt}t}}$ is differentiable at $t=0$ and 
		\begin{align*}
			\nicefrac{\dif}{\dif t}\,\restr_{t=0} \int_{\hspace{-1pt}_{\X}} f \,\dif\upmu_{\substack{\\\hspace{0.5pt}t}}= \shortminus \int_{\hspace{-1pt}_{\X}} \langle \nabla f, \nabla \Upphi \rangle \, \dif\upmu_{\hspace{1pt}_0}.
		\end{align*}
\par In this setting, we also recall the identity $\left|\nabla \Upphi \right|\left( \upgamma\left(0\right) \right) = \dist\left(\upgamma\left(0\right),\upgamma\left(1\right)\right)$ for $\PP$-a.e. $\upgamma$. This is the metric Brenier theorem for $\RCD$ spaces; for the full statement with references, see \cite[Theorem 2.4]{GT-2}. 
	\end{remark}
\begin{theorem}[Second variation formula~\cite{GT}]\label{thm:secvar}
	Let $\left(\X,\dist,\meas\right)$ be a metric measure space that satisfies the $\RCD\left(\K,\N\right)$ condition for $\N<\infty$. Take $\upmu_{\hspace{1pt}_0} \,, \upmu_{\hspace{1pt}_1} \in \Prob_{\hspace{-5pt}_1}\left(\X\right)$ such that $\upmu_{\substack{\\[1pt]\hspace{1pt}i}}=\uprho_{\substack{\\[3pt]\hspace{-2pt}i}} \meas\leq C\meas$ for $C>0$ and $i=0,1$; suppose they are  joined by the unique $\Ltwo$-Wasserstein geodesic $\left( \upmu_{\substack{\\\hspace{0.5pt}t}} \right)_{t\in \left[0,1\right]}$. Let $f\in \mathcal{H}^{^{2,2}}\left(\X\right)$. Then the  map $t\in \left[0,1\right]\mapsto \int f\, \dif \upmu_{\substack{\\\hspace{0.5pt}t}}$ belongs to $\mathcal{C}^{\hspace{0.5pt}^2}\left(\left[0,1\right]\right)$ and for every $t\in \left[0,1\right]$, it holds 
	\begin{align*}
		\nicefrac{\dif^{^{\,2}}}{\dif t^{^2}} \int_{\hspace{-1pt}_{\X}} f \, \dif \upmu_{\substack{\\\hspace{0.5pt}t}}= \int_{\hspace{-1pt}_{\X}} \Hess \left(f\right)\left(\nabla \Upphi_{\substack{\\\hspace{-1pt}t}}, \nabla \Upphi_{\substack{\\\hspace{-1pt}t}}\right)\, \dif\upmu_{\substack{\\\hspace{0.5pt}t}},
	\end{align*}
	where $\Upphi_{\substack{\\\hspace{-1pt}t}}$ is the function with the property that for some $t_{_1}\neq t_{_2} \in \left[0,1\right]$ the function $\shortminus \left(t_{_2} \shortminus t_{_1}\right) \Upphi_{\substack{\\\hspace{-1pt}t}}$ is a Kantorovich potential between $\upmu_{\substack{\\\hspace{0.5pt}t_{_1}}} $ and $\upmu_{\substack{\\\hspace{0.5pt}t_{_2}}} $.
\end{theorem}
\subsubsection*{\small \bf \textit{Maximum principles in metric measure spaces}}
\begin{proposition}[{Strong maximum principle~\cite{BB}}]\label{prop:strong-max-1}
	Suppose $\left( \X,\dist,\meas\right)$ is a doubling metric measure space that supports a $1$-$2$ weak local Poincar\'{e} inequality. Let $f \in \mathcal{S}^{\hspace{1pt}^2}_{\mathsf{loc}}\left(\X\right)$ be continuous and such that
	\begin{align*}
	\int_{\hspace{-1pt}_{\Omega}}  \left| \nabla f \right|^{^2} \, \dif\meas \le \int_{\hspace{-1pt}_{\Omega}}  \left| \nabla \left(f + g\right) \right|^{^2} \, \dif\meas ,
	\end{align*}
	holds for any $g \in \mathsf{Test}\left(\Omega\right) := \Lip_{\sf c}\left( \Omega  \right)$ (${\sf c}$ for compactly supported) and bounded open set $\Omega$ with $\supp \left( g \right) \subset \Omega \subset \X$. If $f$ attains its maximum inside $\Omega$, then $f$ is constant throughout $\Omega$.
	\par In the setting of $\RCD\hspace{-1pt}\left(\K,\N\right)$ spaces, the strong maximum principle yields the following statement which is the maximum principle we will resort to, later on in \textsection\thinspace\ref{sec:rcd-specanal}.
\end{proposition}
\begin{theorem}[{Strong maximum principle in RCD setting~\cite{Gigli-Rigoni}}]\label{thm:Gig-Rig}
	Let $\left(\X, \dist, \meas  \right)$ be an $\RCD^*\hspace{-2pt}\left(\K,\N\right)$ space with $1 \le \N \le \infty$ and $\Omega \subset \X$ an open and connected subset. Suppose $f \in \Sobol\left(\Omega \right) \, \cap \,  \mathcal{C} \left( \ovs{\Omega}  \right)$ be a sub-minimizer of $\CHE$ such that $\max\limits_{x \in \ovs{\Omega}} f\left( x \right) = f\left( x_{\hspace{-0.5pt}_0} \right)$ for some $x_{\hspace{-0.5pt}_0} \in \Omega$. Then $f$ is constant. 
\end{theorem}
\subsection{Gradient comparison and sharp spectral gap for $\RCD^*\hspace{-2pt}\left(\K,\N\right)$ spaces}
\subsubsection*{\small \bf \textit{From Rayleigh quotient to the lowest positive eigenvalue}}
\par The spectral gap $\feig$ for a compact metric measure space $\left( \X,\dist,\meas  \right)$, is defined by
\begin{align*}
		\feig\left(\X\right) =\inf\left\{ \left\| u \right\|^{^{-2}}_{\substack{\\\Ltwo}}\, \left\| \nabla u \right\|^{^{2}}_{\substack{\\\Ltwo}} \;\; \text{\textbrokenbar} \;\; u\in \Lip\left( \X \right) \smallsetminus \left\{0\right\}\; \text{and}\; \int_{\hspace{-1pt}_{\X}}  u\, \dif\meas = 0\right\}\notag.
\end{align*}
\par It turns out that the spectral theorem for the Laplacian for Riemannian manifolds generalizes to the RCD setting without a hitch. 
Indeed, since $\X$ is compact and in virtue of the integration by parts property which is incorporated in the definition of $\Delta$, the local Poincar\'e inequality implies that the resolvent $(\mathrm{I} \shortminus \tau\Delta)^{^{-1}}$ for any $\tau>0$ is a compact self-adjoint operator which obviously does not include $0$ in its spectrum; see~\cite{Raj-3}. 
\par Thus the standard spectral theorem implies $\shortminus \Delta$ has a non-negative discrete real spectrum and there exists an orthonormal basis for the underlying Hilbert space, in this case $\Sobol$, formed by eigenfunctions. Now standard spectral arguments (e.g. see~\cite{Davies}) imply 
$\feig\left(\X\right)$ is the smallest positive eigenvalue of $\shortminus \Delta$; in particular, any minimizer $u$ of the Rayleigh quotient satisfies 
\begin{align*}
\Delta \left(u\right) = \shortminus \feig u.
\end{align*}
Conversely, any such eigenfunction must minimize the Rayleigh quotient. 
\subsubsection*{\small \bf \textit{Gradient comparison and sharp spectral gap}}
\par We recall a crucial gradient comparison for $u$ which will be used extensively in the later sections. This gradient comparison has been proven in~\cite{KeObata} for $\K>0$ and in~\cite{Jiang-Zhang} for all $\K \in \mathbb{R}$; also see~\cite{CM-1,CM-2} for sharp spectral gaps.
\begin{proposition}[\cite{KeObata,Jiang-Zhang}]\label{th:gradestimate}
	Let $\left(\X,\dist,\meas\right)$ be a compact $\RCD^*\hspace{-2pt}\left(\K,\N\right)$ space with $\diam \left(  \X \right) = \mathcal{D}$. Let $\feig$ be the first nonzero eigenvalue of $\X$ and suppose $\feig > \max \left\{   0 , \left(\N^{\hspace{0.5pt}^{'}} \shortminus 1\right)^{\hspace{-4pt}^{-1}}\N^{\hspace{0.5pt}^{'}}\K \right\}$ for some $\N^{^{'}} \ge \N$. Let $u$ be an eigenfunction with respect to $\feig$ and let $v$ be a Neumann eigenfunction for the ${\nicelambda}\hspace{-1pt}\left(\N,\K,\mathcal{D}\right)$ in the $1\D$ model problem. If 
$
		\left[\min u , \max u\right] \subset \left[\min v , \max v\right],
$
then
	\begin{align*}
	\left| \nabla u \right|^{^2} \left(x\right) \le \left( v^{'} \circ v^{-1}  \right)^{\hspace{-3pt}^2}  \left( u\left(x\right) \right)  \quad \text{for} \quad \meas\mbox{-a.e.} \; x \in \X .
	\end{align*}
	Moreover, $\min u= \shortminus \max u$.
\end{proposition}
\subsection{$1$D Localization schemes}\label{subsec:1Dlocschemes}
\par Another crucial tool that we will utilize in our analysis, is a one dimensional localization of the curvature dimension conditions; this localization is given by disintegration of the measure with respect to the geodesics which are in a sense gradient flow curves of a globally locally-$1$-Lipschitz potential function. Here, globally locally-$1$-Lipschitz function means the special case where the function $f$ has local Lipschitz constant equal to $1$ everywhere. We will simply refer to such functions as $1$-Lipschitz potential functions. 
\par In convex geometry (convexity usually implies lower curvature bounds automatically), such a localization technique is called a needle decomposition; see \cite{CM-1} and references therein. The smooth Riemannian version of needle decomposition under lower Ricci curvature bound is proven in~\cite{klartagneedle}. For {$\RCD\hspace{-1pt}\left(\K,\N\right)$ spaces,} such a localization result is shown in~\cite{CM-1, CM-2}. 
\subsubsection*{\small \bf \textit{Setup}}
\par  We assume familiarity with basic concepts in optimal transport.
Let $\left(\X,\dist,\meas\right)$ be a locally compact metric measure space that is essentially non-branching. We  assume that $\supp \left(\meas\right) =\X$.
\par Let $f:\X\rightarrow \mathbb{R}$ be a $1$-Lipschitz potential function. Then 
\begin{align*}
\fgamma_{\substack{\\[2pt]\hspace{-2pt}f}}:=\left\{\left(x,y\right)\in  \X^{^2}  \;\; \text{\textbrokenbar} \;\;  f\left(x\right) \shortminus f\left(y\right)=\dist\left(x,y\right)\right\},
\end{align*}
is a  $\dist$-cyclically monotone set; set
\begin{align*}
	\fgamma_{\substack{\\[2pt]\hspace{-2pt}f}}^{^{-1}}=\left\{\left(x,y\right)\in \X^{^2} \;\; \text{\textbrokenbar} \;\; \left(y,x\right)\in \fgamma_{\substack{\\[2pt]\hspace{-2pt}f}}\right\}.
\end{align*}
Denote by $\mathcal  G^{^{\left[a,b\right]}}\left(\X\right)$, the set of all unit speed geodesics with domain $\left[a,b\right]\subset\R$. If $\upgamma\in \mathcal G^{^{\left[a,b\right]}}\left(\X\right)$ for some $a\le b$ with $\left(\upgamma\left(a\right),\upgamma\left(b\right)\right)\in \fgamma_{\substack{\\[2pt]\hspace{-2pt}f}}$, then $\left(\upgamma\left(t\right),\upgamma\left(s\right)\right)\in \fgamma_{\substack{\\[2pt]\hspace{-2pt}f}}$ for all $a<t\leq s<b$.  It is therefore natural to consider the subset $\mathscr{G}^{^{\left[a,b\right]}}$ of \emph{unit speed transport geodesics} $\upgamma:\left[a,b\right]\rightarrow \R$ for which $\left(\upgamma\left(t\right),\upgamma\left(s\right)\right)\in \fgamma_{\substack{\\[2pt]\hspace{-2pt}f}}$ for all $a\leq t\leq s\leq b$. 
\par The union $\fgamma_{\substack{\\[2pt]\hspace{-2pt}f}} \cup 	\fgamma_{\substack{\\[2pt]\hspace{-2pt}f}}^{^{-1}}$ defines a relation ${\frak R}_{\substack{\\[1pt]\hspace{-0.5pt}f}} $ on $\X\times \X$; the relation ${\frak R}_{\substack{\\[1pt]\hspace{-0.5pt}f}} $ induces a {\it transport set with endpoints} 
\begin{align*}
	\T_{\substack{\\[2pt]\hspace{-2pt}f}} := \mathrm{proj}_{_1}\left( {\frak R}_{\substack{\\[1pt]\hspace{-0.5pt}f}} \smallsetminus \left\{\left(x,y\right) \; \text{\textbrokenbar} \; x=y\right\}\right)\subset \X,
\end{align*}
where $\mathrm{proj}_{_1} \left(x,y\right)=x$. For $x\in \T_{\substack{\\[2pt]\hspace{-2pt}f}}$ one defines
\begin{align*}
\fgamma_{\substack{\\[2pt]\hspace{-2pt}f}} \left(x\right):=\left\{y\in \X \;\; \text{\textbrokenbar} \;\; \left(x,y\right)\in \fgamma_{\substack{\\[2pt]\hspace{-2pt}f}} \right\}, 
\end{align*}
and similarly, $\fgamma_{\substack{\\[2pt]\hspace{-2pt}f}}^{^{-1}}\left(x\right)$ as well as ${\frak R}_{\substack{\\[1pt]\hspace{-0.5pt}f}}\left(x\right)=\fgamma_{\substack{\\[2pt]\hspace{-2pt}f}}\left(x\right)\cup \fgamma_{\substack{\\[2pt]\hspace{-2pt}f}}^{^{-1}}\left(x\right)$. Since $f$ is $1$-Lipschitz, 
$\fgamma_{\substack{\\[2pt]\hspace{-2pt}f}}, \, \fgamma_{\substack{\\[2pt]\hspace{-2pt}f}}^{^{-1}}$ and ${\frak R}_{\substack{\\[1pt]\hspace{-0.5pt}f}}$ are all closed, so are $\fgamma_{\substack{\\[2pt]\hspace{-2pt}f}}\left(x\right),\, \fgamma_{\substack{\\[2pt]\hspace{-2pt}f}}^{^{-1}}\left(x\right)$ and ${\frak R}_{\substack{\\[1pt]\hspace{-0.5pt}f}} \left(x\right)$.
\par The {\it transport set without forward and/or backward branching} (i.e. the non-branching transport set) associated to $f$ is then defined as 
\begin{align*}
\T_{\substack{\\[2pt]\hspace{-2pt}f}}\strut^{\,\textsf{nb}}=\left\{ x\in  \T_{\substack{\\[2pt]\hspace{-2pt}f}} \;\;\text{\textbrokenbar}\;\;    \left(y,z\right)\in {\frak R}_{\substack{\\[1pt]\hspace{-0.5pt}f}}\quad  \forall  y,z\in {\frak R}_{\substack{\\[1pt]\hspace{-0.5pt}f}} \left(x\right)  \right\}.
\end{align*}
\par Notice $\T_{\substack{\\[2pt]\hspace{-2pt}f}}$ and $\T_{\substack{\\[2pt]\hspace{-2pt}f}} \smallsetminus \T_{\substack{\\[2pt]\hspace{-2pt}f}}\strut^{\,\textsf{nb}}$ are $\upsigma$-compact and the sets $\T_{\substack{\\[2pt]\hspace{-2pt}f}}\strut^{\,\textsf{nb}}$ and ${\frak R}_{\substack{\\[1pt]\hspace{-0.5pt}f}} \cap\left( \T_{\substack{\\[2pt]\hspace{-2pt}f}}\strut^{\,\textsf{nb}} \times \T_{\substack{\\[2pt]\hspace{-2pt}f}}\strut^{\,\textsf{nb}} \right)$ are Borel sets.
\par It is shown in~\cite{Cav} that ${\frak R}_{\substack{\\[1pt]\hspace{-0.5pt}f}}$ restricted to $\T_{\substack{\\[2pt]\hspace{-2pt}f}}\strut^{\,\textsf{nb}} \times \T_{\substack{\\[2pt]\hspace{-2pt}f}}\strut^{\,\textsf{nb}}$ is an equivalence relation; denote this relation by ${\frak R}\strut^{\,\textsf{nb}}_{\substack{\\[1pt]\hspace{-0.5pt}f}}$.
Hence, from ${\frak R}_{\substack{\\[1pt]\hspace{-0.5pt}f}}$ one obtains a partition of $\T_{\substack{\\[2pt]\hspace{-2pt}f}}\strut^{\,\textsf{nb}}$ into a disjoint family of equivalence classes $\left\{\X_{\hspace{1pt}q}\right\}_{q\in \mathrm{Q}}$. Moreover, $\T_{\substack{\\[2pt]\hspace{-2pt}f}}\strut^{\,\textsf{nb}}$ is also $\upsigma$-compact. We also write  ${\frak R}\strut^{\,\textsf{nb}}_{\substack{\\[1pt]\hspace{-0.5pt}f}}(x)$ for  
the equivalence class of $x\in  \T_{\substack{\\[2pt]\hspace{-2pt}f}}\strut^{\,\textsf{nb}}$ in  $\T_{\substack{\\[2pt]\hspace{-2pt}f}}\strut^{\,\textsf{nb}}$; in particular, ${\frak R}\strut^{\,\textsf{nb}}_{\substack{\\[1pt]\hspace{-0.5pt}f}}(q)=\X_q$ for $q\in \mathrm{Q}$.

\par Every $\X_{\hspace{1pt}q}$ is isometric to some interval  $\mathrm{I}_{\hspace{1pt}q}\subset\mathbb{R}$ via a distance-preserving map $\upgamma_{\substack{\\[1pt]\hspace{-2pt}q}}:\mathrm{I}_{\hspace{1pt}q} \rightarrow \X_{\hspace{1pt}q}$. The map $\upgamma_{\substack{\\[1pt]\hspace{-2pt}q}}:\mathrm{I}_{\hspace{1pt}q}\rightarrow \X_{\hspace{1pt}q}$ extends to an arclength parametrized geodesic defined on the closure $\ovs{\mathrm{I}}_{\hspace{1pt}q}$ of $\mathrm{I}_{\hspace{1pt}q}$. Let $\ovs{\mathrm{I}}_{\hspace{1pt}q}=\left[\ai_{\hspace{1pt}q},\bi_{\hspace{1pt}q}\right]$.
The set of equivalence classes $\mathrm{Q}$ has a measurable structure such that the quotient map $\mathfrak Q: \T_{\substack{\\[2pt]\hspace{-2pt}f}}\strut^{\,\textsf{nb}} \rightarrow \mathrm{Q}$ onto equivalence classes, is a measurable map. 
We set ${\bf q}:= \mathfrak Q_{\substack{\\ *}} \left(\meas\mres_{\T^{^{\,\textsf{nb}}}_{\substack{\\\hspace{-2.5pt}f}}}\right)$.
\subsubsection*{\small \bf \textit{Disintegration}}
\par It is shown in \cite[Proposition 5.2]{Cav} that there exists a measurable section $\frak{s}$ of $\mathfrak{R}_{\substack{\\[1pt]\hspace{-0.5pt}f}}$ on $\T_{\substack{\\[2pt]\hspace{-2pt}f}}\strut^{\,\textsf{nb}} $. So, one can identify the measurable space $\Q$ with the fixed points of $\frak{s}$ i.e.
\begin{align*}
	\Q \cong \left\{x\in \T_{\substack{\\[2pt]\hspace{-2pt}f}}\strut^{\,\textsf{nb}} \; \text{\textbrokenbar} \; \; \mathfrak{s}\left(x\right) = x\right\},
\end{align*}
where the set of fixed points is equipped with the induced measurable structure. This way, $\mathfrak Q$ is identified with $\mathfrak{s}$. By inner regularity, one can find a $\upsigma$-compact set, $\Q^{'}\subset \X$ with ${\bf q}\left(\Q \smallsetminus \Q^{'}\right)=0$ and we can see ${\bf q}$ as a Borel measure on $\X$.  In what follows, we will replace $\Q$ with $\Q^{'}$ without further notice while (with a slight abuse of notations) keeping the same notation $\Q$. {We parametrize $\upgamma_{\substack{\\[1pt]\hspace{-2pt}q}}$ for $q\in \Q$ such that $\upgamma_{\substack{\\[1pt]\hspace{-2pt}q}}\left(0\right)=\mathfrak{s}\left(q\right)$. In particular, $0\in \mathrm{I}_{\hspace{1pt}q}$. }
Now, let us assume $\left(\X,\dist,\meas\right)$ is an essentially non-branching $\CD^*\hspace{-2pt}\left(\K,\N\right)$ space for $\K\in \R$ and $\N\geq 1$. 
\par The following result on exhaustion by non-branching transport set, is proven in~\cite{CM-1}.
\begin{lemma}[{\cite[Theorem 3.4]{CM-1}}]\label{somelemma}
	Let $\left(\X,\dist,\meas\right)$ be an essentially non-branching $\CD^*\hspace{-2pt}\left(\K,\N\right)$ space for $\K\in \R$ and $\N\in \left(1,\infty\right)$ with $\supp \left(\meas\right) = \X$ and $\meas\left(\X\right)<\infty$.
	Then, for any $1$-Lipschitz potential function $f:\X\rightarrow \R$, it holds $\meas\left(\T_{\substack{\\[2pt]\hspace{-2pt}f}}  \smallsetminus \T_{\substack{\\[2pt]\hspace{-2pt}f}}\strut^{\,\textsf{nb}}  \right)=0$.
\end{lemma}
\par For ${\bf q}$-a.e. $q\in \Q$ it was shown in \cite[Theorem 7.10]{cavmil} that 
\begin{align*}
{\frak R}_{\substack{\\[1pt]\hspace{-0.5pt}f}}\left(x\right)=\ovs{{\frak R}_{\substack{\\[1pt]\hspace{-0.5pt}f}}\strut^{\,\textsf{nb}}\left(x\right)}=\ovs{\X}_{\hspace{1pt}q}\supset \X_{\hspace{1pt}q} =  {\frak R}\strut^{\,\textsf{nb}}_{\substack{\\[1pt]\hspace{-0.5pt}f}}\left(x\right)
\quad \forall x\in \mathfrak Q^{\hspace{-1pt}^{-1}}\left(q\right).
\end{align*}
Consequently, one gets a needle decomposition as follows;
	\begin{theorem}[{\cite{CM-1, CM-2}}]\label{thm:1dlocalization}
	Let $\left(\X,\dist,\meas\right)$ be a compact geodesic metric measure space with $\supp \left(\meas\right) = \X$ and $\meas$ finite. Let $f:\X\rightarrow \mathbb{R}$ be a $1$-Lipschitz potential function;  let $\left(\X_{\hspace{1pt}q}\right)_{q\in \Q}$ be the induced partition of $\T_{\substack{\\[2pt]\hspace{-2pt}f}}\strut^{\,\textsf{nb}}$ via ${\frak R}_{\substack{\\[1pt]\hspace{-0.5pt}f}}$ and $\mathfrak Q: \T_{\substack{\\[2pt]\hspace{-2pt}f}}\strut^{\,\textsf{nb}} \rightarrow \Q$, the induced quotient map as above.
	Then, there exists a unique strongly consistent disintegration $\left\{\meas_{\hspace{1pt}q}\right\}_{q\in \Q}$ of $\meas\mres_{\T^{^{\,\textsf{nb}}}_{\substack{\\\hspace{-2.5pt}f}}}$ w.r.t. $\mathfrak Q$.
\end{theorem}
\subsubsection*{\small \bf \textit{Localization}}
Define the ray map ${\fancy{r}}:  \mathcal V\subset \Q\times \R\rightarrow \X$ via
\begin{align*}
	\mathsf{graph}\left(\ray\right)=\left\{ \left(q, t,x\right) \in \Q\times \R\times \X \; \; \text{\textbrokenbar} \;  \; \upgamma_{\substack{\\[1pt]\hspace{-2pt}q}}\left(t\right)=x\right\},
\end{align*}
where $\mathcal V= \ray^{-1}\left( \T_{\substack{\\[2pt]\hspace{-2pt}f}}\strut^{\,\textsf{nb}} \right)$. The map $\ray$ is Borel measurable; furthermore,  $\ray: \mathcal V\rightarrow  \T_{\substack{\\[2pt]\hspace{-2pt}f}}\strut^{\,\textsf{nb}}$ is bijective and its inverse is given by
\begin{align*}
\ray^{-1}\left(x\right)=\left( \mathfrak Q\left(x\right),  \dist_{\textsf{signed}}\left(x,\mathfrak Q\left(x\right)\right) \right),
\end{align*}
where the signed distance to $\mathfrak Q$ is characterized via the identity
\begin{align*}
	 \upgamma_{\substack{\\[1pt]\hspace{-2pt}q}} \left(  \mathfrak{s}\left(x \right) + \dist_{\textsf{signed}}(x) \right) = x.
\end{align*}
\begin{theorem}[{Localization of CD conditions~\cite{CM-1}}]\label{th:1dlocalisation}
	Let $\left(\X,\dist,\meas\right)$ be an essentially non-branching $\CD^*\hspace{-2pt}\left(\K,\N\right)$  space with $\supp\left(\meas\right)=\X$, $\meas\left(\X\right)<\infty$, $\K\in \R$ and $\N\in \left(1,\infty\right)$.
	\par Then, for any $1$-Lipschitz potential function $f:\X\rightarrow \R$, there exists a disintegration $\left\{\meas_{\hspace{1pt}q}\right\}_{q\in \Q}$ of $\meas$ that is strongly consistent with ${\frak R}\strut^{\,\textsf{nb}}_{\substack{\\[1pt]\hspace{-0.5pt}f}}$.
	\par Moreover, for ${\bf q}$-a.e. $q\in \Q$, $\meas_{\hspace{1pt}q}$ is a Radon measure with $\meas_{\hspace{1pt}q}=\textsf{h}_{\substack{\\[-1pt]q}}\,\mathscr{H}^{^1}\mres_{\X_{\hspace{1pt}q} }$ and $\left(\X_{\hspace{1pt}q} , \dist_{\substack{\\[1pt]\X_{\hspace{1pt}q} }}, \meas_{\hspace{1pt}q} \right)$ verifies the $\CD\hspace{-1pt}\left(\K,\N\right)$ conditions.
	\par More precisely, for ${\bf q}$-a.e. $q\in \Q$,
	\begin{align}\label{kuconcave}
		\textsf{h}_{\substack{\\[-1pt]q}}\left(\upiota_{\substack{\\\hspace{0.5pt}t}}\right)\strut^{\hspace{-2pt}\nicefrac{1}{\left(\N-1\right)}}\geq  \upsigma^{^{\left(1-t\right)}}_{\hspace{-1pt}_{\K,\,\N-1}} \left(\left|\bdot \upiota\right|\right)\textsf{h}_{\substack{\\[-1pt]q}}\left(\upiota_{\hspace{1pt}_0}\right)^{\hspace{-2pt}^{\nicefrac{1}{\left(\N-1\right)}}}+\upsigma^{^{\left(t\right)}}_{\hspace{-1pt}_{\K,\,\N-1}}\left(\left|\bdot \upiota\right|\right)\textsf{h}_{\substack{\\[-1pt]q}}\left(\upiota_{\hspace{1pt}_1}\right)^{\hspace{-2pt}^{\nicefrac{1}{\left(\N-1\right)}}},
	\end{align}
holds along every constant speed geodesic $\upiota:\left[0,1\right]\rightarrow \left(\ai_{\hspace{1pt}q},\bi_{\hspace{1pt}q}\right)$.
\end{theorem}
\begin{remark} From the property
	\eqref{kuconcave}, one deduces $\textsf{h}_{\substack{\\[-1pt]q}}$ is locally Lipschitz continuous and geodesically $\left( \K, \N \right)$-concave on $ \left(\ai_{\hspace{1pt}q},\bi_{\hspace{1pt}q}\right)$. Therefore, $\textsf{h}_{\substack{\\[-1pt]q}} :\R \rightarrow \left(0,\infty\right)$ satisfies
	\begin{align*}
		\nicefrac{\dif^{^{2}}}{\dif r^{^2}}\, \textsf{h}_{\substack{\\[-1pt]q}}^{\hspace{-2pt}^{\nicefrac{1}{\left(\N-1\right)}}}+ \nicefrac{\K}{\left(\N-1\right)}\, \textsf{h}_{\substack{\\[-1pt]q}}^{\hspace{-2pt}^{\nicefrac{1}{\left(\N-1\right)}}}\leq 0,
	\end{align*}
 on $\left(\ai_{\hspace{1pt}q},\bi_{\hspace{1pt}q}\right)$ in distributional sense; see~\cite[Section 4]{CM-1} for further details. 
\end{remark}
\par In \cite{CM-1, CM-2}, the Theorem~\ref{th:1dlocalisation} is used to derive the following localization scheme that is used in particular to derive sharp estimates for the first positive eigenvalue of Laplacian.
\begin{theorem}[Singular 1D localization~\cite{CM-1,CM-2}]\label{th:1Dlocscheme}
	Let $\left(\X,\dist,\meas\right)$ be an essentially non-branching metric measure space with $\meas\left(\X\right)=1$ that satisfies $\CD\left(\K,\N\right)$ for $\K\in \mathbb{R}$ and $\N\geq 1$. Let $u:\X\rightarrow \mathbb{R}$ be $\meas$-integrable such that $\int_{_{\X}} u \, \dif\meas=0$ and assume there exists $x_{\hspace{-0.5pt}_0} \in \X$ such that $\int_{_{\X}} \left|u\left(x\right)\right| \dist\left(x, x_{\hspace{-0.5pt}_0}\right) \, \dif\meas<\infty$. 
\par Then, $\X$ can be written as the disjoint union of two sets $\mathcal Z$ (the zero set) and $\mathcal T$ (the transport set). The latter admits a partition of $\left\{\X_{\hspace{1pt}q}\right\}_{q\in \Q}$ where each $\X_{\hspace{1pt}q}$ is the image of a geodesic $\upgamma_{\substack{\\[1pt]\hspace{-2pt}q}}: \left[\ai_{\hspace{1pt}q},\bi_{\hspace{1pt}q}\right]\rightarrow \X$. Moreover, there exists a family of probability measures $\left\{\meas_{\hspace{1pt}q}\right\}_{q\in \Q}\subset \Prob\left(\X\right)$ with the following properties:
	\begin{enumerate}
		\item For any $\meas$-measurable set $\mathcal{B}\subset \T$ it holds 
		\begin{align*}
		\meas\left(\mathcal{B} \right)=\int_\Q \meas_{\hspace{1pt}q}\left( \mathcal{B} \right) \, \dif\mathbf{q}\left(q\right),
		\end{align*}
		where $\mathbf{q} := \mathfrak Q_{\substack{\\ *}} \left(\meas\mres_{\T^{^{\,\textsf{nb}}}_{\substack{\\\hspace{-2.5pt}f}}}\right) \in \Prob\left(\Q\right)$ and $\Q\subset \X$ is the quotient space.
		\item For $\mathbf{q}$-almost every $q\in \Q$, the set $\X_{\hspace{1pt}q}$ is the image of a geodesic $\upgamma_{\substack{\\[1pt]\hspace{-2pt}q}}$ with positive length and parametrized by arclength which supports the measure  $\meas_{\hspace{1pt}q}$.
		Moreover, $q\mapsto\meas_{\hspace{1pt}q}$ is a $\CD\left(\K,\N\right)$ disintegration, that is $\meas_{\hspace{1pt}q}= \left(\upgamma_{\substack{\\[1pt]\hspace{-2pt}q}}\right)_{\hspace{-1pt}*}\left(\uprho_{\substack{\\[2pt]\hspace{-2pt}q}} \, \mathscr{L}^{^1}\right)$ with 
		\begin{align*}
		\nicefrac{\dif^{^{\,2}}}{\dif r^{^2}}\, \uprho_{\substack{\\[2pt]\hspace{-2pt}q}}^{\hspace{-2pt}^{\nicefrac{1}{\left(\N-1\right)}}}\leq \shortminus \, \nicefrac{\K}{\left(\N-1\right)}\, \uprho_{\substack{\\[2pt]\hspace{-2pt}q}}^{\hspace{-2pt}{\nicefrac{1}{\left(\N-1\right)}}}\quad \text{ on } \quad \left(\ai_{\hspace{1pt}q},\bi_{\hspace{1pt}q}\right).
		\end{align*}
	\item For $\mathbf{q}$-almost every $q\in \Q$, it holds $\int_{_{\X_{_q}}} u \, \dif \meas_{\hspace{1pt}q}=0$ and $u=0$ $\meas$-a.e. in $\mathcal{Z}$.
	\end{enumerate}
\end{theorem}
\section{Spectral analysis in non-negatively curved metric measure spaces}\label{sec:rcd-specanal}
\par This section prepares and presents the key tools that we need in order to characterize the equality case of the sharp spectral gap estimate for $\RCD\left(0,\N\right)$ spaces with $1\le \N < \infty$.  Throughout, we take $\left(\X,\dist,\meas\right)$ to be a compact $\RCD\left(0,\N\right)$ space.
\par For the entirety of this section, we are assuming the equality holds in the spectral gap $\feig \ge \nicefrac{\uppi^2}{\diam^2}$ i.e. we assume $\X$ is spectrally-extremal. By rescaling the metric by $\nicefrac{\uppi}{\diam}$ if necessary, we can assume $\diam = \uppi$ and since the Lipschitz constants scale by $\nicefrac{\diam}{\uppi}$, one deduces $\feig$ scales by $\nicefrac{\diam^2}{\uppi^2}$ hence  we are assuming $\feig=1$. 
\par This section serves the main purpose of proving the following theorem.
\begin{theorem}\label{thm:hessian}
	Let  $u$ satisfy $\Delta \left(u\right)=\shortminus \, u$. Then,
	\begin{enumerate}
		\item $\left|\nabla u\right|^{^2} + u^2 \equiv \left(\max u\right)^{^2}$ \quad $\meas$-a.e.;
		\smallskip
		\item[] and \smallskip
	\item $\Hess\left(u\right)\left(\nabla g, \nabla g \right)= \shortminus \, u \langle \nabla f, \nabla g\rangle^{^2} $  \quad $\meas$-a.e.\  $\forall g\in \mathbb{D}_{\hspace{-1pt}_{\infty}}$,
	\end{enumerate}
	where $f=\sin^{\hspace{-1pt}^{-1}}\, \circ \, \left(\nicefrac{1}{\max u} \ u\right)$.
\end{theorem}
\subsection{Application of $1\D$-localisation}
\begin{lemma}
	Let $u$ be an eigenfunction of $\Delta$ with $\Delta \left( u\right)=\shortminus \, u$, then $u\in \mathbb{D}_{\hspace{-1pt}_{\infty}} $. In particular, $\left|\nabla u\right|^{^2}\in \Sobol\left(\X\right)$. 
	Moreover if $\left|u\right| \le 1$, then $\left|\nabla u\right| \le 1$ and $u$ has a $1$-Lipschitz regular representation.
\end{lemma}
\begin{proof}[\footnotesize \textbf{Proof}] Since $\Delta \left( u \right)=\shortminus \, u$, it holds $u=e^t \, \markov u\in \mathbb{D}_{\substack{\\\hspace{-1pt}\infty}}$. If $\left|u\right|\leq 1$, $\left|\nabla u\right|\leq 1$ follows from the gradient estimates of~\cite{Jiang-Zhang}. By the Sobolev-to-Lipschitz property for $\RCD$ spaces (e.g. see~\cite{AGSRiem, AGMR}) $u$ has a $1$-Lipschitz representation.
\end{proof}
Recall that $\int_{_{\X}} u \, \dif\meas =0$. Since $\X$ is compact, we can apply Theorem \ref{th:1Dlocscheme}.
\begin{lemma}\label{lem:ida}
	Let $\left( \X_{\hspace{1pt}q} \right)_{q\in \Q}$ be the decomposition of the transport set $\mathcal T$ in Theorem \ref{th:1Dlocscheme} and let $\left( \meas_{\hspace{1pt}q}\right)_{q\in \Q}$ be the corresponding disintegration of $\meas$. Then for ${\bf q}$-a.e. $q\in \Q$  we have 
	\begin{align}\label{eq:pi-length}
	\bi_{\substack{\\q}} \shortminus \ai_{\substack{\\q}}=\uppi,
	\end{align}
	\begin{align}\label{id1} 
		\nicefrac{\dif}{\dif r}\ u\circ \upgamma_{\substack{\\[1pt]\hspace{-2pt}q}} \left(r\right) = \left|\nabla u \right|\circ \upgamma_{\substack{\\[1pt]\hspace{-2pt}q}} \left(r\right) \geq 0 \quad \mathscr L^{^1}\mbox{-a.e.},
	\end{align}
	and
	\begin{align}\label{id2} 
		\left( \left( u\circ \upgamma_{\substack{\\[1pt]\hspace{-2pt}q}} \right)^{'}\left(r\right)\right)^{\hspace{-3pt}^2} + \left( u\circ \upgamma_{\substack{\\[1pt]\hspace{-2pt}q}} \left(r\right) \right)^{\hspace{-3pt}^2} \equiv \mathsf{const} = {\mathsf{max}_{\substack{\\q}}\strut^{\hspace{-5pt}2}},
	\end{align}
	where
	\begin{align*}
	\displaystyle \mx_q=\max_{t\in \left[ \ai_q , \bi_q \right]} u\circ \upgamma_{\substack{\\[1pt]\hspace{-2pt}q}} \left(t\right) = \shortminus \min_{t\in \left[ \ai_q, \bi_q \right]} u\circ \upgamma_{\substack{\\[1pt]\hspace{-2pt}q}} \left(t\right).
	\end{align*}
\end{lemma}
\begin{proof}[\footnotesize \textbf{Proof}]
	We have $u=0$ $\meas$-a.e. on $\mathcal Z$ with $\mathcal Z\cup \mathcal T=\X$ and 
\begin{align*}
	\int_{_{\X_{_q}}} u \, \dif\meas_{\hspace{1pt}q}=0,
\end{align*}	
for ${\bf q}$-a.e. $q\in \Q$.
	We identify $u$ with its Lipschitz representative. Then by definition of the local Lipschitz constant $\left|\nabla u\right|$ of $u$ it holds that 
\begin{align*}
		\left| \left( u\circ \upgamma_{\substack{\\[1pt]\hspace{-2pt}q}}  \right)^{\hspace{-2pt}'}\right|\left(r\right)\leq \left|\nabla u\right|\circ \upgamma_{\substack{\\[1pt]\hspace{-2pt}q}} \left(r\right) \quad \ \mathscr L^{^1}\mbox{-a.e. $r \in \left[ \ai_{\substack{\\q}},\bi_{\substack{\\q}} \right]$  and ${\bf q}$-a.e. $q\in \Q$}.
	\end{align*}
	For ${\bf q}$-a.e. $q\in \Q$, by Theorem \ref{th:1Dlocscheme} we have that $\left(  \X_{\hspace{1pt}q}, \dist_{\substack{\\[1pt]\X_{\hspace{1pt}q} }}, \meas_{\hspace{1pt}q} \right)$ satisfies $\CD\left(0,\N\right)$. 
	Revisiting the proof of the sharp spectral inequality in \cite{CM-2} yields
	\begin{align}\label{ineq:1d}
		\int_{\hspace{-1pt}_{\X}} u^2 \, \dif\meas& = \int_{\hspace{-1pt}_{\Q}}\int \left(u\circ \upgamma_{\substack{\\[1pt]\hspace{-2pt}q}}\right)^{^{\hspace{-2pt}2}}\left(r\right)\, \uprho_{\substack{\\[1pt]\hspace{-2pt}q}}\left(r\right)  \, \dif\mathscr{L}^{^1}\left(r\right) \, \dif\mathbf{q}\left(q\right)\nonumber\\
		&\leq\int_{\hspace{-1pt}_{\Q}} \int \left| \left( u\circ \upgamma_{\substack{\\[1pt]\hspace{-2pt}q}}  \right)^{\hspace{-2pt}'} \right|^{^2}\left(r\right) \, \uprho_{\substack{\\[1pt]\hspace{-2pt}q}}\left(r\right)  \, \dif\mathscr{L}^{^1}\left(r\right) \, \dif\mathbf{q}\left(q\right) \\  &\le \int_{\hspace{-1pt}_{\Q}} \int \left|\nabla u\right|^{^2}\circ  \upgamma_{\substack{\\[1pt]\hspace{-2pt}q}}\left(r\right) \,  \uprho_{\substack{\\[1pt]\hspace{-2pt}q}}\left(r\right)  \, \dif\mathscr{L}^{^1}\left(r\right) \, \dif\mathbf{q}\left(q\right) \nonumber \\ & = \int_{\hspace{-1pt}_{\X}} \left|\nabla u\right|^{^2} \, \dif\meas.\nonumber
	\end{align}
	Then the spectral equality $\int_{\hspace{-1pt}_{\X}} \left|\nabla u\right|^{^2} \, \dif\meas= \int_{\hspace{-1pt}_{\X}} u^2 \,  \dif\meas$ yields equalities in all the lines of \eqref{ineq:1d}. Hence $u\circ \upgamma_{\substack{\\[1pt]\hspace{-2pt}q}} $ saturates the spectral inequality for $\mathbf{q}$-almost every $q\in \Q$. Since $\left| \left( u\circ \upgamma_{\substack{\\[1pt]\hspace{-2pt}q}}  \right)^{\hspace{-2pt}'}  \right|\leq \left|\nabla u\right|\circ  \upgamma_{\substack{\\[1pt]\hspace{-2pt}q}} $, this implies 
	\begin{align}\label{id3}
	\left| \left( u\circ \upgamma_{\substack{\\[1pt]\hspace{-2pt}q}}  \right)^{\hspace{-2pt}'}  \right| =  \left|\nabla u\right|\circ  \upgamma_{\substack{\\[1pt]\hspace{-2pt}q}} \quad   \mathscr L^{^1}\mbox{-a.e. on } \left[ \ai_{\substack{\\q}},\bi_{\substack{\\q}} \right].
	\end{align}
	\par We can apply the same gradient estimate arguments as in the proof of Proposition~ \ref{th:gradestimate} (see \cite{Jiang-Zhang}) in order to obtain
	\begin{align*}
			\left| \left( u\circ \upgamma_{\substack{\\[1pt]\hspace{-2pt}q}}  \right)^{\hspace{-2pt}'}  \right|^{^2} + \left( u\circ \upgamma_q \right)^{^{\hspace{-1pt}2}} \leq  {\mathsf{max}_{\substack{\\q}}\strut^{\hspace{-5pt}2}} \quad  \mathscr L^{^1}\mbox{-a.e. on } \left[ \ai_{\substack{\\q}},\bi_{\substack{\\q}} \right].
	\end{align*}
	\par Consider $s_{_0}$ and $s_{_1}$ with  $u\circ \upgamma_{\substack{\\[1pt]\hspace{-2pt}q}}  \left( s_{_0} \right)= \shortminus \,\mx_{\substack{\\q}}$ and $u\circ \upgamma_{\substack{\\[1pt]\hspace{-2pt}q}}  \left( s_{_1} \right)=\max u\circ \upgamma_{\substack{\\[1pt]\hspace{-2pt}q}}  =\mx_{\substack{\\q}}$.  By replacing $u$ with $ \shortminus \, u$ if necessary, we can assume $s_{_0}\leq s_{_1}$.
	\par Up to ${\mathscr L}^{^1}$-null sets, one can define the \emph{extremal set}
\begin{align*}
\mathcal{A} := 
\left\{ s\in \left[ s_{_0}, s_{_1} \right] \;:\;  |u\circ \upgamma_{\substack{\\[1pt]\hspace{-2pt}q}} \left(s\right)| = \mx_{\substack{\\q}}   \right\}\subset 
 \left\{ s\in  \left[ s_{_0}, s_{_1} \right] \;:\; \left|\left( u\circ \upgamma_{\substack{\\[1pt]\hspace{-2pt}q}}  \right)^{\hspace{-2pt}'}\right| = 0 \right\},
\end{align*}
and the \emph{bad set}
\begin{align*}
	\mathcal{B}:=\left\{ s\in  \left[ s_{_0}, s_{_1} \right] \;:\; \left( u\circ \upgamma_{\substack{\\[1pt]\hspace{-2pt}q}}  \right)^{\hspace{-2pt}'} < 0 \right\}.
\end{align*}
Let the auxiliary function $\uppsi : \left[ s_{_0}, s_{_1} \right] \to \R$ be given by
\begin{align*}
\uppsi\left(s\right) := \begin{cases} \nicefrac{|\left( \nicefrac{u\circ\upgamma_{_q} \,}{\,\mx_q} \right)^{'}|}{ \left(1-\left(\nicefrac{u\circ\upgamma_{_q}}{\mx_q} \right)^{^{\hspace{-1pt}2}}\right)^{\nicefrac{1}{2}}} &  s \not\in \mathcal{A},\\ 
0 & s \in \mathcal{A}, \end{cases}
\end{align*}
which is a measurable function since $\mathcal{A}$ is measurable. Notice $\uppsi\left(s\right) \le 1$ for $\mathscr{L}^{^1}$-a.e.
\par Now, we can write	
	\begin{align}\label{EQN:cruciel-computation0}
		\uppi = \diam \geq \int_{_{\left[s_{_0}, s_{_1}  \right]}} \, \dif \mathscr{L}^{^1}\left(s\right)&\ge 
 \int_{_{\left[s_{_0}, s_{_1}  \right]\smallsetminus {\mathcal{A}}}}\uppsi\left(s\right)  \, \dif \mathscr{L}^{^1}\left(s\right)\\
&=
 \int_{_{\left[s_{_0}, s_{_1}  \right]}} \uppsi\left(s\right)  \, \dif \mathscr{L}^{^1}\left(s\right) \notag\\
&\ge \int_{_{[-1,1]}}
\left( 1 \shortminus t^{^2}  \,\right)^{\shortminus \nicefrac{1}{2}} \, \dif t = \uppi. \notag
	\end{align}
So, $\uppsi \in \LL^{\hspace{-3pt}^{1}}$. Notice, the passage from the second line to the third line in (\ref{EQN:cruciel-computation0}) follows by the 1D area formula. Indeed, consider the Lipschitz change of variables function
$\mathsf{t}(s) := \nicefrac{1}{\mx_{\substack{\\q}}} \ u\circ \upgamma_{\substack{\\[1pt]\hspace{-2pt}q}}(s)$;  we have
\begin{align*}
	\textsf{t} \left( \left[ s_{_0}, s_{_1} \right] \smallsetminus \mathcal{A} \right) = \left( \shortminus 1, 1\right),
\end{align*}
and by the 1D area formula (e.g. see \cite{Mag} for a very general such area formula), we get
\begin{align*}
	\int_{_{\left[s_{_0}, s_{_1}  \right]}} \uppsi\left(s\right)  \, \dif \mathscr{L}^{^1}\left(s\right) & = \int_{[-1,1]} \# \left( \mathsf{t}^{-1}(t) \right) \left( 1 \shortminus t^{^2}  \,\right)^{\shortminus \nicefrac{1}{2}} \, \dif t  \ge \int_{[-1,1]} \left( 1 \shortminus t^{^2}  \,\right)^{\shortminus \nicefrac{1}{2}} \, \dif t.
\end{align*}
\begin{remark}
\par Alternatively, for the last inequality in \eqref{EQN:cruciel-computation0}, we could have used a subtle version (notice we only have the Lipschitz regularity of $u$ at our disposal) of change of variables theorem in Lebesgue integration (see~\cite{Serrin2}) in combination with
\begin{align*}
\dif\mathscr{L}^{^1}\left(\mathsf{t}(s)\right) =  \nicefrac{1}{\mx_{\substack{\\q}}} \left| \left( u\circ \upgamma_{\substack{\\[1pt]\hspace{-2pt}q}} \right)^{\hspace{-2pt}'}\left(s\right) \right| \,\dif\mathscr{L}^{^1}\left(s\right).
\end{align*}
\end{remark}
Hence, all the inequalities in (\ref{EQN:cruciel-computation0}) are indeed equalities. In particular, it follows
\begin{align*}
		\left| \left( u\circ \upgamma_{\substack{\\[1pt]\hspace{-2pt}q}} \right)^{\hspace{-2pt}'} \right|^{^2} + \left( u\circ \upgamma_{\substack{\\[1pt]\hspace{-2pt}q}}  \right)^{^{\hspace{-2pt}2}}= {\mathsf{max}_{\substack{\\q}}\strut^{\hspace{-5pt}2}} \quad   \mathscr L^{^1}\mbox{-a.e. on $\left[s_{_0}, s_{_1} \right]$} \notag,
	\end{align*}
	which was claimed in~\eqref{id2}. 
In particular, this yields that $\left( u\circ \upgamma_{\substack{\\[1pt]\hspace{-2pt}q}} \right)^{\hspace{-2pt}'}$ is continuous up to an a.e. modification (and hence, up to an a.e. modification $u\circ \gamma_{\substack{\\[1pt]\hspace{-2pt}q}} $ is continuously differentiable on $[s_0,s_1]$). This in particular implies the extremal  set $\mathcal{A} $ is  $ \mathscr L^{^1}$-null. 
\par Moreover,
	\begin{align*}
	\uppi \geq \bi_{\substack{\\q}} \shortminus \ai_{\substack{\\q}} \geq  \left| s_{_1} \shortminus  s_{_0} \right|=\uppi \quad {\bf q}\mbox{-a.e.} \; q\in \Q,
	\end{align*}
hence, \eqref{eq:pi-length}, and $\ai_{\substack{\\q}}=s_0$ and $\bi_{\substack{\\q}}=s_1$.; recall, the above  argument holds true for $\mathbf{q}$-a.e. $q$. 
\par  Now, we set the auxiliary function (which aposteriori, the same as $\uppsi$)
\begin{align*}
\upphi\left(s\right) := \nicefrac{\left( \nicefrac{u\circ\upgamma_{_q} \,}{\,\mx_q} \right)^{'}}{ \left(1-\left(\nicefrac{u\circ\upgamma_{_q}}{\mx_q} \right)^{^{\hspace{-1pt}2}}\right)^{\nicefrac{1}{2}}}
\end{align*}
with the convention that $\upphi(s) = 0$ whenever $u\circ\upgamma_{_q}(s) =  \plmi \mx_q$.
 \par Thus, one again writes
\begin{align}\label{EQN:cruciel-computation}
		\uppi = \diam \geq \int_{_{\left[s_{_0}, s_{_1}  \right]}} \, \dif \mathscr{L}^{^1}\left(s\right)&\ge 
 \int_{_{\left[s_{_0}, s_{_1}  \right]\smallsetminus {\mathcal{B}}}}\upphi\left(s\right)  \, \dif \mathscr{L}^{^1}\left(s\right)\\
&\ge 
 \int_{_{\left[s_{_0}, s_{_1}  \right]\smallsetminus \mathcal{B}}} \upphi\left(s\right)  \, \dif \mathscr{L}^{^1}\left(s\right)+
 \int_{\mathcal{B}} \upphi\left(s\right)  \, \dif \mathscr{L}^{^1}\left(s\right)\notag\\
&= 
 \int_{_{\left[s_{_0}, s_{_1}  \right]}} \upphi\left(s\right)  \, \dif \mathscr{L}^{^1}\left(s\right)\notag\\
&=
 \int_{s_{_0}}^{s_{_1}} \upphi\left(s\right)  \, \dif \mathscr{L}^{^1}\left(s\right)\notag\\
&=\int_{\nicefrac{u\circ \upgamma_{\hspace{-1pt}_q}(s_{_0})}{\mx_{_q}}}^{\nicefrac{u\circ \upgamma_{\hspace{-1pt}_q}(s_{_1})}{\mx_{_q}}}  
\left( 1 \shortminus t^{^2}  \,\right)^{\shortminus \nicefrac{1}{2}} \, \dif t\notag\\
&=\int_{_{{-1}}}^{{_1}}
\left( 1 \shortminus t^{^2}  \,\right)^{\shortminus \nicefrac{1}{2}} \, \dif t = \uppi; \notag
\end{align}
in \eqref{EQN:cruciel-computation}, the third to fourth line is the passage from Lebesgue integral to Riemann's which is now justified by the $\mathcal{C}^{^1}$ regularity of $u\circ\upgamma_{\hspace{-2pt}_q}$;  the continuous differentiability of $u\circ \upgamma$ also justifies the fourth equality where we have used the  usual substitution rule for one dimensional Lebesgue (also Riemann in this case) integrals.
\par All the inequalities in the chain of inequalities \eqref{EQN:cruciel-computation} are equalities; thus, they result in
\begin{align*}
	\left( u\circ \upgamma_{\substack{\\[1pt]\hspace{-2pt}q}} \right)^{\hspace{-2pt}'}\geq 0 \quad   \mathscr L^{^1}\mbox{-a.e. on $\left[ \ai_{\substack{\\q}},\bi_{\substack{\\q}} \right]$};
	\end{align*}
notice, in particular, the bad set $\mathcal{B}$ is $\mathscr L^{^1}$-null. Now, upon combining with \eqref{id3}, one deduces \eqref{id1}. 
\end{proof}
\begin{corollary}\label{cor:cos}
	For ${\bf q}$-a.e. $q\in \Q$,
	\begin{align*}
	u\circ \upgamma_{\substack{\\[1pt]\hspace{-2pt}q}}  \left(r\right)= \shortminus \, \mx_q \, \cos\left( r \shortminus \ai_{\substack{\\q}}\right),	
	\end{align*}
	holds for all $ r \in \left[ \ai_{\substack{\\q}} ,\bi_{\substack{\\q}} \right]$. In particular, $\left( u\circ \upgamma_{\substack{\\[1pt]\hspace{-2pt}q}} \right)^{\hspace{-2pt}'} \left(r\right)= \mx_q \, \sin\left(r \shortminus \ai_{\substack{\\q}} \right)$ 
	for ${\bf q}$-a.e. $q\in \Q$. 
\end{corollary}
\begin{proof}[\footnotesize \textbf{Proof}]
	This follows immediately from \eqref{id2} and $\mx_q=\max_{\substack{\\\,t\in \left[ \ai_q , \bi_q \right]}} u \,\scalebox{0.8}{$\circ$}\, \upgamma_{\substack{\\[1pt]\hspace{-2pt}q}} \left(r\right)$. Since $\left( u\circ \upgamma_{\substack{\\[1pt]\hspace{-2pt}q}} \right)^{\hspace{-2pt}'}\geq 0$ in $\mathscr{L}^{^1}$-a.e. sense, the continuous function $\left[ \ai_{\substack{\\q}},\bi_{\substack{\\q}} \right] \ni r\mapsto u\circ \upgamma_{\substack{\\[1pt]\hspace{-2pt}q}}$ must exactly be given by $u\circ \upgamma_{\substack{\\[1pt]\hspace{-2pt}q}}= \shortminus \, \mx_q \, \cos\left( r \shortminus \ai_{\substack{\\q}} \right)$ thus it also admits a continuous derivative
	\begin{align*}
\left( u\circ \upgamma_{\substack{\\[1pt]\hspace{-2pt}q}} \right)^{\hspace{-2pt}'} = \mx_q \, \sin\left( r \shortminus \ai_{\substack{\\q}} \right).
	\end{align*}
\end{proof}
\begin{remark}\label{rem:cle}
	Based on Corollary \ref{cor:cos}, we will replace $\left|\nabla u\right|$ with a representative $\widetilde{\left|\nabla u\right|} \in \Ltwo\left( \X,\meas \right)$ which coincides $\meas$-a.e. with $\left|\nabla u\right|$ and in addition, $\left| \left( u\circ \upgamma_{\substack{\\[1pt]\hspace{-2pt}q}} \right)^{\hspace{-2pt}'} \right|=\widetilde{\left|\nabla u\right|}\circ\upgamma_{\substack{\\[1pt]\hspace{-2pt}q}} $ \emph{everywhere} in $\left[ \ai_{\substack{\\q}},\bi_{\substack{\\q}} \right]$ for ${\bf q}$-a.e. $q\in \Q$. 
\end{remark}
\begin{lemma}\label{lem:lemD}
	Let 
	\begin{align*}
	\left( \left[ \shortminus \nicefrac{\uppi}{2}, \nicefrac{\uppi}{2} \right],\dist_{_\mathsf{Euc}} ,\uprho \ \dif \mathscr{L}^{^1} \right),
	\end{align*}	
	be a one dimensional metric measure space that satisfies the $\RCD\left(0,\N\right)$ conditions. Suppose $\uprho$ is locally Lipschitz. Assume that the first nonzero eigenvalue $\feig=1$. Then $\uprho$ is constant. In particular, the disintegration appearing in the Lemma~\ref{lem:ida} now becomes
	\begin{align*}
		\meas\mres_{\T}= \int_{\hspace{-1pt}_{\Q}} {\mathscr H}^{^1}\mres_{\X_{_q}} \, \dif{\bf q}\left(q\right) = \int_{\hspace{-1pt}_{\Q}} \ray\left(q,\cdot\right)_{*}  {\mathscr L}^{^1}  \, \dif{\bf q}\left(q\right),
	\end{align*}
where $\ray$ is the ray map. 
\end{lemma}
\begin{proof}[\footnotesize \textbf{Proof}]
	Let $u$ be an eigenfunction for $\feig$. After rescaling, if necessary, we already know
	\begin{align*}
		\left( u^{'} \right)^{\hspace{-2pt}^2}+u^2 \equiv 1 \quad \text{and} \quad \shortminus u\left(\shortminus \nicefrac{\uppi}{2}\right)= u\left(\nicefrac{\uppi}{2}\right)=1.
	\end{align*}
Therefore, by Corollary~\ref{cor:cos}, we must have $u\left(x\right)=\sin x$. Then for any smooth function $v$ with 
	\begin{align*}
	\mathrm{supp}\,v \; \text{\scalebox{0.7}{$\subset$}} \; \left( \shortminus \nicefrac{\uppi}{2} , \nicefrac{\uppi}{2} \right),
	\end{align*}
	one can calculate
	\begin{align}
		&\int_{\hspace{-1pt}_{\shortminus \nicefrac{\uppi}{2}}}^{^{\nicefrac{\uppi}{2}}}u\left(x\right)v\left(x\right)\uprho\left(x\right)\,\dif x=\int_{\hspace{-1pt}_{\shortminus \nicefrac{\uppi}{2}}}^{^{\nicefrac{\uppi}{2}}}u^{'}\left(x\right)v^{'}\left(x\right)\uprho\left(x\right)\,\dif x\notag\\
		\Leftrightarrow &\int_{\hspace{-1pt}_{\shortminus \nicefrac{\uppi}{2}}}^{^{\nicefrac{\uppi}{2}}}v\left(x\right)\uprho\left(x\right)\sin x\,\dif x=\int_{\hspace{-1pt}_{\shortminus \nicefrac{\uppi}{2}}}^{^{\nicefrac{\uppi}{2}}}v^{'}\left(x\right)\uprho\left(x\right)\cos x\,\dif x.\notag
	\end{align}
	By hypothesis, $\uprho$ is a locally Lipschitz and hence $\uprho$ is differentiable almost everywhere and also included in the $\mathcal{W}^{^{\hspace{0.7pt}1, 1}}$-Sobolev space. Therefore by \emph{integration by parts formula}, we have 
	\begin{align}
		&\int_{\hspace{-1pt}_{\shortminus\nicefrac{\uppi}{2}}}^{^{\nicefrac{\uppi}{2}}} v\left(x\right)\uprho\left(x\right)\sin x\,\dif x \\ &=	\int_{\hspace{-1pt}_{\shortminus\nicefrac{\uppi}{2}}}^{^{\nicefrac{\uppi}{2}}} v^{'}\left(x\right)\uprho\left(x\right)\cos x\,\dif x\notag\\
		&=\Big[v\left(x\right)\uprho\left(x\right)\cos x\Big]^{^{\nicefrac{\uppi}{2}}}_{\substack{\\[2pt]\shortminus \nicefrac{\uppi}{2}}} \shortminus 	\int_{\hspace{-1pt}_{\shortminus\nicefrac{\uppi}{2}}}^{^{\nicefrac{\uppi}{2}}} v\left(x\right)\left(\uprho^{'}\left(x\right)\cos x \shortminus \uprho\left(x\right)\sin x\right)\,\dif x\notag\\
		&\Leftrightarrow 	\int_{\hspace{-1pt}_{\shortminus\nicefrac{\uppi}{2}}}^{^{\nicefrac{\uppi}{2}}} v\left(x\right)\uprho^{'}\left(x\right)\cos x\,\dif x=0.\notag
	\end{align} 
	By the fundamental lemma of calculus of variations, we obtain $\uprho^{'}\left(x\right)=0$. Thus $\uprho$ is a constant. 
\end{proof}
\begin{lemma}\label{lem:repa}
	For ${\bf q}$-a.e. $q\in \Q$, there exists a monotone bijection (reparametrization)
	\begin{align*}
	\upvarphi_{\substack{\\[1pt]\hspace{-2pt}q}}: \mathbb R \rightarrow \left(\ai_{\substack{\\q}}, \bi_{\substack{\\q}} \right),
	\end{align*}
	such that 
	\begin{align*}
		\upvarphi^{'}_{\substack{\\[1pt]\hspace{-2pt}q}}\left(r\right)=\widetilde{\left|\nabla u\right|} \circ \upgamma_{\substack{\\[1pt]\hspace{-2pt}q}} \circ 	\upvarphi_{\substack{\\[1pt]\hspace{-2pt}q}} \left(r\right) \quad  \mbox{$\forall r\in \mathbb R$.}
	\end{align*} 
	Setting $\bar \upgamma_{\substack{\\[1pt]\hspace{-2pt}q}}=\upgamma_{\substack{\\[1pt]\hspace{-2pt}q}} \circ 	\upvarphi_{\substack{\\[1pt]\hspace{-2pt}q}}$, in particular it holds 
	\begin{align*}
	\nicefrac{\dif}{\dif r} \, u\circ \bar \upgamma_{\substack{\\[1pt]\hspace{-2pt}q}} \left(r\right) = \widetilde{ \left|\nabla u\right|}^{^2}\circ \bar \upgamma_{\substack{\\[1pt]\hspace{-2pt}q}}\left(r\right)
	\quad  \forall r\in \mathbb R,
\end{align*}
	and
	\begin{align*}
	\left| \bdot{ \ovs{ \upgamma}_{\substack{\\[1pt]\hspace{-2pt}q}}} \right|\left(t\right) = 	\upvarphi^{'}_{\substack{\\[1pt]\hspace{-2pt}q}}\left(t\right) \quad  \forall r\in \mathbb R,
	\end{align*}
	for ${\bf q}$-a.e. $q\in \Q$.
\end{lemma}
\begin{proof}[\footnotesize \textbf{Proof}] 
	Recall that for ${\bf q}\text{-a.e.}\, q $, one has $\bi_{\substack{\\q}}-\ai_{\substack{\\q}}=\uppi$ and $0\in \left( \ai_{\substack{\\q}},\bi_{\substack{\\q}} \right)$ by the parametrization that we chose for $\upgamma_{\substack{\\[1pt]\hspace{-2pt}q}}$ according to 1D localization. Define 
	\begin{align*}
	 \left( \ai_{\substack{\\q}},\bi_{\substack{\\q}} \right) \ni	r \mapsto \upbeta_{\substack{\\\hspace{-1pt}q}}\left(r\right) :=\int_{\hspace{-1pt}_0}^r \nicefrac{1}{\;
		\text{\raisebox{0pt}[1.2\height]{$\widetilde{\left|\nabla u\right|}\circ \gamma_{_q}\left(\uptau\right)$}}} \, \dif\uptau=\int_{\hspace{-1pt}_0}^r \nicefrac{1}{\;	\text{\raisebox{0pt}[1.2\height]{$	\mx_q \	\sin\left( \uptau- \ai_{\substack{\\q}} \right)$}}} \, \dif\uptau .
	\end{align*}	
	Since $\widetilde{\left|\nabla u\right|}\circ\upgamma_{\substack{\\[1pt]\hspace{-2pt}q}}$ is positive in the interior of $\upgamma_{\substack{\\[1pt]\hspace{-2pt}q}}$,
	$
	 \upbeta_{\substack{\\\hspace{-1pt}q}}: \left( \ai_{\substack{\\q}},\bi_{\substack{\\q}} \right) \rightarrow \mathbb{R}
$
	is well-defined, continuous, monotonically increasing and since the $\nicefrac{1}{\sin x}$ has a simple pole at $x=0$,
		\begin{align*}
		\upbeta_{\substack{\\\hspace{-1pt}q}}\left(r\right) \rightarrow + \infty \; \left( \shortminus \infty \; \text{resp.} \right) \quad \text{as} \quad r\rightarrow \bi_{\substack{\\q}} \; \left( \ai_{\substack{\\q}} \; \text{resp.} \right).
		\end{align*}
	Note that $	\upbeta_{\substack{\\\hspace{-1pt}q}}\left(r\right) \leq 0$ for $r\leq 0$ since we are implicitly adopting the convention $\int_{_0}^{r} \ray \left(r\right)\, \dif r = \shortminus \int_{r}^{^0} \ray\left(r\right)\, \dif r$.
	\par 	The inverse function
	\begin{align*}
	\upbeta_{\substack{\\\hspace{-1pt}q}}^{\hspace{-2pt}^{-1}} \; =: \; \upvarphi_{\substack{\\[1pt]\hspace{-2pt}q}}: \mathbb R \rightarrow  \left( \ai_{\substack{\\q}},\bi_{\substack{\\q}} \right),
	\end{align*}
	is thus well-defined and we compute
	\begin{align*}
		1= \nicefrac{d}{dr}\ \int_{\hspace{-1pt}_0}^{\upvarphi_{_{\hspace{-1pt}q}}\left(r\right)} \nicefrac{1}{\;
			\text{\raisebox{0pt}[1.2\height]{$\widetilde{\left|\nabla u\right|}\circ \gamma_{_q}\left(\tau\right)$}}} \, \dif\tau= \nicefrac{\upvarphi_{_{\hspace{-1pt}q}}'\left(r\right)}{\;
			\text{\raisebox{0pt}[1.2\height]{$\widetilde{\left|\nabla u\right|}\circ \upgamma_{_q} \circ\upvarphi_{_{\hspace{-1pt}q}}\left(r\right)$}}},
	\end{align*}
	which gives the desired conclusion.
\end{proof}
\subsection{The gradient flow of $u$}\label{subsec:gradu} 
Recall from \textsection\thinspace\ref{subsec:1Dlocschemes} that we have a ray map 
\begin{align*}
	\ray : { \mathcal V}\subset \mathbb R\times \Q \rightarrow \mathcal T, \quad \ray\left(t,q\right)=\upgamma_{\substack{\\[1pt]\hspace{-2pt}q}} \left(t\right),
\end{align*}
where $ \left( \ai_{\substack{\\q}},\bi_{\substack{\\q}} \right)\ni t  \mapsto \upgamma_{\substack{\\[1pt]\hspace{-2pt}q}} \left(t\right)$ is a transport geodesic of full range. 
\par By Lemma~\ref{lem:repa}, for ${\bf q}$-a.e. $q\in \Q$ the reparametrization $\upvarphi_{\substack{\\[1pt]\hspace{-2pt}q}}: \mathbb R\rightarrow \left( \ai_{\substack{\\q}},\bi_{\substack{\\q}} \right)$ exists and we set $\bar  \upgamma_{\substack{\\[1pt]\hspace{-2pt}q}}=  \upgamma_{\substack{\\[1pt]\hspace{-2pt}q}} \circ \upvarphi_{\substack{\\[1pt]\hspace{-2pt}q}}$ for those $q \in \Q$. 
\par The construction of $\upvarphi_{\substack{\\[1pt]\hspace{-2pt}q}}$ also gives that $\upvarphi:\left(t,q\right)\in \mathbb R\times \Q \mapsto \upvarphi_{\substack{\\[1pt]\hspace{-2pt}q}}\left(t\right)\in \mathcal V$ is a measurable bijection. We can introduce a reparametrization of the ray map $\ray$ via
\begin{align*}
	\bar \ray:=\ray \circ \upvarphi: \mathbb R \times \Q\rightarrow \mathcal T.
\end{align*}
 Once $\mathbb R\times \Q$ is equipped with the measure given by the density $\upvarphi_{\substack{\\[1pt]\hspace{-2pt}q}}^{'}(t) dt \otimes d{\bf q}(q)$, the map $\bar \ray$ again becomes a measure space isomorphism (measure-preserving bijection), and sets of the form $\bar \ray\left( \left(a,b\right)\times S \right)$ for Borel sets $S\subset \Q$,  generate the class of Borel sets in $\mathcal T$. In simple terms, we are changing the parameter from $r \in \left( \ai_{\substack{\\q}},\bi_{\substack{\\q}} \right)$ to $t \in \R$ in a measure-preserving way.  
\par We also define the shift operator $\ell_{t}: \mathbb R\times \Q\rightarrow \mathbb R\times \Q$ via $\left(\cdot,q\right)\mapsto \left(\cdot + t, q\right)$. Then, we can define a flow map $\F$ as follows
\begin{align*}
	\F:\mathbb R\times \mathcal T \rightarrow \X, \quad \F\left(t, x\right)= \bar \ray \circ \ell_t \circ \bar \ray^{^{-1}}\left(x\right),
\end{align*}
which is a measurable map. One thus gets $\F_{\substack{\\[1pt]\hspace{-3pt}t}}=\F\left(t,\cdot\right): \mathcal T\rightarrow X$. Since $\ell_{t_{_1}}  \circ \ell_{t_{_2}} = \ell_{t_{_1} + t_{_2}}$, $\F$ is indeed a flow (semi-group in $t$). 
\begin{remark}\label{rem:grad}
	Note that 
	\begin{align*}
	\F\left(t,x\right)= \bar\upgamma_{\substack{\\[1pt]\hspace{-2pt}\mathfrak{s}\left(x\right)}} \left( \bar\upgamma^{^{-1}}_{\substack{\\[1pt]\hspace{-2pt}\mathfrak{s}\left(x\right)}}\left(x\right)+ t \right).
	\end{align*}
	Hence by Lemma \ref{lem:repa} we have
	\begin{align*}
	\nicefrac{\dif}{\dif t} \ u\circ \F\left(t,x\right)= \widetilde{ \left|\nabla u\right|}^{^2} \circ \F\left(t,x\right),
	\end{align*} 
	and the disintegration of $\meas$ transforms to 
	\begin{align*}
		\meas\mres_{\mathcal T}= \int_{\hspace{-1pt}_{\Q}} {\mathscr H}^{^1}\mres_{\X_{\hspace{1pt}q}}  \, \dif {\bf q}\left(q\right)= \int_{\hspace{-1pt}_{\Q}} \ray\left(\cdot,q\right)_{*}\left(\dif r \right) \, \dif{\bf q}\left(q\right)= \int_{\hspace{-1pt}_{\Q}} \bar \ray\left(\cdot,q\right)_{*}\left(\upvarphi_q'\left(t\right) \, \dif t\right)\, \dif{\bf q}\left(q\right).
	\end{align*}
	We set $\bar \uprho_{\substack{\\[1pt]\hspace{-2pt}q}}\left(t\right)=\upvarphi_{\substack{\\[1pt]\hspace{-2pt}q}}'\left(t\right)$.
\end{remark}
\par Pick a measurable subset $\EE\subset \T$ with $\meas\left(\EE\right)>0$ and
define a map
\begin{align*}
\Uplambda: \X\rightarrow \mathcal{C}\left(\left[0,1\right], \X\right), \quad \Uplambda\left(x\right) = \upsigma_{\substack{\\x}},
\end{align*}
where
\begin{align*}
	\left[0,1\right] \ni t  \mapsto   \upsigma_{\substack{\\x}} \left(t\right):= \F_{\substack{\\[1pt]\hspace{-3pt}t}}\left(x\right),
	\end{align*}
and a probability measure 
\begin{align*} 
	\ppi_{\substack{\\[1pt]\EE}}= \nicefrac{1}{\meas\left(\EE\right)} \, \Uplambda_{\,*} \meas  \mres_{\EE} \in \mathcal \Prob\left( \mathcal{C}\left(\left[0,1\right],\X\right) \right).
\end{align*}
\begin{proposition}\label{prop:gradu}
Let $\EE\subset \T$ and $\ppi_{\substack{\\[1pt]\EE}}$ be as in above. Assume $\left|\widetilde {\nabla u}\right|\restr_{\EE} \geq \eta>0$ in $\meas$-a.e. sense. Then, $\ppi_{\substack{\\[1pt]\EE}}$ is a test plan that represents the gradient of $u$.
\end{proposition}
\begin{proof}[\footnotesize \textbf{Proof}]
Let $t\in \left[0,1\right]$. For a subset $A \subset \X$, we write 
\begin{align*}
	\meas\left(\EE\right)\left( \e_{\substack{\\t}} \right)_{*}{\ppi}_{\substack{\\[1pt]\EE}}\left(A\right) &=  \meas \left( \EE\cap \F_{\substack{\\[1pt]\hspace{-3pt}-t}}\left(A\right) \right) \\ 
	&=\int_{\hspace{-1pt}_{\Q}} \meas_{\hspace{1pt}q} \left( \EE\cap \F_{\substack{\\[1pt]\hspace{-3pt}-t}}\left(A\right) \right) \, \dif{\bf q}\left(q\right) \\ 
	&= \int_{\hspace{-1pt}_{\Q}} \left( \F_{\substack{\\[1pt]\hspace{-3pt}t}} \right)_{\hspace{-1pt}*} \meas_{\hspace{1pt}q} \mres_{\EE}\left(A\right)\, \dif{\bf q}\left(q\right).
\end{align*}
Let us recall
\begin{align*}
 \meas_{\hspace{1pt}q} \mres_{\EE}=\left( \bar \upgamma_{\substack{\\[1pt]\hspace{-2pt}q}} \right)_{\hspace{-1pt}*} \left(  \bar \uprho_{\substack{\\[1pt]\hspace{-2pt}q}} \  \mathscr L^{^1}\mres_{ \bar \upgamma^{^{-1}}_{_q}\left(\EE\right)} \right).
\end{align*}
It follows
	\begin{align*}
	\left( \F_{\substack{\\[1pt]\hspace{-3pt}t}} \right)_{\hspace{-1pt}*} \meas_{\hspace{1pt}q} \mres_{\EE}&=\left(\bar \upgamma_{\substack{\\[1pt]\hspace{-2pt}q}}  \circ \ell_t \right)_{\hspace{-1pt}*}\left(\bar \uprho_{\substack{\\[1pt]\hspace{-2pt}q}} \, \mathscr L^{^1}\mres_{{\bar{\upgamma}^{^{-1}}_{_q}\left(\EE\right)}}\right) \\  & =  \left(\bar \upgamma_{\substack{\\[1pt]\hspace{-2pt}q}} \right)_{\hspace{-1pt}*} \left(\ell_t\right)_{*} \left( \bar \uprho_{\substack{\\[1pt]\hspace{-2pt}q}} \, \mathscr L^{^1}\mres_{{\bar{\upgamma}^{^{-1}}_{_q} \left(\EE\right)}}\right) \\ &= \left(\bar \upgamma_{\substack{\\[1pt]\hspace{-2pt}q}}\right)_{\hspace{-1pt}*}\left( \bar \uprho_{\substack{\\[1pt]\hspace{-2pt}q}}  \circ \ell_{-t} \; \mathscr L^{^1}\mres_{{\left( \bar{\upgamma}_q^{^{-1}}\left(\EE\right)-t\right) }}\right)\\& = \left(\bar \upgamma_{\substack{\\[1pt]\hspace{-2pt}q}}  \right)_{\hspace{-1pt}*}\left( \bar \uprho_{\substack{\\[1pt]\hspace{-2pt}q}}  \left(\cdot \shortminus t\right) \; \mathscr L^{^1}\mres_{{ \left( \bar{\upgamma}_{_q}^{^{-1}}\left(\EE\right)-t\right)}}\right).
\end{align*}
Moreover, 
 it holds
\begin{align*}
1_{\left( \bar{\upgamma}_{_q}^{^{-1}}\left(\EE\right)-t\right)}(t^{'}) \nicefrac{\bar \uprho_{_q} \left(t^{'} - t\right)}{\bar \uprho_{_q}\left(t^{'} \right)}= 
{1_{\left( \bar{\upgamma}_{_q}^{^{-1}}\left(\EE\right)-t \right)}(t^{'})}\nicefrac{\left|\widetilde{\nabla u}\right|\circ \bar{\upgamma}_{_q}\left( t^{'}-t \right)}{\left|\widetilde{\nabla u}\right|\circ \bar{\upgamma}_{_q}\left( t^{'} \right)} \leq \nicefrac{1}{\upeta},
\end{align*}
for $\mathscr L^{^1}\mbox{-a.e.} \ t^{'}$, where $1_{\hspace{-2pt}_{A}}$ is the indicator function of $A$. 
\par Therefore, again by the translation invariance of $\mathscr L^{^1}$, we have
\begin{align*}
	\left(\F_{\substack{\\[1pt]\hspace{-3pt}t}}\right)_{\hspace{-1pt}*} \meas_{\hspace{1pt}q}\mres_{\EE}= \left(\bar \upgamma_{\substack{\\[1pt]\hspace{-2pt}q}} \right)_{\hspace{-1pt}*}\left( \bar \uprho_{\substack{\\[1pt]\hspace{-2pt}q}}  \circ \ell_{-t}  \; \mathscr L^{^1}\mres_{{\left( \bar{\upgamma}_{_q}^{^{-1}}\left(\EE\right)-t \right)}}\right) \leq \nicefrac{1}{\upeta} \, \left( \bar \upgamma_{\substack{\\[1pt]\hspace{-2pt}q}}  \right)_{\hspace{-1pt}*}\left(\bar \uprho_{\substack{\\[1pt]\hspace{-2pt}q}}  \, \mathscr L^{^1}\mres_{\left( \bar{\upgamma}_{_q}^{^{-1}}\left(\EE\right)-t \right) }\right),
\end{align*} 
and consequently 
\begin{align*}
	\left( \e_{\substack{\\t}}\circ \Uplambda \right)_{*}\meas\mres_{\EE} = \left(\F_{\substack{\\[1pt]\hspace{-3pt}t}}\right)_{\hspace{-1pt}*}\meas\mres_{\EE} &= \int_{\hspace{-1pt}_{\Q}} \left(\bar \upgamma_{\substack{\\[1pt]\hspace{-2pt}q}} \right)_{\hspace{-1pt}*} \left(\bar \uprho_{\substack{\\[1pt]\hspace{-2pt}q}} \circ \ell_{-t}   \, \mathscr L^{^1}\mres_{{ \left( \bar{\upgamma}_{_q}^{^{-1}}\left(\EE\right)-t \right)}}\right) \ \dif{\bf q}\left(q\right)  \\ & \leq \nicefrac{1}{\eta} \int_{\hspace{-1pt}_{\Q}} \left(\bar \upgamma_{\substack{\\[1pt]\hspace{-2pt}q}} \right)_{\hspace{-1pt}*} \left(\bar \uprho_{\substack{\\[1pt]\hspace{-2pt}q}} \mathscr L^{^1}
\right) \, \dif{\bf q}\left(q\right) \\& = \nicefrac{1}{\upeta} \, \meas.
\end{align*}
Then, it follows
\begin{align*}
\left( \e_{\substack{\\t}} \right)_{*}\ppi_{\substack{\\[1pt]\EE}} \leq \left( \upeta\meas\left(\EE\right) \right)^{^{-1}}\; \meas \quad \forall t\in [0,1].
\end{align*}	
Hence,  $\ppi_{\substack{\\[1pt]\EE}}$ has bounded compression. Moreover by Lemma \ref{lem:repa},
\begin{align*}
\left| \bdot{\ovs{\upgamma}}_{\substack{\\[1pt]\hspace{-2pt}q}}\right|^{^{\hspace{-1pt}2}} = \left|\widetilde{\nabla u}\right|^{^2}\circ \bar \upgamma_{\substack{\\[1pt]\hspace{-2pt}q}} \leq 1 \ \Rightarrow \ \int_{\hspace{-1pt}_{\mathsf{Curves}}}  \left|\bdot{\upxi}\right|^2 \, \dif{\ppi}_{\substack{\\[1pt]\EE}}\left(\upxi\right) \leq 1,
\end{align*}
therefore $\ppi_{\substack{\\[1pt]\EE}}$ is a test plan.
Finally we have to show inequality \eqref{ineq:grad} for $\ppi_{\substack{\\[1pt]\EE}}$ and $u$. First, by Fatou's lemma and by Remark \ref{rem:grad} we observe that
\begin{align*}
	&\lowlim_{\tau \downarrow 0} \int_{\hspace{-1pt}_{\mathsf{Curves}}} \ \nicefrac{1}{\tau} \left( u\circ \upxi \left(\tau\right) \shortminus u\circ \upxi \left(0\right) \right) \, \dif \ppi_{\substack{\\[1pt]\EE}}\left(\upxi\right)
	 \\ & =
	\lowlim_{\tau \downarrow 0} \nicefrac{1}{\meas\left(\EE\right)} \, \int_{\hspace{-1pt}_{\EE}} \nicefrac{1}{\tau} \, \left( u\left(\F_{\substack{\\[1pt]\hspace{-5pt}\tau}}\left(x\right)\right) \shortminus u\left(x\right) \right)\, \dif\meas\left(x\right)\\
	&\geq \nicefrac{1}{\meas\left(\EE\right)} \int_{\hspace{-1pt}_{\EE}} \nicefrac{\dif}{\dif t}\restr_{t=0} \left(u\circ \F\right)\left(t,x\right) \, \dif \meas\left(x\right)\\
	&=\nicefrac{1}{\meas\left(\EE\right)}\, \int_{\hspace{-1pt}_{\EE}}\widetilde{\left|\nabla u\right|}^{^2} \, \dif\meas \\
	&=\nicefrac{1}{\meas\left(\EE\right)}\, \int_{\hspace{-1pt}_{\EE}}{\left|\nabla u\right|}^{^2}\, \dif\meas. 
\end{align*}
Before we proceed, we note that  for $\uptau>0$ it holds (by Fubini's theorem)
\begin{align*}
	\int_{\hspace{-1pt}_{\EE}}\int_{\hspace{-1pt}_0}^\uptau\widetilde{ \left|\nabla u\right|}^{^2} \left( \F_{\substack{\\[1pt]\hspace{-3pt}t}} \left(x\right) \right)  \, \dif t \, \dif\meas(x) &=  
	\int_{\hspace{-1pt}_0}^{\uptau}\int_{\hspace{-1pt}_{\EE}} \widetilde{\left|\nabla u\right|}^{^2}\left( \e_{\substack{\\t}}\circ \Uplambda\left(x\right) \right) \, \dif\meas(x) \, \dif t\\
	&=  \meas\left(\EE\right)\int_{\hspace{-1pt}_0}^\uptau\int_{\hspace{-1pt}_{\mathsf{Curves}}}  \widetilde{\left|\nabla u\right|}^{^2}\left(\upxi\left(t\right) \right) \, \dif \ppi_{\substack{\\[1pt]\EE}}\left(\upxi\right) \, \dif t \\ &= \meas\left(\EE\right)  \int_{\hspace{-1pt}_{\mathsf{Curves}}} \int_{\hspace{-1pt}_0}^\uptau \left|\bdot \upxi\left(t\right)\right|^{^2} \,  \dif t \, \dif \ppi_{\substack{\\[1pt]\EE}}\left(\upxi\right).
\end{align*}
where the last equality holds by the definition of of the measure $\ppi_{\substack{\\[1pt]\EE}}\left(\upxi\right)$ and by using
$
\left| \bdot{\ovs{\upgamma}}_{\substack{\\[1pt]\hspace{-2pt}q}}\right|^{^{\hspace{-1pt}2}} = \left|\widetilde{\nabla u}\right|^{^2}\circ \bar \upgamma_{\substack{\\[1pt]\hspace{-2pt}q}} 
$
. Moreover,
\begin{align*}
\mathsf{g}_{\substack{\\[2pt]\hspace{0.5pt}\uptau}}\left(x\right):=\nicefrac{1}{\uptau}\, \int_{\hspace{-1pt}_0}^\uptau \widetilde{\left|\nabla u\right|}^{^2}\left(\F_{\substack{\\[1pt]\hspace{-3pt}t}} \left(x\right)\right)\, \dif t\leq 1 \quad \mbox{$\meas$-a.e.} \;x \quad \forall \uptau>0. 
\end{align*}
Hence by Lebesgue differentiation theorem, $\lim\limits^{}_{\uptau \downarrow 0} \mathsf{g}_{\substack{\\[2pt]\hspace{0.5pt}\uptau}}\left(x\right)$ exists for $\meas$-a.e. $x$ and again by Fatou's lemma,
\begin{align*}\int_{\hspace{-1pt}_{\EE}} \widetilde{\left|\nabla u\right|}^{^2}\, \dif\meas &= \int_{\hspace{-1pt}_{\EE}}	\uplim_{\uptau \downarrow 0}\nicefrac{1}{\uptau} \int_{\hspace{-1pt}_0}^\uptau \widetilde{\left|\nabla u\right|}^{^2}\left( \F_{\substack{\\[1pt]\hspace{-3pt}t}} \left(x\right)\right) \, \dif t \, \dif \meas
	\\ &\ge  \uplim_{\uptau \downarrow 0}\nicefrac{1}{\uptau} 
	\int_{\hspace{-1pt}_{\EE}} \int_{\hspace{-1pt}_0}^\uptau \widetilde{\left|\nabla u\right|}^{^2}\left( \F_{\substack{\\[1pt]\hspace{-3pt}t}}  \left(x\right)\right) \, \dif t \, \dif\meas \\ & = \meas\left(\EE\right) \uplim_{\uptau \downarrow 0}\nicefrac{1}{\uptau}\int_{\hspace{-1pt}_{\mathsf{Curves}}} \int_{\hspace{-1pt}_0}^\uptau \left|\bdot \upxi\left(t\right)\right|^{^2} \, \dif t \, \dif\ppi_{\substack{\\[1pt]\EE}}\left(\upxi\right).
\end{align*}
Therefore, it follows
\begin{align*}
	&	\lowlim_{\tau \downarrow 0} \int_{\hspace{-1pt}_{\mathsf{Curves}}}  \nicefrac{1}{\uptau}\, \left( u\circ \upxi\left(\tau\right)  \shortminus u\circ \upxi \left(0\right) \right) \, \dif\ppi_{\substack{\\[1pt]\EE}} \left(\upxi\right) \ge \\  & \phantom{csycsy} \nicefrac{1}{2\meas\left(\EE\right)} \int_{\hspace{-1pt}_{\X}}  \left|\nabla u\right|^{^2}\left(x\right) \, \dif\meas\left(x\right)+\nicefrac{1}{2}\, \uplim_{\uptau\downarrow 0}\nicefrac{1}{\uptau}\int_{\hspace{-1pt}_{\mathsf{Curves}}} \int_0^\tau \left|\bdot \upxi\left(t\right)\right|^{^2} \, \dif t \, \dif \ppi_{\substack{\\[1pt]\EE}}\left(\upxi\right).
\end{align*}
Hence, $\ppi_{\substack{\\[1pt]\EE}}$ is indeed a test plan that represents the gradient of $u$.
\end{proof}
\subsection{Proof of Theorem~\ref{thm:hessian}}
We start with the following key lemma.
\begin{lemma}\label{LEM:key1}
	Let $\upalpha: \X\rightarrow \left[0,\infty\right)$ be given by
	\begin{align*}
	\upalpha := \left|\nabla u\right|^{^2}+ u^2.
	\end{align*}	
	Then, $\upalpha\in \Sobol\left(\X\right)$ and 
	\begin{align*}
		\langle \nabla \upalpha, \nabla u\rangle =0 \quad \text{$\meas$-a.e. in }\quad \left\{\widetilde{\left|\nabla u\right|}>0\right\} \cap \T. 
	\end{align*}
\end{lemma}
\begin{proof}[\footnotesize \textbf{Proof}]
Recall that $u\in \mathbb{D}_{\hspace{-1pt}_{\infty}}$ and therefore $\left|\nabla u\right|^{^2}\in \Sobol\left(\X\right)$ and $\upalpha = \left|\nabla u\right|^{^2} + u^2 \in \Sobol\left(\X\right)$. 
\par	We define $\mathsf{R}_{_n} :=\left\{\widetilde{\left|\nabla u\right|} > \nicefrac{1}{n}\right\}$. For any measurable set $\EE \subset \mathsf{R}_{_n}\cap \T$, according to Proposition~\ref{prop:gradu}, we define $\ppi_{\substack{\\[1pt]\EE}}$ which represents the gradient of $u$. Hence, we can apply Theorem~\ref{th:firstvariation} -- the first variation formula -- with $f=\upalpha$ to get
	\begin{align}\label{eq:alpha-1}
		\int_{\hspace{-1pt}_{\EE}}  \langle \nabla \upalpha, \nabla u\rangle  \, \dif\meas &= \meas\left(\EE\right) \int \langle \nabla \upalpha,\nabla u\rangle \, \dif \left(\e_{\substack{\\[1pt]\hspace{0.5pt}0}}\right)_{\hspace{-1pt}*}\ppi_{\substack{\\[1pt]\EE}} \\ &= \meas\left(\EE\right)\lim_{t\downarrow 0} \nicefrac{1}{t}\int_{\hspace{-1pt}_{\mathsf{Curves}}}  \Big(\upalpha\left(\upxi\left(t\right)\right) \shortminus \upalpha\left(\upxi\left(0\right)\right)\Big) \, \dif\ppi_{\substack{\\[1pt]\EE}}\left(\upxi\right). \notag
	\end{align}
	The right-hand side of the identity \eqref{eq:alpha-1} can be rewritten as
	\begin{align*}
		\lim_{t\downarrow 0} \int_{\hspace{-1pt}_{\EE}}  \nicefrac{1}{t}\Big(\upalpha\circ \Uplambda\left(x\right)\left(t\right) \shortminus \upalpha\circ \Uplambda\left(x\right)\left(0\right)\Big) \, \dif\meas\mres_{\EE}\left(x\right).
	\end{align*}
We apply the disintegration formula, $\meas\mres_{\T} =\int \mathscr H^{^1}\mres_{\X_{_q}} \, \dif{\bf q}\left(q\right)$ (see Theorem \ref{th:1Dlocscheme}, Corollary~\ref{cor:cos} and the conclusion of Lemma \ref{lem:lemD}) to obtain
	\begin{align*}
&	\lim_{t\downarrow 0} \nicefrac{1}{t}\int_{\hspace{-1pt}_{\Q}}  \int_{\uptau \,\in\, \ray\left(\cdot,q\right)^{-1}\left(\EE\right)} \Big(\upalpha\circ \Uplambda\left( \ray\left(\uptau,q\right) \right)\left(t\right) \shortminus \upalpha\circ \Uplambda\left( \ray\left(\uptau,q\right) \right)\left(0\right)\Big) \bar \uprho_q\left(\uptau\right)\, \dif\mathscr L^{^1}\left(\uptau\right) \, \dif{\bf q}\left(q\right)\\ &= 0,
	\end{align*}
where the equality ($=0$) holds since
	\begin{align}\label{EQ:constancy}
		\upalpha\circ \upgamma_q=\upalpha\circ \ray\left(\cdot, q\right)\equiv \mx_q \quad {\bf q}\mbox{-a.e.} \; q\in \Q,
	\end{align}
	and as a result for each fixed $q$,
	\begin{align*}
	\int_{\uptau \,\in\, \ray\left(\cdot,q\right)^{-1}\left(\EE\right)} \Big(\upalpha\circ \Uplambda\left( \ray\left(\uptau,q\right) \right)\left(t\right)\shortminus \upalpha\circ \Uplambda\left( \ray\left(\uptau,q\right) \right)\left(0\right)\Big) \bar \uprho_q\left(\uptau\right)\,  \dif \mathscr L^{^1}\left(\uptau\right) = 0.
	\end{align*}
	Since $\EE$ is arbitrary in $\mathsf{R}_{_n}\cap\T$, it follows $\langle \nabla \upalpha,\nabla u\rangle =0$ $\meas$-a.e. on $\mathsf{R}_{_n}\cap \T$. Since $n\in \mathbb N$ was arbitrary, it follows $\langle\nabla \upalpha,\nabla u\rangle =0$ $\meas$-a.e. on $\left\{\widetilde{\left|\nabla u\right|} > 0\right\}\cap \T$.
\end{proof}
\begin{corollary}\label{cor:alphabound}
	Set $\upalpha := \left|\nabla u\right|^{^2}+ u^2$ and $\upbeta := \left|\nabla u\right|^{^2} \shortminus u^2$. Then
	\begin{align*}
		{\langle \nabla \upalpha, \nabla \upbeta\rangle}= 
		{\langle  \nabla \upalpha, \nabla \left(\upalpha \shortminus 2u^2\right) \rangle}\geq 0  \quad \meas\mbox{-a.e. in }\quad \{\widetilde{\left|\nabla u\right|}>0\} \cap \T.
	\end{align*}
	In other terms,
	\begin{align*}
	\left|   \nabla \left|\nabla u\right|^{^2}  \right|^{^2} \ge 4 u^2 \left|\nabla u\right|^{^2}  \quad  \ \meas\mbox{-a.e. in }\quad \{\widetilde{\left|\nabla u\right|}>0\} \cap \T .
	\end{align*}
\end{corollary}
\begin{proof}[\footnotesize \textbf{Proof}]
	Applying the calculus rules for the inner product $\langle \nabla \cdot, \nabla \cdot\rangle$, we compute
	\begin{align*}
		{\langle \nabla \upalpha, \nabla \upbeta\rangle} &= 
		{\langle  \nabla \upalpha, \nabla \left(\upalpha \shortminus 2u^2\right) \rangle}= 
		\left|\nabla \upalpha\right|^{^2} \shortminus
		\, {4u\langle \nabla u, \nabla \upalpha \rangle} \\ &\geq  \shortminus {4u\langle \nabla u, \nabla \upalpha \rangle}=0,
	\end{align*}
$\meas \mbox{-a.e. in }\left\{\widetilde{\left|\nabla u\right|}>0\right\}  \cap \T$, where the last equality is obtained in Lemma~\ref{LEM:key1}.
\end{proof}
\begin{proposition}\label{PROP:key}
	$\upalpha\in \Lip\left(\X\right)$.
\end{proposition}
\begin{proof}[\footnotesize \textbf{Proof}]
	The proof is done in two steps.
	\par {\bf 1$\centerdot$}	Recall that in the 1D localization scheme, up to a set of measure $0$, it holds $\X= \mathcal T  \, \dot{\sqcup}  \, \mathcal Z$ where 
	\begin{align*}
	\meas\left( \mathcal Z \smallsetminus \left\{x\in \X \; \text{\textbrokenbar} \; u\left(x\right)=0\right\} \right)=0,
	\end{align*}
	and $\mathcal T$ is the transport set in Theorem \ref{th:1Dlocscheme}.
\par	{\small \bf Claim.} $
	\meas \left( \left\{x\in \X \; \text{\textbrokenbar} \; \widetilde{\left|\nabla u\right|}\left(x\right)=0\right\} \cap \mathcal T \right)=0.$
\par	{\it \small Proof of the Claim.} Assume  
	\begin{align*}
	\meas\left( \left\{x\in \X \; \text{\textbrokenbar} \; \widetilde{\left|\nabla u\right|}\left(x\right)=0\right\}\cap \mathcal T \right)>0,
	\end{align*}
	and let the transported critical set to be given by
	\begin{align*}
	\mathsf{C}_{\substack{\\\mathcal{T}}} :=  \mathcal T \cap \left\{x\in \X \; \text{\textbrokenbar} \; \widetilde{\left|\nabla u\right|}\left(x\right)=0\right\}.
	\end{align*}
	Then, by the disintegration of $\meas\mres_{\T}$ as given by Theorem \ref{th:1Dlocscheme},
	\begin{align*}
		0<\meas \left( \mathsf{C}_{\substack{\\\mathcal{T}}} \right) = \int_{\hspace{-1pt}_{\Q}}  \mathscr L^{^1} \left( \ray\left(\cdot, q\right)^{\hspace{-2pt}^{-1}}\left( \mathsf{C}_{\substack{\\\mathcal{T}}} \right) \right) \, \dif{\bf q}\left(q\right),
	\end{align*} 
	and therefore there exists $q\in \Q$ with $ \mathscr L^{^1} \left( \ray\left(\cdot, q\right)^{\hspace{-2pt}^{-1}}\left(  \mathsf{C}_{\substack{\\\mathcal{T}}}  \right) \right)>0$. Moreover,  by Lemma \ref{lem:ida} and Remark \ref{rem:cle}, we can assume that for this $q$, we have
	\begin{align*}
	\widetilde{\left|\nabla u\right|} \circ \upgamma_{\substack{\\[1pt]\hspace{-2pt}q}} \left(t\right)=  \left( u\circ \upgamma_{\substack{\\[1pt]\hspace{-2pt}q}} \right)^{\hspace{-2pt}'}(t) \quad
\mbox{ for } \mathscr L^{^1}\mbox{-a.e. } t. 
	\end{align*}
	But $u\circ \upgamma_{\substack{\\[1pt]\hspace{-2pt}q}} \left(r\right)=\mx_q\, \cos\left(r \shortminus \ai_q\right)$ by Corollary \ref{cor:cos}. This is a contradiction and proves the claim. \hfill \scalebox{0.6}{$\blacksquare$}
	\par Consequently it also follows that 
	\begin{align}\label{EQ:zero-set}
			\meas\left( \left\{\left|\nabla u\right|=0\right\} \smallsetminus { \left\{u=0\right\} }\right)=0;
	\end{align}
	recall that $\widetilde{\left|\nabla u\right|} = \left|\nabla u\right|$ $\meas$-a.e. 
\par 	{\bf 2$\centerdot$}
Now the \emph{improved Bochner inequality} (Theorem \ref{thm:improved})  for $u\in \mathbb{D}_{\hspace{-1pt}_{\infty}} $ is
	\begin{align*}
		{4\left|\nabla u\right|^{^2}}\left(\nicefrac{1}{2}\, { \Delta}_{\text{abs}} \left|\nabla u\right|^{^2} + \left|\nabla u\right|^{^2} \right)\geq \left| \nabla \left|\nabla u\right|^{^2} \right|^{^2} \quad \meas\mbox{-a.e.}\ .
	\end{align*}
	Therefore for $\upalpha$, it follows
	\begin{align*}
		{\left|\nabla u\right|^{^2}}{ \Delta}_{\text{abs}} \upalpha 
		& \geq \shortminus \, {2\left|\nabla u\right|^{^4}}+ \nicefrac{1}{2} 
		\left| \nabla \left|\nabla u\right|^{^2} \right|^{^2}+ {2\left|\nabla u\right|^{^2}} u \Delta u + 2 \left|\nabla u\right|^{^4} 
		\\ &= \nicefrac{1}{2}  \left| \nabla \left|\nabla u\right|^{^2} \right|^{^2} \shortminus 2{\left|\nabla u\right|^{^2}} u^2 \quad  \meas\mbox{-a.e.}\ .
	\end{align*}
	Hence, taking the above claim into consideration, we deduce
	\begin{align}\label{first}
		{ \Delta}_{\text{abs}} \upalpha \geq \nicefrac{1}{2}\ \left(\nicefrac{\left| \nabla \left|\nabla u\right|^{^2} \right|^{^2}}{\left|\nabla u\right|^{^2}}\right) \shortminus  2 u^2\ \ \meas\mbox{-a.e. on } \T.
	\end{align}
	Moreover, the Bochner inequality already yields that ${ \Delta}_{\text{abs}}\upalpha\geq 2 u \Delta u$ $\meas$-a.e. Hence, by~(\ref{EQ:zero-set}) we have
	\begin{align}\label{second}
		{\Delta}_{\text{abs}}\alpha\geq 0\quad \meas\mbox{-a.e. on }  \mathcal Z\subset \left\{u=0\right\}. 
	\end{align}
	Also recall that ${\bf \Delta}_{\text{sing}} \left|\nabla u\right|^{^2}\geq 0$ in measure sense. Upon integrating \eqref{first} and \eqref{second}, one arrives at
	\begin{align}\label{EQ:chain}
		0&=\int_{\hspace{-1pt}_{\X}}  \dif {\bf \Delta}\upalpha\geq \int_{\hspace{-1pt}_{\X}} {\Delta}_{\text{abs}} \upalpha \, \dif\meas \geq \int_{\hspace{-1pt}_{\T}} {\Delta}_{\text{abs}}\upalpha \, \dif\meas \notag\\
		&\geq  \int_{\hspace{-1pt}_{\T}}  \ \nicefrac{\left( \left| \nabla \left|\nabla u\right|^{^2}\right|^{^2} - {4u^2\left|\nabla u\right|^{^2}}\right)}{{2}\left|\nabla u\right|^{^2} } \, \dif\meas \\ &=\int_{\hspace{-1pt}_{\T}}  \nicefrac{
			\big\langle {\nabla\left(\left|\nabla u\right|^{^2} + u^2 \right)} , \nabla \left( \left|\nabla u\right|^{^2} - u^2 \right) \big\rangle}{{2}|\nabla u|^{^2}} \, \dif\meas \geq 0 \notag,
	\end{align}
	where the last inequality is manifested
	 in Corollary \ref{cor:alphabound}. It follows that the first and second inequalities in the chain of inequalities (\ref{EQ:chain}) are actually equalities. In particular, ${\bf \Delta}\upalpha=\left( {\Delta}_{\text{abs}}\upalpha \right) \meas$ and
	\begin{align*}
		\nicefrac{ \left( \left| \nabla \left|\nabla u\right|^{^2} \right|^{^2}- {4u^2\left|\nabla u\right|^{^2}} \right)}{\left|\nabla u\right|^{^2}}=0 \quad  \meas\mbox{-a.e. on }\T. 
	\end{align*}
	Hence
	\begin{align*}
		\left| \nabla \left|\nabla u\right|^{^2} \right|={2u\left|\nabla u\right|} \ \ \meas\mbox{-a.e. on }\ \T. 
	\end{align*}
Moreover, by a standard property of the minimal weak upper gradient, one gets
	\begin{align*}
	\left| \nabla \left|\nabla u\right|^{^2} \right|^{^2}=0 \quad  \meas\mbox{-a.e. on }\mathcal Z\subset \left\{u=0\right\}. 
	\end{align*}
To see this, one  applies~\cite[Corollary 2.25]{Cheeg1} twice. First, it holds that $\left|\nabla u\right|=0$ $\meas$-a.e. on $\left\{u=0\right\}$. Then, since $\left|\nabla u\right|^{^2}\in \Sobol\left(\X\right)$ (because $u\in \mathbb{D}_{\hspace{-1pt}_{\infty}} $), it holds $\left|\nabla \left|\nabla u\right|^{^2}\right|=0$ $\meas$-a.e. in $\left\{\left|\nabla u\right|^{^2}=0\right\}$; note that according to \cite{agslipschitz},  minimal weak upper gradients are relaxed gradients in the sense of \cite{Cheeg1}. Consequently
	\begin{align*}
		\left|\nabla \left|\nabla u\right|^{^2}\right| = {2u\left|\nabla u\right|} \quad \meas\mbox{-a.e.}
	\end{align*}
	Since $u, \left|\nabla u\right|\in \LL^{\hspace{-2pt}^{\infty}}\left(\meas\right)$, it follows $\left|\nabla \left|\nabla u\right|^{^2}\right|\in \LL^{\hspace{-2pt}^{\infty}}$; by the Sobolev-to-Lipschitz property (see e.g. \cite[Theorem 6.2]{AGSRiem}), $\left|\nabla u\right|^{^2}$ has a Lipschitz representative and the thesis of this proposition is proven.
\end{proof}
\begin{corollary}
	$\meas \left( \left\{\left|\nabla u\right|=0\right\} \right)=0$ and $\upalpha\equiv \left(\max u\right)\strut^{\hspace{-2pt}2}$ $\meas$-a.e. on $\X$.
\end{corollary}
\begin{proof}[\footnotesize \textbf{Proof}]
	By the proof of the previous proposition, we also obtain ${\bf \Delta}\upalpha\geq 0$ on $\X$. Now since $\upalpha$ is Lipschitz continuous and $\X$ is compact, the strong maximum principle given in \cite[Theorem 2.8]{Gigli-Rigoni} (see Theorem~\ref{thm:Gig-Rig} in our \textsection\thinspace\ref{sec:prelim}; also see \cite{BB, GM}) yields constancy of $\upalpha$ throughout $\X$. Now by (\ref{EQ:constancy}), 
	\begin{align*}
	\upalpha \equiv \text{const} \ge \; \mx_{\substack{\\q}}\strut^{\hspace{-4pt}2} > 0 \ \ {\bf q}\mbox{-a.e.} \;q \in \Q,
	\end{align*}
	which combined with the denseness of the transport set and continuity of $\upalpha$ and $u$, yields
	$
	\upalpha\equiv \left(\max u\right)\strut^{\hspace{-2pt}2}.
	$
\par	It then follows that 
	$\meas\left( \left\{\left|\nabla u\right|=0\right\} \right)=0$ since by~(\ref{EQ:zero-set}), $\upalpha=0$ $\meas$-a.e. on $\left\{\left|\nabla u\right|=0\right\}$.
\end{proof}
Let us normalize $u$ such that $\upalpha=1$ $\meas$-a.e.\ . In the sequel, $u$, $\left|\nabla u\right|$ and $\upalpha$ are the Lipschitz continuous representatives of the corresponding functions involved. 
\begin{definition}(Singular set)\label{defn:sing-set}
	Define the extremal sets $\Sing_{_\pm} :=u^{-1}\left(\plmi 1  \right)$ and set  the singular set to be $\Sing :=\Sing_{_-} \, \dot{\sqcup} \,  \Sing_{_+}$. 
\end{definition}
It is clear that $\X \smallsetminus \Sing$ is open in $\X$. We also note that $\left\{\left|\nabla u\right|=0\right\}= \Sing$ (of course in $\meas$-a.e. sense) and hence $\meas\left(\Sing\right)=0$.
\begin{corollary}
		$\meas\left(\left\{ u=0\right\}\right)=0$. In particular $\left(\X \smallsetminus \Sing\right)\cap \T$ has full $\meas$-measure.
\end{corollary}
	\begin{proof}[\footnotesize \textbf{Proof}]
		By \cite[Corollary 2.25]{Cheeg1} (that we have already used a couple of times) $\left|\nabla u\right|=0$ $\meas$-a.e. on $\left\{u=0\right\}$. On the other hand $\upalpha= \left|\nabla u\right|^{^2} + u^2 = \left(\max u\right)\strut^{\hspace{-2pt}2}=\textsf{const}>0$. Therefore, $\meas\left(\left\{u=0\right\}\right)=0$.
\end{proof}
\begin{lemma}[Distributional harmonicity of $f = \sin^{\hspace{-1pt}^{-1}} \circ\ u$]
	The function $f = \sin^{\hspace{-1pt}^{-1}} \circ\ u$ is locally $1$-Lipschitz in $\X\smallsetminus \Sing$ and $f$ is in the domain of the distributional Laplacian on $\X \smallsetminus \Sing$, $\Dom \left({\bf \Delta}, \X \smallsetminus \Sing\right)$. furthermore, $f$ satisfies ${\bf \Delta} f = 0$ on $\X \smallsetminus \Sing$.
\end{lemma}
\begin{proof}[\footnotesize \textbf{Proof}]
	Using the  chain rule for the minimal weak upper gradient (see e.g. \cite[Proposition 3.15]{Gigli-1}), we first compute
	\begin{align}\label{EQ:arcsin-der}
		\left|\nabla f\right|= \left| \nabla \left( \sin^{\hspace{-1pt}^{-1}} \circ \ u \right) \right|= \left(1 \shortminus u^2\right)^{\shortminus \nicefrac{1}{2}}\left|\nabla u\right|= 1 \quad \mbox{ $\meas$-a.e. on }\quad  \X \smallsetminus \Sing.
	\end{align}
	Hence, by the Sobolev-to-Lipschitz property, $f$ is locally $1$-Lipschitz in $\X\smallsetminus \Sing$. 
	\par Using the chain rule for the distributional Laplacian (see Proposition~\ref{prop:chainrulelaplacian}), we have
	\begin{align*}
		{\bf \Delta}\left(f\right) &= {\bf \Delta} \left( \sin^{\hspace{-1pt}^{-1}} \circ\ u  \right) = \left( \sin^{\hspace{-1pt}^{-1}} \right)^{\hspace{-2pt}'}\left(u\right) {\bf \Delta} u + \left( \sin^{\hspace{-1pt}^{-1}} \right)^{\hspace{-2pt}''} \left(u\right) \left| \nabla u \right|^{^2} {\meas} \notag \\ &= \shortminus u \left(  1 \shortminus u^2   \right)^{\shortminus \nicefrac{1}{2}} {\meas} +   u \left( 1 \shortminus u^2 \right)^{\shortminus \nicefrac{3}{2} } \left| \nabla u \right|^{^2} {\meas}  \notag  = 0 \quad \mbox{ on } \quad \X \smallsetminus \Sing.
	\end{align*}
\end{proof}
\begin{corollary} \label{cor:hess_id}
	For any $\upvarphi \in \mathbb{D}_{\hspace{-1pt}_{\infty}}$ supported in $\ovs{V}\subset \X \setminus \Sing$ for $V\subset \X \smallsetminus \Sing$ open, it holds $\upvarphi \cdot f\in \Dom\left(\Delta\right)$; in particular, $\upvarphi \cdot f\in \mathcal{H}^{\hspace{0.7pt}^{2, 2}}\left(\X\right)$ (see the exposition following Theorem \ref{thm:improved}). Moreover, it holds
	\begin{align*}
	\Hess\left(\upvarphi \cdot f\right)\restr_U=0 \quad \meas\mbox{-a.e. on }\quad  U\subset \X \smallsetminus \Sing,
	\end{align*}	
	whenever $\upvarphi\equiv 1$ on $U$ for some $U$ open in $\X$ and $U\subset \X \smallsetminus \Sing$.
\end{corollary}
\begin{proof}[\footnotesize \textbf{Proof}]
	From \eqref{EQ:arcsin-der} we deduce $\left|\nabla f\right|=1$ $\meas$-a.e.\ on $\X \smallsetminus \Sing$. Let $x\in \X \smallsetminus \Sing$ and  pick a function $\upvarphi\in \mathbb{D}_{\hspace{-1pt}_{\infty}} $ with $\supp \left(\upvarphi\right)\subset \ovs{V}\subset \X \smallsetminus \Sing$ for some open set $V$ with $U\subset V$, and $\upvarphi\equiv 1$ on $U$, where $U$ is as in the statement of the corollary. The existence of such cutoff functions is guaranteed by Lemma \ref{lem:cutoff}.
\par	{\bf \small Claim.} $\upvarphi \cdot f\in \Dom\left(\Delta\right)$. 
	\par {\it \small Proof of the Claim.} By \cite[Theorem 4.29]{Gigli-1} (Leibniz rule for ${\bf \Delta}$), we have that $\upvarphi\cdot f\in \Dom\left(\Delta, \X \smallsetminus \Sing \right)$ and 
	\begin{align}\label{id:uv}
		{\bf \Delta}\left(\upvarphi\cdot f\right) &= \upvarphi {\bf \Delta} f + f{\bf  \Delta} \upvarphi + 2\langle\nabla \upvarphi, \nabla f\rangle {\meas} \notag \\ &= f{\bf \Delta } \upvarphi + 2\langle \nabla \upvarphi, \nabla f\rangle {\meas} \\ &= \left( f\Delta \upvarphi +  2\langle \nabla \upvarphi , \nabla f\rangle\right) {\meas} \notag,
	\end{align}
	since $\upvarphi$ and $f$ are Lipschitz on $\X \smallsetminus \Sing$.  
	\par It holds $\Delta \upvarphi, \left|\nabla \upvarphi\right| \in \LL^{\hspace{-2pt}^\infty}$ since $\upvarphi\in \mathbb{D}_{\hspace{-1pt}_{\infty}}$. Therefore, one has ${\bf \Delta}\left(\upvarphi \cdot f\right) = \Delta\left(\upvarphi \cdot f\right){\meas}$ and {$\Delta \left(\upvarphi \cdot f\right) \in \Ltwo\left(\meas\right)$}. Also the support of $\upvarphi\cdot f$ is contained in $\ovs{V}\subset \X \smallsetminus \Sing$. Therefore $\upvarphi\cdot f \in \Dom\left(\Delta\right)$. This proves the claim. \hfill \scalebox{0.6}{$\blacksquare$}
	\par Since $\upvarphi\cdot f\in \Dom\left(\Delta\right)$, we can pick a sequence $k_{_n}\in \mathbb{D}_{\hspace{-1pt}_{\infty}}$ such that (eventually after passing to a subsequence), $ k_{_n}, \left| \nabla k_{_n} \right|$ and $\Delta \left( k_{_n} \right)$ converge in $\Ltwo\left(\meas\right)$ sense to $\upvarphi\cdot f, \left|\nabla \left(\upvarphi \cdot  f\right)\right|$ and $\Delta\left(\upvarphi\cdot f\right)$ respectively.  Let $\ball{\varepsilon}\left(x\right)\subset U$. 
	We pick another non-negative cutoff test function $\uppsi\in \mathbb{D}_{\hspace{-1pt}_{\infty}} $ supported in $\ball{\varepsilon}\left(x\right)$ such that $\uppsi\equiv 1$ on $\ball{\nicefrac{\varepsilon}{2}}\left(x\right)$. Then, we can write the improved Bochner formula for $k_{_n}$ in the form of
	\begin{align*}
		\uphat{\BG_{\hspace{-2pt}_2}} \left( k_{_n} ;\uppsi \right)\geq \int_{\hspace{-1pt}_\X} \left\| \Hess\left(k_{_n}\right)\right\|^{^2}_{_{\textsf{HS}}} \uppsi  \, \dif\meas.
	\end{align*}
Notice that the above version of  Bochner inequality follows from \eqref{eq:bochner-modified} in the remark directly after Theorem \ref{thm:improved} together with the exposition before Definition~\ref{def:be} where $\uphat{\BG}_{\substack{\\[1pt]\hspace{-2pt}2}}$ is introduced. Indeed, we know $\uphat{\BG_{\hspace{-2pt}_2}} $ and ${\BG}_{\hspace{-2pt}_2} $ coincide on the intersection of their domains. Then, exactly as in the first part of the proof of \cite[Proposition 3.2]{GKK}, it follows
	\begin{align}\label{id:vvww}
		\uphat{\BG_{\hspace{-2pt}_2}} \left(\upvarphi \cdot f;\uppsi\right)\geq  \int_{\hspace{-1pt}_\X} \left\| \Hess\left(\upvarphi \cdot f \right)\right\|^{^2}_{_{\textsf{HS}}} \uppsi \, \dif\meas\geq \int_{\hspace{-1pt}_{\ball{\nicefrac{\varepsilon}{2}}\left(x\right)}} \left\| \Hess\left( \upvarphi \cdot f \right)\right\|^{^2}_{_{\textsf{HS}}} \, \dif\meas.
	\end{align}
	Let us analyze the term $\uphat{\BG}_{\hspace{-2pt}_2}  \left(f\upvarphi;\uppsi\right)$ which is well-defined because $\upvarphi \cdot f \in \Dom\left(\Delta\right)$ and $\uppsi\in \mathbb{D}_{\hspace{-1pt}_{\infty}}$. It holds
	\begin{align}\label{id:vw}
		&\uphat{\BG}_{\hspace{-2pt}_2} \left(\upvarphi \cdot f;\uppsi\right) = \\ &\nicefrac{1}{2} \int_{\hspace{-1pt}_\X} \left| \nabla \left(\upvarphi \cdot f \right) \right|^{^2} \Delta \uppsi \, \dif\meas + \int_{\hspace{-1pt}_\X} \left( \Delta \left( \upvarphi \cdot f \right) \right)^{^2} \uppsi \, \dif\meas + \int_{\hspace{-1pt}_\X}  \langle \nabla \uppsi, \nabla \left( \upvarphi \cdot f \right)\rangle \Delta \left( \upvarphi \cdot f \right) \, \dif\meas. \notag
	\end{align}
	For the first term on the right-hand side of \eqref{id:vw}, we have
	\begin{align}\label{EQ:analyze1}
		\int_{\hspace{-1pt}_\X} \left| \nabla \left( \upvarphi \cdot f \right) \right|^{^2} \Delta \uppsi \, \dif\meas&= \int_{\hspace{-1pt}_{\B_{_\varepsilon}\left(x\right)}} \left| \nabla \left( \upvarphi \cdot f \right) \right|^{^2} \Delta \uppsi \, \dif\meas \notag \\ &= \int_{\hspace{-1pt}_{\B_{_\varepsilon}\left(x\right)}} \left|\nabla f\right|^{^2} \Delta \uppsi \, \dif\meas  \\
		&= \int_{\hspace{-1pt}_\X} \Delta\uppsi \, \dif\meas=0 \notag.
	\end{align}
	In the first equality in (\ref{EQ:analyze1}), we have used the fact that $\uppsi$ is compactly supported in $\ball{\varepsilon}\left(x\right)$ and hence $\Delta \uppsi =0$ on $\X \smallsetminus \ovs{\ball{\varepsilon}\left(x\right)}$ by locality of $\Delta$; indeed, the latter can be checked by testing $\Delta \uppsi$ against Lipschitz functions compactly supported in $\X \smallsetminus \ovs{\ball{\varepsilon}\left(x\right)}$. The second equality uses  $f=\upvarphi \cdot f$ on $\ball{\varepsilon}\left(x\right) \subset U$ and the locality property of $\left|\nabla \,\cdot\right|$ while in the third one, the fact that  $\left|\nabla f\right|=1$ $\meas$-a.e. in $\X \smallsetminus \Sing$ and again $\uppsi$ being compactly supported in $\ball{\varepsilon}\left(x\right)$ are used.
\par Similarly -- using \eqref{id:uv} -- one can check that the other two terms on the right hand side in \eqref{id:vw} must vanish as well. Indeed, in the last integral in \eqref{id:vw}, we observe that 
	\begin{align*}
		\int_{\hspace{-1pt}_\X} \langle \nabla \uppsi , \nabla \left(\upvarphi \cdot f \right)\rangle \Delta \left( \upvarphi \cdot f \right) \, \dif\meas &= \int_{\hspace{-1pt}_{\B_{_\varepsilon}\left(x\right)}} \langle \nabla \uppsi , \nabla \left( \upvarphi \cdot f \right)\rangle \Delta \left( \upvarphi \cdot f \right) \, \dif\meas \\
		&= \int_{\hspace{-1pt}_{\B_{_\varepsilon}\left(x\right)}}\langle \nabla \uppsi , \nabla \left( \upvarphi \cdot f \right)\rangle {{\dif\bf \Delta} \left(\upvarphi \cdot  f\right) } \\ &= \int_{\hspace{-1pt}_{\B_{_\varepsilon}\left(x\right)}} \langle \nabla \uppsi, \nabla f\rangle {{\dif\bf \Delta} f}  =0.
	\end{align*}
	Let us elaborate. For the first equality, we use that $\uppsi$ is compactly supported in $\ball{\varepsilon}\left(x\right)$ and hence $\left|\nabla \uppsi\right|$ vanishes on $\X \smallsetminus \ovs{\ball{\varepsilon}\left(x\right)}$. The second equality uses the coincidence of distributional Laplacian and the Cheeger Laplacian for $\upvarphi f\in \Dom\left(\Delta\right)$ and in the third equality, $f= \upvarphi \cdot f$ on $\ball{\varepsilon}\left(x\right)\subset U$ and the locality property of ${\bf \Delta}$ and $\left|\nabla\, \cdot \right|$ are used. Finally, the last equality uses ${\bf \Delta} f=0$ on  $\X \smallsetminus \Sing$ and hence on $\ball{\varepsilon}\left(x\right)$ by the locality property of ${\bf \Delta}$ again.
	\par In the same way, one can verify
	\begin{align*}
		\int_{\hspace{-1pt}_{\X}} \left( \Delta \left(\upvarphi \cdot f\right) \right)^{^2} \uppsi \, \dif\meas=0.
	\end{align*}
	\par Consequently, together with \eqref{id:vvww} it follows
	\begin{align*}
		\int_{\hspace{-1pt}_{\B_{\varepsilon/2}\left(x\right)}} \left\| \Hess\left( \upvarphi \cdot f \right)\right\|^{^2}_{_{\textsf{HS}}} \, \dif\meas=0,
	\end{align*}
	thus, it follows that $\left\| \Hess\left( \upvarphi \cdot f \right)\right\|^{^2}_{_{\textsf{HS}}} =0$ $\meas$-a.e.  on $\B_{\nicefrac{\varepsilon}{2}}\left(x\right)$. 
\end{proof}
\subsubsection*{\small \bf \textit{Proof of theorem \ref{thm:hessian} (Second part)}}
Let $\upvarphi\in \mathbb{D}_{\hspace{-1pt}_{\infty}}$ be as in Corollary \ref{cor:hess_id}. The functions
\begin{align*}
\upvarphi \cdot f \quad \text{and} \quad \sin\circ\, \left(\upvarphi \cdot f\right)=:\tilde u,
\end{align*}
both belong to the function space $\mathbb{D}_{\hspace{-1pt}_{\infty}}$; as a result, they have well-defined Hessians on $\X$; so, we can apply the Hessian chain rule (Proposition \ref{prop:chainrulehessian}) to $\tilde u$ on $\X \smallsetminus \Sing$.  Thus, we get
\begin{align*}
	\Hess \left(\tilde{u}\right) &= \sin^{''} \circ\, \left(\upvarphi \cdot f\right) \, \dif\left(\upvarphi \cdot f\right)  \otimes \dif\left(\upvarphi \cdot f\right)   + \sin^{'} \circ\, \left(\upvarphi \cdot f\right) \, \Hess\left(\upvarphi \cdot f\right) \\ &= \shortminus \sin \circ\, \left(\upvarphi \cdot f\right) \, \dif\left(\upvarphi \cdot f\right)  \otimes \dif\left(\upvarphi \cdot f\right) \quad \mbox{ on } \quad U.
\end{align*}
From the fact that 
\begin{align*}
	\Hess \left(\upvarphi\cdot f\right) = 0 \quad \meas\mbox{-a.e. on}\quad  U\subset\left\{\upvarphi\equiv1\right\},
\end{align*}
we deduce
%
\begin{align*}
	\Hess \left(\tilde{u}\right)\left(  \nabla g, \nabla g \right) &= \shortminus \tilde{u} \left< \nabla \left(\upvarphi \cdot f\right), \nabla g \right>^{\hspace{-1pt}^2} \quad \mbox{ on } \quad  U, \quad \forall g \in  \mathbb{D}_{\hspace{-1pt}_{\infty}}.
\end{align*}
We also know that the eigenfunction $u$ belongs to $ \mathbb{D}_{\hspace{-1pt}_{\infty}}$. Moreover, for any $U\subset \X \smallsetminus \Sing$ with $\meas\left(\partial U\right)=0$ and such that $\upvarphi\restr_U\equiv 1$, it follows $\tilde u\restr_U= u\restr_U$.  We also recall that since $\tilde u, u\in \mathbb{D}_{\hspace{-1pt}_{\infty}}$, the Hessian for both $u$ and $\tilde u$ is computed according to the Riemannian formula
\begin{align*}
	&\Hess \left(u\right)\left(\nabla g, \nabla h\right)= \\ &\phantom{cys}\nicefrac{1}{2} \, \Big(\big\langle \nabla g, \nabla \langle \nabla u, \nabla h\rangle \big\rangle + \big\langle \nabla h, \nabla \langle \nabla g, \nabla u\rangle \big\rangle  - \big\langle \nabla u, \nabla \langle \nabla g, \nabla h\rangle \big\rangle \Big) \quad \forall g,h \in\mathbb{D}_{\hspace{-1pt}_{\infty}}.  \notag
\end{align*}
Now based on higher-order locality (see Lemma~\ref{lem:ho-locality}) and without ambiguity, we can write
\begin{align*}
	\Hess \left(u\right)\restr_{\ovs{U}} &=\Hess\left(u\restr_{\ovs{U}}\right)= \Hess\left(\tilde u\restr_{\ovs{U}}\right)=\Hess \left(\tilde u\right)\restr_{\ovs{U}} \\ &= \shortminus \, \tilde u\restr_{\ovs{U}}\langle \nabla \left(\upvarphi \cdot f\right)\restr_{\ovs{U}}, \cdot \rangle\otimes \langle \nabla \left(\upvarphi \cdot f \right)\restr_{\ovs{U}}, \cdot\rangle \\  &= \shortminus \, u\restr_{\ovs{U}}\langle \nabla f\restr_{\ovs{U}}, \cdot \rangle\otimes \langle \nabla f\restr_{\ovs{U}}, \cdot \rangle.
\end{align*}
Since $u\in \mathbb{D}_{\hspace{-1pt}_{\infty}}$ and since $\meas\left(\Sing\right)=0$, it follows that $\Hess \left(u\right)= \shortminus \, u\langle \nabla f, \cdot \rangle\otimes \langle \nabla f, \cdot \rangle$.
This concludes the proof of Theorem~\ref{thm:hessian}. 
\qed
\section{The gradient flow that provides the splitting}\label{sec:grad-flow}
\subsection{Setup}
Recall the defnition of the singular set $\Sing$ in Definition~\ref{defn:sing-set}.
\begin{remark}
	Since $\left|\nabla f\right|= 1$ on $\X \smallsetminus \Sing$ and $\meas\left(\Sing\right)=0$,  one deduces $\left|\nabla f\right|= 1$ $\meas$-a.e. Furthermore, by the Sobolev-to-Lipschitz property (see Definition \ref{def_sobtolip}) $f$ is $1$-Lipschitz on $\X$.
\end{remark}
\begin{remark}\label{rem:useful}
	We note that, by $\upalpha\equiv 1$ and  \eqref{EQ:arcsin-der} (see also Lemma \ref{lem:repa}), one deduces
	\begin{align*}
		\nicefrac{\dif}{\dif t} \ f\circ \upgamma_{\substack{\\[1pt]\hspace{-2pt}q}} = 1  \quad \mbox{ for } \quad {\bf q}\mbox{-a.e. } q\in \Q;
	\end{align*}
which upon being integrated from $s$ to $t (\ge s)$ for $s, t \in \left( \ai_{\substack{\\q}},\bi_{\substack{\\q}} \right) $, yields
	\begin{align*}
	f\left(  \upgamma_{\substack{\\[1pt]\hspace{-2pt}q}} \left(t\right) \right) - f\left( \upgamma_{\substack{\\[1pt]\hspace{-2pt}q}} \left(s\right) \right)= \dist \left( \upgamma_{\substack{\\[1pt]\hspace{-2pt}q}}\left(s\right),  \upgamma_{\substack{\\[1pt]\hspace{-2pt}q}}\left(t\right) \right).
	\end{align*}
Hence, $ \upgamma_{\substack{\\[1pt]\hspace{-2pt}q}}$ is a transport geodesic for the $1$-Lipschitz potential function $f$ and for ${\bf q}$-a.e. $q\in \Q$, we find $t_{_0}\in \left( \ai_{\substack{\\q}},\bi_{\substack{\\q}} \right)$ such that $f\circ  \upgamma_{\substack{\\[1pt]\hspace{-2pt}q}} \left( t_{_0} \right)=0$. 
	In particular, we can replace $\T$ with $\widetilde \T$ where $\meas\left( \T \smallsetminus  \widetilde \T \right)=0$ and $\widetilde \T$ is a union of transport geodesics w.r.t. $f$. In the following, we denote $\widetilde \T$ by $\T$ again with a slight abuse of notation. 
\par So, we can choose a section $\mathfrak{s}: \T\rightarrow \T$ with $\mathfrak{s}\left(\T\right)\subset f^{^{-1}}\left(0\right)$. Consequently, up to a $\bf q$-null subset, we can consider $\Q$ as a subset of $f^{^{-1}}\left(0\right)$ and therefore $\ai_{\substack{\\q}}= \shortminus \nicefrac{\uppi}{2}$ and $\bi_{\substack{\\q}}=\nicefrac{\uppi}{2}$ for ${\bf q}$-a.e. $q\in \Q$.
	\par With arclength parametrization, the section $\mathfrak{s}:\T\rightarrow \T$ is given by $x\mapsto \upgamma_{\substack{\\[1pt]\hspace{-2pt}q}}\left(0\right)$  when $x= \upgamma_{\substack{\\[1pt]\hspace{-2pt}q}} \left(t\right)$ for some $t\in \left( \shortminus \nicefrac{\uppi}{2},\nicefrac{\uppi}{2} \right)$ and $q\in \Q$; for why such a measurable selection is possible, see e.g.~\cite{cavmil} or~\cite[Lemma 3.9]{CM-1}.
\end{remark}
\begin{definition}[Horizontal geodesic]\label{defn:horiz-geo}
	We refer to a geodesic with endpoints in $\Sing_{_+}$ and $\Sing_{_-}$ as a horizontal geodesic.
\end{definition}	
\begin{remark}\label{rem:horiz-foliation}
	Since $\T$ has full measure in $\X \smallsetminus \Sing$,
	$\meas$-almost every $x\in \X \smallsetminus \Sing$ lies in the interior of a horizontal geodesic. In particular, such a horizontal geodesic is a transport geodesic of $f$.
\end{remark}
\subsubsection*{\small \bf \emph{ Revisiting the ray map}}
	Recall the ray map from \textsection\thinspace\ref{subsec:gradu}, given by
	\begin{align*}
		\ray: \mathcal V\subset \mathbb R\times \Q\rightarrow \T , \ \ray\left(t,q\right)=\upgamma_{\substack{\\[1pt]\hspace{-2pt}q}} \left(t\right),
	\end{align*}
	that is a measurable map. By the Remark \ref{rem:useful}, we have $\mathcal V\supset \left( \shortminus \nicefrac{\uppi}{2},\nicefrac{\uppi}{2} \right) \times \Q$ up to an $\meas$-null set.
\begin{remark}\label{rem:q-on-Y}
	The previous considerations yield the following:  The $1\D$ localization w.r.t. $f$ is a measurable decomposition of $\X$ in the sense that $\meas$-almost every point $x\in \X \smallsetminus \Sing$ lies on  a transport geodesic $\upgamma: \left[ \shortminus \nicefrac{\uppi}{2}, \nicefrac{\uppi}{2} \right] \rightarrow \X$ w.r.t. $f$ and of length $\uppi$, with endpoints in $\Sing_{_-}$ and $\Sing_{_+}$ respectively  i.e. on a horizontal geodesic.  As in the Remark \ref{rem:useful}, we can identify this geodesic with the point $\upgamma\left(0\right)\in f^{^{-1}}\left(0\right)=:\Y$. Hence, we can also consider the measure ${\bf q}$ as a Borel measure concentrated on $\Y$; these considerations lead to the disintegration 
	\begin{align}\label{eq:disintegration2}
		\meas=\int_{\hspace{-1pt}_\Y} \left(\upgamma_{\substack{\\[1pt]\hspace{-2pt}q}}\right)_{\hspace{-1pt}*} \mathscr L^{^1}\mres_{\left[\shortminus \nicefrac{\uppi}{2},\nicefrac{\uppi}{2}\right]} \, \dif{\bf q}\left(q\right),
	\end{align}
	that -- up to a ${\bf q}$-null set -- coincides with the previous disintegration given in Theorem~\ref{th:1Dlocscheme}. Indeed, the subtle advantage compared to the previous disintegration in~Theorem~\ref{th:1Dlocscheme}  is that, now, we can find for $\meas$-a.e. $x\in \X \smallsetminus S$, a geodesic $\upgamma_{\substack{\\[1pt]\hspace{-2pt}q}}$ with $q\in \Y$, and some $t\in \left(\shortminus \nicefrac{\uppi}{2},\nicefrac{\uppi}{2}\right)$ such that $\upgamma_{\substack{\\[1pt]\hspace{-2pt}q}}\left(t\right) = x$. However, with a slight abuse of notation, we keep the previous notations. 
\end{remark}
\subsubsection*{\small \bf \textit{Local geodesic linearity of $f$ on the regular set}}
The next proposition is an immediate consequence of the second variation formula~\cite{GT, GT-2} and the Hessian identity for {$f := \sin^{\hspace{-2pt}^{-1}} \circ\, u$} in Corollary \ref{cor:hess_id}.
\begin{proposition}
	Let $U\subset \X \smallsetminus \Sing$ be an open set and let $\left( \upmu_{\substack{\\[1pt]\hspace{1pt}t}} \right)_{t\in \left[0,1\right]}$ be an $\Ltwo$-Wasserstein geodesic with $\upmu_{\substack{\\[1pt]\hspace{1pt}t}} \in \Prob_{\hspace{-5pt}_2} \left(\X,\meas\right)$ and $\upmu_{\substack{\\[1pt]\hspace{1pt}t}} \left(U\right)=1$ for every $t\in \left[0,1\right]$. Then 
	\begin{align*}
		{\nicefrac{\dif^{^{\,2}}}{\dif t^2}} \int_{\hspace{-1pt}_U} f \, \dif \upmu_{\substack{\\[1pt]\hspace{1pt}t}}=0. 
	\end{align*}
	In particular, we have the linearity
	\begin{align*}
		\left[0,1\right]	\ni t \mapsto \int_{\hspace{-1pt}_U} f \, \dif\upmu_{\substack{\\[1pt]\hspace{1pt}t}}= \left(1 \shortminus t\right)\int_{\hspace{-1pt}_U} f \, \dif\upmu_{_0} + t\int_{\hspace{-1pt}_U} f \, \dif \upmu_{_1}.
	\end{align*}
\end{proposition}
\begin{corollary}[Local geodesic linearity of $f$]\label{cor:loccon}
Suppose $\ball{2\updelta}\left( x_{\hspace{-0.5pt}_0} \right) \subset \X \smallsetminus \Sing$ and $x\in \ball{\updelta}\left( x_{\hspace{-0.5pt}_0} \right)$ are given. Then, for $\meas$-a.e.  $y$ in $\ball{\updelta}\left( x_{\hspace{-0.5pt}_0} \right)$ there exists a geodesic $\upeta:\left[0,1\right]\rightarrow \ball{2\updelta}\left( x_{\hspace{-0.5pt}_0} \right)$ between $x$ and $y$ such that 
	\begin{align*}
		f\circ \upeta \left(t\right) = \left(1 \shortminus t\right) f\left(x\right) + t f\left(y\right). 
	\end{align*}
\end{corollary}
\begin{remark}
	Note that the assumption that $\upmu_{\substack{\\[1pt]\hspace{1pt}t}}\left(U\right)=1$ for every $t\in \left[0,1\right]$ and some open set $U \subset \X \smallsetminus \Sing$, is paramount in the statement of Proposition~\ref{cor:loccon}. Indeed, consider the  space $\X=\mathbb S^{^1}$ equipped with the standard metric and the standard measure. One can easily observe that the function $f=\sin^{\hspace{-2pt}^{-1}}\circ\, u$ will be given by a distance function $\dist\left(p,\cdot\right)$ for some $p\in \mathbb S^{^1}$.  Let $q$ be the antipodal point to $p$; since $q$ is the conjugate point of $p$, $f$ is not affine on  any Wasserstein geodesics $\upmu_{\substack{\\[1pt]\hspace{1pt}t}}$ with $\upmu_{\substack{\\[1pt]\hspace{1pt}t}}\left( \ball{\updelta}(q) \right)>0$ for some $\updelta>0$.
\end{remark}
\subsection{Gradient flow of $f$.}\label{subsec:gf-f}
\subsubsection*{\small \bf \emph{Definition of the flow}}
We can define $ \F_{\substack{\\[1pt]\hspace{-3pt}t}}  \left(x\right):=\F\left(t,x\right):= \upgamma_{\substack{\\[1pt]\hspace{-2pt}q}} \left(s+t\right)$ whenever $x= \upgamma_{\substack{\\[1pt]\hspace{-2pt}q}} \left(s\right)$ for $q\in \Y$ and $s+t\in \left( \ai_{\substack{\\q}}, \bi_{\substack{\\q}} \right)$.  Note that 
\begin{align*}
\mathrm{Domain} \left(\F_{\substack{\\[1pt]\hspace{-3pt}t}} \right)= \left\{x\in \T \; \text{\textbrokenbar} \; \exists y\in {\frak R}\left(x\right) \mbox{ s.t. } f\left(y\right) = f\left(x\right) + t\right\}\subset \T.
\end{align*}
The intuitive picture is that $\F_{\substack{\\[1pt]\hspace{-3pt}t}}$, whenever defined, is the gradient flow of $f$; in Proposition~\ref{prop:evi}, we will attach a precise meaning to this intuition.
\subsubsection*{\small \bf \emph{Push-forward of measures along the flow}}
Set ${\rm U}_{\substack{\\\!\updelta}} := f^{^{-1}}\left( \left(\shortminus 1+\updelta, 1 \shortminus \updelta \right) \right)$. Note that $\rm U_{\substack{\\\!\updelta}}$ is an open subset of $\X \smallsetminus \Sing$ and that $\ball{\updelta}\left(x\right)\subset \X \smallsetminus \Sing$ for all $x\in \rm U_{\substack{\\\!\updelta}}$.
Let $\upmu \in \Prob \left(\X,\meas\right)$ such that $\upmu\left( \rm U_{\substack{\\\updelta}} \right)=1$.  Since $\upmu\left( \rm U_{\substack{\\\updelta}} \right)=1$, it follows $\upmu\left( \mathrm{Domain}  \left( \F_{\substack{\\[1pt]\hspace{-3pt}t}} \right)\right)=1$ for $t\in \left[ \shortminus \updelta, \updelta \right]$. Hence, we can define $\left( \F_{\substack{\\[1pt]\hspace{-3pt}t}} \right)_{\hspace{-1pt}*}\upmu=: \upmu_{\substack{\\[1pt]\hspace{1pt}t}} $ for $t\in \left[ \shortminus \updelta,\updelta \right]$.  
\par As the following lemma shows, the pushforward of a measure along this flow is a Wasserstein geodesic. We do not use this fact in an essential way; however, it justifies sticking to the same notation $\upmu_{\substack{\\[1pt]\hspace{1pt}t}}$ which we had previously used for Wasserstein geodesics. 
\begin{lemma}
The curve $ \upmu_{\substack{\\[1pt]\hspace{1pt}t}} := \left( \F_{\substack{\\[1pt]\hspace{-3pt}t}} \right)_{\hspace{-1pt}*}\upmu$ as defined above, is an $\Ltwo$-Wasserstein geodesic. 	
\end{lemma}
\begin{proof}[\footnotesize \textbf{Proof}]
In order to show this curve is an $\Ltwo$-Wasserstein geodesic, it suffices to show for each $s \le t$, the graph of the transport map $\F_{\substack{\\[1pt]\hspace{-3pt}s-t}}$ -- which is the support of the transport plan $\left(\F_{\substack{\\[1pt]\hspace{-3pt}s}} ,  \F_{\substack{\\[1pt]\hspace{-3pt}t}} \right)_{\hspace{-1pt}*} \upmu$ -- is concentrated on a $\dist^{^2}$-cyclically monotone set; see~\cite{Pr}. The $\dist^{^2}$-cyclical monotonicity follows from~\cite[Lemma 4.1]{CM-1} since for any $s \le t$, we have
\begin{align*}
	\left(  f\left( \bar  \upgamma_{\substack{\\[1pt]\hspace{-2pt}q_{_1}}}(t)\right) -  f\left( \bar  \upgamma_{\substack{\\[1pt]\hspace{-2pt}q_{_1}}}(s)\right) \right) \cdot \left( f\left( \bar  \upgamma_{\substack{\\[1pt]\hspace{-2pt}q_{_2}}}(t)\right) -  f\left( \bar  \upgamma_{\substack{\\[1pt]\hspace{-2pt}q_{_2}}}(s)\right) \right) = \left(  t \shortminus s \right)^{^2}  \ge 0,
\end{align*}
for any two reparametrized transport geodesics $\bar  \upgamma_{\substack{\\[1pt]\hspace{-2pt}q_{_i}}}$, $i= 1,2$.
\end{proof}
\par Now we want to construct a test plan which represents the gradient of $f$. We are going to proceed along the same lines as in~\textsection\thinspace\ref{subsec:gradu}. Define
\begin{align*}
	\Uplambda: \X\rightarrow \mathcal{C}\left(\left[ \shortminus \updelta, \updelta \right],\X\right), \ \ \Uplambda\left(x\right)=  \upsigma_{\substack{\\x}} \ \mbox{ where }  \upsigma_{\substack{\\x}} \left(t\right)=\F_{\substack{\\[1pt]\hspace{-3pt}t}} \left(x\right). 
\end{align*}
and set $\Uplambda_{\, *}\ \upmu=: \ppi_{\substack{\\[1pt]\rm U_{_\updelta}}}\in \Prob\left(\mathcal{C} \left( \left(\shortminus \updelta, \updelta\right), \X \right) \right)$ with $\upmu$ having the specific properties as in the above discussion. Clearly, $\left(\e_{\substack{\\t}}\right)_{*} \ppi_{\substack{\\[1pt]\rm U_{_\updelta}}} = \upmu_{\substack{\\[1pt]\hspace{1pt}t}}$. Now, similar to what we did in \textsection\thinspace\ref{subsec:gradu} -- with $f$ in place of $u$ -- or directly from \cite[Proposition 4.4 and Theorem 4.5]{cavmonlap}, one gets the following first variation formula. 
\begin{proposition}[Test plan representatives of $\nabla f$]\label{prop:gradf} $\ppi_{\substack{\\[1pt]\rm U_{_\updelta}}}$ is a test plan that represents $\nabla f$. In particular, for any $\upvarphi\in \Sobol\left(\X\right)$,
		\begin{align*}
			\int_{\hspace{-1pt}_\X} \langle \nabla f, \nabla \upvarphi\rangle \,\dif \upmu_{\substack{\\[1pt]\hspace{1pt}t}}= \lim_{h\downarrow 0}\nicefrac{1}{h}\int_{\hspace{-1pt}_\X} \left( \upvarphi\left(\upxi_{\substack{\\[0.5pt]\,t+h}}\right) - \upvarphi\left(\upxi_{\substack{\\\,t}}\right)\right) \, \dif\ppi_{\substack{\\[1pt]\rm U_{_\updelta}}},
		\end{align*}
		holds by Theorem \ref{th:firstvariation}.
\end{proposition}
\par Set
\begin{align*}
\mathcal F\left(\upmu\right):=\int_{\hspace{-1pt}_\X} f \, \dif\upmu, \ \ \upmu\in  \Prob_{\hspace{-5pt}_2}\left(\X,\meas\right).
\end{align*}
\begin{lemma}\label{lem:2prop} Let $t,s\in \left[\shortminus \updelta, \updelta\right]$ and $\left( \F_{\substack{\\[1pt]\hspace{-3pt}t}} \right)_{\hspace{-1pt}*}\upmu= \upmu_{\substack{\\[1pt]\hspace{1pt}t}}$ as in above.
Then, the following hold
\smallskip
\begin{enumerate}
	\item semi-group isomorphism: {$\mathcal F\left(  \upmu_{\substack{\\[1pt]\hspace{1pt}t}} \right)=\mathcal F\left(\upmu\right)+t$};
	\medskip
	\item  contraction property: $\Was_{\hspace{-1pt}_2}\left(  \upmu_{\substack{\\[1pt]\hspace{1pt}s}},  \upmu_{\substack{\\[1pt]\hspace{1pt}t}}\right)\leq  \left|s \shortminus t\right|$.
\end{enumerate}
\smallskip
In particular, $ \left[\shortminus \updelta,\updelta\right] \ni t\mapsto \upmu_{\substack{\\[1pt]\hspace{1pt}t}}$ is $\Was_{\hspace{-1pt}_2}$-absolutely continuous (we actually know it is a geodesic).
\end{lemma}
\begin{proof}[\footnotesize \textbf{Proof}]
The semi-group isomorphism follows from the computation
\begin{align*}
	\nicefrac{d}{dt} \restr_{t=t_0} \int_{\hspace{-1pt}_\X} f \, \dif\upmu_{\substack{\\[1pt]\hspace{1pt}t}} &= \lim_{h\rightarrow 0} \nicefrac{1}{h} \left( \int_{\hspace{-1pt}_\X}  f \, \dif\upmu_{\substack{\\[1pt]\hspace{1pt}t_{_0}+h}} -  \int_{\hspace{-1pt}_\X} f \, \dif\upmu_{\substack{\\[1pt]\hspace{1pt}t_{_0}}}   \right) \\ &=  \lim_{h \rightarrow 0} \nicefrac{1}{h} \left( \int_{\hspace{-1pt}_\X}  f \, \dif \left( \F_{\substack{\\[1pt]\hspace{-3pt}t_{_0}+h}} \right)_{\hspace{-1pt} *}\upmu -  \int_{\hspace{-1pt}_\X} f \, \dif \left( \F_{\substack{\\[1pt]\hspace{-3pt}t_{_0}}} \right)_{\hspace{-1pt}*}\upmu  \right) \\ &= \lim_{h \rightarrow 0} \nicefrac{1}{h}  \int_{\hspace{-1pt}_\X} \left( f \circ \F_{\substack{\\[1pt]\hspace{-3pt}t_{_0}+h}}\left(x\right) -  f\circ \F_{\substack{\\[1pt]\hspace{-3pt}t_{_0}}} \left(x\right)  \right) \, \dif \upmu \\ &= \lim_{h \rightarrow 0} \nicefrac{1}{h} \int_{\hspace{-1pt}_\X} h \, \dif \upmu = 1.
\end{align*}
The contraction property is obvious since for the dynamic plan $\ppi_{\substack{\\[1pt]U_{_\updelta}}}$, we have
\begin{align*}
 \int_{\hspace{-1pt}_\X} \dist^{^2}  \,   \dif \left( \e_{\substack{\\[1pt]\hspace{0.5pt}s}} \,, \e_{\substack{\\[1pt]\hspace{0.5pt}t}} \right)_{\!*}\ppi_{\substack{\\[1pt]U_{_\updelta}}} = \int_{\hspace{-1pt}_\X}  \dist\left( \F_{\substack{\\[1pt]\hspace{-3pt}s}} , \F_{\substack{\\[1pt]\hspace{-3pt}t}} \right)^{\hspace{-2pt}^2} \, \dif \upmu = \left| s \shortminus t \right|^{^2}.
\end{align*}
\end{proof}
\subsubsection*{\small \bf \emph{Verifying the EVI formulation}}
In order to obtain a localized gradient flow, we need to determine some suitable domains on which the  EVI formulation of the gradient flow will be verified. Let $\bar \updelta>0$ be fixed and sufficiently small; let $\updelta\in \left(0,\bar \updelta\right)$. Pick $x_{\hspace{-0.5pt}_0} \in \rm U_{\substack{\\3\updelta}}$ with $\ball{\updelta}\left( x_{\hspace{-0.5pt}_0} \right)\subset \rm U_{\substack{\\3\updelta}}$ and  consider $\upmu\in \Prob\left(\X,\meas\right)$ such that $\upmu\left(\ball{\updelta}\left( x_{\hspace{-0.5pt}_0} \right)\right)=1$. Since $f$ is a $1$-Lipschitz function (and as a result, it locally does not increase distances), one verifies that $t\in \left( \shortminus \updelta, \updelta \right)\mapsto \upmu_{\substack{\\[1pt]\hspace{1pt}t}}$ is concentrated in $\ball{2\updelta}\left( x_{\hspace{-0.5pt}_0} \right)\subset \rm U_{\substack{\\2\updelta}}$. For any $\upnu\in \Prob\left(\X,\meas\right)$ with $\upnu\left(\ball{2\updelta}\left( x_{\hspace{-0.5pt}_0} \right)\right)=1$,  let $\left(\upnu^s_{\substack{\\\!t}}\right)_{s\in \left[0,1\right]}$ be the Wasserstein geodesic between $\upnu$ and $\upmu_{\substack{\\[1pt]\hspace{1pt}t}}$. Then $\upnu^s_{\substack{\\\!t}}$ is concentrated in $\ball{4\updelta}\left( x_{\hspace{-0.5pt}_0} \right)\subset \rm U_{\substack{\\\updelta}} \subset \X \smallsetminus \Sing$ for all $s\in \left[0,1\right]$ and $t\in \left(\shortminus \updelta, \updelta\right)$. 
\begin{proposition}[Evolutional variational {equality} for the gradient flow]\label{prop:evi} 
Let $\upnu$ and 
\begin{align*}
\left(\shortminus\updelta, \updelta \right) \ni t \mapsto \left( \F_{\substack{\\[1pt]\hspace{-3pt}t}} \right)_{\hspace{-1pt}*}\upmu = \upmu_{\substack{\\[1pt]\hspace{1pt}t}},
\end{align*}
be as before. Then,
\begin{align}\label{id:evi}
	\nicefrac{\dif}{\dif t} \, \nicefrac{1}{2}\, \Was_{\hspace{-1pt}_2}\left(\upmu_{\substack{\\[1pt]\hspace{1pt}t}},\upnu\right)^{\hspace{-2pt}^2}= \mathcal F\left(\upnu\right) - \mathcal F\left(\upmu_{\substack{\\[1pt]\hspace{1pt}t}}\right)  \quad \mbox{ for $\mathscr{L}^{^1}$-a.e. }  t\in \left(\shortminus\updelta, \updelta\right) .
\end{align}
\end{proposition}
\begin{proof}[\footnotesize \textbf{Proof}]
Pick a pair of Kantorovich potentials $\Upphi_{\substack{\\\!t}}, \Uppsi_{\substack{\\\!t}}$ satisfying
\begin{align*}
	\Upphi_{\substack{\\\!t}}\left(x\right)+\Uppsi_{\substack{\\\!t}}\left(y\right)\leq {\nicefrac{1}{2}}\, \dist\left(x,y\right)^{\hspace{-1pt}^2} \quad \forall x,y\in \X,
\end{align*}
and
\begin{align*}
	\nicefrac{1}{2}\, \Was_{\hspace{-1pt}_2} \left(\upmu_{\substack{\\[1pt]\hspace{1pt}t}},\upnu\right)^{\hspace{-2pt}^2}= \int_{\hspace{-1pt}_\X} \Upphi_{\substack{\\\!t}} \, \dif \upmu_{\substack{\\[1pt]\hspace{1pt}t}} + \int_{\hspace{-1pt}_\X} \Uppsi_{\substack{\\\!t}} \, \dif\upnu \quad \forall t \in \left(\shortminus \updelta, \updelta\right).
\end{align*}	
Notice that this can be done based on the, by now, standard developments in optimal transport theory; e.g. see~\cite[Theorem 7.35]{V}.
Since $\upmu_{\substack{\\[1pt]\hspace{1pt}t}} $ and $\upnu$ have bounded supports, we can assume $\Upphi_{\substack{\\\!t}} $ is Lipschitz continuous (see~\cite[Lemma 1]{mccannpolar}) and consequently, $\Upphi_{\substack{\\\!t}} \in \Sobol\left(\X\right)$. 	
By Kantorovich's dual formulation of Wasserstein distance (also called Kantorovich duality), it readily follows
\begin{align*}
	\nicefrac{1}{2}\, \Was_{\hspace{-1pt}_2} \left(\upmu_{\substack{\\[1pt]\hspace{1pt}t+s}},\upnu\right)^{\hspace{-2pt}^2}\geq \int_{\hspace{-1pt}_\X} \Upphi_{\substack{\\\!t}} \, \dif\upmu_{\substack{\\[1pt]\hspace{1pt}t+s}} + \int_{\hspace{-1pt}_\X} \Uppsi_{\substack{\\\!t}} \, \dif\upnu \quad \mbox{ for } \quad s+t> \shortminus\updelta, \quad  s>0.
\end{align*}
Since $\left(\shortminus\updelta,\updelta\right) \ni t \mapsto\upmu_{\substack{\\[1pt]\hspace{1pt}t}} $ is $\Was_{\hspace{-1pt}_2}$-absolutely continuous, $ \left(\shortminus\updelta,\updelta\right) \ni t \mapsto \Was_{\hspace{-1pt}_2}\left(\upmu_{\substack{\\[1pt]\hspace{1pt}t}} ,\upnu\right)^{\hspace{-2pt}^2}$ is differentiable $\mathscr L^{^1}$-a.e. So, \emph{at a point of differentiability} $t\in \left(\shortminus\updelta,\updelta\right)$, we deduce
\begin{align*}
	\nicefrac{\dif}{\dif t} \, \nicefrac{1}{2} \, \Was_{\hspace{-1pt}_2} \left( \upmu_{\substack{\\[1pt]\hspace{1pt}t}} ,\upnu\right)^{\hspace{-2pt}^2}&=\lim_{\varepsilon\rightarrow 0} \, \nicefrac{1}{2\varepsilon}\, \left(\Was_{\hspace{-1pt}_2} \left( \upmu_{\substack{\\[1pt]\hspace{1pt}t+\varepsilon}} ,\upnu\right)^{\hspace{-2pt}^2} -  \Was_{\hspace{-1pt}_2} \left( \upmu_{\substack{\\[1pt]\hspace{1pt}t}} , \upnu\right)^{\hspace{-2pt}^2} \right)\\
	& \geq \uplim_{\varepsilon\downarrow 0}\, \nicefrac{1}{\varepsilon} \, \left(\int_{\hspace{-1pt}_\X} \Upphi_{\substack{\\\!t}} \, \dif \upmu_{\substack{\\[1pt]\hspace{1pt}t+\varepsilon}} - \int_{\hspace{-1pt}_\X} \Upphi_{\substack{\\\!t}} \, \dif \upmu_{\substack{\\[1pt]\hspace{1pt}t}}\right) \\&= \uplim_{\varepsilon\downarrow 0} \, \nicefrac{1}{\varepsilon}\int_{\hspace{-1pt}_\X} \Upphi_{\substack{\\\!t}} \circ  \F_{\substack{\\[1pt]\hspace{-3pt}t+\varepsilon}}  \left(x\right) - \Upphi_{\substack{\\\!t}} \circ  \F_{\substack{\\[1pt]\hspace{-3pt}t}}  \left(x\right)  \, \dif\upmu.
\end{align*}
Now, it follows from Proposition~\ref{prop:gradf} and  Fatou's Lemma that
\begin{align}\label{ineq:Aa}
	\nicefrac{\dif}{\dif t}\, \nicefrac{1}{2}\, \Was_{\hspace{-1pt}_2} \left( \upmu_{\substack{\\[1pt]\hspace{1pt}t}} ,\upnu\right)^{\hspace{-2pt}^2} &\geq \int_{\hspace{-1pt}_\X}  \uplim_{\varepsilon\downarrow 0}\nicefrac{1}{\varepsilon} \; \Upphi_{\substack{\\\!t}} \left( \F_{\substack{\\[1pt]\hspace{-3pt}t+\varepsilon}}  \left(x\right) \right) -   \Upphi_{\substack{\\\!t}} \left( \F_{\substack{\\[1pt]\hspace{-3pt}t}}  \left(x\right) \right) \, \dif\upmu \notag \\ &\ge \int_{\hspace{-1pt}_\X}	\nicefrac{\dif}{\dif\uptau} \ \Upphi_{\substack{\\\!t}} \circ \F_{\substack{\\[1pt]\hspace{-3pt}\uptau}} \restr_{\uptau=t}\dif\upmu \\&= \int_{\hspace{-1pt}_\X}  \langle \nabla \Upphi_{\substack{\\\!t}}, \nabla f\rangle \circ \F_{\substack{\\[1pt]\hspace{-3pt}t}}\ \dif\upmu  \notag\\&= \int_{\hspace{-1pt}_\X} \langle \nabla \Upphi_{\substack{\\\!t}} , \nabla f\rangle\, \dif\upmu_{\substack{\\[1pt]\hspace{1pt}t}} \notag,
\end{align}
where the passage from the first to second line is due to the fact that $\Upphi_{\substack{\\\!t}}$ and $\F_{\substack{\\[1pt]\hspace{-3pt}t}}$ are both Lipschitz so $\uplim$ (after it has moved inside the integral) is replaced with $\lim$. 
\par Recall $\left(\upnu^s_{\substack{\\\!t}}\right)_{s\in[0,1]}$ is the $\Was_{\hspace{-1pt}_2}$-geodesic between $\upmu_{\substack{\\[1pt]\hspace{1pt}t}}$ and $\upnu$. The first variation formula (see Remark~\ref{rem:metbre}) implies
\begin{align}\label{ineq:Bb}
	\nicefrac{\dif}{\dif s}\restr_{s=0} \int_{\hspace{-1pt}_\X} f\,\dif \upnu^s_{\substack{\\\!t}} = \shortminus \int_{\hspace{-1pt}_\X} \langle \nabla f,\nabla \Upphi_{\substack{\\\!t}} \rangle\, \dif\upmu_{\substack{\\[1pt]\hspace{1pt}t}}.
\end{align}
On the other hand, along $\left( \upnu^s_{\substack{\\\!t}} \right)_{s\in [0,1]}$, we have the linearity
\begin{align*}
	\int_{\hspace{-1pt}_\X} f \,\dif \upnu^s_{\substack{\\\!t}}= \left(1 \shortminus s\right)\int_{\hspace{-1pt}_\X} f \, \dif\upnu^{^0}_{\substack{\\\!t}} + s\int_{\hspace{-1pt}_\X} f \, \dif\upnu^{^1}_{\substack{\\\!t}},
\end{align*}
therefore, upon combining \eqref{ineq:Aa} and \eqref{ineq:Bb} we obtain 
	\begin{align*}
		\nicefrac{\dif}{\dif t}\, \nicefrac{1}{2}\, \Was_{\hspace{-1pt}_2} \left( \upmu_{\substack{\\[1pt]\hspace{1pt}t}} ,\upnu\right)^{\hspace{-2pt}^2} \ge \mathcal F\left(\upnu\right) - \mathcal F\left( \upmu_{\substack{\\[1pt]\hspace{1pt}t}} \right)  \mbox{ for $\mathscr{L}^{^1}$-a.e. }  t\in \left(\shortminus \updelta, \updelta\right).
\end{align*}
The other direction can be shown by writing
\begin{align*}
	\nicefrac{\dif}{\dif t}\, \nicefrac{1}{2}\, \Was_{\hspace{-1pt}_2}\left( \upmu_{\substack{\\[1pt]\hspace{1pt}t}}  ,\upnu\right)^{\hspace{-2pt}^2}&=\lim_{\varepsilon\rightarrow 0} \nicefrac{1}{2s}\, \left(\Was_{\hspace{-1pt}_2}\left( \upmu_{\substack{\\[1pt]\hspace{1pt}t}}  ,\upnu\right)^{\hspace{-2pt}^2} -  \Was_{\hspace{-1pt}_2} \left(\upmu_{\substack{\\[1pt]\hspace{1pt}t-\varepsilon}} , \upnu\right)^{\hspace{-2pt}^2}\right),
\end{align*}
and arguing in a similar fashion. 
Consequently,
\begin{align}\label{ineq:Abb}
	\nicefrac{\dif}{\dif t}\, \nicefrac{1}{2}\, \Was_{\hspace{-1pt}_2} \left( \upmu_{\substack{\\[1pt]\hspace{1pt}t}} ,\upnu\right)^{\hspace{-2pt}^2} = \int_{\hspace{-1pt}_\X} \langle \nabla \Upphi_{\substack{\\\!t}} , \nabla f\rangle\, \dif\upmu_{\substack{\\[1pt]\hspace{1pt}t}},
\end{align}
holds at any point $t\in \left( \shortminus \updelta, \updelta \right)$ at which, $
\nicefrac{1}{2}\, \Was_{\hspace{-1pt}_2} \left( \upmu_{\substack{\\[1pt]\hspace{1pt}t}} ,\upnu\right)^{\hspace{-2pt}^2}$ is differentiable.
\end{proof}
\begin{remark}\label{rem:add} A posteriori, by bootstrapping, the equality \eqref{id:evi} in fact holds for every $t\in \left( \shortminus\updelta, \updelta \right)$; this is due to the fact that, $t\in \left( \shortminus\updelta, \updelta \right)\mapsto \Was_{\hspace{-1pt}_2} \left(\upmu_{\substack{\\[1pt]\hspace{1pt}t}},\upnu\right)^{\hspace{-2pt}^2}$ is locally Lipschitz and since its almost everywhere derivative, $t\in \left( \shortminus\updelta, \updelta \right)\mapsto \mathcal F(\upmu_{\substack{\\[1pt]\hspace{1pt}t}}) $, is continuous (compare with Lemma \ref{lem:2prop}). Hence, $ \Was_{\hspace{-1pt}_2} \left(\upmu_{\substack{\\[1pt]\hspace{1pt}t}},\upnu\right)^{\hspace{-2pt}^2}$ is differentiable for every $t\in \left( \shortminus\updelta, \updelta \right)$ and \eqref{ineq:Abb} holds for every $t\in \left( \shortminus \updelta, \updelta \right)$.
\end{remark}
\subsection{Distance-preserving property of the flow}
In order to establish the distance-preserving property of the flow $\F_{\substack{\\[1pt]\hspace{-3pt}t}}$ (or rather an extension thereof), we start with showing the Wasserstein-distance preserving property. 
\begin{proposition}[$\Was_{\hspace{-1pt}_2}$-distance-preserving property of the flow]\label{Prop:dist-preserving}
Let $\upmu,\upnu\in \Prob_{\hspace{-5pt}_2} \left(\X,\meas\right)$ be concentrated in $\ball{\updelta}\left( x_{\hspace{-0.5pt}_0}  \right)$ and let $\upmu_{\substack{\\[1pt]\hspace{1pt}t}}$ and $\upnu_{\substack{\\[1pt]\hspace{-1pt}t}}$ denote their pushforward under $\F_{\substack{\\[1pt]\hspace{-3pt}t}}$ respectively. Then
\begin{align*}
\Was_{\hspace{-1pt}_2} \left(\upmu,\upnu\right)= \Was_{\hspace{-1pt}_2} \left(\upmu_{\substack{\\[1pt]\hspace{1pt}t}} ,\upnu_{\substack{\\[1pt]\hspace{-1pt}t}} \right) \quad \forall t\in \left[ \shortminus \updelta,\updelta \right].
\end{align*}
\end{proposition}
\begin{proof}[\footnotesize \textbf{Proof}]
By Lemma \ref{lem:2prop} and applying the triangle inequality, we deduce $t\mapsto\Was_{\hspace{-1pt}_2} \left(\upmu_{\substack{\\[1pt]\hspace{1pt}t}} ,\upnu_{\substack{\\[1pt]\hspace{-1pt}t}}\right)$ is $2$-Lipschitz. To show what is claimed, it suffices to show that the derivative of the map
\begin{align*}
\left(\shortminus \updelta,\updelta\right) \ni t \mapsto \Was_{\hspace{-1pt}_2} \left(\upmu_{\substack{\\[1pt]\hspace{1pt}t}} ,\upnu_{\substack{\\[1pt]\hspace{-1pt}t}}\right)^{\hspace{-2pt}^2},
\end{align*}
vanishes almost everywhere. Fix $t\in \left(\shortminus \updelta, \updelta\right)$ and let $\upomega_{\substack{\\t}} \in \Prob\left(\X, \meas\right)$ be a midpoint of $\upmu_{\substack{\\[1pt]\hspace{1pt}t}}$ and $\upnu_{\substack{\\[1pt]\hspace{-1pt}t}}$. $\upomega_{\substack{\\t}}$ is supported in $\ball{2\updelta}\left(x_{\hspace{-0.5pt}_0} \right)$. By Young's Inequality and the Proposition~\ref{prop:evi} along with a computation similar to the one in ~\cite[Lemma 4.8]{GKK}, we get 
\begin{align*}
	&\uplim_{\varepsilon\downarrow 0} \,\nicefrac{1}{2\varepsilon} \left( \Was_{\hspace{-1pt}_2} \left(\upmu_{\substack{\\[1pt]\hspace{1pt}t+ \varepsilon}} ,\upnu_{\substack{\\[1pt]\hspace{-1pt}t+ \varepsilon}}\right)^{\hspace{-2pt}^2} -  \Was_{\hspace{-1pt}_2} \left(\upmu_{\substack{\\[1pt]\hspace{1pt}t}} ,\upnu_{\substack{\\[1pt]\hspace{-1pt}t}}\right)^{\hspace{-2pt}^2} \right) \le 0,
\end{align*}
where in the proof, one uses the affinity of $\upmu\mapsto \int_{\hspace{-1pt}_\X} {f} \dif\upmu$ along $\Was_{\hspace{-1pt}_2}$-geodesic in $\rm U_{\substack{\\\updelta}}$.
\par Similarly it follows that 
\begin{align*}
	\lowlim_{\varepsilon\downarrow0} \, \nicefrac{1}{2\varepsilon}\, \left( \ \Was_{\hspace{-1pt}_2} \left(\upmu_{\substack{\\[1pt]\hspace{1pt}t+ \varepsilon}} ,\upnu_{\substack{\\[1pt]\hspace{-1pt}t+ \varepsilon}}\right)^{\hspace{-2pt}^2} - \Was_{\hspace{-1pt}_2} \left(\upmu_{\substack{\\[1pt]\hspace{1pt}t}} ,\upnu_{\substack{\\[1pt]\hspace{-1pt}t}}\right)^{\hspace{-2pt}^2}  \right) \geq  0.
\end{align*}
and one deduces the right derivative is zero. The left derivative is also zero using similar computations. Hence, $\nicefrac{\dif}{\dif t}\, \Was_{\hspace{-1pt}_2} \left( \upmu_{\substack{\\[1pt]\hspace{1pt}t}} ,\upnu_{\substack{\\[1pt]\hspace{-1pt}t}} \right)^{\hspace{-2pt}^2}  =0$ whenever $t\in \left( \shortminus \updelta,\updelta\right)$ is a point of differentiability. 
\end{proof}
\begin{remark}
The $\Was_{\hspace{-1pt}_2}$-distance-preserving property and measure-preserving property along the flow can also be obtained by using the heat flow and the fact that $f$ is Hessian-free, similar to the proof of distance-preserving property of the gradient flow of the Busemann function; see~\cite{Gsplit}.  
\end{remark}
\begin{proposition}[Local isometries]\label{prop:local-sio}
There exists a map $\widetilde \F^\updelta: \left[ \shortminus \updelta,\updelta\right]\times \rm U_{\substack{\\3\updelta}} \rightarrow \X$ satisfying the following properties. 
\begin{enumerate}
	\item $\widetilde \F^\updelta$ almost everywhere coincides with $\F$ namely,
	\begin{align*}
	\meas \left( \left\{x\in \rm U_{\substack{\\3\updelta}}  \; \text{\emph{\textbrokenbar}} \; \exists t\in \left[ \shortminus \updelta,\updelta\right] \quad \F_{\substack{\\[1pt]\hspace{-3pt}t}}\left(x\right)\neq \widetilde{\F}^\updelta_{\substack{\\[1pt]\hspace{-3pt}t}}\left(x\right) \right\} \right)=0;
	\end{align*}
	\smallskip
	\item $\widetilde{\F}^\updelta_{\substack{\\[1pt]\hspace{-3pt}t}}\restr_{\B_{_\updelta}\left(x_{\hspace{-0.5pt}_0} \right)}$ is distance-preserving for $t\in \left[  \shortminus\updelta, \updelta \right]$ and $\ball{\updelta}\left( x_{\hspace{-0.5pt}_0} \right)\subset \rm U_{\substack{\\3\updelta}} $;
	\smallskip
	\item $\widetilde{\F}^{\updelta^{'}}_{\substack{\\[1pt]\hspace{-3pt}t}}= \widetilde{\F}^\updelta_{\substack{\\[1pt]\hspace{-3pt}t}} \restr_{\rm U_{3\updelta^{'}}}$ for $0<\updelta<\updelta^{'}$ and $t\in \left[ \shortminus \updelta, \updelta\right]$.
\end{enumerate}
Note that by item (2),  the map $\widetilde{\F}^\updelta_{\substack{\\[1pt]\hspace{-3pt}t}}$ is continuous on $\rm U_{\substack{\\3\updelta}}$.
\end{proposition}
\begin{proof}[\footnotesize \textbf{Proof}]
First, we fix $t\in \left[ \shortminus \updelta, \updelta\right]\cap \mathbb Q$. Consider $\upvarphi \in \mathcal{C}_{\substack{\\\textsf{bs}}}\left(\X\right)$. By the Lebesgue differentiation theorem, there exists $\mathrm{A}_{\upvarphi,t}\subset \ball{\updelta}\left( x_{\hspace{-0.5pt}_0} \right)$ of full $\meas$-measure in $\ball{\updelta}\left( x_{\hspace{-0.5pt}_0} \right)$ such that for every $z\in \mathrm{A}_{\upvarphi,t}$ and for the uniform measures $ \upmu^{\upeta}:= \nicefrac{1}{\meas\left(\B_{_\upeta}\left(z\right)\right)} \, \meas\mres_{\B_{_\upeta}\left(z\right)}$ and their pushforwards $\upmu^\upeta_{\substack{\\t}} := \left( \F_{\substack{\\[1pt]\hspace{-3pt}t}} \right)_{\hspace{-1pt}*} \; \upmu^\upeta$, it holds
	\begin{enumerate}
		\item $\displaystyle\int_{\hspace{-1pt}_\X} \upvarphi \, \dif\upmu^\upeta$  converges to $\displaystyle \upvarphi\left(z\right)= \int_{\hspace{-1pt}_\X} \upvarphi \, \dif\updelta_z$ as $\upeta \to 0$;
		\smallskip
		\item  the integral
		\begin{align*}
		\int_{\hspace{-1pt}_\X} \upvarphi \, \dif \upmu^\upeta_{\substack{\\t}} = \int_{\hspace{-1pt}_\X} \upvarphi\circ \F_{\substack{\\[1pt]\hspace{-3pt}t}} \; \dif\upmu^\upeta =  \nicefrac{1}{\meas\left(\B_{_\upeta}\left(z\right)\right)} \int_{\hspace{-1pt}_\X} \upvarphi \circ \F_{\substack{\\[1pt]\hspace{-3pt}t}} \; \dif\meas\mres_{\B_{_\upeta}\left(z\right)},
		\end{align*}
		converges to 
		$ \displaystyle
		\upvarphi\circ \F_{\substack{\\[1pt]\hspace{-3pt}t}} \left(z\right)=  \int_{\hspace{-1pt}_\X} \upvarphi\, \dif\updelta_{\substack{\\[2pt]\hspace{-2pt}\F_{\substack{\\\hspace{-3pt}t}} \left(z\right)}}
		$
		as $\upeta \rightarrow 0$.
\end{enumerate}
\par Let $\Upxi\subset \mathcal{C}_{\substack{\\\textsf{bs}}}\left(\X\right)$ be a countable subset that is dense in $\mathcal{C}_{\substack{\\\textsf{bs}}}\left(\X\right)$ w.r.t. uniform convergence (recall $\mathcal{C}_{\substack{\\\textsf{bs}}}\left(\X\right)$ is separable; e.g. see \cite{Boga}). Set $\mathrm{A}_{\substack{\\t}}:=\bigcap_{\upvarphi\in \Upxi} \mathrm{A}_{\upvarphi, t}$; notice $\mathrm{A}_{\substack{\\t}}$ has full $\meas$-measure in $\ball{\updelta}\left( x_{\hspace{-0.5pt}_0} \right)$.  By items (1) and (2), for every $z\in \mathrm{A}_{\substack{\\t}}$,
it follows that $\upmu^\upeta \rightharpoonup \updelta_z$ and  $\upmu^\upeta_t \rightharpoonup \updelta_{\substack{\\\F_{\substack{\\\hspace{-3pt}t}}\left(z\right)}}$ (weakly) as $\upeta\rightarrow 0$.
\par Now for $z_{_1},z_{_2}\in \mathrm{A}_{\substack{\\t}}$, the continuity of $\Was_{\hspace{-1pt}_2}$ w.r..t weak convergence and the $\Was_{\hspace{-1pt}_2}$-distance-preserving property of the flow, in combination with Proposition~\ref{Prop:dist-preserving}, yield
	\begin{align*}
		\dist\left( \F_{\substack{\\\hspace{-3pt}t}}\left( z_{_1} \right), \F_{\substack{\\\hspace{-3pt}t}}\left( z_{_2} \right) \right)=\dist\left( z_{_1} , z_{_2} \right).
\end{align*}
\par By continuity of $t\mapsto \F_{\substack{\\\hspace{-3pt}t}}\left(z\right)$ for $z\in \bigcap_{t\in \left[-\updelta,\updelta\right]\cap \mathbb Q}\mathrm{A}_{\substack{\\t}} $, the previous identity holds for every $t\in \left(\shortminus \updelta,\updelta\right)$ and for all $z_{_1} , z_{_2} \in \mathrm{A} := \bigcap_{t\in \left( \shortminus \updelta, \updelta\right)\cap \mathbb{Q}} \mathrm{A}_{\substack{\\t}}$; thus, $\F_{\substack{\\\hspace{-3pt}t}}$ is an isometry on $\mathrm{A}$ for every $t\in \left[ \shortminus \updelta, \updelta\right]$. 
\par Since $\mathrm{A}$ is a set of full measure in $\ovs{\ball{\updelta}\left( x_{\hspace{-0.5pt}_0} \right)}$, there exists a unique isometric extension $\widetilde{\F}^\updelta_{\substack{\\[1pt]\hspace{-3pt}t}}$ of $\F_{\substack{\\\hspace{-3pt}t}}$ on $\ovs{\ball{\updelta}\left( x_{\hspace{-0.5pt}_0} \right)}$ for every $t\in \left[ \shortminus \updelta,\updelta \right]$. Since the extension $\widetilde \F$ of $\F$ is unique on $\ball\updelta\left( x_{\hspace{-0.5pt}_0} \right)$ whenever $\ball{\updelta}\left( x_{\hspace{-0.5pt}_0} \right)\subset \rm U_{\substack{\\\updelta}}$, we see that there exists a well-defined extension $\widetilde{\F}^\updelta_{\substack{\\[1pt]\hspace{-3pt}t}}$ of $\F_{\substack{\\[1pt]\hspace{-3pt}t}}$ on $\rm U_{\substack{\\3\updelta}}$ for every $t\in \left[ \shortminus \updelta, \updelta \right]$
such that the properties  {(1) - (3)} are  satisfied. 
\end{proof}
\begin{remark}\label{rem:F-map}
Given $\widetilde{\F}_{\hspace{-4pt}s}^\updelta$, $s\in [\shortminus \updelta, \updelta]$ (as in Proposition~\ref{prop:local-sio}) and for $\updelta>0$ sufficiently small,  we consider $t\in \left[ \shortminus \nicefrac{\uppi}{2}, \nicefrac{\uppi}{2} \right] $. Then, there exists $k\in \mathbb N$ and $\upeta\in \left(0,\updelta\right)$ such that $\mbox{sgn}\left(t\right) t= k\updelta + \upeta$. 
We define
\begin{align*}
	\rm U\strut^{\updelta, t}:=  \rm U_{\substack{\\3\updelta}} \cap \left(\widetilde \F^{\updelta}_{\hspace{-2pt}\sgn\left(t\right)\updelta}\right)^{\hspace{-2pt}^{-1}} \circ \dots \circ \left(\widetilde \F^{\updelta}_{\hspace{-2pt}\sgn\left(t\right)\updelta}\right)^{\hspace{-2pt}^{-1}} \circ \left(\widetilde{\F} ^\updelta_{\hspace{-2pt}\sgn\left(t\right)\upeta}\right)^{\hspace{-2pt}^{-1}} \,\left(\X \smallsetminus \Sing\right),
	\end{align*}
and the map $\widetilde{\F}^\updelta_{\substack{\\[1pt]\hspace{-3pt}t}}:	\rm U\strut^{\updelta, t} \rightarrow \X $ by
\begin{align*}
\widetilde{\F}^\updelta_{\substack{\\[1pt]\hspace{-3pt}t}}= \widetilde \F_{\hspace{-2pt}\sgn\left(t\right)\updelta}^{\updelta} \circ\dots\circ \widetilde \F_{\hspace{-2pt}\sgn\left(t\right) \updelta}^{\updelta} \circ \widetilde \F_{\hspace{-2pt}\sgn\left(t\right)\upeta}^{\updelta};
\end{align*}
the map $\widetilde{\F}^\updelta_{\substack{\\[1pt]\hspace{-3pt}t}}$ is continuous. It follows from the construction of the maps $\widetilde{\F}^\updelta_{\substack{\\[1pt]\hspace{-3pt}t}}$, that
\begin{enumerate}
	\item $\meas\left( \left\{x\in \rm U\strut^{\updelta, t}  \; \text{\textbrokenbar} \;  \exists t \in \left[ \shortminus \nicefrac{\uppi}{2}, \nicefrac{\uppi}{2} \right]  \quad  \F_{\substack{\\[1pt]\hspace{-3pt}t}}\left(x\right)\neq \widetilde{\F}^\updelta_{\substack{\\[1pt]\hspace{-3pt}t}}\left(x\right) \right\}\right)=0$;\smallskip
	\item $\widetilde{\F}^\updelta_{\substack{\\[1pt]\hspace{-3pt}t}}\restr_{\B_{_\updelta}\left( x_{\hspace{-0.5pt}_0} \right)}$ is distance-preserving for $t\in \left[ \shortminus \nicefrac{\uppi}{2}, \nicefrac{\uppi}{2} \right]$ provided $\ball{\updelta}\left( x_{\hspace{-0.5pt}_0} \right)\subset	\rm U\strut^{\updelta, t}$; \smallskip
	\item $	\rm U\strut^{\updelta^{'}, t} \subset 	\rm U\strut^{\updelta, t}$ and $\widetilde{\F}^{\updelta^{'}}_{\substack{\\[1pt]\hspace{-3pt}t}}= \widetilde{\F}^\updelta_{\substack{\\[1pt]\hspace{-3pt}t}}\restr_{U^{\updelta^{'}, t}}$ if $\updelta^{'}\in \left(0, \updelta\right)$; \smallskip
	\item For any $x\in \X \smallsetminus \Sing$, there exists $\updelta>0$ small enough such that $x \in \rm U_{\substack{\\3\updelta}}$.
\end{enumerate}
\end{remark}
\section{Isometric Splitting}\label{sec:isom-split}
\subsection{The Projection Map}
Set $\Y := u^{{-1}}\left(0\right)$. $\Y$ is the cross section, the isometric propagation of which along the gradient flow of $f$, we wish to show comprises the whole $\X$. For this, define the projection map, $\mathrm{proj}_{\substack{\\[2pt]\Y}}: \X \smallsetminus \Sing \rightarrow \Y$, by
\begin{align*}
	\mathrm{proj}_{\substack{\\[2pt]\Y}}\left(x\right)= \widetilde \F^{\updelta}_{\hspace{-4pt} \shortminus f\left(x\right) }\left(x\right) \quad \mbox{ provided } \quad  x\in {\rm U}\strut^{\updelta, \shortminus f\left(x\right)}.
\end{align*}
By item (3) in Remark~\ref{rem:F-map}, one deduces that the projection map onto $\Y$, is well-defined and does not depend on $\updelta$. 
\begin{lemma}[Local $1$-Lipschitz property of the projection map]\label{Lem:loc-lip}
	Let $x_{\hspace{-0.5pt}_0} \in {\rm U}_{\substack{\\3\updelta}}$. There exists $\upeta>0$ such that 
	$ \mathrm{proj}_{\substack{\\[2pt]\Y}}\restr_{\ovs{\B_{_\upeta}\left(x_{\hspace{-0.5pt}_0} \right)}}: \ovs {\ball{\upeta}\left( x_{\hspace{-0.5pt}_0} \right)} \rightarrow \Y$ is $1$-Lipschitz.
\end{lemma}
\begin{proof}[\footnotesize \textbf{Proof}]
We need to argue
		\begin{align*}
		\dist\left(  \mathrm{proj}_{\substack{\\[2pt]\Y}} \left(y\right),x\right)  \le  \dist\left(y,x\right).
	\end{align*}
\par	Pick $\upeta\in \left( 0,\nicefrac{\updelta}{6} \right)$ with
	\begin{align*}
		\ball{\upeta}\left( x_{_0} \right)\subset {\rm U}\strut^{\updelta, \shortminus f\left( x_{_0}\right)}\cap f^{^{-1}}\big(\left(f\left(x_{_0}\right) \shortminus \nicefrac{\updelta}{4}, f\left(x_{_0}\right)+\nicefrac{\updelta}{4}\right)\big).
	\end{align*}
	One can chose such $\upeta >0$ since ${\rm U}\strut^{\updelta, \shortminus f\left(x_{_0}\right)}$ and $f^{^{-1}}\big(\left(f\left(x_{_0}\right) \shortminus \nicefrac{\updelta}{4}, f\left(x_{_0}\right)+\nicefrac{\updelta}{4}\right)\big)$ are open and both contain $x_{_0}$.  Let us point out that $\ball{\upeta}\left(x_{_0}\right)\subset {\rm U}^{\updelta, \shortminus f\left(y\right)}$ for all $y\in \ball{\upeta}\left(x_{_0}\right)$. Moreover, since $0<\upeta<\nicefrac{\updelta}{6}<\updelta$, we know that $\widetilde \F_{\hspace{-4pt}\shortminus f\left(x\right)}^{\updelta}\restr_{{\B_{_\upeta}\left(x_{_0}\right)}}$ is distance-preserving.  
\par	Now let $x,y \in \ball{\upeta}\left(x_{_0}\right)$. Without loss of generality, we can assume $ \mathrm{proj}_{\substack{\\[2pt]\Y}}\left(x\right)=x$;  indeed, replace $\ball{\upeta}\left(x_{_0}\right)$ with $\ball{\upeta}\left(\widetilde \F_{\hspace{-4pt} \shortminus f\left(x\right)}^\updelta \left(x_{_0}\right)\right)$ if necessary and take into account the  distance-preserving property of $\widetilde \F_{\hspace{-4pt}\shortminus f\left(x\right)}^{\updelta}\restr_{{\B_{_\upeta}\left(x_{_0}\right)}}$. 
\par Consider a $\upmu\in \Prob\left(\X, \meas\right)$ which is concentrated in $\ball{\upeta}\left(x_{_0}\right)$ along with its pushforward, 
	$ \upmu_{\substack{\\[1pt]\hspace{1pt}t}} = \left( \F_{\substack{\\[1pt]\hspace{-3pt}t}} \right)_{\hspace{-1pt}*}\upmu$. Clearly, $\Was_{\hspace{-1pt}_2} \left(\upmu, \updelta_{\substack{\\\hspace{-1pt}x}}\right)\leq 2\upeta$ holds. For $\left|t\right|\geq 4\upeta$, it also holds
	\begin{align*}
	\Was_{\hspace{-1pt}_2} \left( \upmu_{\substack{\\[1pt]\hspace{1pt}t}} , \updelta_{\substack{\\\hspace{-1pt}x}}\right)^{\hspace{-2pt}^2} &=  \int_{\hspace{-1pt}_\X} \dist\left( \F_{\substack{\\[1pt]\hspace{-3pt}t}} \left(z\right),x\right)^{\hspace{-2pt}^2} \, \dif\upmu\left(z\right) \notag \\ &\geq \int_{\hspace{-1pt}_\X} \Big( \dist\left( \F_{\substack{\\[1pt]\hspace{-3pt}t}} \left(z\right),z\right) \shortminus \dist\left(z,x\right)\Big)^{\hspace{-4pt}^2}\, \dif\upmu\left(z\right) \\ & \geq \left(\left|t\right| \shortminus 2\upeta\right)^{\hspace{-2pt}^2} \geq 4\upeta^2 \notag,
	\end{align*}
	where the first equality is a result of the fact that for a.e. $y \in \supp\left( \upmu_{\substack{\\[1pt]\hspace{1pt}t}}  \right)$, there exists a unique geodesic from $x$ to $y$ and these geodesics provide the optimal transport map from $\updelta_{\substack{\\[1pt]\hspace{-0.5pt}x}}$ to $\upmu_{\substack{\\[1pt]\hspace{1pt}t}} $. 
\par Hence, $\Was_{\hspace{-1pt}_2} \left(\upmu_{\substack{\\[1pt]\hspace{1pt}t}} ,\updelta_{\substack{\\\hspace{-1pt}x}}\right)\geq 2\upeta$ if $\left|t\right|\geq 4\upeta$.  Recall, $\upeta<\nicefrac{\updelta}{6}$ and hence $4\upeta\leq \nicefrac{2\updelta}{3}<\updelta$. Since $ \left(\shortminus \updelta, \updelta\right) \ni t \rightarrow \Was_{\hspace{-1pt}_2} \left( \upmu_{\substack{\\[1pt]\hspace{1pt}t}}  ,\updelta_{\substack{\\\hspace{-1pt}x}}\right)$ is continuous, it follows that 
	$
	\left[ \shortminus 4\upeta, 4\upeta\right] \ni  t \mapsto \Was_{\hspace{-1pt}_2} \left(\upmu_{\substack{\\[1pt]\hspace{1pt}t}} ,\updelta_{\substack{\\\hspace{-1pt}x}}\right)
	$
	attains a minimum for some $t_{_0}\in \left( \shortminus 4\upeta, 4\upeta\right)$. 
	Moreover, we compute
	\begin{align*} 
		\dist\left( \F_{\substack{\\[1pt]\hspace{-3pt}t_{_0}}} \left(z\right), x\right)\leq \dist\left( \F_{\substack{\\[1pt]\hspace{-3pt}t_{_0}}} \left(z\right),z\right)+ \dist\left(z,x\right)\leq 4\upeta + 2\upeta=6\upeta \leq  \updelta, \quad \forall z\in \B_{2\upeta}\left(x\right).
	\end{align*}
	As a result, $\upmu_{\substack{\\t_{_0}}}\left(\ball{\updelta}\left(x\right)\right)=1$.
\par	Let $ [0,1] \ni s \mapsto \upnu_{\substack{\\[1pt]\hspace{-2pt}s}}$ be the geodesic from $\upmu_{\substack{\\t_{_0}}}$ to $\updelta_{\substack{\\\hspace{-1pt}x}}$. Then, for every $s\in [0,1]$ and $t\in \R$ with $t+t_{_0}\in \left[ \shortminus4\upeta, 4\upeta \right]$, it follows 
	\begin{align*}
		\Was_{\hspace{-1pt}_2} \left( \upnu_{\substack{\\[1pt]\hspace{-2pt}s}} ,\updelta_{\substack{\\\hspace{-1pt}x}} \right)&= \left(1 \shortminus s\right) \Was_{\hspace{-1pt}_2} \left( \upmu_{\substack{\\t_{_0}}} ,\updelta_{\substack{\\\hspace{-1pt}x}}\right)\\
		&\leq \left(1 \shortminus s\right) \Was_{\hspace{-1pt}_2} \left(\upmu_{\substack{\\t + {t_{_0}}}},\updelta_{\substack{\\\hspace{-1pt}x}}\right) \\
		&\leq \left(1 \shortminus s\right) \Big( \Was_{\hspace{-1pt}_2} \left(\upmu_{\substack{\\t + {t_{_0}}}}, \left( \F_{\substack{\\[1pt]\hspace{-3pt}t}} \right)_{\hspace{-1pt}*}\upnu_{\substack{\\[1pt]\hspace{-2pt}s}} \right)+ \Was_{\hspace{-1pt}_2} \left(\left( \F_{\substack{\\[1pt]\hspace{-3pt}t}} \right)_{\hspace{-1pt}*}\upnu_{\substack{\\[1pt]\hspace{-2pt}s}}, \updelta_{\substack{\\\hspace{-1pt}x}}\right)\Big)\\
		&\leq \left(1 \shortminus s\right) \Was_{\hspace{-1pt}_2} \left( \upmu_{\substack{\\t_{_0}}}, \upnu_{\substack{\\[1pt]\hspace{-2pt}s}}\right) + \left(1 \shortminus s\right) \Was_{\hspace{-1pt}_2} \left( \left( \F_{\substack{\\[1pt]\hspace{-3pt}t}} \right)_{\hspace{-1pt}*}\upnu_{\substack{\\[1pt]\hspace{-2pt}s}}, \updelta_{\substack{\\\hspace{-1pt}x}}\right)\\
		&= s \Was_{\hspace{-1pt}_2} \left(  \upnu_{\substack{\\[1pt]\hspace{-2pt}s}} ,\updelta_{\substack{\\\hspace{-1pt}x}}\right)+ \left(1 \shortminus s\right) \Was_{\hspace{-1pt}_2} \left( \left( \F_{\substack{\\[1pt]\hspace{-3pt}t}} \right)_{\hspace{-1pt}*} \upnu_{\substack{\\[1pt]\hspace{-2pt}s}}, \updelta_{\substack{\\\hspace{-1pt}x}}\right).
	\end{align*}
	Therefore, 
	\begin{align*}
	\Was_{\hspace{-1pt}_2} \left( \upnu_{\substack{\\[1pt]\hspace{-2pt}s}},\updelta_{\substack{\\\hspace{-1pt}x}}\right) \le \Was_{\hspace{-1pt}_2} \left( \left( \F_{\substack{\\[1pt]\hspace{-3pt}t}} \right)_{\hspace{-1pt}*} \upnu_{\substack{\\[1pt]\hspace{-2pt}s}} , \updelta_{\substack{\\\hspace{-1pt}x}}\right) \quad \forall t \in \left[ \shortminus 4\updelta \shortminus t_{_0}, 4\updelta \shortminus t_{_0}  \right];
	\end{align*}
	i.e. $\Was_{\hspace{-1pt}_2} \left( \left( \F_{\substack{\\[1pt]\hspace{-3pt}t}} \right)_{\hspace{-1pt}*} \upnu_{\substack{\\[1pt]\hspace{-2pt}s}} , \updelta_{\substack{\\\hspace{-1pt}x}}\right)$ attains a local minimum at $t=0$. 
		\par Since $\upmu_{\substack{\\t_{_0}}}$ is concentrated in $\ball{\updelta}\left(x\right)$, we deduce $\upnu_{\substack{\\[1pt]\hspace{-2pt}s}}$ is concentrated in $\ball{2\updelta}\left(x\right)$. In particular, from Remark \ref{rem:add} we infer $t\mapsto \Was_{\hspace{-1pt}_2} \left( \upmu_{\substack{\\t}}, \updelta_{\substack{\\\hspace{-1pt}x}}\right)^{\hspace{-1pt}^2}$ is differentiable at $t_{_0}$ and $t\mapsto \Was_{\hspace{-1pt}_2} \left( \left( \F_{\substack{\\[1pt]\hspace{-3pt}t}} \right)_{\hspace{-1pt}*} \upnu_{\substack{\\[1pt]\hspace{-2pt}s}}, \updelta_{\substack{\\\hspace{-1pt}x}}\right)^{\hspace{-2pt}^2} $ is differentiable at $t=0$. Since we have extrema, the derivatives vanish. So, by an analogue of \eqref{ineq:Abb} we obtain
		\begin{align*}
			0&= \nicefrac{\dif}{\dif t}\restr_{t=0}\, \Was_{\hspace{-1pt}_2} \left( \left( \F_{\substack{\\[1pt]\hspace{-3pt}t}} \right)_{\hspace{-1pt}*} \upnu_{\substack{\\[1pt]\hspace{-2pt}s}}, \updelta_{\substack{\\\hspace{-1pt}x}}\right)^{\hspace{-2pt}^2} \\ &=  \nicefrac{\dif}{\dif t}\restr_{t=0}  \  \int \dist\left( \F_{\substack{\\[1pt]\hspace{-3pt}t}} \left(z\right),x\right)^{\hspace{-2pt}^2}\, \dif\upnu_{\substack{\\[1pt]\hspace{-2pt}s}}\left(z\right) \\ &= \int_{\hspace{-1pt}_\X} \langle \nabla f, \nabla \, \dist_{x}^{^2}\rangle \, \dif \upnu_{\substack{\\[1pt]\hspace{-2pt}s}}.
		\end{align*}
	Notice here, we have used the simple formula for the Wasserstein distance when one of the measures is a delta measure. 
	\par	With the first variation formula (Theorem \ref{th:firstvariation} together with the subsequent Remark \ref{rem:metbre}) this yields
		\begin{align*}
			\lim_{h\rightarrow 0} \, \nicefrac{1}{h} \, \Big(\int_{\hspace{-1pt}_\X} f \, \dist \upnu_{\substack{\\[1pt]\hspace{-2pt}s+h}} - \int_{\hspace{-1pt}_\X} f \,\dif \upnu_{\substack{\\[1pt]\hspace{-2pt}s}} \Big)=\nicefrac{1}{\left(1-s\right)}\int_{\hspace{-1pt}_\X} \langle \nabla f ,\nabla \, \dist_{x}^{^2} \rangle \, \dif\upnu_{\substack{\\[1pt]\hspace{-2pt}s}}=0,
	\end{align*}
(here, $\dist_{x}$ is the distance to $x$) which in turn, by using the definition of $\F_{\substack{\\[1pt]\hspace{-3pt}t_{_0}}}$, yields
	\begin{align*}
		t_{_0} + 	\int_{\hspace{-1pt}_\X} f \, \dif\upmu =\int_{\hspace{-1pt}_\X} f\circ \F_{\substack{\\[1pt]\hspace{-3pt}t_{_0}}} \, \dif\upmu=\int_{\hspace{-1pt}_\X} f \, \dif\upmu_{\substack{\\t_{_0}}}=\lim_{s\uparrow 1} \int_{\hspace{-1pt}_\X} f\, \dif\upnu_{\substack{\\[1pt]\hspace{-2pt}s}} = f\left(x\right) = 0.
	\end{align*}
	Therefore, one has $t_{_0}\rightarrow \shortminus f\left(y\right)$ as $\upmu \rightharpoonup \updelta_{\substack{\\\hspace{-1pt}y}}$. 
	Consider the chain of inequalities
	\begin{align}\label{eq:chain-ineq}
		\Was_{\hspace{-1pt}_2} \left(\upmu,\updelta_{\substack{\\\hspace{-1pt}x}}\right)^{\hspace{-2pt}^2} &\geq \Was_{\hspace{-1pt}_2} \left( \upmu_{\substack{\\t_{_0}}} ,\updelta_{\substack{\\\hspace{-1pt}x}}\right)^{\hspace{-2pt}^2} \notag \\&=  \int_{\hspace{-1pt}_\X} \dist\left(\widetilde \F^{\updelta}_{t_{_0}}\left(z\right), x\right)^{\hspace{-2pt}^2} \, \dif\upmu\\
		& \geq \int_{\hspace{-1pt}_\X} \Big(  \dist\left( \widetilde \F_{\hspace{-3pt} t_{_0}}^{\updelta}\left(z\right), x \right) \shortminus \dist\left(  \widetilde \F_{\hspace{-3pt} t_{_0}}^{\updelta} \left(z\right), \widetilde \F_{\hspace{-4pt} \shortminus f(y) }^{\updelta}\left(z\right)\right)\Big)^{\hspace{-2pt}^2}\, \dif\upmu. \notag
	\end{align}
	As $\upmu \rightharpoonup \updelta_{\substack{\\\hspace{-1pt}y}}$, the first term in \eqref{eq:chain-ineq} converges to $\dist\left(y,x\right)^{\hspace{-2pt}^2}$ and the last term in \eqref{eq:chain-ineq} converges to 
	\begin{align*}
		\dist\left( \widetilde \F_{\hspace{-4pt} \shortminus f(y) }^\updelta \left(y\right),x\right)^{\hspace{-2pt}^2} = \dist\left( \mathrm{proj}_{\substack{\\[2pt]\Y}} \left(y\right),x\right)^{\hspace{-2pt}^2}.
	\end{align*}
	Hence,
	\begin{align*}
		\dist\left(  \mathrm{proj}_{\substack{\\[2pt]\Y}} \left(y\right),x\right)  \le  \dist\left(y,x\right),
	\end{align*}
	which is what we had claimed. 
\end{proof}	
\subsection{The spaces $\widetilde{\X}$ and $\widetilde{\Y}$:}\label{sec:defn-tilde-spaces} The subset $\Y$ equipped with the distance $\dist$ is locally geodesically convex. Indeed, the Corollary \ref{cor:loccon} implies that, for every $y\in \Y$, there exists $\updelta>0$ with the property that, for every pair $x,z\in \ball{\updelta}\left(y\right)$, one can find a geodesic $\upgamma:\left[0,1\right]\rightarrow \X$ between $x$ and $z$ along which
\begin{align*}
	f\left(\upgamma(t)\right) = \left(1 \shortminus t\right) f\left(x\right) + t f\left(z\right),
\end{align*}
holds. Choosing $x,z\in \ball{\updelta}\left(x_{_0}\right)\cap \Y$ then yields $f\left(\upgamma\left(t\right)\right)=0$ for all $t\in \left(0,1\right)$.
\par  Let $\left(\widetilde \X,  \dist_{\substack{\\\widetilde \X}}\right)=: \widetilde \X$ be the completion of the \emph{extended intrinsic distance}, $\widetilde{\dist}_{\substack{\\[1pt]\X \smallsetminus \Sing}}$, on $\X \smallsetminus  \Sing$ (see Definition~\ref{defn:sing-set} for $\Sing$) where $\widetilde{\dist}_{\substack{\\[1pt]\X \smallsetminus \Sing}}$ is the extended distance which is induced by $\dist_{\substack{\\\X}}$. Here, extended means we set the distance between two points in different path connected components, to be infinity. Similarly, let $\widetilde \dist_{\substack{\\\Y}}$ be the extended induced intrinsic distance of $\Y$. We set $\left(\Y, \widetilde \dist_{\substack{\\\Y}}\right)=:\widetilde \Y$. Note that the topology of $\widetilde \X$ may be different from the one of $\X$. In particular, $\widetilde \X$ and $\widetilde \Y$ may have several path connected components.
\par The projection map $\mathrm{proj}_{\substack{\\[2pt]\Y}}: \X \smallsetminus \Sing \rightarrow \Y$ is continuous hence it can be closed to get $\mathrm{proj}_{\substack{\\[2pt]\widetilde \Y}}: \widetilde \X \rightarrow \widetilde \Y$. 
\begin{corollary}
	The projection  $\mathrm{proj}_{\substack{\\[2pt]\widetilde \Y}}: \widetilde \X \rightarrow \widetilde \Y$ is $1$-Lipschitz (with the convention $\infty \le 1 \cdot \infty$).
\end{corollary}
\begin{proof}[\footnotesize \textbf{Proof}]
	This is straightforward from Lemma~\ref{Lem:loc-lip}. 
\end{proof}
\begin{proposition}\label{prop:geocon}
	Let $x_{_0}\in \Y$ and let $\updelta >0$ be as in~\textsection\thinspace\ref{subsec:gf-f}. For every pair $\upmu_{\hspace{1pt}_0}, \upmu_{\hspace{1pt}_1}\in  \Prob\left(\Y,{\bf q}\right)$ with $\upmu_{\substack{\\[1pt]\hspace{1pt}i}}\left(\ball{\nicefrac{\updelta}{2}}\left(x_{_0}\right)\right)=1$, there exists a \emph{unique geodesic} $\upmu_{\substack{\\t}}\in \Prob\left(\X\right)$ with $\upmu_{\substack{\\t}}\left(\Y\right)=1$. Furthermore, $\upmu_{\substack{\\t}}\in \Prob\left(\Y,{\bf q}\right)$; recall by Remark~\ref{rem:q-on-Y}, $\Q$ can be almost identified with $\Y$.
\end{proposition}
\begin{proof}[\footnotesize \textbf{Proof}] 
	Let $\ray: \Y\times \left[ \shortminus \nicefrac{\uppi}{2},\nicefrac{\uppi}{2} \right]\rightarrow \X$ be the ray map and define the map ${\Upupsilon}\strut^\upeta$ to be
	\begin{align*} \upmu\in \Prob\left(Y\right)\mapsto
		\Upupsilon\strut^\upeta \left(\upmu\right) := \ray_{\substack{\\\hspace{-2pt}{*}}}\left(\upmu\otimes \nicefrac{1}{2\eta}\,\mathscr L^{^1}\mres_{\left[-\upeta,\upeta\right]}\right), \quad \upeta > 0.
	\end{align*}
\par	Then, the measures $\upmu_{\substack{\\i}}^\upeta= \Upupsilon\strut^\upeta\left( \upmu_{\substack{\\i}} \right)\in \Prob\left(\X,\meas\right), i=0,1$, are concentrated in 
\[
\left\{ x\in \X \; \text{\textbrokenbar} \; \shortminus \upeta \leq f\left(x\right) \leq \upeta \right\}.
\]
Thus, for $\upeta>0$ sufficiently small, the measures $\upmu_{\substack{\\i}}^\upeta$ are concentrated in 
$\ball{\updelta}\left(x_{_0}\right)$.
\par	Now, as a result of Corollary \ref{cor:loccon}, there exist \emph{unique geodesics} $\left(\upmu_{\substack{\\t}}^\upeta\right)_{t\in \left[0,1\right]}$ in $\Prob\left(\X,\meas\right)$, concentrated in 
$
	\left\{x\in \X \; \text{\textbrokenbar} \; \shortminus \upeta \leq f\left(x\right) \leq \upeta\right\},
$
 joining $\upmu_{\hspace{1pt}_0}^\upeta$ to $\upmu_{\hspace{1pt}_1}^\upeta$ in the $2$-Wasserstein space. 
\par	By the standard compactness results for Wasserstein geodesics, as $\upeta \to 0$ and after passing to a subsequence, the $\Was_{\hspace{-1pt}_2}$-geodesics $\left( \upmu_{\substack{\\t}}^\upeta \right)_{t\in \left[0,1\right]}$ weakly converges to a $\Was_{\hspace{-1pt}_2}$-geodesic $\upmu_{\substack{\\t}}$ between $\upmu_{\hspace{1pt}_0}$ and $\upmu_{\hspace{1pt}_1}$. Any such geodesic $\upmu_{\substack{\\t}}$ is concentrated in $\Y= \left\{x\in \X \; \text{\textbrokenbar} \; f\left(x\right) =0\right\}$ as a result of the weak convergence of measures.
\par	To show uniqueness, let  $\upnu_{\substack{\\t}}$ be any given geodesic joining $\upmu_{\substack{\\i}}$'s, $i=0,1$, and with $\upnu_{\substack{\\t}}\left(\Y\right)=1$. Since $\F_{\substack{\\[1pt]\hspace{-3pt}t}} : \ball{\updelta}\left(x_{_0}\right)\rightarrow \X$ for $t\in \left[ \shortminus \updelta, \updelta\right]$ is distance-preserving, the map $\Upupsilon\strut^\upeta$ is distance-preserving (compare with the proof of~\cite[Corollary 5.30]{Gsplit}). Hence, $\Upupsilon\strut^\upeta\left( \upnu_{\substack{\\t}} \right)=: \upnu_{\substack{\\t}}^\upeta$ is a $t$-midpoint  between $\upmu_{\hspace{1pt}_0}^\upeta$ and $\upmu_{\hspace{1pt}_1}^\upeta$. Since $\upmu_{\hspace{1pt}_0}^\upeta$ and $\upmu_{\hspace{1pt}_1}^\upeta$ are $\meas$-absolutely continuous,  $ \upnu_{\substack{\\t}}^\upeta$ coincides with the unique $t$-midpoint between  $\upmu_{\hspace{1pt}_0}^\upeta$ and $\upmu_{\hspace{1pt}_1}^\upeta$. But  $\upnu_{\substack{\\t}}^\upeta$ converges back to $\upnu_{\substack{\\t}}$ as $\upeta \to 0 $. Since $\upnu_{\substack{\\t}}^\upeta$ is unique, it must coincide with $\upmu_{\substack{\\t}}^\upeta$ from before and consequently $\upnu_{\substack{\\t}}=\upmu_{\substack{\\t}}$.
\par	Moreover,  $\upnu_{\substack{\\t}}^\upeta=: \ray_{\substack{\\\hspace{-2pt}{*}}}\left(\upnu_t\otimes \nicefrac{1}{2\eta}\,\mathscr L^{^1}\mres_{\left[-\upeta,\upeta\right]}\right)$ is $\meas$-absolutely continuous because $X$ is $\RCD$ and since $\ray$ is a measure space isomorphism  and $\upnu_t$ is ${\bf q}$-absolutely continuous. 
\end{proof}
\subsection{The splitting map}\label{subsec:split}
Define the splitting map $ {\fancy{$\mathcal{S}$}}$ and its a.e. inverse, ${\fancy{$\mathcal T$}}$, by
\begin{align*}
	&{\fancy{$\mathcal{S}$}}: \X\rightarrow \Y\times \left[ \shortminus \nicefrac{\uppi}{2}, \nicefrac{\uppi}{2} \right], \quad  {\fancy{$\mathcal{S}$}}\left(x\right) = \left(  \mathrm{proj}_{\substack{\\[2pt]\Y}}\left(x\right), f\left(x\right)\right),\\
	&{\fancy{$\mathcal T$}}: Y\times \left[ \shortminus \nicefrac{\uppi}{2}, \nicefrac{\uppi}{2} \right] \rightarrow \X, \quad  {\fancy{$\mathcal T$}}\left(y, t\right) = \widetilde{\F}_{\substack{\\\hspace{-3pt} t}}\left(y\right).
\end{align*}
The space $\Y\times \left[ \shortminus \nicefrac{\uppi}{2}, \nicefrac{\uppi}{2} \right]$ is equipped with the measure ${\bf q}\otimes \mathscr L^{^1}\mres_{\left[ \shortminus \nicefrac{\uppi}{2}, \nicefrac{\uppi}{2} \right]}$ where ${\bf q}$ is concentrated in $\Q\underset{\mathrm{almost}}{\subset} Y$.
\par By construction, we have 
$ {\fancy{$\mathcal{S}$}} \circ {\fancy{$\mathcal T$}}\left(y,t\right)= \left(y,t\right)$  for ${\bf q}\otimes \mathscr L^{^1}\mres_{\left[ \shortminus \nicefrac{\uppi}{2}, \nicefrac{\uppi}{2} \right]}$-a.e.  $\left(y,t\right)$ and $ {\fancy{$\mathcal{T}$}} \circ {\fancy{$\mathcal{S}$}} \left(x\right)=x$  for $\meas$-a.e. $x$. 
\par Based on the disintegration formula \eqref{eq:disintegration2}, for sets of the form $S \times \left[ a, b  \right] \subset \Y \times \left[ \shortminus \nicefrac{\uppi}{2}, \nicefrac{\uppi}{2} \right]$, we clearly get 
\begin{align*}
	\meas \left(  {\fancy{$\mathcal T$}}  \left( S \times \left[ a, b  \right] \right) \right) &=\int_{_{{\fancy{$\mathcal T$}} \left( S \times \left[ a, b  \right] \right)}} \left(\upgamma_{\substack{\\[1pt]\hspace{-2pt}q}}\right)_{\hspace{-1pt}*} \mathscr L^{^1}\mres_{\left[\shortminus \nicefrac{\uppi}{2},\nicefrac{\uppi}{2}\right]} \, \dif{\bf q}\left(q\right) \\ &= \int_{\hspace{-1pt}_S}   \left(\upgamma_{\substack{\\[1pt]\hspace{-2pt}q}}\right)_{\hspace{-1pt}*} \mathscr L^{^1}\mres_{\left[\shortminus \nicefrac{\uppi}{2},\nicefrac{\uppi}{2}\right]}\left(   {\fancy{$\mathcal T$}} \big( \{q\} \times \left[ a, b  \right] \big)  \right)       \dif{\bf q}\left(q\right) \\ &= \int_{\hspace{-1pt}_S}    (b \shortminus a)     \dif{\bf q}\left(q\right) \\ &= (b \shortminus a) \,{\bf q}\left( S \right)  \\ &= {\bf q}\otimes \mathscr L^{^1}\mres_{\left[ \shortminus \nicefrac{\uppi}{2}, \nicefrac{\uppi}{2} \right]} \big( S \times \left[ a,b \right]   \big).
\end{align*}
So by Carath\'eodory's extension theorem, the splitting map ${\fancy{$\mathcal T$}}$ is measure-preserving; so is its inverse (in a.e. sense) $ {\fancy{$\mathcal{S}$}}$.
\begin{proposition}\label{prop:lipmap}
	For any pair $\left(y,s\right), \left(z,t\right)\in \Y\times \left[ \shortminus \nicefrac{\uppi}{2}, \nicefrac{\uppi}{2} \right]$, we have 
	\begin{align}\label{bilip}
	2^{\shortminus \nicefrac{1}{2}}\, \uphat \dist\left(\left(y,s\right), \left(z,t\right)\right) \leq \widetilde \dist_{\substack{\\\X}}\left( {\fancy{$\mathcal{T}$}} \left(y,s\right), {\fancy{$\mathcal{T}$}} \left(z,t\right)\right)\leq 2^{\nicefrac{1}{2}}\ \uphat \dist\left(\left(y,s\right),\left(z,t\right)\right),
	\end{align}
	where $\uphat \dist\left(\left(y,s\right), \left(z,t\right)\right)= \Big( \widetilde \dist_{\substack{\\ \Y}}\left(y,z\right)^{\hspace{-2pt}^2} + \left|t \shortminus s\right|^{^2} \Big)^{\hspace{-4pt}^{\nicefrac{1}{2}}}$ is the isometric product of distances. 
\end{proposition}
\begin{proof}[\footnotesize \textbf{Proof}]
	Based on the definition of $\widetilde \dist_{\substack{\\\X}}$ and $\widetilde \dist_{\substack{\\\Y}}$, it is enough to show \eqref{bilip} on balls $\ball{\updelta}\left(x_{_0}\right)\subset \X \smallsetminus \Sing$. Then, we can use the fact that $f$ is Lipschitz and $\mathrm{proj}_{\substack{\\[2pt]\Y}}$ is locally bi-Lipschitz and proceed exactly as in the proof of \cite[Lemma 5.2]{GKK}.
\end{proof}
\subsection{Local Sobolev space}
We recall the definition of the local Sobolev space $\Sobol_{\hspace{-1pt}\textsf{loc}}\left(\Omega\right)$ for an open subset $\Omega\subset \X$ (for instance see \cite{Gsplit}). For $v \in \Ltwo\left(\Omega\right)$, we say $v\in \Sobol_{\hspace{-1pt}\textsf{loc}}\left(\Omega\right)$ if $\upchi\cdot v\in \Sobol\left(\X\right)$ for any test Lipschitz function $\upchi: \X\rightarrow \R$ with $\supp \left(\upchi\right)\subset \Omega$. By locality of the minimal weak upper gradient, one can define $\left|\nabla v\right|_{\substack{\\[1pt]\Omega}}= \left|\nabla \left(\upchi\cdot v\right)\right|_X$ on $\upchi\equiv 1$. If $\left|\nabla v\right|_{\substack{\\[1pt]\Omega}}\in \Ltwo\left(\Omega, \meas\right)$ then we say $v\in \Sobol\left(\Omega\right)$ and we set the local Cheeger-Dirichlet energy to be 
\begin{align*}
	\CHE\strut^{\hspace{-10pt}\Omega}\,\left(v\right) :=\int_{\hspace{-1pt}_{\Omega}} \left|\nabla v\right|^{^2}_{\substack{\\[1pt]\Omega}} \, \dif\meas.
\end{align*}
One should note that in general, $\Sobol\left(\Omega\right)$ does not coincide with $\Sobol\left(\hspace{1pt}\ovs{\Omega}\hspace{1pt}\right)$; e.g. for the radially slit disk, these spaces are not the same. 
\par Let $\Omega =\X \smallsetminus \Sing$. 
We note that $\widetilde \dist_{\substack{\\ \X}}$ and $\dist=\dist_{\substack{\\ \X}}$ coincide on sufficiently small balls that are contained in $\X \smallsetminus \Sing$. In particular,  for $v\in \Sobol_{\hspace{-1pt}\textsf{loc}}\left(\X \smallsetminus \Sing\right)$ it holds that $\left|\nabla v\right|_{\substack{\\[1pt]\X \smallsetminus \Sing}} = \left|\nabla \left(\upchi \cdot v\right)\right|_{\substack{\\[1pt]\widetilde \X}}$ on $\upchi\equiv 1$ whenever $\supp\left(\upchi\right)\subset \X \smallsetminus \Sing$.
\begin{proposition}\label{prop:pullback-Sobolev}
	The Sobolev spaces $\Sobol_{\hspace{-1pt}\mathsf{loc}}\left(\X \smallsetminus \Sing\right)$ and $\Sobol\left(\Y\right)$ interact in the following manner.
	\begin{enumerate}
		\item For $\mathscr L^{^1}\mbox{-a.e.}  \ t\in \left( \shortminus \nicefrac{\uppi}{2}, \nicefrac{\uppi}{2} \right)$, it holds
		\begin{align*} 
			{\fancy{$\mathcal T$}}\left(\cdot, t\right)\strut^{\hspace{-2pt}*} \left(  \W_{\mathsf{loc}}^{^{1,2}}\left(\X \smallsetminus \Sing\right) \right) \subset \Sobol\left(\Y\right)  ;
		\end{align*} 
		\item One has
		\begin{align*}
			\mathrm{proj}_{\substack{\\[2pt]\Y}}\strut^{\hspace{-6pt}*} \left( \Sobol\left(\Y\right)  \right) \subset \W_{\mathsf{loc}}^{^{1,2}}\left(\X \smallsetminus \Sing\right),
		\end{align*}
		and furthermore, for any $\uphat{g} \in \Sobol\left(\Y\right)$, 
		\begin{align*}
			\left|\nabla \left( \uphat{g}  \circ \mathrm{proj}_{\substack{\\[2pt]\Y}}\right) \right|_{\substack{\\[1pt]\X \smallsetminus \Sing}} \left(x\right) = \left|\nabla\, \uphat{g} \right|_{\substack{\\[1pt]\Y}}\circ \mathrm{proj}_{\substack{\\[2pt]\Y}} \left(x\right),
		\end{align*}
		holds for $\meas\mbox{-a.e. } x\in \X$. 
	\end{enumerate}
\end{proposition}
\begin{proof}[\footnotesize \textbf{Proof}]
	The proof is along the same lines of the proof of \cite[Proposition 5.29]{Gsplit} so, it is omitted. 
\end{proof}
\begin{corollary}[$\RCD$ for the vertical fiber]
	Consider $\Y$ equipped with the induced intrinsic distance $\widetilde \dist_{\substack{\\\Y}}$ and the measure ${\bf q}$. Then every path connected component of $\left(\Y, \widetilde \dist_{\substack{\\\Y}}, {\bf q}\right)$ satisfies $\RCD\left(0,\N\right)$.
\end{corollary}
\begin{proof}[\footnotesize \textbf{Proof}]
	The  proof of the corollary is exactly as for in \cite[Corollary 5.30]{Gsplit} using the map $\Upupsilon\strut^\upeta$ that was previously introduced in Proposition~\ref{prop:geocon}. 
	\par Let us just check that $\left(\Y, \widetilde \dist_\Y, {\bf q}\right)$ is infinitesimally Hilbertian. We pick  ${\uphat{g}}_{_1}, {\uphat{g}}_{_2}\in \Sobol\left(\Y\right)$. By Proposition~\ref{prop:pullback-Sobolev}, we know that $w_{_1} := {\uphat{g}}_{_1} \circ \mathrm{proj}_{\substack{\\[2pt]\Y}}$ and $ w_{_2} :=  {\uphat{g}}_{_2} \circ \mathrm{proj}_{\substack{\\[2pt]\Y}}$ belong to $\Sobol\left(\X \smallsetminus \Sing\right)$ with
	\begin{align*}
		\left|\nabla {\uphat{g}}_{_1}  \right|_{\substack{\\[1pt]\Y}}^{^2}\left(y\right) = \left|\nabla w_{_1} \right|_{\substack{\\[1pt]\X \smallsetminus \Sing}}^{^2} \circ 	{\fancy{$\mathcal T$}}\left(y,t\right), \quad \left|\nabla {\uphat{g}}_{_2} \right|_{\substack{\\[1pt]\Y}}^{^2}\left(y\right) = \left|\nabla w_{_2} \right|_{\substack{\\[1pt]\X \smallsetminus \Sing}}^{^2} \circ 	{\fancy{$\mathcal T$}} \left(y,t\right),
	\end{align*}
	and
	\begin{align*}
		\left|\nabla \left({\uphat{g}}_{_1} \plmi \, {\uphat{g}}_{_2}\right)\right|_{\substack{\\[1pt]\Y}}^{^2}\left(y\right) = \left|\nabla \left( w_{_1} \plmi \, w_{_2} \right)\right|_{\substack{\\[1pt]\X \smallsetminus \Sing}}^{^2}\circ {\fancy{$\mathcal T$}} \left(y,t\right),
	\end{align*}	
	for ${\bf q}\otimes \mathscr L^{^1}\mbox{-a.e. } \left(y, t\right)\in \Y\times \left( \shortminus \nicefrac{\uppi}{2}, \nicefrac{\uppi}{2} \right)$; note that $w_{_1} \plmi \, w_{_2}= \left({\uphat{g}}_{_1}\plmi \, {\uphat{g}}_{_2}\right)\circ \mathrm{proj}_{\substack{\\[2pt]\Y}}$. 
	\par Since $\X$ is infinitesimally Hilbertian and by the definition of $\left| \nabla \cdot   \right|_{\substack{\\[1pt]\X \smallsetminus \Sing}}$, it holds 
	\begin{align*}
		\left|\nabla \left(w_{_1} \scalebox{0.7}{+}\, w_{_2}\right)\right|_{\substack{\\[1pt]\X \smallsetminus \Sing}}^{^2}\circ {\fancy{$\mathcal T$}} + \left|\nabla \left(w_{_1} \shortminus w_{_2}\right)\right|_{\substack{\\[1pt]\X \smallsetminus \Sing}}^{^2}\circ {\fancy{$\mathcal T$}} = 2\left|\nabla w_{_1}\right|_{\substack{\\[1pt]\X \smallsetminus \Sing}}^{^2} \circ {\fancy{$\mathcal T$}}+ 2\left|\nabla w_{_2}\right|_{\substack{\\[1pt]\X \smallsetminus \Sing}}^{^2}\circ {\fancy{$\mathcal T$}},
	\end{align*}
	for ${\bf q}\otimes \mathscr L^{^1}\mbox{-a.e.  on } \Y\times \left( \shortminus \nicefrac{\uppi}{2}, \nicefrac{\uppi}{2} \right)$. 
	Consequently, by the Fubini's theorem, we obtain similar polarization identities for ${\uphat{g}}_{_1}$ and $ {\uphat{g}}_{_2}$ in place of $w_{_1}$ and $w_{_2}$ (resp.). 
\end{proof}
\subsubsection*{\small \bf \textit{Modified Sobolev-to-Lipschitz property}}
\begin{proposition}[Modified Sobolev-to-Lipschitz property] \label{prop:modSobLip}
	Let $v\in \Sobol_{\hspace{-1pt}\mathsf{loc}}\left(\X \smallsetminus \Sing\right)$ with $\left|\nabla v\right|\in \LL^{\hspace{-2pt}^{\infty}}\left(\X \smallsetminus \Sing, \meas\right)$ and $\left|\nabla v\right|\leq  \sf L<\infty$. Then there exists $\widetilde v: \widetilde \X\rightarrow \R$ which is $\sf L$-Lipschitz w.r.t. $\widetilde \dist$ on each path connected component of $\widetilde \X$ so that $\widetilde v \restr_{\X \smallsetminus \Sing}$ is a representative for $v$. 
\end{proposition}
\begin{proof}[\footnotesize \textbf{Proof}] Let $v\in \Sobol_{\hspace{-1pt}\textsf{loc}}\left(\X \smallsetminus \Sing\right)$ and pick $x_{_0}\in \X \smallsetminus \Sing$ and $\updelta>0$ such that  $\ball{4\updelta}\left(x_{_0}\right)\subset \X \smallsetminus \Sing$. 
	\par {\bf \small Claim.} $v\restr_{\ovs{\B_\updelta\left(x_{_0}\right)}}$ is Lipschitz w.r.t. $\dist$.  
\par	{\it \small Proof of the claim.} Let $\upchi: \X\rightarrow \R$ a Lipschitz  function w.r.t. $\dist$ with $\supp  \left( \upchi\right)\subset \ball{4\updelta}\left( x_{_0} \right)$ and $\upchi\equiv 1$ on $\ball{2\updelta}\left(x_{_0}\right)$; i.e. $\upchi$ is a cutoff function. From the definition of $\W_{\textsf{loc}}^{1,2}\left(\Omega\right)$, it follows $\upchi \cdot u\in \Sobol\left(\X\right)$. 
\par	Pick $z_{_0}, z_{_1}\in \ball{\updelta}\left(x_{_0}\right)$ and measures $\upmu^n_{_0}, \upmu^n_{_1}  $ with $\upmu^n_{\substack{\\i}} \left(\B_{\updelta}\left(x_{_0}\right)\right)=1$  and $\upmu^n_{\substack{\\i}} \rightharpoonup \updelta_{z_{_i}}$ weakly, $i=0,1$. Let $\PP_{n}$ be the unique dynamical optimal plan between $\upmu^n_{_0}$ and $ \upmu^n_{_1}$. 
	Since $\upmu^n_{\substack{\\i}} \left( \ball{\updelta}\left(x_{_0}\right) \right)=1$, it follows $\left(\e_{\substack{\\t}}\right)_{*}\PP_n\left( \ball{2\updelta}\left(x_{_0}\right) \right)=1$ for all $t\in \left(0,1\right)$. 
	\par Since $\upchi\cdot v\in \Sobol\left(\X\right)$ and $\PP_n$ is a test plan, by applying Cauchy-Schwartz twice, it follows 
	\begin{align}\label{eq:ineq}
		&\int_{\hspace{-1pt}_{\mathsf{Curves}}}   \left|v\left(\upxi\left(1\right)\right)- v\left(\upxi\left(0\right)\right)\right| \, \dif\PP_n\left(\upxi\right) \notag \\ &\leq \int_{\hspace{-1pt}_{\mathsf{Curves}}}  \int_{_0}^1 \left|\nabla \left(\upchi \cdot v\right)\right|\left(\upxi\left(t\right)\right) \left|\bdot\upxi\right| \, \dif t \, \dif \PP_n\left(\upxi\right)  \\
		&\leq  \int_{\hspace{-1pt}_{\mathsf{Curves}}}  \left( \int_{_0}^1 \left|\nabla \left(\upchi \cdot v\right)\right|^{^2}\left(\upxi\left(t\right)\right) \dif t    \right)^{\nicefrac{1}{2}} \left( \int_{_0}^1  \left|\bdot\upxi\right|^{^2}  \dif t    \right)^{\nicefrac{1}{2}} \,  \dif \PP_n\left(\upxi\right) \notag \\
		&\leq \left(\int_{\hspace{-1pt}_{\mathsf{Curves}}}   \int_{_0}^1 \left|\nabla \left(\upchi \cdot v\right)\right|^{^2}\left(\upxi\left(t\right)\right) \,\dif t \,  \dif \PP_n \right)^{\nicefrac{1}{2}} \left( \int_{\hspace{-1pt}_{\mathsf{Curves}}}  \int_{_0}^1   \left|\bdot\upxi\right|^{^2}   \, \dif t   \,  \dif \PP_n \right)^{\nicefrac{1}{2}} \notag \\ &\leq {\sf L} \, \Was_{\hspace{-1pt}_2} \left( \upmu^n_{_0}, \upmu^n_{_1}\right). \notag
	\end{align}
	Now, upon letting $n \rightarrow \infty$, we get that the inequality in (\ref{eq:ineq}) (involving the first and the last term)  converges to $\left| v\left( z_{_0} \right) \shortminus v\left( z_{_1} \right) \right|\leq {\sf L}\,  \dist\left(z_{_0},z_{_1}\right)$.
	\phantom{123} \hfill\scalebox{0.6}{$\blacksquare$}
	\par Now let $z_{_0}, z_{_1}\in \X \smallsetminus \Sing$ be arbitrary in the same path connected component, $\uphat{\X}$ of $\widetilde \X$. Recall $\widetilde \dist$ is the induced intrinsic distance on $\widetilde \X$. Then, for $\varepsilon>0$, there exists a constant speed curve $\uphat \upgamma: [0, 1]\rightarrow \X \smallsetminus \Sing$ between $z_{_0}$ and $z_{_1}$ such that $\widetilde \dist\left(z_{_0}, z_{_1}\right) \geq {\length}\left(\uphat\upgamma\right) \shortminus \varepsilon$.  Then, we can find a partition $\mathbb{P}$ via
	\begin{align*}
		0=t_{_0}\leq t_{_1}\leq \dots \leq t_{_{k-1}}\leq t_{_k}=1 \quad  \text{and} \quad \varkappa_{\substack{\\\mathbb{P}}} := \sup_i \left| t_i \shortminus t_{i-1}   \right|,
	\end{align*}
	with $ \varkappa_{\substack{\\\mathbb{P}}}  $ sufficiently small such that 
	\begin{align*}
		\left| v\left(\uphat \upgamma\left(t_{i-1}\right) \shortminus v\left( \uphat \upgamma \left( t_i \right)\right)\right) \right| \leq {\sf L} \, \dist\left(\uphat{\upgamma}\left(t_{i-1}\right), \uphat{\upgamma}\left(t_i\right)\right) \quad  i=1, \dots, k.
	\end{align*}
	holds by local $\sf L$-Lipschitz property for $v$. 
	\par As a result,
	\begin{align*}
		\left| v\left(z_{_0}\right) \shortminus v\left(z_{_1}\right) \right|\leq {\sf L} \sum_{i=1}^k \dist\left( \uphat{\upgamma}\left(t_{i-1}\right), \uphat{\upgamma}\left(t_i\right) \right).
	\end{align*}
	After taking the supremum over all such partitions $\mathbb{P}$, we obtain 
	\begin{align*}
		\left| v\left(z_{_0}\right) \shortminus v\left(z_{_1}\right) \right|\leq {\sf L} \length\left(\uphat{\upgamma}\right)\leq {\sf L}(\widetilde \dist\left(z_{_0},z_{_1}\right)+ \varepsilon).
	\end{align*}
	Upon letting $\varepsilon\rightarrow 0$, we deduce $v\restr_{\uphat{\X} \smallsetminus \Sing}$ is $\sf L$-Lipschitz w.r.t. $\widetilde \dist$. Then, we can take $\widetilde{v}\restr_{\uphat{X}}$ to be the completion of $v\restr_{\uphat{\X} \smallsetminus \Sing}$ w.r.t. $\widetilde \dist$. 
\end{proof}
\subsubsection*{\small \bf \textit{Energy identities/estimates}}
A corollary of Proposition \ref{prop:lipmap} is the following estimates for the Sobolev norms. 
\begin{corollary}\label{cor:rough}
	It holds $k \in \Sobol(\Y\times  \left[ \shortminus \nicefrac{\uppi}{2}, \nicefrac{\uppi}{2} \right])$ if and only if $k\circ {\fancy{$\mathcal S$}}\in \Sobol_{\hspace{-1pt}\mathsf{loc}}(X\smallsetminus \Sing)$; in such a case, one has
	\begin{align}\label{rough_energy_estimate}
		2^{\shortminus \nicefrac{1}{2}}\, \CHE^{ \Y \times \left[ \shortminus \nicefrac{\uppi}{2}, \nicefrac{\uppi}{2} \right]}(k)\leq
		\CHE^{X\smallsetminus \Sing}(k \circ {\fancy{$\mathcal S$}})
		\leq 2^{\nicefrac{1}{2}} \, \CHE^{ \Y \times \left[ \shortminus \nicefrac{\uppi}{2}, \nicefrac{\uppi}{2} \right]}(k).
	\end{align}
\end{corollary}
\begin{proof}[\footnotesize \textbf{Proof}]
	We  only show that $k \in \Sobol(\Y\times  \left[ \shortminus \nicefrac{\uppi}{2}, \nicefrac{\uppi}{2} \right])$ implies $k\circ {\fancy{$\mathcal S$}}\in \Sobol_{\hspace{-1pt}\textsf{loc}}(\X\smallsetminus \Sing)$ and the second inequality in \eqref{rough_energy_estimate}.
\par	Let $k$ be Lipschitz and let $\upchi: \X\rightarrow \R$ be a Lipschitz function with support in $\X \smallsetminus \Sing$. Then, using Proposition \ref{prop:lipmap} it is not hard to check  that  we have 
	\begin{align}\label{ineq:pres}
	\left|\nabla (\upchi \cdot (k\circ {\fancy{$\mathcal S$}}))\right|_{\substack{\\[2pt]\X}}  \leq 2^{\nicefrac{1}{2}}\, k\circ {\fancy{$\mathcal S$}} \left| \nabla \upchi \right|_{\substack{\\[2pt]\X}} + 2^{\nicefrac{1}{2}}\, \left( \left| \nabla k\right|_{\substack{\\[2pt] \Y \times \left[ \shortminus \nicefrac{\uppi}{2}, \nicefrac{\uppi}{2} \right] }} \circ {\fancy{$\mathcal S$}} \right) \upchi.
	\end{align}
	Note that the minimal weak upper gradients are in fact local Lipschitz constants.
\par	Now, let $k\in \Sobol(\Y \times \left[ \shortminus \nicefrac{\uppi}{2}, \nicefrac{\uppi}{2} \right])$ and let $k_{_n}$ be a sequence of Lipschitz function on $\Y \times \left[ \shortminus \nicefrac{\uppi}{2}, \nicefrac{\uppi}{2} \right]$ that converges in $\Sobol(\Y\times \left[ \shortminus \nicefrac{\uppi}{2}, \nicefrac{\uppi}{2} \right])$ to $k$.
	Then, it also follows that $ k_{_n}\circ{\fancy{$\mathcal S$}}$ converges in $\Ltwo(\meas)$ to $k\circ {\fancy{$\mathcal S$}}$ since ${\fancy{$\mathcal S$}}$ is a measure space isomorphism. Consequently, considering \eqref{ineq:pres} with $k_{_n}$ instead of $k$, it follows that that $\CHE^\X(\upchi\cdot (k_{_n}\circ {\fancy{$\mathcal S$}}))$ is uniformly bounded. 
\par	Hence, one can extract a subsequence of $
	\left| \nabla (\upchi \cdot (k_{_n}\circ {\fancy{$\mathcal S$}})) \right|_{\substack{\\[2pt]\X}}$ that converges weakly to $G\in \Ltwo(\meas)$. By the stability property of weak upper gradients (Theorem \ref{th:stagra}), $G$ is a weak upper gradient of $\upchi \cdot (k\circ{\fancy{$\mathcal S$}})$. 
	Moreover, the inequality \eqref{ineq:pres} is preserved in the limit and hence 
	\begin{align*}
		\left| \nabla  (\chi\cdot (k\circ {\fancy{$\mathcal S$}}))\right|_{\substack{\\[2pt]\X}} \leq G\leq 2^{\nicefrac{1}{2} \,}k\circ {\fancy{$\mathcal S$}} \left| \nabla \upchi \right|_{\substack{\\[2pt]\X}} + 2^{\nicefrac{1}{2}} \, \upchi \left| \nabla k \right|_{\substack{\\[2pt] \Y \times \left[ \shortminus \nicefrac{\uppi}{2}, \nicefrac{\uppi}{2} \right]}} \circ  {\fancy{$\mathcal S$}} \quad \meas\mbox{-a.e.};
	\end{align*}
	therefore,  
	\begin{align*}
		\left| \nabla   (\upchi \cdot (k\circ {\fancy{$\mathcal S$}})) \right|_{\substack{\\[2pt]\X}}  \leq 2^{\nicefrac{1}{2}} \, \left| \nabla k \right|_{\substack{\\[2pt]\Y \times \left[ \shortminus \nicefrac{\uppi}{2}, \nicefrac{\uppi}{2} \right]}} \circ {\fancy{$\mathcal S$}}  \quad \meas\mbox{-a.e.  on } {\chi\equiv 1}.
	\end{align*}
	Since $\upchi$ was arbitrary, it follows
	\begin{align*}
		\left| \nabla   (k\circ {\fancy{$\mathcal S$}}) \right|_{\substack{\\[2pt] \X\smallsetminus \Sing}} \leq 2^{\nicefrac{1}{2}} \, \left| \nabla k \right|_{\substack{\\[2pt]\Y \times \left[ \shortminus \nicefrac{\uppi}{2}, \nicefrac{\uppi}{2} \right]}}\circ {\fancy{$\mathcal S$}}  \quad \meas\mbox{-a.e. in $\X$. }
	\end{align*}
	The claim then follows upon integration. 
\end{proof}
\begin{proposition}\label{prop:isometry} For every $k\in \Ltwo \left( \Y \times \left[ \shortminus \nicefrac{\uppi}{2}, \nicefrac{\uppi}{2} \right] \right)$ we have 
	\begin{align*}
		\CHE^{\X\smallsetminus \Sing}\left( k \circ {\fancy{$\mathcal S$}} \right) = \CHE^{ \Y \times \left[ \shortminus \nicefrac{\uppi}{2}, \nicefrac{\uppi}{2} \right]}(k).
	\end{align*}
	In particular, 
	the (dual) map
	\begin{align*}
		{\fancy{$\mathcal{S}$}}\strut^{*}:  \Sobol\big(\Y\times \left[ \shortminus \nicefrac{\uppi}{2}, \nicefrac{\uppi}{2} \right]\big) \to \Sobol_{\hspace{-1pt}\mathsf{loc}}\left(\X \smallsetminus \Sing\right),
	\end{align*}
	is an energy-preserving isomorphism. 
\end{proposition}
\begin{proof}[\footnotesize \textbf{Proof}]
	By density reasons (c.f. \cite{Gsplit} or \cite{GKK}), it is sufficient to prove the claim for functions of the form
	\begin{align}\label{id:special}
		k = \sum_{i\in \mathcal{I}} g_{\substack{\\[2pt]\hspace{-1pt}i}} h_i,
	\end{align}
	for a finite index set $\mathcal{I}$ and for $g\in \mathcal G$ and $h\in \mathcal H$ where 
	\begin{align*}
		\mathcal G= \left\{g: \Y\times \left[ \shortminus \nicefrac{\uppi}{2}, \nicefrac{\uppi}{2} \right] \rightarrow \R \; \text{\textbrokenbar} \; g\left(y,t\right) = \uphat g\left(y\right)  \ \mbox{ for some } \ \uphat g\in \Sobol\left(\Y\right) \cap  \LL^{\hspace{-2pt}^\infty} \left(\Y, {\bf q}\right)\right\},
	\end{align*}
	and
	\begin{align*}
		\mathcal H = &\left\{ h: \Y\times \left[ \shortminus \nicefrac{\uppi}{2}, \nicefrac{\uppi}{2} \right] \rightarrow \R \; \text{\textbrokenbar} \; h\left(y,t\right)= \uphat h\left(t\right)  \right.  \\ & \phantom{cyscyscys}\left.  \mbox{ for some } \uphat h \in \Sobol\left(  \left[ \shortminus \nicefrac{\uppi}{2}, \nicefrac{\uppi}{2} \right] \right) \cap  \LL^{\hspace{-2pt}^\infty} \left(  \left[ \shortminus \nicefrac{\uppi}{2}, \nicefrac{\uppi}{2} \right] \right) \right\}.
	\end{align*}
\par If $u$ is of the form (\ref{id:special}),  recalling that $\widetilde \Y\times  \left[ \shortminus \nicefrac{\uppi}{2}, \nicefrac{\uppi}{2} \right]$ is $\RCD\left(0,\N+1\right)$, one can expand $\left|\nabla u\right|_{\substack{\\[2pt]\Y\times  \left[ \shortminus \nicefrac{\uppi}{2}, \nicefrac{\uppi}{2} \right]}}^{^2} $ to get
	\begin{align*}
		\sum_{i,j\in \mathcal{I}} \left\{g_{\substack{\\[2pt]\hspace{-1pt}i}} g_{\substack{\\[2pt]\hspace{-1pt}j}} \langle \nabla h_i, \nabla h_j\rangle_{\substack{\\[3pt]\hspace{-1pt}\Y\times  \left[ \shortminus \nicefrac{\uppi}{2}, \nicefrac{\uppi}{2} \right]}} + 2 g_{\substack{\\[2pt]\hspace{-1pt}i}}h_j \langle \nabla g_{\substack{\\[2pt]\hspace{-1pt}j}}, \nabla h_i\rangle_{\substack{\\[3pt]\hspace{-1pt}\Y\times  \left[ \shortminus \nicefrac{\uppi}{2}, \nicefrac{\uppi}{2} \right]}} + h_ih_j \langle \nabla g_{\substack{\\[2pt]\hspace{-1pt}i}}, \nabla g_{\substack{\\[2pt]\hspace{-1pt}j}}\rangle_{\substack{\\[3pt]\hspace{-1pt}\Y\times  \left[ \shortminus \nicefrac{\uppi}{2}, \nicefrac{\uppi}{2} \right]}}\right\}.
	\end{align*}
\par	It is evident form Proposition~\ref{prop:pullback-Sobolev}, that
	\begin{align*}
		\left|\nabla \left(\uphat g_{\substack{\\[2pt]\hspace{-1pt}i}}\circ 	\mathrm{proj}_{\substack{\\[2pt]\Y}} \right)\right|_{\substack{\\\X \smallsetminus \Sing}}= \left|\nabla \uphat g_{\substack{\\[2pt]\hspace{-1pt}i}} \right|_{\substack{\\[1pt]\Y}} \circ 	\mathrm{proj}_{\substack{\\[2pt]\Y}}  \quad \meas\mbox{-a.e. in} \quad \X \smallsetminus \Sing,
	\end{align*}
	which upon polarization, implies 
	\begin{align*}
		\big\langle \nabla \left(\uphat{g}_{\substack{\\[2pt]\hspace{-1pt}i}} \circ 	\mathrm{proj}_{\substack{\\[2pt]\Y}} \right), \nabla \left(\uphat{g}_{\substack{\\[2pt]\hspace{-1pt}j}}\circ	\mathrm{proj}_{\substack{\\[2pt]\Y}} \right) \big\rangle_{\substack{\\[3pt]\hspace{-2pt} \X \smallsetminus \Sing}} &= \langle \nabla \uphat{g}_{\substack{\\[2pt]\hspace{-1pt}i}}, \nabla \uphat{g}_{\substack{\\[2pt]\hspace{-1pt}j}}\rangle_{\substack{\\[3pt]\hspace{-2pt} \Y}}  \circ 	\mathrm{proj}_{\substack{\\[2pt]\Y}} \\ &= \langle \nabla g_{\substack{\\[2pt]\hspace{-1pt}i}}, \nabla g_{\substack{\\[2pt]\hspace{-1pt}j}}\rangle_{\substack{\\[3pt]\hspace{-1pt}\Y\times  \left[ \shortminus \nicefrac{\uppi}{2}, \nicefrac{\uppi}{2} \right]}}  \circ {\fancy{$\mathcal{S}$}},
	\end{align*}
	for $ \meas\mbox{-a.e. } x \ \mbox{in } \X \smallsetminus \Sing$. 
\par	In below, we proceed in a similar fashion as in the proof of the splitting theorem in \cite{Gsplit}. For $g \in 	\mathcal G$ and $h \in \mathcal H$, we set $w := \uphat{g} \circ 	\mathrm{proj}_{\substack{\\[2pt]\Y}}$ and $v = \uphat{h} \circ f$; recall $\uphat{g} \in \Sobol\left( \Y \right)$ and $\uphat{h}\in \Sobol\left( \left[ \shortminus \nicefrac{\uppi}{2},  \nicefrac{\uppi}{2} \right] \right)$. Then, $v \in \Sobol\left( \X \smallsetminus \Sing \right)$
	and 
	\begin{align*}
		\left|\nabla v \right|_{\substack{\\\X \smallsetminus \Sing}} \left(x\right) = \left| \uphat{h}^{'}\right|\circ f\left(x\right)\ \mbox{ for }\meas\mbox{-a.e. }x\in \X. 
	\end{align*}
\par	Thus again by polarization, one obtains
	\begin{align*}
		\big\langle \nabla \left(\uphat{h}_i\circ f \right), \nabla \left(\uphat{h}_j\circ f\right) \big\rangle_{\substack{\\[3pt]\hspace{-2pt} \X \smallsetminus \Sing}}= \langle \nabla h_i, \nabla h_j\rangle_{\substack{\\[3pt]\hspace{-1pt}\Y\times  \left[ \shortminus \nicefrac{\uppi}{2}, \nicefrac{\uppi}{2} \right]}}.
	\end{align*}
\par	At this point, a verbatim argument as in \cite{Gsplit} shows
	\begin{align*}
		\langle \nabla g, \nabla h\rangle_{\substack{\\[3pt]\hspace{-1pt}\Y\times  \left[ \shortminus \nicefrac{\uppi}{2}, \nicefrac{\uppi}{2} \right]}} =0 \quad {\bf q}\otimes \mathscr L^{^1}\mbox{-a.e.},
	\end{align*}
	and consequently,
	\begin{align*}
		\big\langle \nabla \left(g \circ {\fancy{$\mathcal{S}$}} \right), \nabla \left(h\circ {\fancy{$\mathcal{S}$}} \right) \big\rangle_{\substack{\\[3pt]\hspace{-2pt} \X}} =0 \quad \meas\mbox{-a.e.}\ .
	\end{align*}
	Putting these steps together, on the level of miminal weak upper gradients, we obtain the identity
	\begin{align}\label{id:before}
		\left|\nabla \left(gh \circ {\fancy{$\mathcal{S}$}} \right)\right|_{\substack{\\[2pt]\X \smallsetminus \Sing}} = \left|\nabla \left(gh\right) \right|_{\substack{\\ \Y\times \left[ \shortminus \nicefrac{\uppi}{2}, \nicefrac{\uppi}{2} \right]}} \circ {\fancy{$\mathcal{S}$}},
	\end{align}
	for any $g\in \mathcal{G}$ and $h\in \mathcal{H}$; hence, the same holds for $k = \sum_{i\in \mathcal{I}} g_{\substack{\\[2pt]\hspace{-1pt}i}} h_i$. 
	\par For general $k \in \Sobol\left( \Y\times \left[ \shortminus \nicefrac{\uppi}{2}, \nicefrac{\uppi}{2} \right] \right)$, we can find a sequence $k_{_n}$ of functions in the form \eqref{id:special} so that $k_{_n} \rightarrow k$ in $\Sobol\left( \Y\times \left[ \shortminus \nicefrac{\uppi}{2}, \nicefrac{\uppi}{2} \right] \right)$.  Corollary \ref{cor:rough} yields that $k\circ {\fancy{$\mathcal S$}}, k_{_n}\circ {\fancy{$\mathcal S$}}\in \Sobol_{\hspace{-1pt}\textsf{loc}}(X\smallsetminus \Sing)$ and $k_{_n}\circ {\fancy{$\mathcal S$}}\rightarrow k\circ {\fancy{$\mathcal S$}}$ in $\Sobol_{\hspace{-1pt}\textsf{loc}}(X\smallsetminus \Sing)$.  Finally, since the energy identity holds for $k_{_n}\circ {\fancy{$\mathcal S$}}$, it also passes over to the limit. 
\end{proof}
\subsubsection*{\small \bf \textit{Isometric splitting of $\widetilde \X$}}
\par We have, so far, collected all the needed ingredients to prove the splitting phenomenon by following the proof of \cite[Proposition 4.20]{Gsplit}. 
\begin{theorem}\label{thm:split-iso}
	Restricted to path connected components, the maps $	{\fancy{$\mathcal{S}$}}: \widetilde \X\rightarrow \widetilde \Y\times \left[ \shortminus \nicefrac{\uppi}{2}, \nicefrac{\uppi}{2} \right] $ and $	{\fancy{$\mathcal{T}$}}: \widetilde \Y\times \left[\shortminus \nicefrac{\uppi}{2},\nicefrac{\uppi}{2}\right]\rightarrow \widetilde \X$ are isomorphisms between metric measure spaces.
\end{theorem}
\begin{proof}[\footnotesize \textbf{Proof}]
	For the proof of the theorem we can follow verbatim the same argument as for  the direction $(ii)\Rightarrow (i)$ in the proof of \cite[Proposition 4.20]{Gsplit}. 
\par	Let us outline the main steps for the sake of completeness. 
\par	First, we have the following localization of the norm of gradients. 
	\begin{lemma}[Localization]\label{hhh}
		For every $k\in \Sobol(\Y\times \left[ \shortminus \nicefrac{\uppi}{2}, \nicefrac{\uppi}{2} \right] )$ it holds
		\begin{align*}
			\left| \nabla (k\circ {\fancy{$\mathcal{S}$}}) \right|_{\substack{\\[2pt] \X \smallsetminus \Sing}}= \left| \nabla k \right|_{\substack{\\[2pt] \Y\times \left[ \shortminus \nicefrac{\uppi}{2}, \nicefrac{\uppi}{2} \right]}}\circ {\fancy{$\mathcal{S}$}}.
		\end{align*}
	\end{lemma}
\begin{proof}[\footnotesize \textbf{Proof}]
	The proof of the lemma is exactly the proof of Lemma 4.17 in \cite{Gsplit}. It uses that the map ${\fancy{$\mathcal S$}}$ is a measure space isomorphism and that the dual map ${\fancy{$\mathcal S$}}^*$ is an energy isomorphism; this was established in Proposition~\ref{prop:isometry}.
\end{proof}
	\par Another key lemma is 
	\begin{lemma}[Contraction by local duality]\label{hh}
		There exists a  map $\widetilde{\fancy{$\mathcal{S}$}}: \widetilde \X \rightarrow {\widetilde\Y}\times \left[ \shortminus \nicefrac{\uppi}{2}, \nicefrac{\uppi}{2} \right]$ that is $1$-Lipschitz on every path connected component and coincides with ${\fancy{$\mathcal S$}}$ in $\meas$-a.e. sense. 
	\end{lemma}
\begin{proof}[\footnotesize \textbf{Proof}]
The proof of this lemma again goes exactly along the same lines as in the argument provided for the direction $(ii)\Rightarrow (i)$ in the proof of \cite[Lemma 4.19]{Gsplit}; so for this proof, one takes advantage of the previous lemma and  the modified Sobolev-to-Lipschitz property that we established in Proposition \ref{prop:modSobLip}.
\par	Lemmas~\ref{hhh} and~\ref{hh} allow us to finish the proof of Theorem~\ref{thm:split-iso} by following the direction $(ii)\Rightarrow (i)$ in the  proof of \cite[Proposition 4.20]{Gsplit}. Note that  in the beginning of~\textsection\thinspace\ref{subsec:split}, we have already established that ${\fancy{$\mathcal S$}}$ and ${\fancy{$\mathcal T$}}$ are measure space isomorphisms and that ${\fancy{$\mathcal T$}}\circ {\fancy{$\mathcal S$}}(x)=x$ and  ${\fancy{$\mathcal S$}}\circ {\fancy{$\mathcal T$}}(y,t)=(y,t)$ almost everywhere.  In particular, this yields an isomorphism between the corresponding $\Ltwo$-spaces.
\end{proof}
\par	Hence, by the energy identity for ${\fancy{$\mathcal S$}}$ (Proposition \ref{prop:isometry}),  we deduce $h= h\circ {\fancy{$\mathcal T$}}\circ {\fancy{$\mathcal S$}}\in \Ltwo(\meas)$ satisfies $\CHE^{\X\smallsetminus \Sing}(h)<\infty$ if and only if $\CHE^{\Y\times \left[ \shortminus \nicefrac{\uppi}{2}, \nicefrac{\uppi}{2} \right]}(h\circ {\fancy{$\mathcal T$}}) <\infty$.  At this point we can repeat the previous steps with ${\fancy{$\mathcal T$}}$ in place of ${\fancy{$\mathcal S$}}$, while again closely following the arguments in \cite{Gsplit}. We underline the local nature of the arguments.	This finishes the proof.
\end{proof}
\begin{remark}
	Following the line of arguments in~\cite[Section 7]{Gsplit}, one can verify the stronger $\RCD\left(  0, \N  \shortminus1\right)$ conditions for $\Y$. We omit the details.
\end{remark}
\section{Proof of the main theorem}\label{sec:proof-main}
This section concludes the proof of our main result, Theorem~\ref{thm:main-YZ}. Recall $\left(\widetilde \X,  \dist_{\substack{\\\widetilde \X}}\right)=: \widetilde \X$ is the completion of the \emph{extended intrinsic distance} on $\X \smallsetminus  \Sing$; see~\textsection\thinspace\ref{sec:defn-tilde-spaces}. 
\par So far, we have shown $\left( \widetilde \X, \dist_{\widetilde \X}, \meas \right)$  is isomorphic to $\left( \widetilde \Y\times \left[ \shortminus \nicefrac{\uppi}{2}, \nicefrac{\uppi}{2} \right], {\bf q}\otimes \mathscr L^{^1}\mres_{\left[ \shortminus \nicefrac{\uppi}{2}, \nicefrac{\uppi}{2} \right]} \right)$ and 
we want to show $\left(\X, \dist_{\X}, \meas\right)$ must be isomorphic to $\mathbb S^{^1}$ or $\left[0,\uppi\right]$. 
\subsection{Foliation by horizontal geodesics}
Recall the definition of the singular sets $\Sing_{\pm}$ in Definition~\ref{defn:sing-set}.
\begin{definition}[horizontal geodesics]
A geodesic $\upgamma$ with end points $p_{_-}  \in\Sing_{_{-}}$ and $p_{_+}  \in \Sing_{_{+}}$ is called a horizontal geodesic.
\end{definition}
	\begin{lemma}
	A horizontal geodesic $\upgamma$ intersects each of $\Sing_{_{-}}$ and $\Sing_{_{+}}$ at exactly one point (the end points $p_{_-}$ and $p_{_+}$ resp.).
\end{lemma}
\begin{proof}[\footnotesize \textbf{Proof}]
A computation similar to \eqref{EQN:cruciel-computation} (with $\upgamma_{\substack{\\[1pt]\hspace{-2pt}q}}$ replaced by $\upgamma$), shows $\upgamma$ is of length $\uppi$. Since the diameter of $\X$ is also $\uppi$, the geodesic $\upgamma$ can not intersect each of  $\Sing_{_{-}}$ and $\Sing_{_{+}}$ more than one time otherwise it would have excess length. 
\end{proof}
\begin{remark}
Notice that not every geodesic of length $\uppi$ is horizontal; e.g. when $\X$ is a circle.
\end{remark}
\begin{lemma}\label{lem:horiz-char}
	Any curve $\upgamma$ of length at most $\uppi$ intersecting both $\Sing_{_{-}}$ and $\Sing_{_{+}}$, is a horizontal geodesic.
\end{lemma}
\begin{proof}[\footnotesize \textbf{Proof}]
	Suppose $\upgamma(t_{_0}) \in \Sing_{_{-}}$ and $\upgamma(t_{_1})\in \Sing_{_{+}}$ then 
	\[
	\length \left( \upgamma\restr_{\left[ t_{_0},t_{_1} \right]} \right) \ge \uppi
	\]
	hence by the hypothesis we must have $t_{_0} = 0$ and $t_{_1} = 1$ and consequently $
	\length (\upgamma) = \uppi.
	$ and $\upgamma$ is a horizontal geodesic. 
\end{proof}

The properties established in the following proposition provide a foliation of $\X \smallsetminus \Sing$ by (the interior of) horizontal geodesics.
\begin{proposition}\label{prop:foliation}
	The following statements hold true;
	\begin{enumerate}
		\item (weak linearity) every point $x\in \X$ lies on at least one horizontal geodesic along which $f$ is linear with slope $1$;
		\medskip
		\item (distance to singularity) for any point $x$, 
		\[
		\dist\left( x,\Sing_{_{-}} \right) = f(x) + \nicefrac{\uppi}{2}  \quad  \text{and} \quad \dist \left( x,\Sing_{_{+}}\right)  = \nicefrac{\uppi}{2}  - f(x);
		\]\smallskip
		\item (extension property) 	for any point $x \in \X$, any shortest path $\upgamma$ from $x$ to either $\Sing_{_{-}}$ or $\Sing_{_{+}}$ can be extended to a horizontal geodesic i.e. $\upgamma$ is a segment of a horizontal geodesic; in particular, $x$ is equidistant to all points of $\Sing_{_{-}}$ and also $x$ is equidistant to all points of $\Sing_{_{+}}$.
		\medskip
		\item (strong linearity) $f$ is linear with slope $1$ along all horizontal geodesics;
		\medskip
		\item (essential uniqueness) $\meas$-a.e. $x \in \X \smallsetminus \Sing$ lies on a unique horizontal geodesic along which $f$ is linear with slope $1$; consequently, this geodesic is one of the maximal transport geodesics $\upgamma_{\substack{\\[1pt]\hspace{-2pt}q}}$.
	\end{enumerate}
\end{proposition}
\begin{proof}[\footnotesize \textbf{Proof}]
	\hfill
	\begin{enumerate}
		\item 
	\par  By the key localization lemma (Lemma \ref{lem:ida}), we know for ${\bf q}$-a.e. $q\in \Q$, $\upgamma_{\substack{\\[1pt]\hspace{-2pt}q}}$ is of length $\uppi$. Furthermore, along these geodesics, $u$ (and $f$) are monotonically increasing; hence one deduces $\upgamma_{\substack{\\[1pt]\hspace{-2pt}q}}$ must join a point from $u^{-1}(-1)$ to $u^{-1}(1)$ i.e. $\upgamma_{\substack{\\[1pt]\hspace{-2pt}q}}$ is a horizontal geodesic. Consequently, since the transport set if of full measure, we deduce almost every point lies on a horizontal geodesic; see Remark~\ref{rem:horiz-foliation}.
	\par Now let $x\in \X$ be arbitrary. For any $r>0$, $\ball{r}(x)$ intersects a horizontal geodesic. Hence by applying Arzela-Ascoli, as $r\to0$, one gets a sequence of horizontal geodesics that are subsequentially converging to a geodesic (e.g.~\cite{BBI}). This limiting geodesic must pass through $x$. Since $\Sing_{_{-}}$ and $\Sing_{_{+}}$ are compact subsets, one deduces that the limiting geodesic has end points in  $\Sing_{_{-}}$ and $\Sing_{_{+}}$ and by the lower semi-continuity of length, the length of the limiting geodesic is at most $\pi$, hence by Lemma~\ref{lem:horiz-char}, $\upgamma$ is a horizontal geodesic. 
	
	\par The linearity of $f$ follows from the fact that for ${\bf q}$-a.e. $q\in \Q$, $f$ is linear with slope $1$ along $\upgamma_{\substack{\\[1pt]\hspace{-2pt}q}}$ therefore, repeating the above arguments with these horizontal geodesics instead, one gets a limiting geodesic along which $f$ is linear with slope $1$. 
	\medskip
	\item We know there exists a horizontal geodesic passing through $p$ along which $f$ is linear with slope $1$. This implies
	\[
	\dist\left( p,\Sing_{_{-}} \right) \le f(p) + \nicefrac{\uppi}{2}, \quad  \text{and} \quad \dist \left( p,\Sing_{_{+}}\right)  \le \nicefrac{\uppi}{2}  - f(p).
	\]
	On the other hand,
	\[
	\uppi = \dist\left( \Sing_{_{-}}, \Sing_{_{+}}\right) \le \dist\left( p, \Sing_{_{-}}\right)  +  \dist\left( p, \Sing_{_{+}}\right) \le \left(   f(p) + \nicefrac{\uppi}{2} \right) +  \left( \nicefrac{\uppi}{2}  - f(p) \right) = \uppi.
	\]
	hence, all the above inequalities must be identities.
	\medskip
	\item Let $\upgamma$ be a shortest path from $\Sing_{_{-}}$ to $x$ and $\uptheta$ be a shortest path from $x$ to $\Sing_{_{+}}$. By item (2), one has
	\[
	\length (\upgamma) = f(x) + \nicefrac{\uppi}{2}, \quad \text{and} \quad  \length (\uptheta) =   \nicefrac{\uppi}{2} - f(x). 
	\]
	\par Set $\upbeta := \upgamma + \uptheta$ (concatenation of the two geodesics). $\upbeta$ intersects both $\Sing_{_{-}}$ and $\Sing_{_{+}}$  and
	\[
\length (\upbeta) = \length(\upgamma) + \length (\uptheta) = \uppi;
	\]
	hence, by Lemma~\ref{lem:horiz-char}, $\upbeta$ is a horizontal geodesic. 
	\medskip
	\item This directly follows from item (3). 
	\medskip
	\item Let $x\in \X \smallsetminus \Sing$; then by item (1), there exists a horizontal geodesic $\upgamma$ along which $f$ is linear with slope $1$. This implies $\upgamma$ is included in the transport set and hence for $\meas$-a.e. such $x$, the corresponding $\upgamma$ is included in $\upgamma_{\substack{\\[1pt]\hspace{-2pt}q(x)}}$ for a unique value of $q(x)$.  Hence  $\upgamma_{\substack{\\[1pt]\hspace{-2pt}q(x)}}$ is of length $\uppi$ and  intersects both  $\Sing_{_{-}}$ and $\Sing_{_{+}}$; this means $\upgamma_{\substack{\\[1pt]\hspace{-2pt}q(x)}}$ is a horizontal geodesic. 
		\par Since $\upgamma_{\substack{\\[1pt]\hspace{-2pt}q}}$'s are non-branching for ${\bf q}$-a.e. $q\in \Q$ (this follows from dis-integration in Theorem~\ref{thm:1dlocalization} and the fact that the non-branching transport set is of full measure), we established the essential uniqueness of horizontal geodesics. 
%
%
	\end{enumerate}
\end{proof}
\begin{remark}[uniqueness of horizontal geodesics]
	In~\cite{qin}, a proof of non-branching of RCD spaces has been presented. So if one assumes non-branching, then horizontal geodesics are unique and foliate $X$; for our purposes, the essential uniqueness established in item (5) in Proposition~\ref{prop:foliation} would suffice. So we are not assuming non-branching of RCD spaces. 
\end{remark}
\subsection{Conclusion of the proof of the main theorem}
 Recall $\mathcal R_1$ is the set of points in $X$ that admit unique tangent cone isometric to $\R$.
To complete the proof, we will consider two cases:
\subsubsection*{\small \bf Case I} $\mathcal R_1 \neq \emptyset$. \\
\par In this case,  $\X$ is isometric to $\mathbb S^{^1}$ or $[0,\uppi]$ by \cite{KL-2}. Moreover, in this case $\X$ is even isomorphic to $\left(\mathbb S^{^1}, \vol_{\mathbb S^{^1}} \right)$ or $\left( [0,\uppi], \mathscr L^{^1}|_{[0,\uppi]}  \right)$ as a metric measure space; indeed the measure is $\meas = {\mathscr H}^{^1}$ because of Lemma~\ref{lem:lemD}. Hence, we are done.
\subsubsection*{\small \bf Case II} $\mathcal R_1 = \emptyset$. \\
\par By Proposition~\ref{prop:foliation}, we know $\meas$-a.e. $y\in \X \smallsetminus \Sing$ lies on a unique horizontal geodesic.  Corresponding to $y$, we consider  the potential function (distance function) $f_y:= \dist_{\X}(y,\cdot)$ that is a globally $1$-Lipschitz function on $\X$. Theorem~\ref{thm:1dlocalization} implies that there exists a family of non-branching geodesics $\{\uptheta_p:  \left[ \ai_{\substack{\\p}}, \bi_{\substack{\\p}} \right] \rightarrow \X\}_{p\in {\rm P}}$ with disjoint interior images such that 
\[
\bigcup_{p\in {\rm P}}\mbox{Im}(\uptheta_p) \underset{\textsf{full measure}}{\subset} \X 
\]
and $\uptheta_p\left(\bi_p\right)=y$ for every $p\in {\rm P}$; recall that there is also a measure ${\bf p}$ on $\rm P$  that is the cross section measure in the disintegration. 
\par { \bf \small Claim.} ${\bf p}(\{p\})=0$ for all $p\in \rm P$ i.e.  ${\bf p}$ is atom-less.
\par	{\it \small Proof of the claim.} There are two cases to discuss. 
\begin{enumerate}
	\item The first case is when there is no point on $\uptheta_p\left(\left( \ai_{\substack{\\p}}, \bi_{\substack{\\p}} \right) \right)$  with a unique tangent cone. In this case, by~\cite{MN}, it follows that $\uptheta_p\left(\left( \ai_{\substack{\\p}}, \bi_{\substack{\\p}} \right) \right)$ is contained in a subset of $\meas$-measure $0$. Hence, one can see from the disintegration of $\meas$ w.r.t. ${\bf p}$ that ${\bf p}(\{p\})=0$ must hold. 
\item The second case is when $\uptheta_p\left(\left( \ai_{\substack{\\p}}, \bi_{\substack{\\p}} \right) \right)$ contains a point where the tangent cone is unique. 
Then we argue by contradiction. Assume $\uptheta_p:  \left[ \ai_{\substack{\\p}}, \bi_{\substack{\\p}} \right] \rightarrow \X$ is a geodesic with $p\in {\rm P}$ such that ${\bf p}(\{p\})>0$, and let  $x_{_0}\in \uptheta_p\left(\left( \ai_{\substack{\\p}}, \bi_{\substack{\\p}} \right) \right)$ with $\uptheta_p(t_{_0})=x_{_0}$ where the tangent cone is unique. Then, given a sequence $(r_{\!_i})_{_{i\in \mathbb N}}$ with $r_{\!_i}\downarrow 0$, it is true by disintegration that 
$$\meas\left(\uptheta_p\left(\left[ t_{_0}-r_{\!_i}, t_{_0}+r_{\!_i}\right] \right)\right)= \int_{t_0-r_i}^{t_0+r_i} \uprho_p  \; d \mathcal L^1\restr_{\left[t_0-r_i, t_0+r_i\right]} \cdot {\bf p}(\{p\})>0$$ 
for all $i\in \mathbb N$; recall $\uprho_p$ is the density that appears in Theorem \ref{th:1Dlocscheme}.

We define $(\X,  r_i^{-1}\dist, \meas^{x_0}_{r_i})_{i\in \mathbb N}$ 
where 
\[
\meas_{r_{\!_i}}^{x_0}= \left( \int_{B_{r_{\!_i}}(x_0)}1- \frac{1}{r_{\!_i}}\dist(\cdot, x_{_0}) \; d\meas \right)^{\!-1} \meas.
\]
 By  the Gromov's compactness theorem and by the assumption of uniqueness of tangent cones,  $\left(X, r_{\!_i}^{-1} \dist, \meas_{r_i}^{x_0}, x_{_0}\right)$ subsequentially converges in pointed measured Gromov-Hausdorff sense to a pointed metric measure space $(C, \dist_C, \meas_C, o)$ that is the unique tangent cone in $\uptheta_p(t_0)$  and satisfies $\RCD(0,\N)$.  Since $x_0$ lies on a geodesic, $C$ contains a line $\gamma$ with $\gamma(0)=o$ and therefore, by the  generalization of the splitting theorem established in \cite{Gsplit,Goverview}, splits as $C\simeq \R\times {\sf Y}$ for an $\RCD(0, \N-1)$ space $\left( {\sf Y}, \dist_{\sf Y}, \meas_{\sf Y} \right)$ that is not discrete because $\mathcal R_1=\emptyset$.  Moreover, the geodesic segment $\uptheta_p|_{[t_0-r_i, t_0+r_i]}$ has length $1$ w.r.t. $r_i^{-1} \dist$ and  by  standard Gromov-Hausdorff theory, one can embed the sequence in a single complete metric space $(Z,\dist_Z)$  where measure Gromov-Hausdorff convergence is realized; hence, $\uptheta_p|_{[t_0-r_i, t_0+r_i]}$ converges to $\gamma|_{[-1, 1]}$ in $Z$. 
\par  The pointed measure GH convergence now further implies that the (pushforward of) measures also weakly converge in $Z$. Hence 
\begin{align*}
0<\uprho_p(t_0) \cdot {\bf p}(\{p\}) &= \lim_{i\rightarrow \infty} \frac{1}{r_i} \int_{t_0-r_i}^{t_0+r_i} \uprho_p\; d\mathcal L^1\restr_{[t_0-r_i, t_0+r_i]}\cdot {\bf p}(\{p\}) \\ &\leq \meas_C(\gamma([-1,1]).
\end{align*}
Since $C$ splits along $\gamma$, there exists $y_{_0}\in {\sf Y}$ such that $\gamma|_{[-1,1]}= [-1,1]\times \{y_0\}$  and $\meas_C= \mathcal L^1 \otimes \meas_{\sf Y}$. Hence $\meas_C(\gamma([-1,1]))= \mathcal L^1([-1,1]) \cdot \meas_{\sf Y}(\{y_0\})>0$ yielding $\meas_{\sf Y}(\{y_0\})>0$. But ${\sf Y}$ is an $\RCD(0, \N)$ space and consequently has no atoms. That is the contradiction.  \hfill\scalebox{0.6}{$\blacksquare$}
\end{enumerate}
\par Consider the subsets (set-valued maps)
\[
{\sf C}_{_-}(y) := \{p: \exists t\in \left( \ai_{\substack{\\p}}, \bi_{\substack{\\p}} \right)\mbox{ s.t. }\uptheta_p(t)\in \Sing_{_-}\} \quad \text{and} \quad   {\sf C}_{_+} (y):= \{p: \exists t\in \left( \ai_{\substack{\\p}}, \bi_{\substack{\\p}} \right)  \mbox{ s.t. } \uptheta_p(t)\in \Sing_{_+}\}.
\]
of $\rm P$. 
\par {\bf \small Claim.} For a.e. $y\in \X \smallsetminus \Sing$, $ {\sf C}_{_\pm} (y)$ are single-valued; in particular,  
${\bf p} \left( {\sf C}_{_\pm} (y) \right) = 0$ holds for  $\meas$-a.e. $y\in \X \smallsetminus \Sing$. 
\par	{\it \small Proof of the claim.} 
Recall by the horizontal foliation (Proposition~\ref{prop:foliation}), we know that any geodesic $\uptheta$ that connects $\Sing_{_-}$ to a point $x$ in $X \smallsetminus \Sing$ with $x$ lying on a horizontal geodesic $\upgamma$, must be a restriction of this horizontal geodesic $\upgamma$. 
\par Now consider a geodesic $\uptheta_p$ with $p \in {\sf C}_{_-}(y)$ and $\uptheta_p\left( t_{_0}\right) \in \Sing_{_-}$.  On account of the aforementioned extension property and the essential uniqueness of the horizontal geodesics (item (5) in Proposition~\ref{prop:foliation}), for $\meas$-a.e. such $y$,  any such $\uptheta_p\restr_{\left[ t_{_0},  \bi_{\substack{\\p}} \right]}$ is a segment of the \emph{unique horizontal geodesic} (which is also a maximal transport geodesic) $\upgamma_{\substack{\\[1pt]\hspace{-2pt}q}}$ that passes through $y$.  This in conjunction with the fact that the interior of $\uptheta_p$'s are disjoint, yields ${\sf C}_{_-}(y) $ is a singleton for $\meas$-a.e. $y$. Same argument applies to ${\sf C}_{_+}(y)$ as well. \phantom{123} \hfill\scalebox{0.6}{$\blacksquare$}
\par Consequently, $\meas$-a.e. point $x$ in the same connected component of $\X \smallsetminus \Sing$ as $y$, is connected to $y$ via a geodesic that stays in  $\X \smallsetminus \Sing$. In particular, for $\meas\otimes \meas$-almost every pair of points $(x,y)$  in $(\X \smallsetminus \Sing)^2$, there exists a geodesic that stays in $\X \smallsetminus \Sing$.  This implies $\widetilde \X$ is isometric to $\X$ and in particular, $\widetilde \X$ is connected. 
\par $\widetilde{X} = X$ together with the isometric splitting given by Theorem~\ref{thm:split-iso}, $\diam {\X} = \uppi$ and connectedness, imply $\widetilde \Y$  has diameter zero and hence is a point. This means $\X$ is one dimensional that is a contradiction to the hypothesis $\mathcal R_1 = \emptyset$; hence, we showed {\bf Case II} ($\mathcal R_1 = \emptyset$) never occurs. \qed

\appendix
\section{Rigidity under non-negative Bakry-\'Emery Ricci tensor (alternative proof)}\label{app:Drif-Laplacian}
Here, we will present a proof of rigidity under non-negative Bakry-\'Emery Ricci tensor using smooth analysis techniques. We will also discuss that this case is indeed an special case of our main Theorem~\ref{thm:main-YZ}.
\par This proof entails adjustments in the classic proof of the rigidity for the non-drifted Laplacian which has been proven only in the \emph{closed case} in~\cite{Hang-Wang}; our proof is new and covers the drifted Laplacian and the the case where the manifold has a convex boundary. Also there are intricacies when $\mathscr{X}$ is not integrable that we will address. 
\subsection{Setup}
\par Associated to a smooth vector field $\mathscr{X}$ on the Riemannian manifold $\M^{^n}$, one defines the drift Laplacian
\begin{align*}
	\BEL = \Delta \shortminus \mathcal{L}_{\substack{\\\mathscr{X}}},
\end{align*}
\par Suppose $n<\N$. The $\N$-Bakry-\'{E}mery Ricci tensor associated to the vector field $\mathscr{X}$ is defined as
\begin{align*}
	\BERic = \Ric + \, \nicefrac{1}{2} \, \mathcal{L}_{\substack{\\\mathscr{X}}} \g \, \shortminus \, \nicefrac{1}{\left(\N-n\right)} \, \mathscr{X}\strut^{^\flat} \otimes \mathscr{X}\strut^{^\flat},
\end{align*}
in which $\mathscr{X}\strut^{^\flat}$ is the $1$-form dual to $\mathscr{X}$ and
\begin{align*}
	\mathcal{L}_{\substack{\\\mathscr{X}}}  \g \left(V,W\right)= \g\left(V,\nabla_{\substack{\\\hspace{-2pt}\mathscr{X}}}   W\right)+ \g\left(\nabla_{\substack{\\\hspace{-2pt}\mathscr{X}}}   V, W\right).
\end{align*}
\begin{remark}
	The definition of $\BERic$ for $n=\N$ only makes sense for the trivial vector field $\vX = 0$. Said another way, for a constant weight $\upvarphi = \mathsf{const}$, one has $\Ric^{\,n}_{ \, \upvarphi} = \Ric$. The other extreme is when $\N = \infty$ where we have $\BERic = \Ric + \nicefrac{1}{2} \, \mathcal{L}_{\substack{\\\mathscr{X}}} \g$. 
\end{remark}
\subsubsection*{\small \bf \textit{On curvature conditions in presence of boundaries}}
$\BE\left(\K,\N\right)$ conditions can also unambiguously be defined for smooth functions that are supported in the interior. We will denote this function space by $\mathcal{D}^{^\infty}\left(\M^{^n} \right)$; the $\BE\left(\K,\N\right)$  conditions restricted to  $\mathcal{D}^{^\infty}\left(\M^{^n} \right)$ are called 
 \emph{smooth $\BE\left(\K,\N\right)$ conditions}.  
\begin{proposition}
	Suppose $\M^{^n}$ is a manifold possibly with boundary. Then, the \textsf{smooth} $\BE\left(\K,\N\right)$ ($n<\N$) conditions is equivalent to $\BERic \ge \K$ (both happening in the interior). 
\end{proposition}
\begin{proof}[\footnotesize \textbf{Proof}]
	Both these bounds are local.  Any interior point admits a simply connected neighborhood in the interior on which the vector field $\vX$ is integrable. Thus, the equivalence follows from~\cite{Ledoux}. 
\end{proof}
Since the proof of comparison theorem in its fullest is contingent upon convexity assumptions  on the boundary (in some cases, only mean convexity and in some cases the positivity of the second fundamental form) and since the proof of rigidity will also require geodesic convexity of the interior (as we will soon see), we will set the following convexity assumption on the boundary for the rest of this section. 
\begin{definition}[Boundary convexity]
	In these notes, we say the boundary is convex whenever the second fundamental form of the boundary (as an embedded manifold in $\M$) is positive definite. 
\end{definition}
\par It follows that when  the boundary of $\M^{^n}$ is convex, then the interior is geodesically convex i.e. every geodesic joining two interior points lies entirely in the interior. Furthermore, any geodesic connecting two boundary points must have its interior included in the interior of $\M^{^n}$. 
\begin{remark}\label{rem:subcase}
	The fact that the interior is geodesically convex combined with the fact that boundary is a null set and dynamical optimal plans are concentrated on geodesics enables us to invoke the exact same proof of equivalence of ``$\Ric \ge \K$ and $n \le \N$'' with ``$\CD^*\hspace{-1pt}\left(\K, \N \right)$ conditions'' to assert that $\M^{^n}$ is an $\RCD\hspace{-1pt}\left( \K, \N\right)$ space. Said another way, when the boundary is convex, the smooth $\BE\hspace{-1pt}\left(\K, \N \right)$, the $\RCD\hspace{-1pt}\left(\K, \N \right)$ and the $\BE\hspace{-1pt}\left( \K, \N  \right)$ are all equivalent. We however will not go into further details of the proof and leave it as a remark. 
\end{remark}
\par Based on the Remark~\ref{rem:subcase}, the rigidity result for Bakry-\'Emery manifolds is actually a special case of the $\RCD$ case; however, this does not diminish the importance of the proof in the smooth case since the smooth proof clarifies the main strategy and puts much less machinery to use than the nonsmooth case does. It is worth mentioning that a more detailed investigation of the RCD condition for manifolds with boundary has been carried out in~\cite{Han}.
\subsection{Isometric Splitting}
\par The Bakry-\'{E}mery curvature dimension conditions $\BE\left(\K,\N\right)$ for the drift Laplacian $\BEL$ are equivalent to
\begin{align*}
	\BERic \ge \K\g;
\end{align*} 
e.g. see~\cite{sturm-tensor}.
\begin{remark} 
	When $\mathscr{X}$ is globally integrable i.e. when $\mathscr{X}  = \nabla \upvarphi$ for some smooth potential $\upvarphi$, the triple $\left(\M^{^n} , \g , e^{-\upvarphi} \dif\vol
	_{\g}   \right)$ is also called a Bakry-\'{E}mery manifold or a smooth measure space, a manifold with density or a weighted manifold. In this case, the drift Laplacian takes the form
	\begin{align*}
		\Delta_{\, \upvarphi} = \Delta + \nabla \upvarphi \boldsymbol{\cdot} \nabla,
	\end{align*}
	and by a straightforward computation the $\N$-Bakry-\'{E}mery Ricci tensor becomes
	\begin{align*}
		\Ric^{\,\N}_{ \, \upvarphi} = \Ric + \Hess \left( \upvarphi \right)\, \shortminus \, \nicefrac{1}{\left(\N-n\right)} \ \dif\upvarphi \otimes \dif\upvarphi.
	\end{align*}
\end{remark}
\par The Bochner formula, for the drift Laplacian, takes the form
\begin{align}\label{id:bo}
	\nicefrac{1}{2} \, \BEL \left| \nabla f \right|^{^2} = \left\|\Hess \left( f \right)\right\|^{^2}_{_{\textsf{HS}}} + \Ric \left( \nabla f , \nabla f \right) + \nicefrac{1}{2}\ \mathcal{L}_{\substack{\\\mathscr{X}}}  \g \left( \nabla f , \nabla f \right) +  \nabla \BEL f \boldsymbol{\cdot}  \nabla f.  
\end{align}
One can prove the following Bochner inequality \cite[Chapter 14]{V}
\begin{align*}
	\nicefrac{1}{2}\ \BEL \left| \nabla f \right|^{^2}  \ge \nicefrac{1}{\N}\ \left(\BEL \left(f\right) \right)^{\hspace{-2pt}^2}  + \K \left| \nabla f \right|^{^2}   +  \nabla \BEL f \boldsymbol{\cdot}  \nabla f.
\end{align*}
\par When $\mathscr{X} = \nabla \upvarphi$, from the work in~\cite{BQ} or~\cite{AC}, one has the usual eigenvalue comparison estimates (\ref{eq:eig-compar}). In particular, we know that $\Ric^{\,\N}_{ \, \upvarphi} \ge 0$ implies that $\feig \ge \nicefrac{\uppi^{\hspace{0.5pt}^2}}{\diam^{^2}}$.   
\begin{remark}
\par For a general vector field $\mathscr{X}$ and the corresponding drift Laplacian,~\cite[Theorem 14]{BQ} yields the gradient comparison and eigenvalue comparison for the \emph{real spectrum}. It should be noted that for $\mathscr{X}$ that is not a gradient vector field $\BEL$ is not symmetric and might not admit a completely real spectrum however the same lower bound still holds for the bottom (positive) of the \emph{real spectrum} (of course, when the real spectrum is nonempty). Hence, we can still prove the rigidity result stated in Theorem~\ref{thm:main-2}. 
\end{remark}
Now we can state and prove the rigidity in this case. 
\begin{theorem}\label{thm:main-2}
	Suppose $\M^{^n}$ is compact (without boundary or with a convex one) and $\BERic \ge 0$ for $n \le \N$ and a smooth (complete) vector field $\vX$.
	\par  Then, there exists $u$ with $\BEL u=\shortminus \feig u$ for  $\feig = \nicefrac{\uppi^{\hspace{0.5pt}2}}{\diam^2}$ if and only if $n=1$ and $\M^{^1}$ is isometric to either a circle or a line segment with $\vX$ vanishing everywhere. 
\end{theorem}


\subsubsection*{ \bf \textit{Proof of Theorem~\ref{thm:main-2}}}
\par To show $n=1$, we show there exists an open region $\Omega$ within $\M^{^n}$ which consists of one or two open intervals.  We note here that the $\mathcal{C}^\infty$-regularity for the Neumann eigenfunctions of the drift Laplacian follows from the standard elliptic regularity theory; e.g. see~\cite{Gil-Tru}.
\par As before, by rescaling, we can assume $\diam \left( \M^{^n} \right) = \uppi$ and $\feig = 1$. Suppose $u$ is an eigenfunction associated to $\feig = 1$. From the proof of the comparison theorem in~\cite{BQ}, we can assume
\begin{align*}
	\shortminus \min v = \max v = \max u = \shortminus \min u = 1
\end{align*}
Also we have the gradient comparison  
\begin{align*}
	\left| \nabla u \right|^{^2}  \le \left( v^{'} \circ v^{-1}  \right)^{\hspace{-2pt}^2} \circ u,
\end{align*}
where, $v$ is the eigenfunction of the corresponding $1$-dimensional model space; see (\ref{eq:eig-compar})-(\ref{eq:model-2}).
\par In the $1$-dimensional model space, $v(t) = \sin (t) $ hence, 
\begin{align}\label{eq:grad-ineq}
	\left| \nabla u \right|^{^2} \le \left( v^{'} \circ v^{-1}  \right)^{\hspace{-2pt}^2} \circ u = 1 \shortminus u^2.
\end{align}
Let $\alpha = \left| \nabla u \right|^{^2} + u^2$. From above we know that $\alpha \le 1$. 
s before, define the singular sets $ \Sing_{_{-}}:=u^{-1}\left( \shortminus1  \right)$ and $ \Sing_{_{+}}:=u^{-1}\left(1 \right)$ and let 
\begin{align*}
	\mathsf{C}_{\hspace{1pt}u} := \left\{ x\; \text{ \textbrokenbar} \; \;   \nabla u\left(x\right) =0    \right\},
\end{align*}
be the set of critical points of $u$. Then it is straightforward to see that \eqref{eq:grad-ineq} implies
\begin{align*}
	\Sing_{_{-}} \, \sqcup \, \Sing_{_{+}} \subset 	\mathsf{C}_{\hspace{1pt}u}.
\end{align*}
\begin{remark}
	By e.g.~\cite{Hardt-Nad}, we know $\mathsf{C}_{\hspace{1pt}u}$, being the critical set of a nontrivial solution of an elliptic partial differential equation in a compact manifold, has Hausdorff codimension $\ge 2$, unless the manifold is $1$-dimensional where the codimension of the critical set can be $1$. In any case, the critical set is of measure zero.   
\end{remark}
Recall the definition of a horizontal geodesic (Definition~\ref{defn:horiz-geo}).
\begin{lemma}\label{lem:key-BE}
	Let $\upgamma$ be a horizontal geodesic. Then, $\length \left( \upgamma \right) = \uppi$ and the identity
	\begin{align*}
		\left| \nabla u \right|^{^2}  \left( \upgamma\left(t\right) \right) = \left( v^{'} \circ v^{-1} \left(u\right)  \right)^{\hspace{-2pt}^2}  = 1 \shortminus u^{2} \left( \gamma\left(t\right) \right),
	\end{align*}
	holds along $\upgamma$.  Furthermore, 
	\begin{align*}
		\left| \nicefrac{\dif}{\dif t}  \; u \circ \gamma \right| = \left( \left| \nabla u \right| \circ \upgamma\right)    \left|  \bdot{\upgamma} \right| \quad \forall t. 
	\end{align*}
\end{lemma}
\begin{proof}[\footnotesize \textbf{Proof}]
	Along such horizontal geodesic $\upgamma:\left[ 0,l \right] \to \M^{^n}$, one has
	\begin{align*}
		\left| \nicefrac{\dif}{\dif t}  \; u \circ \gamma \right| \le \left( \left| \nabla u \right| \circ \upgamma\right)    \left|  \bdot{\upgamma} \right|.
	\end{align*}	
	Hence,
	\begin{align}\label{eq:horiz-geo-comp}
		\uppi = \diam \ge \int_{\hspace{-1pt}_0}^l \; \left| \bdot{\upgamma} \right| \, \dif t & \ge \int_{ \left| \nabla u \right| \circ \upgamma \neq 0}  \; \nicefrac{\left|\bdot {u \circ \gamma}\right|}{\left| \nabla u \right|\left( \upgamma(t)  \right) }  \;   \dif t \notag \\ & =\lim_{\varepsilon \downarrow 0} \int_{- 1 + \varepsilon}^{1-\varepsilon} \; \nicefrac{1}{| \nabla u |} \; \di u  \\ & \ge \lim_{\varepsilon \downarrow 0} \int_{- 1 + \varepsilon}^{1-\varepsilon}  \; \nicefrac{1}{\left( v^{'} \circ v^{-1} (u)  \right)} \; \di u \notag \\ & =  \lim_{\varepsilon \downarrow 0} \int_{- 1 + \varepsilon}^{1-\varepsilon} \left( 1 \shortminus u^2 \right)^{\hspace{-1pt}^{\shortminus\nicefrac{1}{2}}} = \uppi. \notag
	\end{align}
	This means any such horizontal geodesic $\upgamma$, has length equal to $\uppi$ and furthermore, all the inequalities in (\ref{eq:horiz-geo-comp}) are indeed equalities. Therefore, along $\upgamma$, one has
	\begin{align}\label{eq:id-hor-geo}
		\left| \nabla u \right|\left( \upgamma\left(t\right) \right) = \left( v^{'} \circ v^{-1} (u)  \right) = 1 \shortminus u^2 \left( \upgamma\left(t\right) \right)  \quad \mbox{$\mathscr{L}^{^1}$-a.e.
			$t$}.  
	\end{align}
	Since the functions involved are continuous, the identity in \eqref{eq:id-hor-geo} indeed holds for all $t$. From the chain of equalities in (\ref{eq:horiz-geo-comp}) and the continuity, it follows $	\left| \nicefrac{\dif}{\dif t}  \; u \circ \gamma \right| = \left( \left| \nabla u \right| \circ \upgamma\right)    \left|  \bdot{\upgamma} \right|$ for all $t$.
\end{proof}
\begin{lemma}
	$u$ is strictly monotonically increasing along any horizontal geodesic $\upgamma$.
\end{lemma}
\begin{proof}[\footnotesize \textbf{Proof}]
	Let $\upgamma: \left[ 0,\uppi \right] \to \M^{^n}$ be a unit speed horizontal geodesic.  Suppose $u$ is not strictly increasing along $\upgamma$, then there exists $t \in \left( 0, \uppi \right)$ such that $\left| \nicefrac{\dif}{\dif t}  \; u \circ \gamma \right| = \left( \left| \nabla u \right| \circ \upgamma\right)    \left|  \bdot{\upgamma} \right| = 0 $ therefore $\left| \nabla u \right| \circ \left(\upgamma (t)\right) = 0$ therefore $\gamma(t) \in 	\Sing_{_{-}} \, \sqcup \,  \Sing_{_{+}}$ by Lemma~\ref{lem:key-BE}. If $\gamma(t) \in 	\Sing_{_{-}}$, then $\upgamma\restr_{\left[t,\uppi\right]}$ is a horizontal geodesic hence has length $\pi$; also $\upgamma$ has length $\pi$ which is a contradiction. The case $\upgamma(t) \in \Sing_{_{+}}$ would lead to a contradiction in a similar fashion. 
\end{proof}
\begin{remark}
	Lemma~\ref{lem:key-BE} is indeed stating that the horizontal geodesics in an spectrally-extremal Bakry-\'Emery manifold are themselves spectrally-extremal. This phenomenon which is key for the proof of the rigidity can be shown in spaces with much less regularity; see the key Lemma~\ref{lem:ida} in the $\RCD$ context. 
\end{remark}
\par  In below, wherever there is a risk of confusion, we have used $\nabla$ for gradient and $\nabla_{\bullet}$ for the covariant derivative to differentiate between the two. 
\par The following key lemma provides us with a suitable partial differential inequality, to solutions of which, we can apply the maximum principle. In what follows, $\ring{\M}^{^n}$ denotes the interior of $\M^{^n}$. 
\begin{lemma}\label{lem:PDI}
	Let $u$ be as in above and set
	\begin{align*}
		\upalpha:= \left| \nabla u \right|^{^2} + u^2 \quad \text{and} \quad \upbeta:= \left| \nabla u \right|^{^2} \shortminus u^2.
	\end{align*}
	Then
	\begin{align}\label{eq:pdi}
		\BEL \upalpha \, \shortminus \, \nicefrac{1}{2} \, \left| \nabla u \right|^{^{-2}} \,  \nabla \upalpha \boldsymbol{\cdot} \nabla \upbeta  \ge 2\BERic,
	\end{align}
	holds on $\ring{\M}^{^n} \smallsetminus \mathsf{C}_{\hspace{1pt}u}$ which is an open subset of $\M^{^n}$ and $\ring{\M}^{^n}$. 
\end{lemma}
\begin{proof}[\footnotesize \textbf{Proof}]
	In terms of $1$-forms, we have
	\begin{align*}
		\nabla_{\hspace{-2pt}_\bullet} \left| \nabla u \right|^{^2} = 2  \left( \nabla_{\hspace{-2pt}_\bullet} \nabla u \right) \boldsymbol{\cdot} \nabla u.
	\end{align*}
	Therefore,
	\begin{align*}
		\nabla \upalpha \boldsymbol{\cdot} \nabla \upbeta &= \nabla \left| \nabla u \right|^{^2}  \boldsymbol{\cdot} \nabla \left| \nabla u \right|^{^2} \shortminus \left| \nabla u^2   \right|^{^2}\\ &\le 4 \left\| \Hess \left( u \right) \right\|^{^2}_{_\textsf{HS}}  \left| \nabla u \right|^{^2} \shortminus 4 u^2  \left| \nabla u \right|^{^2},
	\end{align*}
	so,
	\begin{align*}
		\nicefrac{1}{4}\ \left| \nabla u \right|^{^{-2}} \, \nabla \upalpha \boldsymbol{\cdot} \nabla \upbeta \le \left\| \Hess \left( u \right) \right\|^{^2}_{_\textsf{HS}}  \shortminus u^2.
	\end{align*}
	\par Now with the aid of the Bochner formula \eqref{id:bo} and chain rule for $\BEL$, we compute
	\begin{align}\label{eq:const-hess}
		\BEL \upalpha &= \BEL  \left( \left| \nabla u \right|^{^2} \right) + \BEL \left( u^2 \right) \notag \\ &=  2\left\| \Hess \left( u \right) \right\|^{^2}_{_{\textsf{HS}}}  + 2\Ric \left( \nabla u , \nabla u \right) + \mathcal{L}_{\substack{\\\mathscr{X}}}  \g \left( \nabla u , \nabla u \right) + 2 \nabla \BEL u \boldsymbol{\cdot}  \nabla u  + \BEL \left( u^2 \right) \notag \\ &= 2\left\| \Hess \left( u \right) \right\|^{^2}_{_{\textsf{HS}}} \shortminus 2u^2 + 2\Ric \left( \nabla u , \nabla u \right) + \mathcal{L}_{\substack{\\\mathscr{X}}} \g \left( \nabla u , \nabla u \right)  \notag \\ &\ge \nicefrac{1}{2} \ \left| \nabla u \right|^{^{-2}}\, \nabla \alpha \boldsymbol{\cdot} \nabla \beta + 2\BERic \left( \nabla u , \nabla u \right),
	\end{align}
	which finishes the proof of the lemma.
\end{proof}
Now we wish to use the elliptic maximum principle to the elliptic differential inequality (\ref{eq:const-hess}) to assert constancy of $\upalpha$ on some suitable domains. 
\begin{definition}
	A connected component $\Omega$ of  $\ring{\M}^{^n} \smallsetminus \mathsf{C}_{\hspace{1pt}u}$ is said to be admissible whenever there exists a horizontal geodesic $\upgamma$ such that $\Omega \cap \upgamma \neq \emptyset$.
\end{definition}
Let
\begin{align*}
	\mathcal{M} := \left\{ x \; \text{ \textbrokenbar} \; \left| \nabla u \right|^{^2} + \left| u \right|^{^2} =1  \right\} \ne \emptyset. 
\end{align*}
\begin{lemma}\label{lem:con-comp}
	Let $\Omega$ be an admissible connected component of the open domain $\ring{\M}^{^n} \smallsetminus \mathsf{C}_{\hspace{1pt}u}$; then, $\upalpha\restr_{\Omega} \equiv 1$. 
\end{lemma}
\begin{proof}[\footnotesize \textbf{Proof}]
	Since $\Omega$ is admissible, there exists a point $p \in \Omega$ which lies on a horizontal geodesic hence by~\eqref{eq:id-hor-geo}, we have $\alpha(p) = 1$. Therefore, $\mathcal{M} \cap \Omega \neq \emptyset$. 
	\par Being a level set of a Lipschitz continuous function, $\mathcal{M} \cap \Omega$ is a closed subset of $\Omega$. Applying the strong maximum principle (smooth version) to the inequality (\ref{eq:pdi}) (restricted to the open domain $\Omega$) and upon using standard arguments, one deduces $	\mathcal{M} \cap \Omega$ is also open. Therefore, $\mathcal{M} \cap \Omega = \Omega$.  
\end{proof}
\par By Lemma~\ref{lem:con-comp}, we deduce $\mathcal{M} \cap \left( \ring{\M}^{^n} \smallsetminus \mathsf{C}_{\hspace{1pt}u}\right)$ coincides with the union of a number of connected components of the open domain $\ring{\M}^{^n} \smallsetminus \mathsf{C}_{\hspace{1pt}u}$. 
\begin{lemma}\label{lem:component-interior}
	Let $\Omega$ be an admissible connected component of the open domain $\ring{\M}^{^n} \smallsetminus \mathsf{C}_{\hspace{1pt}u}$; then, for any horizontal geodesic $\upgamma$ with $\upgamma \cap \Omega \neq \emptyset$,  the interior of $\upgamma$ is entirely included in $\Omega$. 
\end{lemma}
\begin{proof}[\footnotesize \textbf{Proof}]
	By~\eqref{eq:id-hor-geo}, we deduce the interior of $\upgamma$ is included in~$\mathcal{M} \cap \left( \ring{\M}^{^n} \smallsetminus \mathsf{C}_{\hspace{1pt}u}\right)$. Since the interior of $\upgamma$ is a connected set, it has to be included in exactly one connected component; hence the conclusion follows. 
\end{proof}
\par The Lemma~\ref{lem:component-interior} implies $\mathcal{M} \cap \left( \ring{\M}^{^n} \smallsetminus \mathsf{C}_{\hspace{1pt}u}\right)$ contains the interior of all horizontal geodesics. 
\par Now let $\Upomega$ be a connected component of $\ring{\M}^{^n} \smallsetminus \mathsf{C}_{\hspace{1pt}u}$ with $\alpha\restr_{\Upomega} = 1$ (which is not necessarily admissible!). Note that the vector field $\upzeta := \nicefrac{\nabla u}{\left| \nabla u \right|}$ is well-defined on $\Upomega$. 
\begin{lemma}\label{lem:hess-smooth}
	On $\Upomega$, $\Hess \left( u \right)$ is given by the identities
	\begin{align*}
		\Hess \left( u \right) = \shortminus u \upzeta^{^\flat} \otimes \upzeta^{^\flat} = \shortminus u\left|  \nabla u \right|^{^{-2}} \dif u \otimes \dif u = \shortminus u\left(1 \shortminus u^2  \right)^{\hspace{-1pt}^{-1}} \nabla_{\hspace{-2pt}_{\bullet}} u \otimes \nabla_{\hspace{-2pt}_{\bullet}} u .
	\end{align*}
	In particular, these Hessian identities hold on any admissible connected component $\Omega$.
\end{lemma}
\begin{proof}[\footnotesize \textbf{Proof}]
	Using the inequality~\eqref{eq:const-hess}, constancy of $\upalpha$ and $\BERic \ge 0$, we deduce
	\begin{align}\label{eq:norm-hes}
		\left\| \Hess \left( u \right)\right\|^{^2}_{_{\textsf{HS}}}   = u^2. 
	\end{align}
	\par By the Hessian formula \eqref{eq:hessian}, one computes 
	\begin{align*}
		\Hess \left( u \right) \left( \nabla u , \nabla g \right) = \nicefrac{1}{2} \ \nabla g \boldsymbol{\cdot} \nabla \left| \nabla u \right|^{^2}  = \nicefrac{1}{2} \ \nabla g  \boldsymbol{\cdot} \nabla \left( 1 \shortminus u^2 \right) = \shortminus u  \nabla u  \boldsymbol{\cdot} \nabla g.
	\end{align*}
	In particular, for the vector field $\upzeta$, it holds
	\begin{align}\label{eq:hes}
		\Hess \left( u \right)  \left( \upzeta, \upzeta \right) = \shortminus u.
	\end{align}
	\par Now, (\ref{eq:norm-hes}) and (\ref{eq:hes}) imply for an \emph{orthonormal} frame $\left\{\upzeta = e_{_1}, e_{_2} , \dots ,e_{_n}\right\}$,
	\begin{align*}
		\Hess \left( u \right)  \left(e_{_i}, e_{_i}\right) = 0 \quad  i = 2, \dots , n,
	\end{align*}
	holds which by virtue of \eqref{eq:hes}, gives the desired conclusion. 
\end{proof}
\par Recall on $\ovs{\Upomega}$, we have $\alpha \equiv 1$. Let $f:= \sin^{^{-1}} \circ \, u$. Let the local singular sets $\Sing^{^\pm}_{\ovs{\Upomega}}$ be given by
\begin{align*}
	\Sing^{^\pm}_{\ovs{\Upomega}} :=  \mathsf{C}_{\hspace{1pt}u} \cap \ovs{\Upomega} = u^{-1}\left( \left\{ \plmi 1\right\} \right) \cap \ovs{\Upomega} = f^{^{-1}}\left( \left\{ \plmi \nicefrac{\pi}{2} \right\} \right) \cap \ovs{\Upomega}.
\end{align*}
Set $\Sing_{\ovs{\Upomega}} := \Sing^{^-}_{\ovs{\Upomega}} \ \dot{\sqcup} \ \Sing^{^+}_{\ovs{\Upomega}}$. Consider the local regular set
$
\mathsf{R}_{\Upomega} :=  \Upomega \smallsetminus 	\Sing_{\ovs{\Upomega}}.
$
Since $\Upomega$ is disjoint from $\mathsf{C}_{\hspace{1pt}u}$, we get $\Sing_{\ovs{\Upomega}} \subset \partial \Upomega$.
\begin{lemma}
	$f$ is Hessian-free on $\mathsf{R}_{\Upomega}$; in particular, it is geodesically affine and harmonic within $\mathsf{R}_{\Upomega}$.
\end{lemma}
\begin{proof}[\footnotesize \textbf{Proof}]
	On $\mathsf{R}_{\Upomega}$, we have
	\begin{align*}
		\nabla_{\hspace{-2pt}_{\bullet}} \, f =  \left( 1 \shortminus u^2\right)^{\hspace{-1pt}^{\shortminus \nicefrac{1}{2}}} \nabla_{\hspace{-2pt}_{\bullet}}  u =  \left| \nabla u \right|^{^{-1}} \nabla_{\hspace{-2pt}_{\bullet}}  u = \upzeta^{^{\flat}}.
	\end{align*}
	Therefore, by the aid of Lemma~\ref{lem:hess-smooth}, we get
	\begin{align*}
		\Hess \left( f \right)  =  \nabla_{\hspace{-2pt}_{\bullet}}^{^2} f = u\left( 1 \shortminus u^2\right)^{\hspace{-1pt}^{\shortminus \nicefrac{3}{2}}} \nabla_{\hspace{-2pt}_{\bullet}} u \otimes \nabla_{\hspace{-2pt}_{\bullet}} u + \left( 1 \shortminus u^2 \right)^{\hspace{-1pt}^{\shortminus \nicefrac{1}{2}}} \nabla_{\hspace{-2pt}_{\bullet}}^{^2} u =0.
	\end{align*}
\end{proof}
\begin{lemma}
	Let $\upgamma:[0,\uppi] \to \M^{^n}$ be a (unit speed) horizontal geodesic whose interior is included in $\Upomega$. Then, from the above discussion, we know $\Upomega$ is joining a point $p_{_-}\in \Sing^{^-}_{\ovs{\Upomega}}$ to a point $p_{_+} \in \Sing^{^+}_{\ovs{\Upomega}}$.
	\par Then, the relations
	\begin{align*}
		\Hess\left( u \right)_{p_{_-}} = \left( \bdot{\upgamma}(0) \right)^{^\flat} \otimes \left( \bdot{\upgamma}(0) \right)^{^\flat} \quad \text{and} \quad \Hess\left(u\right)_{p_{_+}} =  \left( \bdot{\upgamma}(\uppi) \right)^{^\flat} \otimes \left( \bdot{\upgamma}(\uppi) \right)^{^\flat},
	\end{align*}
	hold true. In particular, 
	\begin{align*}
		\bdot{\upgamma}\left(0\right) =  \plmi \, \upzeta\left( p_{_-} \right) \quad \text{and} \quad \bdot{\upgamma}\left(\uppi\right) = \plmi \, \upzeta\left( p_{_+} \right).
	\end{align*}
\end{lemma}
\begin{proof}[\footnotesize \textbf{Proof}]
	The same proof of~\cite{Hang-Wang} can be applied here hence the proof is omitted. 
\end{proof}
Let $\left\{  \Omega_{i} \right\}_{i \in \Lambda}$ be the set of all admissible connected components of $\ring{\M}^{^n} \smallsetminus \mathsf{C}_{\hspace{1pt}u}$. Let $\mathsf{Sing}^{^\pm} := \bigcup_{i \in \Lambda} \Sing^{^\pm}_{\ovs{\Omega}_i}$  and set $\mathsf{Sing} = \bigcup_{i \in \Lambda} \Sing^{^\pm}_{\ovs{\Omega}_i} = \mathsf{Sing}^{^-} \, \dot{\sqcup} \, \mathsf{Sing}^{^+}$. 
\begin{lemma}
	$\mathsf{Sing}$ is a finite set and $2 \le \# \left( \mathsf{Sing} \right) \le 4$.
\end{lemma}
\begin{proof}[\footnotesize \textbf{Proof}]
	The proof follows the same argument as in~\cite{Hang-Wang}, hence omitted. 
\end{proof}
\begin{lemma}
	Any admissible connected component $\Omega$ isometrically splits off an interval and consequently, so does $\bigcup_{i \in \Lambda} \Omega_i$. Furthermore, any such component consists of a finite number of horizontal geodesics.
\end{lemma}
\begin{proof}[\footnotesize \textbf{Proof}]
	Let $\Omega$ be an admissible connected component. To show the splitting phenomenon, we argue as follows. The fact that $f$ is Hessian-free means that $\upzeta = \nicefrac{\nabla u}{\left| \nabla u \right|}$ is a parallel vector field in $\Omega$. Thus,
	the integral curves of $\upzeta$ are geodesics that do not intersect one another. 
	\par By invoking standard extension theorem from the theory of ordinary differential equations (e.g. see~\cite{Arnold}), one deduces the following; besides closed orbits in $\Omega$, the (local) geodesics which are maximal integral curves of $\upzeta$, must either have singular points as either of their endpoints, or alternatively must run into the boundary of $\Omega$. Since the vector field is a gradient field with no singularities within $\Omega$, again the standard ode theory implies that there are no closed orbits inside $\Omega$ and furthermore, $f$ is increasing along the maximal flow lines (which are exactly the horizontal geodesics) and the maximal flow line go from absolute minima of $f$ in $\ovs{\Omega}$ to absolute maxima of $f$ in $\ovs{\Omega}$. Also notice that since the flow lines do not intersect each other, the maximal flow lines are indeed minimizing geodesics. 
	\par Putting all these facts together, one deduces the maximal flow lines of $\upzeta$ are minimizing geodesics that connect points form $\Sing^{^-}_{\ovs{\Upomega}} $ to points form $\Sing^{^+}_{\ovs{\Upomega}}$ hence they are horizontal geodesics whose interior is included in $\Omega$. Since the the local singular sets $\Sing^{^\pm}_{\ovs{\Upomega}}$ are finite, we deduce the desired splitting phenomenon and also the fact that $\Omega$ is $1$-dimensional. 
\end{proof}
To finish the proof of Theorem~\ref{thm:main-2}, since $\M^{^1}$ is one dimensional, it is topologically either a circle or an interval and metrically, it has diameter $= \uppi$. 
\par To show the drift term must vanish, first notice that if $\M^{^1}$ is an interval then, $\vX$ is integrable i.e. $\vX = \nabla \upvarphi$ due to simple connectivity. Therefore, $\feig = 1$ is equivalent to 
\begin{align*}
	\feig	\left( \left[ \shortminus \nicefrac{\uppi}{2}, \nicefrac{\uppi}{2} \right],\dist_{_\textsf{Euc}} ,e^{-\upvarphi} \ \dif \mathscr{H}^{^1} \right) = 1.
\end{align*}
At this point, we can apply Lemma~\ref{lem:lemD} to deduce the density is constant. 
\par If $\M^{^1}$ is a circle, from $u^{2} + \left(u^{'}\right)^{^2} = 1$ and the gradient comparison, we deduce $u\left(\uptheta\right) = \sin \left(\theta\right)$ (here $\shortminus \nicefrac{\uppi}{2} \le \uptheta < \nicefrac{\uppi}{2}$ is the angle coordinate measured from $u^{-1}\left(0 \right)$) therefore, in addition to $\BEL u = \shortminus \, u$, we get $\Delta u =  \shortminus  \, u$ hence, $\mathcal{L}_{\substack{\\\mathscr{X}}} \sin \left(\theta\right) \equiv 0$ which implies the vector field $\vX$ must vanish.
\qed

	\bibliographystyle{amsplain}
	\bibliography{RSGNR}
\end{document}